\documentclass[titlepage,leqno,12pt]{report}
\pdfoutput=1

\usepackage{amsmath} 
\usepackage{amsthm}
\usepackage{amsfonts}
\usepackage{amssymb}
\usepackage{mathtools}
\usepackage[scr=rsfs]{mathalpha}
\usepackage[english]{babel}
\usepackage{epigraph}
\usepackage[utf8]{inputenc} 
\usepackage{csquotes}
\usepackage{lmodern} 
\usepackage{graphicx} 
\graphicspath{{Figures/}}
\usepackage{xcolor}
\usepackage{float}
\usepackage{lastpage}
\usepackage{caption}
\usepackage{tensor}
\usepackage{tikz}
\usepackage{forest,textcomp}
\usepackage{cancel}
\usepackage[noadjust]{cite}
\usepackage{subcaption}
\usepackage{fancyref}
\usepackage{subfiles}
\usepackage{gensymb}
\usepackage{enumerate}
\usepackage[bottom]{footmisc}
\usepackage{url}
\usepackage{setspace}
\usepackage[colorlinks]{hyperref}
\hypersetup{linkcolor=blue, urlcolor=blue, citecolor=red}

\usepackage[a4paper,width=150mm,top=30mm,bottom=30mm]{geometry}

\usepackage{fancyhdr}
\newcommand{\navn}{Martin Hampenberg Christensen}

\newcommand{\dato}{January, 2025}

\newcommand{\set}[1]{\relax\ifmmode{\{\,#1\,\}}\else{$\{\,#1\,\}$}\fi}
\newcommand{\given}{\mid}
\newcommand*\from{\colon}

\newcommand{\h}{\frac{1}{2}}

\newcommand{\paa}[1]{\left(#1\right)}

\let\oldhat\hat
\renewcommand{\hat}[1]{\oldhat{\mathbf{#1}}}

\renewcommand{\i}{\text{i}}

\newcommand{\identity}{\relax\ifmmode{\mathbb{1}}\else{$\mathbb{1}$}\fi}

\newcommand{\bigO}[1]{\relax\ifmmode{\mathcal{O}(#1)}\else{$\mathcal{O}(#1)$}\fi}
\newcommand{\inv}[1]{\relax\ifmmode{{#1}^{-1}}\else{${#1}^{-1}$}\fi}
\newcommand{\dom}[1]{\relax\ifmmode\operatorname{dom}(#1)\else\text{dom}$(#1)$\fi}
\newcommand{\im}[1]{\relax\ifmmode\operatorname{im}(#1)\else\text{im}$(#1)$\fi}
\newcommand{\fix}[1]{\relax\ifmmode\operatorname{fix}(#1)\else\text{fix}$(#1)$\fi}
\newcommand{\closure}[1]{\relax\ifmmode\operatorname{cl}(#1)\else\text{cl}$(#1)$\fi}
\newcommand{\verset}[1]{\relax\ifmmode\operatorname{V}(#1)\else\text{V}$(#1)$\fi}
\newcommand{\ary}[1]{\relax\ifmmode\operatorname{dim}(#1)\else\text{dim}$(#1)$\fi}
\newcommand{\Sym}[1]{\relax\ifmmode\operatorname{Sym}(#1)\else\text{Sym}$(#1)$\fi}
\newcommand{\Stab}[1]{\relax\ifmmode\operatorname{Stab}(#1)\else\text{Stab}$(#1)$\fi}
\newcommand{\AStab}[1]{\relax\ifmmode\operatorname{AStab}(#1)\else\text{AStab}$(#1)$\fi}
\newcommand{\PStab}[1]{\relax\ifmmode\operatorname{PStab}(#1)\else\text{PStab}$(#1)$\fi}
\newcommand{\setStab}[2]{\relax\ifmmode{#1}_{\lbrace #2\rbrace}\else{$#1$}\textsubscript{$\lbrace #2\rbrace$}\fi}
\newcommand{\pointStab}[2]{\relax\ifmmode{#1}_{(#2)}\else{$#1$}\textsubscript{$(#2)$}\fi}
\newcommand{\Iso}[1]{\relax\ifmmode\operatorname{Iso}(#1)\else\text{Iso}$(#1)$\fi}
\newcommand{\Endo}[1]{\relax\ifmmode\operatorname{End}(#1)\else\text{End}$(#1)$\fi}
\newcommand{\Aut}[1]{\relax\ifmmode\operatorname{Aut}(#1)\else\text{Aut}$(#1)$\fi}
\newcommand{\pAut}[1]{\relax\ifmmode\operatorname{pAut}(#1)\else\text{pAut}$(#1)$\fi}
\newcommand{\ipEnd}[1]{\relax\ifmmode\operatorname{ipEnd}(#1)\else\text{IpEnd}$(#1)$\fi}
\newcommand{\orbit}[1]{\relax\ifmmode\operatorname{Orb}(#1)\else\text{Orb}$(#1)$\fi}
\newcommand{\order}[1]{\relax\ifmmode\mathcal{O}_{#1}\else$\mathcal{O}_{#1}$\fi}
\newcommand{\n}{\{\,0,1,\ldots,n-1\,\}}
\newcommand{\N}{\relax\ifmmode{\mathbb{N}}\else{$\mathbb{N}$}\fi}
\newcommand{\Z}{\relax\ifmmode{\mathbb{Z}}\else{$\mathbb{Z}$}\fi}
\newcommand{\Q}{\relax\ifmmode{\mathbb{Q}}\else{$\mathbb{Q}$}\fi}
\newcommand{\R}{\relax\ifmmode{\mathbb{R}}\else{$\mathbb{R}$}\fi}
\newcommand{\GD}{\relax\ifmmode{\mathcal{D}}\else{$\mathcal{D}$}\fi}
\newcommand{\GH}{\relax\ifmmode{\mathcal{H}}\else{$\mathcal{H}$}\fi}
\newcommand*{\zf}{\relax\ifmmode\mathbf{ZF}\else\textbf{ZF}\fi}
\newcommand*{\zfc}{\relax\ifmmode\mathbf{ZFC}\else\textbf{ZFC}\fi}
\newcommand*{\ac}{\relax\ifmmode\mathbf{AC}\else\textbf{AC}\fi}
\newcommand{\strt}{\rule{0em}{.6em}} 
\newcommand*{\bsleq}[1][]{\preccurlyeq_{\strt #1}}
\newcommand*{\bsless}[1][]{\prec_{\strt #1}}
\newcommand*{\bsgeq}[1][]{\succcurlyeq_{\strt #1}}
\newcommand*{\bsgreat}[1][]{\succ_{\strt #1}}
\newcommand*{\bsequal}[1][]{\approx_{\strt #1}}
\newcommand{\fin}{\relax\ifmmode{\mathfrak{F}}\else{$\mathfrak{F}$}\fi}
\newcommand{\partition}{\relax\ifmmode{\mathcal{P}}\else{$\mathcal{P}$}\fi}
\newcommand{\filter}{\relax\ifmmode{\mathcal{F}}\else{$\mathcal{F}$}\fi}
\newcommand{\ideal}{\relax\ifmmode{\mathcal{I}}\else{$\mathcal{I}$}\fi}
\newcommand{\tree}[1]{\relax\ifmmode{\mathcal{S}(#1)}\else{$\mathcal{S}(#1)$}\fi}
\newcommand{\regtree}[1]{\relax\ifmmode{T^{(#1)}}\else{$T^{(#1)}$}\fi}
\newcommand{\prune}{\relax\ifmmode{F}\else{$F$}\fi}
\newcommand{\relStructure}{\relax\ifmmode{\mathcal{R}}\else{$\mathcal{R}$}\fi}
\newcommand{\relations}{\relax\ifmmode{\mathbf{R}}\else{$\mathbf{R}$}\fi}
\newcommand{\topI}{\relax\ifmmode{\mathcal{I}}\else{$\mathcal{I}$}\fi}
\newcommand{\topology}{\relax\ifmmode{\mathcal{T}}\else{$\mathcal{T}$}\fi}
\newcommand{\card}[1]{\relax\ifmmode\lvert{#1}\rvert\else$\lvert{#1}\rvert$\fi}
\newcommand{\rank}[1]{\relax\ifmmode\operatorname{rank}(#1)\else\text{rank}$(#1)$\fi}
\newcommand{\genset}[1]{\relax\ifmmode\langle#1\rangle\else$\langle#1\rangle$\fi}
\newcommand{\eval}[1]{\left.\vphantom{\h}#1\right\rvert}
\newcommand{\arrows}[3]{\rightarrow (#1)_{#3}^{#2}}

\renewcommand{\iff}{\Leftrightarrow}
\renewcommand{\implies}{\Rightarrow}

\makeatletter
\newcommand\niton{\mathrel{\m@th\mathpalette\canc@l\owns}}
\newcommand\canc@l[2]{{\ooalign{$\hfil#1/\mkern1mu\hfil$\crcr$#1#2$}}}
\makeatother

\theoremstyle{plain}
\newtheorem{theorem}{Theorem}[chapter]
\newtheorem{lemma}[theorem]{Lemma}
\newtheorem{proposition}[theorem]{Proposition}
\newtheorem{conjecture}[theorem]{Conjecture}
\newtheorem{corollary}[theorem]{Corollary}

\theoremstyle{definition}
\newtheorem{definition}[theorem]{Definition}

\theoremstyle{remark}

\DeclareFontFamily{U}{MnSymbolD}{}
\DeclareFontShape{U}{MnSymbolD}{m}{n}{
	<-6>  s*[1.2] MnSymbolD5
	<6-7>  s*[1.2] MnSymbolD6
	<7-8>  s*[1.2] MnSymbolD7
	<8-9>  s*[1.2] MnSymbolD8
	<9-10> s*[1.2] MnSymbolD9
	<10-12> s*[1.2] MnSymbolD10
	<12->   s*[1.2] MnSymbolD12}{}
\DeclareFontShape{U}{MnSymbolD}{b}{n}{
	<-6>  s*[1.2] MnSymbolD-Bold5
	<6-7>  s*[1.2] MnSymbolD-Bold6
	<7-8>  s*[1.2] MnSymbolD-Bold7
	<8-9>  s*[1.2] MnSymbolD-Bold8
	<9-10> s*[1.2] MnSymbolD-Bold9
	<10-12> s*[1.2] MnSymbolD-Bold10
	<12->   s*[1.2] MnSymbolD-Bold12}{}
\DeclareSymbolFont{MnSyD}{U}{MnSymbolD}{m}{n}
\SetSymbolFont{MnSyD}{bold}{U}{MnSymbolD}{b}{n}

\DeclareMathSymbol\prec{\mathrel}{MnSyD}{"68}
\DeclareMathSymbol\succ{\mathrel}{MnSyD}{"69}
\DeclareMathSymbol\preceq{\mathrel}{MnSyD}{"6A}
\DeclareMathSymbol\succeq{\mathrel}{MnSyD}{"6B}
\DeclareMathSymbol\preccurlyeq{\mathrel}{MnSyD}{"6C}
\DeclareMathSymbol\succcurlyeq{\mathrel}{MnSyD}{"6D}
\DeclareMathSymbol\precsim{\mathrel}{MnSyD}{"6E}
\DeclareMathSymbol\succsim{\mathrel}{MnSyD}{"6F}
\DeclareMathSymbol\precapprox{\mathrel}{MnSyD}{"70}
\DeclareMathSymbol\succapprox{\mathrel}{MnSyD}{"71}
\DeclareMathSymbol\nprec{\mathrel}{MnSyD}{"E0}
\DeclareMathSymbol\nsucc{\mathrel}{MnSyD}{"E1}
\DeclareMathSymbol\npreceq{\mathrel}{MnSyD}{"E2}
\DeclareMathSymbol\nsucceq{\mathrel}{MnSyD}{"E3}
\DeclareMathSymbol\npreccurlyeq{\mathrel}{MnSyD}{"E4}
\DeclareMathSymbol\nsucccurlyeq{\mathrel}{MnSyD}{"E5}
\DeclareMathSymbol\nprecsim{\mathrel}{MnSyD}{"E6}
\DeclareMathSymbol\nsuccsim{\mathrel}{MnSyD}{"E7}
\DeclareMathSymbol\nprecapprox{\mathrel}{MnSyD}{"E8}
\DeclareMathSymbol\nsuccapprox{\mathrel}{MnSyD}{"E9}

\numberwithin{equation}{section}

\newenvironment{nalign}{
	\begin{equation}
		\begin{aligned}
		}{
		\end{aligned}
	\end{equation}
	\ignorespacesafterend
}

\begin{document}
	
\begin{titlepage}
	\begin{center}
		{\scshape\LARGE \vphantom{University of Hertfordshire} \par}
		\vspace{0.5cm}
		{\huge\bfseries Submonoids of Infinite \\ Symmetric Inverse Monoids \par}
		\vspace{1cm}
		{\Large \textit{Generating Uncountable Semigroups of \\ Partial Permutations} \par}
		\vspace{1.5cm}
		{\Large \navn \par}
		\vspace{3.5cm}

        {
        \textit{A dissertation submitted to the \\ University of Hertfordshire \\ in partial fulfilment of the requirements for the degree of \\ Doctor of Philosophy}
        }
        
        \vspace{2cm}

        {
        Mathematics and Theoretical Physics \\
        School of Physics, Engineering and Computer Science \\
        University of Hertfordshire
        }

		\vspace{2cm}
        
		{\large \dato \par}
	\end{center}
\end{titlepage}

\chapter*{Abstract}

In this thesis we study the subsemigroup structure of the symmetric inverse monoid $I_X$, the inverse semigroup of bijections between subsets of the set $X$, when $X$ is an infinite set. We explore three different approaches to this task.

First, we classify the maximal subsemigroups of $I_X$ containing some of the most well-known subgroups of the symmetric group $\Sym{X}$.
The subgroups in question are the symmetric group itself, the pointwise stabiliser of a finite non-empty subset of $X$, the stabiliser of an ultrafilter on $X$, and the stabiliser of a finite partition of $X$.
As part of this process, we develop methods for determining whether a given list of subsets of $I_X$ is the set of all maximal subsemigroups of $I_X$ containing a fixed and `sufficiently large' inverse subsemigroup of $I_X$.

Next, we study subsemigroups of $I_X$ which are closed in semigroup topologies on $I_X$ introduced by Elliot et al. in 2023.
We discover that the closed subsemigroups in these topologies that contain all the idempotents of $I_X$ coincide exactly with semigroups of partial endomorphisms and partial automorphisms of relational structures defined on $X$.
Furthermore, we show that if a relational structure $\relStructure$ on a countable set $X$ only contains a finite number of relations, then there exists a finite subset $U$ of $I_X$ such that the union of the partial automorphisms of $\relStructure$ with $U$ generates all of $I_X$.

Finally, we study the subsemigroup structure of $I_X$ under a preorder introduced by George Bergman and Saharon Shelah in 2006 for $\Sym{X}$.
Extending the preorder to $I_X$, if $S_1, S_2$ are subsemigroups of $I_X$, we say that $S_1 \bsleq[I] S_2$ if there exists a finite subset $U \subseteq I_X$ such that $S_1$ is contained in the semigroup generated by $S_2 \cup U$.
We classify certain types of subsemigroups of $I_X$ according the Bergman-Shelah preorder, and we formulate a conjecture analogous to the main result by Bergman and Shelah about the closed subgroups of $\Sym{X}$.

\chapter*{Acknowledgements}
I would like to thank the many people who made this thesis possible.

Thank you to my primary supervisor Yann Péresse, whose time and expertise I have made extensive use of.
Thanks to your guidance, I can now comfortably call myself a mathematician.

Thank you to my secondary supervisors Catarina Carvalho and Jim Collett, whom I have relied on greatly for help with all things administrative.

Thank you to Vicky Gould and Nik Ru\v{s}kuc for your work as external progress reviewers for my thesis and the interest that both of you have shown in my work.

Thank you to James East, whose {\LaTeX} preambles permeate this thesis.

Thank you to Charles Young and Tomasz Łukowski for helping me and everyone else at the Mathematics and Theoretical Physics group settle into our duties and life at the university in general.

Thank you to my fellow PhD students Jonah Stalknecht, Tommaso Franzini, and Simon Jonsson.
Though we work within different fields, yours being physics and mine mathematics, we were still able to support each other and talk about our individual projects.

Thank you to my parents, who always support me unconditionally in everything I do.
Even if you do not understand what I am doing.

Thank you to my aunt and uncle, Ulla and Wolfgang Kirschning, for the moral and financial support throughout my studies.

Thank you to Emil Vyff, Holger Heebøll, Jacob Kirketerp, Anders Frederiksen, and Jens Phannoi for for all of the support and good times, even when we live in different parts of the world.

Thank you to the people of the Smol Ame Speedrunning community, thanks to whom I was able to complete this thesis with most of my sanity intact.

Again, thank you to Nik Ru\v{s}kuc and Charles Young for agreeing to be examiners for my viva.

Finally, thank you to anyone who reads this thesis.

\tableofcontents

\chapter{Introduction}
\epigraph{
    \textit{``In modern mathematics the investigation of the symmetries of a given mathematical structure has always yielded the most powerful results.''
    }}{
    Emil Artin \cite{artin2016geometric}
}
Symmetry is one of the most powerful tools available in the natural sciences.
In Mathematics, symmetry plays an equally important role as emphasised in the epigraph above.
Artin follows the quote up with a definition; \textit{``Symmetries are maps which preserve certain properties''} \cite{artin2016geometric}.
This understanding of symmetry as property preserving maps has lead most mathematicians to effectively equate symmetry with Group Theory.
Groups are algebraic objects that are, in some sense, concerned with \emph{permutations} and \emph{automorphisms}.
Group Theory is perhaps one of the most successful topics within the field of Algebra, with applications found in virtually every scientific field.
Examples of such applications range from cryptography and error-correction in Information Theory to polymers and elementary particles in Physics.
In fact, the Standard Model of modern fundamental Physics is expressed in terms of a symmetry group.

However, while groups are certainly useful and do a good job of describing many symmetries, it has been strongly argued that they do not encapsulate the entire concept of symmetry.
This concern regarding the failings of group theory to fully describe what people mean to by `symmetry' is expressed by Alan Weinstein in \cite{weinstein1996groupoids}:

\begin{displayquote}
    \textit{
        ``Mathematicians tend to think of symmetry as being virtually synonymous with the theory of groups{\dots} 
        In fact, though groups are indeed sufficient to characterize homogeneous structures, there are plenty of objects which exhibit what we clearly recognize as symmetry, but which admit few or no non-trivial automorphisms.''
    }
\end{displayquote}

So, what is a suitable replacement for groups in the study of symmetry?
While there is no definitive answer to this question, a strong candidate is a type of algebraic objects known as \emph{inverse semigroups}.
These are a type of semigroups, which manage to stay very `group-like' by insisting that all elements have an inverse.
Inverse semigroups were first introduced independently by Viktor Wagner and Gordon Preston in 1952 \cite{wagner1952groups} and 1954 \cite{preston1954inverse, preston1954ideals, preston1954representations} respectively. 
They have since found numerous applications in fields such as quasicrystals, tilings, and $C^*$-algebras \cite{paterson2012groupoids}.
The symmetries captured by inverse semigroups are called \emph{partial symmetries} (a generalisation of the \emph{global symmetries} described by groups).
In \cite{lawson_inverse_semigroups}, Mark Lawson describes partial symmetries as those that consider \textit{``the relationship between the parts and the whole.''}
He further defines a partial symmetry of a structure as an isomorphism between two substructures of the structure, which is also called a \emph{partial automorphism}.

Chief among these inverse semigroups are the \emph{symmetric inverse monoids}.
Given a set $X$, the symmetric inverse monoid on $X$ is the set of all bijections between subsets of $X$ and we denote it by $I_X$.
This semigroup is to inverse semigroups what the \emph{symmetric group} $\Sym{X}$, the group of all permutations on the set $X$, is to groups.
As such, there is an argument to be made that the concept of partial symmetries can be studied by classifying symmetric inverse monoids and their subsemigroups, which is what we will do in this thesis.
However, when $X$ is an infinite set, $I_X$ is an uncountably infinite inverse semigroup, in which case the space of subsemigroups is far too large to explore fully.
In cases such as this, there are two typical approaches one can take.
The first approach is to find examples of particularly interesting substructures that satisfy a given condition.
The second approach is to `squint your eyes' and treat classes of substructures as if they are the same structure.
We will be taking both of these approaches in this thesis.
For the first approach, we will be classifying maximal subsemigroups of $I_X$ containing certain interesting subgroups, when $X$ is a countably infinite set.
This is similar to research done in works such as \cite{ball1966maximal, ball1968indices, baumgartner1993maximal, biryukov2000set, brazil1994maximal, covington1996some, macpherson1993large, subgroups_macpherson_neumann, maximal_macpherson_preager, richman1967maximal}, where the maximal subgroups of $\Sym{X}$ have been studied.
This also extends previous work by Xiuliang in \cite{xiuliang1999classification}, in which all the maximal subsemigroups of $I_X$ are classified for $X$ any finite set.
For the second approach, we study the subsemigroup structure of $I_X$ under the so-called \emph{Bergman-Shelah equivalence relation}.
Given two subsemigroups $S_1$ and $S_2$ of a semigroup $T$, we say that $S_1$ and $S_1$ are \emph{equivalent} if there exits a finite subset $U$ of $T$ such that $S_1 \cup U$ and $S_2 \cup U$ generate the same subsemigroup of $T$.
This equivalence relation was first introduced by George Bergman and Saharon Shelah in their 2006 paper \cite{Bergman_2006}, where they studied the equivalence classes of subgroups of the symmetric group under this relation.
What they discovered was that, when $X$ is countably infinite, the subgroups of $\Sym{X}$ that are closed in the \emph{pointwise topology} fall into four distinct equivalence classes.
This pointwise topology is in a sense the `canonical' or `natural' choice for a topology on $\Sym{X}$.
Similarly, two particularly `natural' topologies for $I_X$ have recently been discovered in \cite{part1topological}, called $\topI_1$ and $\topI_4$.
We will also explore the closed subsemigroups of $I_X$ in these topologies and study their relation to relational structures as well as their equivalence classes under the Bergman-Shelah equivalence relation.

\section{Outline}

This dissertation consists of three main chapters and a chapter on background material.
We outline the contents of these chapters here.

\begin{itemize}
    \item[] \textbf{Chapter \ref{chap:background}}: 
    In this chapter we cover the background material needed to read this thesis.
    This mostly consists of giving definitions of some of the concepts from Set Theory and Semigroup Theory that we will be using.
    We also give a brief introduction to Topological Algebra and the Bergman-Shelah preorder.

    \item[] \textbf{Chapter \ref{chap:maximal_subsemigroups}}:
    In this chapter we classify the maximal subsemigroups and maximal inverse subsemigroups of $I_X$ containing certain interesting subgroups of $\Sym{X}$.
    The subgroups in question are the symmetric group $\Sym{X}$ itself, the pointwise stabiliser of a finite non-empty subset $\Sigma$ of $X$, the stabiliser of an ultrafilter $\filter$ on $X$, and the stabiliser of a finite partition $\partition$ of $X$.
    The chapter is structured so that we first develop a general method for classifying maximal subsemigroups of $I_X$ containing `large' inverse subsemigroups, and we then spend the rest of the chapter applying this method to the examples given above.

    \item[] \textbf{Chapter \ref{chap:relational_structures}}:
    In this chapter we study partial automorphisms of relational structures and their connection to the topologies $\topI_1$ and $\topI_4$.
    In the first section of the chapter we introduce semigroups of partial endomorphisms and partial automorphisms on relational structures, and we prove theorems relating these to the set of closed subsemigroups of $I_X$ in the topologies $\topI_1$ and $\topI_4$.
    In the second section we classify a large portion of these semigroups that are equivalent to $I_X$ under the Bergman-Shelah equivalence relation.

    \item[] \textbf{Chapter \ref{chap:BS_IX}}:
    In this final chapter we study the Bergman-Shelah equivalence classes of $I_X$ more broadly.
    First, we extend results for $\Sym{X}$ from \cite{Bergman_2006} to $I_X$, mainly by looking at how orbit and cone sizes affect how `big' a subsemigroup of $I_X$ can at most be.
    Next, we introduce a type of inverse subsemigroup of $I_X$ defined on rooted trees, and we fully classify all possible such semigroups under the Bergman-Shelah equivalence relation when $X$ is countably infinite.
    Finally, we establish a list of conditions that describe the, in some sense, `well-behaved' subsemigroups if $I_X$ and formulate a conjecture about $I_X$ analogous to the main result of Bergman and Shelah for $\Sym{X}$ in \cite{Bergman_2006}.
\end{itemize}

\chapter{Background} \label{chap:background}

In this chapter, we outline the necessary preliminaries for reading and understanding this thesis.
This includes both introductory background material and the research which served as motivation for this project.
The chapter is structured so that the earlier sections serve mostly as background material, and the later sections delve deeper into more specialised topics.

\section{Set Theory}

This thesis operates within the framework known as Zermelo-Fraenkel Set Theory with the Axiom of Choice (abbreviated as \zfc).
For a brief and simple introduction to \zfc, see \cite{halmos1960naive}. For a more thorough coverage of the field of Set Theory, see e.g. \cite{Jech2003}.

\subsection{Conventions}
A proper introduction to the fundamental concepts of Axiomatic Set Theory is beyond the scope of this text.
As such, we will assume that the reader is familiar with everything presented in \cite{halmos1960naive}.
However, we will still briefly summarise some notation and simple definitions:

The \emph{cardinality} of a set $X$ will be denoted $\card{X}$, and a \emph{moiety} is a subset $Y$ of $X$ such that $\card{Y} = \card{X} = \card{X \setminus Y}$. 
A \emph{cardinal} is understood to be the least ordinal of its cardinality, and an \emph{ordinal} to be the well-ordered (see Definition \ref{def:well_order}) set of all smaller ordinals. 
We will denote the set of natural numbers by $\omega$, the smallest infinite ordinal, and its associated cardinal by $\aleph_0$ (also known as \emph{countably infinite}).
The first \emph{uncountable} ordinal will be denoted by $\omega_1$ (the set of all countable ordinals) and its associated cardinal by $\aleph_1$.
Each natural number itself will also be treated as an ordinal (as well as a cardinal), and we will therefore regularly use notation such as $i \in 5 = \set{0,1,2,3,4}$ to indicate that $i$ is a natural number strictly less than 5. 
Given an ordinal $\alpha$ and a cardinal $\kappa$, we denote the \emph{successor ordinal} of $\alpha$ by $\alpha+1 = \set{\beta \given \beta \leq \alpha}$ (the least ordinal which is strictly greater than $\alpha$) and the \emph{successor cardinal} of $\kappa$ by $\kappa^+$ (the least cardinal which is strictly greater than $\kappa$).
An infinite ordinal is called a \emph{limit ordinal} if it not the successor of any ordinal.
An infinite cardinal $\kappa$ is called \emph{regular} if every unbounded subset $A \subseteq \kappa$ has cardinality $\kappa$, and it is called \emph{singular} otherwise.
Finally, we will be composing functions from left to right, so function arguments will be written to the left of the function (same as in e.g. \cite{Bergman_2006, part0topological, part1topological}).

\subsection{Induction and the Axiom of Choice}

As mentioned above, we will be assuming the \emph{Axiom of Choice} (\ac) throughout this text.
One way one stating {\ac} is as follows (as found in \cite{halmos1960naive}).
\begin{itemize}
    \item[] \textbf{Axiom of Choice}: The Cartesian product of a non-empty family of non-empty sets is non-empty.
\end{itemize}

{\ac} has been a somewhat controversial axiom, as it only postulates the existence of a set but does not define how to construct it \cite{Jech2003}.
However, the axiom is immensely useful when working with infinite sets and we will make thorough use of it throughout this text.
Or, to be more specific, we will explicitly make use of some results that are equivalent to assuming {\ac}.
These results are the famous Well-Ordering Theorem and Zorn's Lemma.

\begin{theorem}[Zermelo's Well-Ordering Theorem {\cite[Theorem 5.1]{Jech2003}}] \label{well-ordering theorem}
    Every set can be well-ordered.
\end{theorem}

\begin{theorem}[Zorn's Lemma {\cite[Theorem 5.4]{Jech2003}}] \label{zorn's Lemma}
    If $(X,\leq)$ is a non-empty partially ordered set such that every chain in $(X,\leq)$ has an upper bound, then $(X,\leq)$ has a maximal element.
\end{theorem}

Theorems \ref{well-ordering theorem} and \ref{zorn's Lemma} make use of the concept of an ordering, which we will define later (see Definitions \ref{def:partial_order}, \ref{def:total_order}, and \ref{def:well_order}).
For now we will list some of the consequences of {\ac} and the above theorems, which will be useful in this thesis.

\begin{itemize}
    \item Since every set can be well-ordered, our definition of cardinals as being the smallest ordinal of their cardinality is well-defined.

    \item Every cardinal has a successor cardinal.

    \item If $\kappa$ and $\lambda$ are infinite cardinals, then $\kappa + \lambda = \kappa \lambda = \max(\kappa,\lambda)$.

    \item If $X$ is an infinite set and $\kappa$ a cardinal such that $1 \leq \kappa \leq \card{X}$, then $X$ can be partitioned (see Definition \ref{def:partition}) into $\kappa$ moieties. 

    \item If $X$ is an infinite set, then there exist $2^{\card{X}}$ distinct moieties of $X$.
\end{itemize}

With the full power of the {\zfc} axioms behind us, we can also extend the usual notion of induction on the natural number.

\begin{theorem}[Transfinite Induction {\cite[Theorem 2.14]{Jech2003}}] \label{def:transfinite_induction}
    Let $P$ be a property such that $P(\alpha)$ is well-defined for all ordinals $\alpha$. If $P(\beta)$ being true for all $\beta <\alpha$ implies that $P(\alpha)$ is also true, then $P$ is true for all ordinals.
\end{theorem}

Theorem \ref{def:transfinite_induction} is used for what is called \emph{proof by transfinite induction}, which can be used to prove that a given statement holds for all ordinals (similar to how usual induction can be used to show that something holds true for all natural numbers).
Usually such proofs are divided into three parts, though technically this is not necessary.
\begin{itemize}
    \item[] \textbf{Zero case}: Show that $P(0)$ holds.

    \item[] \textbf{Successor case}: Show that if $P(\alpha)$ holds for any ordinal $\alpha$, then $P(\alpha+1)$ holds.

    \item[] \textbf{Limit case}: Show that if $P(\beta)$ holds for all $\beta < \alpha$, where $\alpha$ is any limit ordinal, then $P(\alpha)$ holds.
\end{itemize}

\subsection{Relations}

Relations are important objects in modern mathematics and they play a central role in this thesis.

\begin{definition}[Relation]
    Let $n$ be a natural number. Then an \emph{$n$-ary relation} $R$ is any subset of a Cartesian product $\bigtimes_{i \in n} X_i$, where $X_i$ is a set for all $i \in n$. 
    We denote the arity of a relation $R$ by dim$(R)=n$ and if $X_i = X$ for all $i \in n$ we say that $R$ is a \emph{relation on $X$}.
\end{definition}

Among the many different types of relations, the most prevalent are certainly binary relations (here binary is simply another term for 2-ary).
There are many different properties and types of binary relations, some of which we shall cover shortly.
One such property is the notion of the \emph{domain} and \emph{image} of a binary relation.

\begin{definition}[Domain] \label{def:domain}
    Let $R \subseteq X \times Y$ be a binary relation.
    Then the set of points $\set{x \in X \given (\exists y \in Y)~ (x,y) \in R}$ is called the \emph{domain of $R$}.
    We denote the domain of $R$ by $\dom{R}$.
\end{definition}

\begin{definition}[Image] \label{def:image}
    Let $R \subseteq X \times Y$ be a binary relation.
    Then the set of points $\set{y \in Y \given (\exists x \in X)~ (x,y) \in R}$ is called the \emph{image of $R$}.
    We denote the image of $R$ by $\im{R}$.
\end{definition}

For any binary relation we can also define the notion of an \emph{inverse}.

\begin{definition}[Inverse] \label{def:inverse}
    Let $R \subseteq X \times Y$ be a binary relation.
    Then the relation $\set{(y,x) \in Y \times X \given (x,y) \in R}$ is called the \emph{inverse} of $R$ and is usually denoted by $\inv{R}$.
\end{definition}

It follows from Definition \ref{def:inverse} that $\inv{(\inv{R})} = R$ for any binary relation $R$.
Given two binary relations we can define an operation called \emph{composition}.

\begin{definition}[Composition] \label{def:composition}
    Let $R \subseteq X \times Y$ and $S \subseteq V \times W$ be binary relations.
    Then we can \emph{compose} $R$ and $S$ to get the \emph{composite} binary relation $R \circ S$ defined as follows:
    \begin{equation*}
        R \circ S = \set{(x,w) \in X \times W \given (\exists u \in Y \cap V)~ (x,u) \in R ~\land~ (u,w) \in S}
    \end{equation*}
\end{definition}

Note that composing any binary relation $R$ with its inverse results in a relation for which each element in the domain of $R$ is related to itself.
That is, $\set{(x,x) \given x \in \dom{R}} \subseteq R \circ \inv{R}$.
A relation purely of the form $\set{(x,x) \given x \in X}$ is known as the \emph{identity map} on $X$.
Further note, that Definition \ref{def:composition} is consistent with our notion of function composition going from left to right.
On the topic of functions, we can (and will) define functions as a type of binary relations.

\begin{definition}[Function] \label{def:function}
    A \emph{function} $f \subseteq X \times Y$ is a binary relation with the following properties:
    \begin{itemize}
        \item[] \textbf{Functional}: For all $x \in X$ and all $y,z \in Y$, if $(x,y) \in f$ and $(x,z) \in f$ then $y=z$.
        
        \item[] \textbf{Total}: The domain of $f$ is all of $X$. That is, $\dom{f} = X$.
    \end{itemize}
    We say that $f$ is a \emph{function from $X$ to $Y$}, which we write as $f \from X \to Y$, and if $(x,y) \in f$ we can instead express this fact as $(x)f = y$ (or simply $xf = y$).
    If in addition $X=Y$, then $f$ is sometimes referred to as a \emph{transformation}.
\end{definition}

From Definition \ref{def:function} we see that whether a given binary relation is considered to be a function depends on the context of what it is a relation on due to the condition of totality.
Two binary relations $f$ and $g$ might be identical as sets, but only one is considered a function if they are relations on different sets.
A binary relation which is functional, but not necessarily total, is called a \emph{partial function} (whereas the usual functions are sometimes referred to as \emph{total functions}).
Given two sets $X$ and $Y$ we write the set of all total functions from $X$ to $Y$ as $Y^X$ (though we will also use the exponentiation notation $y^x$ in its usual meaning of taking the $x$-fold product of $y$ with itself).
Composition of partial functions is the same as composition of binary relations as described in Definition \ref{def:composition}.
We will however often forego writing the composition symbol $\circ$ and instead just express the composition of partial functions $f \circ g$ as $fg$.

\begin{definition}[Bijection] \label{def:bijection}
    A total function $f \in Y^X$ is called a bijection if it satisfies the following conditions:
    \begin{itemize}
        \item[] \textbf{Injective}: For all $x,z \in X$ and all $y \in Y$, if $(x)f = y$ and $(z)f = y$ then $x=z$.
        
        \item[] \textbf{Surjective}: The image of $f$ is all of $Y$. That is, $\im{f} = Y$.
    \end{itemize}
    If in addition $X=Y$, then $f$ is often referred to as a \emph{permutation}.
\end{definition}

We are ready to define the main type of relational object that we will be studying in this thesis.

\begin{definition}[Partial bijection] \label{def:partial_bijection}
    A \emph{partial bijection} $f \subseteq X \times Y$ is a bijection from a subset of $X$ to a subset of $Y$.
    That is, $f$ is an injective partial function.
    If in addition $X=Y$, then $f$ is sometimes referred to as a \emph{chart} or \emph{partial permutation}.
\end{definition}

Spaces of partial bijections will be the main objects of study in this thesis.
Note that the inverse of any (partial) bijection is also a (partial) bijection (unlike functions, whose inverses are simply binary relations in general).
There is one final type of function that we wish to introduce.

\begin{definition}[Sequence] \label{def:sequence}
    Let $X$ be a set and $\alpha$ an ordinal.
    Then an \emph{$\alpha$-sequence on $X$} (or simply a \emph{sequence}) is a function from $\alpha$ to $X$ and it will often be denoted by $(x_i)_{i \in \alpha} \in X^\alpha$, where $x_i \in X$ is the image of $i \in \alpha$ in the sequence.
    If the ordinal $\alpha$ is not explicitly stated, then `sequence' typically refers to an $\omega$-sequence.
\end{definition}

We will sometimes refer to the elements in the image of a sequence as the \emph{elements} or \emph{entries} of the sequence.
In a slight abuse of notation we will even write things like $x \in p$ to refer to an element $x$ in the image of a sequence $p$.

\begin{definition}[Restriction] \label{def:restriction}
    Let $f$ be a partial function from a set $X$ to a set $Y$.
    Then a \emph{restriction} of $f$ to a subset $A$ of $X$ is the function;
    \begin{equation*}
        \eval{f}_A = \set{(x,y) \in f \given x \in A}
    \end{equation*}
    which has domain $A \cap \dom{f}$ and image $Af = \set{y \in Y \given (\exists x \in A)~ xf = y}$.
    Furthermore, we denote the restriction of $f$ to the preimage of $B \subseteq Y$ by $\eval{f}^B = \set{(x,y) \in f \given y \in B}$.
\end{definition}

We have covered various types of partial functions, but there are also other kinds of binary relations that will play a major role in this thesis.
We will now list these objects and some of their properties.

\begin{definition}[Preorder] \label{def:pre_order}
    Let $X$ be a set and $\lesssim$ a binary relation on $X$. Then we will refer to $\lesssim$ as a \emph{preorder} (or \emph{quasiorder}) on $X$, if it satisfies the following conditions:
    \begin{itemize}
        \item[] \textbf{Reflexive}: For all $x \in X$, $(x,x) \in~\lesssim$.
        
        \item[] \textbf{Transitive}: For all $x,y,z \in X$, if $(x,y) \in~\lesssim$ and $(y,z) \in~\lesssim$ then $(x,z) \in~\lesssim$.
    \end{itemize}
    Given a preorder $\lesssim$, if $(x,y) \in~\lesssim$ we will usually express this fact as $x \lesssim y$ or $y \gtrsim x$ and we say that \emph{$x$ is less than or equivalent to $y$} (or that \emph{$y$ is greater than or equivalent to $x$}).
\end{definition}

Preorders are a sort of precursor to orders and equivalence relations, which we will define now.

\begin{definition}[Partial order] \label{def:partial_order}
    Let $\leq$ be a preorder on a set $X$.
    We will then refer to $\leq$ as a \emph{partial order} (or an \emph{order relation}) on $X$, if it satisfies the following additional condition:
    \begin{itemize}
        \item[] \textbf{Antisymmetric}: For all $x,y \in X$, if $x \leq y$ and $y \leq x$ then $x=y$.
    \end{itemize}
    Given a partial order $\leq$ on a set $X$, we will often refer to the ordered pair $(X,\leq)$ as an \emph{ordered set}.
    Due to the requirement of antisymmetry, we can express the statement $x \leq y$ as \emph{$x$ is less than or equal to $y$} whenever $\leq$ is a partial order.
\end{definition}

\begin{definition}[Total order] \label{def:total_order}
    Let $(X,\leq)$ be an ordered set.
    We will then refer to $\leq$ as a \emph{total order} (or \emph{linear order}), if it satisfies the following additional condition:
    \begin{itemize}
        \item[] \textbf{Strongly Connected}: For all $x,y \in X$, $x \leq y$ or $y \leq x$.
    \end{itemize}
    Total orders also give rise to the notion of a \emph{chain}, which is a totally ordered subset of any partially ordered set.
\end{definition}

\begin{definition}[Well-order] \label{def:well_order}
    A totally ordered set $(X,\leq)$ is called \emph{well-ordered}, if it has the property that every non-empty subset of $X$ has a \emph{least element} under $\leq$.
    That is, for all non-empty $Y \subseteq X$ there exists $y \in Y$ such that $(\forall z \in Y)~ y \leq z$.
\end{definition}

Any ordered set has a corresponding \emph{strict ordering}.
Given a partial order $\leq$ we will denote the corresponding strict order by $<$ and it is obtained from $\leq$ by removing all identity pairs $(x,x) \in X \times X$.
\begin{itemize}
    \item[] \textbf{Irreflexive}: For all $x \in X$, $(x,x) \notin~<$.
\end{itemize}
If $<$ is a strict partial order and $x < y$, we then say that \emph{$x$ is (strictly) less than $y$}.
Definitions \ref{def:partial_order}, \ref{def:total_order}, and \ref{def:well_order} describe the three types of order relations that will be relevant to this thesis.
They are all examples of preorders, but there is also a different direction one can take starting from preorders.

\begin{definition}[Equivalence relation] \label{def:equivalence_relation}
    Let $\approx$ be a preorder on a set $X$.
    We will then refer to $\approx$ as an \emph{equivalence relation}, if it satisfies the following additional condition:
    \begin{itemize}
        \item[] \textbf{Symmetric}: For all $x,y \in X$, if $x \approx y$ then $y \approx x$.
    \end{itemize}
    We say that the $x$ and $y$ above are \emph{equivalent} under $\approx$.
\end{definition}

It is a well-known fact that any equivalence relation defines a \emph{partition} of the set that it is a relation on.

\begin{definition}[Partition] \label{def:partition}
    Let $X$ be a set.
    Then a \emph{partition} $\partition$ of $X$ is a collection of disjoint non-empty subsets of $X$ such that the union over $\partition$ is all of $X$.
    That is, 
    \begin{enumerate}[~\normalfont(i)]
        \item $\emptyset \notin \partition$;
        \item $(\forall A \in \partition)$ $A \subseteq X$;
        \item $(\forall A,B \in \partition)$ $A \cap B = \emptyset$ whenever $A \neq B$; and
        \item $\bigcup_{A \in \partition} A = X$.
    \end{enumerate}
\end{definition}

The elements of a partition are often referred to as \emph{parts} or \emph{blocks}.
A \emph{refinement} of a partition $\partition$ is another partition $\mathcal{Q}$ on the same set $X$ such that for all points $x,y \in X$, if $x$ and $y$ belong to distinct parts in $\partition$ (that is, there exist distinct parts $A,B \in \partition$ such that $x \in A$ and $y \in B$), then $x$ and $y$ also belong to distinct parts in $\mathcal{Q}$.
That is, $Q$ can at most `break' the parts in $\partition$ into smaller parts.
Let $\partition$ be a partition on a set $X$ and $\mathcal{Q}$ a partition on a set $Y$.
Then $\partition$ and $\mathcal{Q}$ are \emph{isomorphic as partitions} if there exists a bijection $f \from X \to Y$ such that for all $x,y \in X$, there exists a part $A \in \partition$ such that $x,y \in A$ if and only if there exists a part $B \in \mathcal{Q}$ such that $xf,yf \in B$.
We say that such a function $f$ \emph{acts as a bijection} between $\partition$ and $\mathcal{Q}$.

It is a well-known fact that any equivalence relation defines a partition of the set that its a relation on and vice versa \cite{halmos1960naive}.
If the parts of a partition are explicitly obtained from an equivalence relation, then these parts are often referred to as \emph{equivalence classes}.
Looking at the definitions one can also conclude that given any preorder $\lesssim$, the preorder can be split into an equivalence relation $\approx~ = \set{(x,y) \in ~\lesssim~ \given x \lesssim y ~\land~ y \lesssim x}$ and a strict partial order $<~ = \{\, (x,y) \in  ~\lesssim~ \given \linebreak x \lesssim y ~\land~ y \not\lesssim x \,\}$ which in some sense orders the equivalence classes of $\approx$ (if $x \approx a$ and $y \approx b$, then $x < y$ implies that $a < b$).
This fact plays a major role in this thesis, as we will be studying an equivalence relation resulting from a preorder.

\subsection{Relational Structures}

Extending the notion of a relation, one can define what is known as a \emph{relational structure}.

\begin{definition}[Relational structure]
    Let $X$ be a set and $\relations$ a set of relations on $X$. 
    Then the ordered pair $\relStructure = (X,\relations)$ is called a \emph{relational structure} on $X$.
\end{definition}

In a slight abuse of notation we will often not distinguish between a relation $R$ on a set $X$ and the relational structure $(X,\set{R})$.
Taking this abuse of notation a step further, we will also omit the set brackets around the relation $R$ and simple write a relational structure on the single relation $R$ as $(X,R)$.
In this notation, the previously defined ordered sets can also be considered as relational structures.
Given a relation or a relational structure on a set $X$ we can also define an \emph{induced subrelation} or \emph{induced substructure} on a subset $Y$ of $X$.

\begin{definition}[Induced subrelation]
    Let $X$ be a set and $R$ a relation on $X$. 
    Then the \emph{induced subrelation} $R'$ of $R$ on a subset $Y$ of $X$ is the set of tuples in $R$ which only include elements from $Y$.
    That is, $R' = R \cap Y^{\ary{R}}$.

    Similarly, if $\relStructure = (X,\relations)$ is a relational structure on $X$, then the \emph{induced substructure} of $\relStructure$ on a subset $Y$ of $X$ is a relational structure $\relStructure' = (Y,\relations')$, where $\relations' = \set{R \cap Y^{\ary{R}} \given R \in \relations}$.
\end{definition}

Relational structures are regularly studied in Universal Algebra and Model Theory.
In this thesis, the main objects of interest will typically be spaces of certain functions between relational structures or from a relational structure to itself.

\begin{definition}[Endomorphism of relational structure] \label{def:end_rel}
    Let $R$ be an $n$-ary relation on a set $X$.
    Then an \emph{endomorphism} $f$ on $R$ is a function from $X$ to $X$ such that; 
    \begin{equation*}
        (\forall \alpha \in R)~ \alpha f \in R
    \end{equation*}
    where $\alpha = (a_0, \dots, a_{n-1})$ is an \emph{n-tuple} in $X$ and $\alpha f = (a_0f, \dots, a_{n-1}f)$.
    Similarly, if $\relStructure = (X, \relations)$ is a relational structure, then an endomorphism $f$ on $\relStructure$ is a function from $X$ to $X$ such that $f$ is an endomorphism for all $R \in \relations$.
\end{definition}

\begin{definition}[Automorphism of relational structure] \label{def:aut_rel}
    Let $R$ be a relation on a set $X$.
    Then an \emph{automorphism} $f$ on $R$ is a bijection from $X$ to $X$ such that: 
    \begin{equation*}
        (\forall \alpha \in R)~ \alpha f \in R ~\land~ \alpha\inv{f} \in R.
    \end{equation*}
    Similarly, if $\relStructure = (X, \relations)$ then an automorphism $f$ on $\relStructure$ is a bijection from $X$ to $X$ such that $f$ is an automorphism for all $R \in \relations$.
\end{definition}

A simple example of a type of relational structure is a \emph{graph}, which is a relational structure with only a single binary relation (similar to the ordered sets discussed in the previous section).

\begin{definition}[Graph] \label{def:graph}
    A \emph{graph} $G = (V,E)$ is a relational structure where $E$ is a \emph{symmetric irreflexive binary relation}.
    The elements of $V$ are often referred to as \emph{vertices}. 
    Two vertices $x,y \in V$ are said to be \emph{adjacent} if they are related in $E$ and the unordered pair $\set{x,y}$ is then called an \emph{edge} in $G$.
\end{definition}

\begin{definition}[Path] \label{def:path}
    Let $G = (V,E)$ be a graph.
    Then a \emph{path} $(x_i)_{i \in I}$ in $G$ is a sequence of vertices in $V$ such that each entry is distinct ($x_i \neq x_j$ whenever $i \neq j$) and $x_{i+1}$ is adjacent to to $x_i$ for all $i \in I$.
    A \emph{cycle} refers to a path with 3 or more entries in which the first and last entries are adjacent.
\end{definition}

We say that two vertices are \emph{connected} if there exists a path between them (that is, a path that starts at one vertex and ends at the other).
A graph is then called a \emph{connected graph}, if all its vertices are connected to each other.
Given two connected vertices we say that the \emph{distance} between them is the length (i.e. the cardinality) of the minimal path connecting them.
From this definition we see that the distance between two vertices in a connected graph is always finite (due the requirement of adjacency in paths), but not all paths in a connected graph are necessarily finite (some graphs allow for paths that are $\omega$-sequences, but such have no `last entry').

\section{Semigroup Theory}

Considering that the word `monoid' appears twice in the title of this dissertation, surprisingly little Semigroup Theory is needed to read and understand the thesis.
For a general introduction to Semigroup Theory see e.g. \cite{howie1995fundamentals} and for an introduction more specifically on the topic of inverse semigroups see e.g. \cite{lawson_inverse_semigroups}.
We will assume that the reader is familiar with basic Abstract Algebra and Group Theory, otherwise see e.g. \cite{hungerford2013abstract} or \cite{aluffi2021algebra} for an introduction to the field of Abstract Algebra.

\subsection{Semigroup Basics}
Semigroup Theory starts with the notion of a semigroup.
This is a generalisation of the concept of a group and as such is simply defined by fewer axioms.

\begin{definition}[Semigroup] \label{def:semigroup}
    Let $S$ be a set and $\mu$ a \emph{binary operation} on $S$ (a function from $S \times S$ to $S$).
    We then call the ordered pair $(S,\mu)$ a \emph{semigroup}, if the binary operation satisfies the following condition:
    \begin{itemize}
        \item[] \textbf{Associative}: For all $x,y,z \in S$, $\paa{(x,y)\mu,z}\mu = \paa{x,(y,z)\mu}\mu$.
    \end{itemize}
\end{definition}

It is customary to forego explicitly writing the binary operation $\mu$ and instead just refer to a semigroup $(S,\mu)$ by the underlying set $S$.
Likewise, we will not write out the binary operation $\mu$ when applying the operation but instead just use concatenation of the elements (that is, we express the binary operation by the usual notation for \emph{multiplication}).
In this notation the property of associativity can be written as follows:
\begin{itemize}
    \item[] \textbf{Associative}: For all $x,y,z \in S$, $(xy)z = x(yz)$.
\end{itemize}
We will now cover some special types of semigroups.

\begin{definition}[Monoid] \label{def:monoid}
    Let $S$ be a semigroup.
    We then refer to $S$ as a \emph{monoid}, if it satisfies the following condition:
    \begin{itemize}
        \item[] \textbf{Identity}: There exists an element $e \in S$ with the property that for all $x \in S$, $ex=xe=x$.
    \end{itemize}
    Such an element $e$ is referred to as the \emph{identity} of $S$.
    It is a standard exercise to show that if a semigroup contains such an identity element, then it is necessarily unique.
    If $S$ is a semigroup, then we will use the notation $S^1$ to refer to the monoid obtained from adjoining an identity to $S$ if it does not already have one (that is, if $S$ is already a monoid then $S^1 = S$).
\end{definition}

\begin{definition}[Inverse semigroup] \label{def:inverse_semigroup}
    Let $M$ be a semigroup.
    We then refer to $M$ as an \emph{inverse semigroup}, if it satisfies the following condition:
    \begin{itemize}
        \item[] \textbf{Inverses}: For all $x \in M$ there exists a unique $y \in M$ such that $xyx = x$ and $yxy = y$.
    \end{itemize}
    This element $y$ is referred to as the \emph{inverse of $x$} and we typically denote it by $\inv{x}$.
\end{definition}

It again follows from the definition of inverses that $\inv{(\inv{x})} = x$ for any element $x$ of an inverse semigroup $M$.
Furthermore, given any inverse semigroup $M$ and a subset $A \subseteq M$ we can define the set $\inv{A} = \set{x \in M \given \inv{x} \in A}$ and a unary operation $\iota$ on $M$ (a function from $M$ to $M$) which maps each element to its inverse.
That is, $(\forall x \in M)~ (x)\iota = \inv{x}$.
Such a unary operation is referred to as \emph{inversion}.

The definition of an \emph{inverse monoid} can be obtained by combining the above definitions of monoids and inverse semigroups.
Finally, we can now define the usual groups as follows.

\begin{definition}[Group] \label{def:group}
    Let $G$ be an inverse monoid.
    We will then refer to $G$ as a \emph{group}, if it satisfies the following condition:
    \begin{itemize}
        \item[] \textbf{Units}: For all $x \in G$, $x\inv{x} = e$ where $e$ is the identity of $G$.
    \end{itemize}
\end{definition}

We have now covered the different types of semigroups that appear in this thesis.
An important concept when working with semigroups and other algebraic structures is the notion of a \emph{substructure}.

\begin{definition}(Subsemigroup) \label{def:subsemigroup}
    Let $(S,\mu)$ be a semigroup.
    A \emph{subsemigroup} $T$ of $S$ is then a subset of $S$ which happens to also be a semigroup under the same binary operation $\mu$ (when restricted to the subset $T$).
    That is, $\mu$ maps $T \times T$ to $T$ which can also be expressed as;
    \begin{equation*}
        (\forall x,y \in T)~ xy \in T.
    \end{equation*}
\end{definition}

Using Definition \ref{def:subsemigroup} above, one can derive the notion of other substructures.
Notably one can also talk about the subsemigroups of a group or vice versa (and we will talk about such objects).
Not every subset of a group or semigroup is necessarily a semigroup, but we can define the notion of \emph{generating} a semigroup from any subset of a given semigroup.

\begin{definition}(Generated semigroup) \label{def:generate_semigroup}
    Let $S$ be a semigroup and $U$ a subset of $S$.
    We then denote by $\genset{U}$ the intersection of all subsemigroups of $S$ containing $U$.
    This $\genset{U}$ is itself a subsemigroup of $S$ \cite{howie1995fundamentals} and we call it the semigroup \emph{generated by $U$}.
    An equivalent definition would be that $\genset{U}$ is the set of all finite products of elements in $U$ under the binary operation of $S$.
\end{definition}

In a slight abuse of notation, we will often simply list any number of subsets and elements of a semigroup $S$ within the brackets $\genset{{\dots}}$ to denote the union of the included sets and the sets containing the elements listed (e.g. $\genset{U,V,f,g}$ will be used to denote the properly written $\genset{U \cup V \cup \set{f,g}}$, where $U,V \subseteq S$ and $f,g \in S$).
In this thesis, we will always use $\genset{{\dots}}$ to mean the \emph{semigroup} generated by the content within the brackets.
In other contexts, it could also mean the generated group or inverse semigroup, but we will consistently use it to refer to only semigroup generation.
There is another type of semigroups substructure that we will make use of in this thesis.

\begin{definition}(Semigroup ideal) \label{def:semigroup_ideal}
    Let $S$ be a semigroup and $I$ a subset of $S$.
    Then $I$ is called a \emph{left ideal} if $SI = \set{xy \given x \in S ~\land~ y \in I} \subseteq I$, a \emph{right ideal} if $IS = \set{yx \given y \in I ~\land~ x \in S} \subseteq I$, and a (two-sided) \emph{ideal} if it is both a left and right ideal.
\end{definition}

This notation of taking left or right `acts' by subsets of a semigroup can be extended further to products of any length. 
Furthermore, we can also multiply by single elements of the semigroup rather than a subset: $xSz = \set{xyz \in S \given y \in S}$, where $x$ and $z$ are fixed elements in $S$.
Next, we define a certain type of elements that a semigroup might posses, which are of particular interest.

\begin{definition}[Idempotent] \label{def:idempotent}
    Let $S$ be a semigroup and $e$ an element of $S$.
    If $e^2 = ee = e$, then $e$ is called an \emph{idempotent} of $S$.
    If a subsemigroup $T$ contains all the idempotents of $S$, then $T$ is referred to as a \emph{full} subsemigroup of $S$.
\end{definition}

It is a well-known fact that groups contain exactly one idempotent, namely the identity, and clearly every idempotent is its own inverse.
Finally, we note that the concept of a homomorphism (which the reader should be familiar with from Group Theory) can also be defined for semigroups.

\begin{definition}[Homomorphism] \label{def:homomorphism}
    Let $S$ and $T$ be semigroups and $f$ a function from $S$ to $T$.
    Then $f$ is called a \emph{homomorphism}, if it satisfies the following condition:
    \begin{equation*}
        (\forall x,y \in S)~ (xy)f = (x)f (y)f
    \end{equation*}
\end{definition}

Similarly, we still have that an \emph{isomorphism} is a bijective homomorphism, an \emph{endomorphism} is a homomorphism from a semigroup to itself, and an \emph{automorphism} is a bijective endomorphism.

\subsection{The Symmetric Inverse Monoid}

Finally, we get to the main star of the show, the titular \emph{symmetric inverse monoid}.

\begin{figure}
    \centering
    \includegraphics[width=0.9\linewidth]{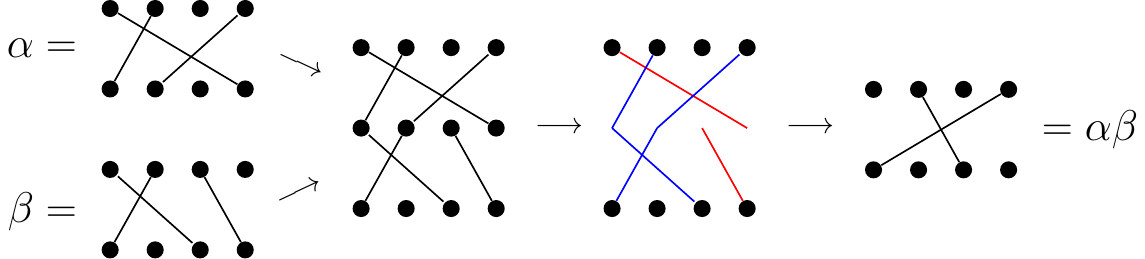}
    \caption{Composition of charts $\alpha$ and $\beta$ in $I_4$.}
    \label{fig:IX_composition}
\end{figure}

\begin{definition}[Symmetric inverse monoid] \label{def:IX}
    Let $X$ be a set.
    Then the \emph{symmetric inverse monoid} on $X$, denoted $I_X$, is the semigroup of all charts (partial bijections) on $X$ under the usual composition of binary relations.
    For an example of composition in a finite symmetric inverse monoid, see Figure \ref{fig:IX_composition}.
\end{definition}

Since the elements of $I_X$ are all bijections (and the identity function of $X$ is also a partial bijection), it follows that $I_X$ is an inverse monoid.
The symmetric inverse monoid belongs to a family of semigroups usually referred to as \emph{binary relation semigroups}.
These semigroups are of particular interest due to a set of representation theorems.
Especially famous among these binary relation semigroups are the \emph{symmetric group} \Sym{X} of all permutations on a set $X$ and the \emph{full transformation monoid} $T_X$ of all functions on $X$ (another notation for $T_X$ would be the previously mentioned $X^X$).
The permutation group \Sym{X} is probably the most well-known due to Cayley's Theorem, which states that every group is isomorphic to a subgroup of a symmetric group.
A similar results exists for $T_X$, stating that every semigroup embeds into some transformation semigroup.
Likewise, there is a representation theorem for inverse semigroups.

\begin{theorem}[Wagner-Preston Theorem {\cite[Theorem 1]{lawson_inverse_semigroups}}] \label{wagner-preston}
    Let $M$ be an inverse semigroup.
    Then there exists a set $X$ and an injective homomorphism $\phi \from M \to I_X$.
\end{theorem}

In this sense, we can consider $I_X$ to be a `canonical' example of an inverse semigroup.
For further reading on finite symmetric inverse monoids see e.g. \cite{lipscomb1996symmetric}.
Two other important binary relation monoids are the \emph{partial transformation monoid} $P_X$ of all partial functions on a set $X$ and the \emph{binary relation monoid} $B_X$ of all binary relations on $X$.
It is well-known that all the semigroups mentioned in this section have cardinality $2^{\card{X}}$ whenever $X$ is an infinite set.
In order to prove this fact, we will use that that $\Sym{X}$ is contained in $I_X$, $T_X$, $P_X$, and $B_X$ (that is, all the other monoids mentioned in this section) and that $B_X$ contains all other binary relation semigroups on $X$.

\begin{proposition}
    Let $X$ be an infinite set.
    Then $2^{\card{X}} = \card{\Sym{X}} = \card{I_X} = \card{T_X} = \card{P_X} = \card{B_X}$.
\end{proposition}

\begin{proof}
    It suffices to show that $2^{\card{X}} \leq \card{\Sym{X}}$ and $2^{\card{X}} \geq \card{B_X}$.

    To show that $2^{\card{X}} \leq \card{\Sym{X}}$, we invoke the Axiom of Choice and use that $X$ contains $2^{\card{X}}$ many moieties and that each moiety defines a partition of $X$ into two equal sized subsets, namely the moiety itself and its complement in $X$ (this technically results in a double counting, but that does not matter for infinite sets as $2 \cdot 2^{\card{X}} = 2^{\card{X}}$).
    For each such partition into two subsets we then pick a permutation in $\Sym{X}$, which swaps the two parts.
    These permutations then make up a proper subset of $\Sym{X}$ with cardinality $2^{\card{X}}$.

    To show that $2^{\card{X}} \geq \card{B_X}$, we use the fact that $2^{\card{X}}$ is the cardinality of the powerset of $X$ \cite{halmos1960naive}.
    Meanwhile, the elements of $B_X$ by definition make up the powerset of $X \times X$, and since $\card{X \times X} = \card{X}$ when $X$ is infinite, we get that $B_X$ has cardinality equal to the powerset of $X$.    
\end{proof}

For any of the monoids introduced in this section, a useful tool is to consider the image of points in $X$ under the elements of a subsemigroup of these binary relation monoids.

\begin{definition}[Cones and orbits] \label{def:cones_orbits}
    Let $X$ be a set, $x$ an element of $X$, and $S$ a subset of $B_X$.
    Then the \emph{cone} of $x$ under $S$ is the set containing the image of $x$ under all elements of $S$;
    \begin{equation*}
        xS = \set{xs \in X \given s \in S}.
    \end{equation*}
    If $S$ is an inverse semigroup, then the cones under $S$ are called \emph{orbits} and they partition a subset of $X$ \cite{howie1995fundamentals}.
    It is important to note, that this partitioned subset can potentially be a proper subset of $X$ and an element $x \in X$ is always contained in its own orbit unless that orbit is empty.
\end{definition}

Before moving on to the next topic, we have one final piece of background information about the symmetric inverse monoid to cover.

\begin{definition}[Semilattice of idempotents] \label{def:EX}
    Let $X$ be a set.
    Then the idempotents of $I_X$ are exactly those charts which act as the identity on a subset of $X$ \cite{lawson_inverse_semigroups}.
    We call these charts \emph{partial identities} and they form a submonoid of $I_X$ called the \emph{semilattice of idempotents}, which we denote by $E_X$.
    \begin{equation*}
        E_X = \set{e \in I_X \given (\forall x \in \dom{e})~ xe = x}
    \end{equation*}
\end{definition}

\subsection{Sierpiński rank and relative rank}

Generating sets and semigroup generation will be a central concept in this thesis.
A typical endeavour on the topic of generating sets is to find a \emph{minimal} generating set for a given semigroup.

\begin{definition}[Rank] \label{def:rank}
    Let $S$ be a semigroup.
    Then the \emph{rank} of $S$, denoted by \rank{S}, is the minimum cardinality of a generating set for $S$.
    That is, $\rank{S} = \min(\card{A} \given \genset{A} = S)$.
\end{definition}

In the study of countable semigroups, determining ranks is a natural consideration.
However, when a semigroup $S$ is uncountable the rank of $S$ will always simply be the cardinality of $S$ (this is a consequence of cardinal arithmetic).
So the rank does not tell us much when studying uncountable semigroups.
Instead we can consider the \emph{relative size} of a given subset of a semigroup.

\begin{definition}[Relative rank] \label{def:relative_rank}
    Let $S$ be a semigroup and $A$ a subset of $S$.
    Then the \emph{relative rank of $A$ in $S$} is the minimum cardinality of a set $U \subseteq S$ such that $\genset{A,U} = S$ and we denote it by \rank{S:A}.
\end{definition}

This notion of relative rank tells us whether a given subsemigroup is `large' or `small' relative to the supersemigroup it sits within.
It is important to note, that one could also consider the relative rank as a \emph{group} or \emph{inverse semigroup} instead, if the supersemigroup is either of those type of objects.
In this text however, we will mostly consider the relative rank as a semigroup unless explicitly stating otherwise.
There exists a rich study of the relative rank of subsemigroups of various semigroups, especially for the full transformation semigroup $T_X$ (see e.g. \cite{howie1998relative, higgins2003countable, higgins2003generating}).

A separate question one could ask on the topic of generating sets is whether all countable subsets of a given semigroup can be generated by subsets of cardinality less than some number $n \in \omega$.
This is the very question which Sierpiński considered for the full transformation monoid $T_X$ in 1935.

\begin{theorem}[{\cite[Théorème I]{sierpinski1935suites}}] \label{sierpinksi_theorem}
    Let $X$ be an infinite set.
    Then every countable subset of $T_X$ is contained in a 2-generated subsemigroup of $T_X$.
\end{theorem}

A consequence of Theorem \ref{sierpinksi_theorem} above is that the relative rank of any subset of $T_X$ must either be uncountable or at most 2 when $X$ is infinite.
Theorem \ref{sierpinksi_theorem} is a significant result in Semigroup Theory and in light of it the following definition has been introduced.

\begin{definition}[Sierpiński rank] \label{def:sierpinski_rank}
    Let $S$ be a semigroup.
    Then the \emph{Sierpiński rank} of $S$ is the least number $n \in \omega$ such that every countable subset of $S$ is contained in an $n$-generated subsemigroup of $S$. If no such $n$ exists, we say that $S$ has \emph{infinite Sierpiński rank}.
\end{definition}

Again we note that Definition \ref{def:sierpinski_rank} refers to the Sierpiński rank of a structure \emph{as a semigroup}.
In the case that $S$ is a group or inverse semigroup one could just as well consider the Sierpiński rank \emph{as a group} or \emph{as an inverse semigroup}.
As usual, we will only consider the Sierpiński rank of objects as semigroups unless otherwise stated.
An interesting question now is what the Sierpiński ranks of various semigroups are.
Of interest to us are the following results.

\begin{theorem}[{\cite[Theorem 3.5]{galvin1995generating}}] \label{thm:galvin_group_rank}
    Let $X$ be an infinite set. 
    Then every countable subset of \Sym{X} is contained in a 2-generator subsemigroup of \Sym{X}.
\end{theorem}

\begin{theorem}[{\cite[Theorem 1.4]{hyde2012sierpi}}] \label{thm:sierpinksi_rank_IX}
    Let $X$ be an infinite set.
    Then every countable subset of $I_X$ is contained in a 2-generated subsemigroup of $I_X$.
\end{theorem}

Theorems \ref{sierpinksi_theorem}, \ref{thm:galvin_group_rank}, and \ref{thm:sierpinksi_rank_IX} can thus be summarised as ``the semigroups $T_X$, \Sym{X}, and $I_X$ all have Sierpiński rank 2.''
In \cite{galvin1995generating} Galvin also shows that $\Sym{X}$ has Sierpiński rank 2 as a group, and in \cite{higgins2003countable} Higgins, Howie, Mitchell, and Ru\v{s}kuc show that $I_X$ has Sierpiński rank 2 as an inverse semigroup.

\section{Topological Algebra}

The field of Topology is concerned with the study of continuous functions through the definition of open and closed sets.
We will assume that the reader is familiar with basic Point Set Topology and metric spaces, but we include a short glossary here for convenience. For further information see e.g. \cite{munkrestopology}.

\begin{definition}[Topological space]
    Let $X$ be a set and $\topology$ a collection of subsets of $X$.
    Then $(X,\topology)$ is called a \emph{topological space} and $\topology$ a \emph{topology} on $X$, if the following hold:
    \begin{enumerate}[\normalfont~(i)]
        \item $\emptyset \in \topology$ and $X \in \topology$.
        \item If $\mathcal{A}$ is a subset of $\topology$, then $\bigcup_{A \in \mathcal{A}} A \in \topology$.
        \item If $\mathcal{A}$ is a finite subset of $\topology$, then $\bigcap_{A \in \mathcal{A}} A \in \topology$.
    \end{enumerate}
    The elements of the topology $\topology$ are then typically referred to as \emph{open sets} in $X$.
\end{definition}

\begin{definition}[Closed set]
    Let $(X,\topology)$ be a topological space.
    Then a subset $U \subseteq X$ is called \emph{closed} if the complement $X \setminus U$ is open in $(X,\topology)$.
\end{definition}

\begin{theorem}[{\cite[Theorem 17.1]{munkrestopology}}]
    Let $(X,\topology)$ be a topological space.
    Then the following statements about closed sets hold:
    \begin{enumerate}[\normalfont~(i)]
        \item $\emptyset$ and $X$ are closed.
        \item If $\mathcal{A}$ is a finite collection of closed sets, then $\bigcup_{A \in \mathcal{A}} A$ is closed.
        \item If $\mathcal{A}$ is a collection of closed sets, then $\bigcap_{A \in \mathcal{A}} A$ is closed.
    \end{enumerate}
\end{theorem}

\begin{definition}[Subspace topology]
    Let $(X,\topology)$ be a topological space and $Y$ a subset of $X$.
    Then the set $\topology_Y$ of all intersections of open sets in $\topology$ with $Y$ is a topology on $Y$.
    That is,
    \begin{equation*}
        \topology_Y = \set{U \cap Y \given U \in \topology}
    \end{equation*}
    is a topology on $Y$ and we call it the \emph{subspace topology}.
    Furthermore, we call the topological space $(Y,\topology_Y)$ a \emph{subspace} of $(X,\topology)$.
\end{definition}

\begin{theorem}[{\cite[Theorem 17.2]{munkrestopology}}]
    Let $(X,\topology)$ be a topological space with a subspace $(Y,\topology_Y)$.
    Then a subset $A \subseteq Y$ is closed in $(Y,\topology_Y)$ if and only if there exists a closed set $C$ of $(X,\topology)$ such that $A = C \cap Y$.
\end{theorem}

\begin{definition}[Topological neighbourhood]
    Let $(X,\topology)$ be a topological space and $x$ a point in $X$.
    Then a \emph{neighbourhood} of $x$ is any subset $V \subseteq X$ such that there exists an open set $U \in \topology$ with the property that $x \in U$ and $U \subseteq V$)
\end{definition}

\begin{definition}[Topological basis]
    Let $(X,\topology)$ be a topological space.
    Then a \emph{basis} $\mathcal{B}$ for the topology $\topology$ is a subset of $\topology$ such that every open set can be written as the union of basis elements.
    That is,
    \begin{equation*}
        \topology = \left\{\, \bigcup_{A \in \mathcal{A}} A \mathrel{\Big|} \mathcal{A} \subseteq \mathcal{B} \,\right\}.
    \end{equation*}
    The elements of topological basis are typically referred to as \emph{basic open sets}.
\end{definition}

\begin{definition}[Subbasis]
    Let $(X,\topology)$ be a topological space.
    Then a \emph{subbasis} $\mathcal{S}$ of the topology $\topology$ is a subset of $\topology$ such that $\topology$ is the intersection of all topologies on $X$ that contain $\mathcal{S}$.
    We say that $\topology$ is the topology \emph{generated} by the collection of subsets $\mathcal{S}$.
    Furthermore, the set of all finite intersections of elements of $\mathcal{S}$ form a basis for $\topology$.
    That is,
    \begin{equation*}
        \mathcal{B} = \left\{\, \bigcap_{A \in \mathcal{A}} A \mathrel{\Big|} \mathcal{A} \subseteq \mathcal{S} ~\land~ \card{A} < \aleph_0 \,\right\}
    \end{equation*}
    is a a basis for $\topology$.
\end{definition}

\begin{definition}[Topological closure]
    Let $(X,\topology)$ be a topological space and $A$ a subset of $X$.
    Then the \emph{closure} of $A$, denoted by $\closure{A}$, is the intersection of all closed sets in $X$ that contain $A$.
    It should be clear that $\closure{A}$ is closed and that if $A$ is closed then $\closure{A} = A$.
\end{definition}

\begin{theorem}[{\cite[Theorem 17.5]{munkrestopology}}] \label{thm:closure_alt}
    Let $(X,\topology)$ be a topological space and $A$ a subset of $X$.
    Then the following hold:
    \begin{enumerate}[\normalfont~~(a)]
        \item \label{thm:closure_alt:a}
        For all $x \in X$, $x \in \closure{A}$ if and only if every open neighbourhood of $x$ intersects $A$.
        That is,
        \begin{equation*}
            \closure{A} = \set{x \in X \given (\forall U \in \topology)~ x \in U \implies U \cap A \neq \emptyset}.
        \end{equation*}
        \item \label{thm:closure_alt:b}
        Suppose that $\topology$ is generated by a basis $\mathcal{B}$. 
        Then for all $x \in X$, $x \in \closure{A}$ if and only if every basic open set containing $x$ intersects $A$. That is,
        \begin{equation*}
            \closure{A} = \set{x \in X \given (\forall B \in \mathcal{B})~ x \in B \implies B \cap A \neq \emptyset}.
        \end{equation*}
    \end{enumerate}
\end{theorem}

\begin{definition}[Compactness]
    Let $(X,\topology)$ be a topological space.
    Then $\topology$ is called \emph{compact} if every open cover of $X$ has a finite subcover.
    That is, whenever there exists a collection of open sets $\mathcal{C} \subseteq \topology$ such that $X = \bigcup_{U \in \mathcal{C}} U$, then there exists a finite subset $\mathcal{A} \subseteq \mathcal{C}$ such that $X = \bigcup_{U \in \mathcal{A}} U$.
\end{definition}

\begin{definition}[$T_1$]
    Let $(X,\topology)$ be a topological space.
    Then $\topology$ is called $T_1$ if every singleton in $X$ is closed in $\topology$.
\end{definition}

\begin{definition}[Hausdorff]
    Let $(X,\topology)$ be a topological space.
    Then $\topology$ is called \emph{Hausdorff} or $T_2$ if for all distinct points $x,y \in X$ there exist open sets $U,V \in \topology$ such that $x \in U$, $y \in V$, and $U \cap V = \emptyset$.
    We say that distinct points in $X$ are \emph{separated by neighbourhoods}.
    Every Hausdorff space is also $T_1$.
\end{definition}

\begin{definition}[Second countable]
    Let $(X,\topology)$ be a topological space.
    Then $\topology$ is called \emph{second countable} or \emph{completely separable} if $\topology$ has a countable basis.
\end{definition}

\begin{definition}[Product topology]
    Let $I$ be an index set and $(X_i,\topology_i)$ a topological space for each $i \in I$.
    Then the \emph{product topology} on the Cartesian product $\bigtimes_{i \in I} X_i$ is the collection of Cartesian products $\bigtimes_{i \in I} U_i$, where $U_i \in \topology_i$ for each $i \in I$ with the additional condition that $U_i = X_i$ for all but finitely many $i \in I$.
\end{definition}

\begin{definition}[Continuous function]
    Let $(X,\topology_X)$ and $(Y,\topology_Y)$ be topological spaces.
    Then a function $f \from X \to Y$ is called \emph{continuous} if for all $U \in \topology_Y$, $U\inv{f} \in \topology_X$.
\end{definition}

\begin{definition}[Limit of a sequence]
    Let $(X,\topology)$ be a topological space and $(x_i)_{i\in\omega}$ a sequence in $X$.
    We then say that a point $x \in X$ is a \emph{limit} of the sequence $(x_i)_{i\in\omega}$, or that the sequence $(x_i)_{i\in\omega}$ \emph{converges} to the point $x$, if for every open neighbourhood $U$ of $x$ there exists a natural number $N \in \omega$ such that $x_n \in U$ for all $n \geq N$.
\end{definition}

\begin{definition}[Metric space]
    Let $X$ be a set and $\mu$ a function from $X \times X$ to $\R$.
    Then $(X,\mu)$ is called a \emph{metric space} and $\mu$ a \emph{metric} on $X$, if the following hold:
    \begin{enumerate}[\normalfont~(i)]
        \item For all $x,y \in X$, $(x,y)\mu \geq 0$ and $(x,y)\mu = 0$ if and only if $x = y$.
        \item For all $x,y \in X$, $(x,y)\mu = (y,x)\mu$.
        \item  For all $x,y,z \in X$, $(x,z)\mu \leq (x,y)\mu + (y,z)\mu$.
    \end{enumerate}
\end{definition}

\begin{definition}[Metric topology]
    Let $(X,\mu)$ be a metric space.
    Then the collection of all \emph{$\epsilon$-balls} $B_\mu(x,\epsilon) = \set{y \in X \given (x,y)\mu < \epsilon}$, for all $x \in X$ and positive real number $\epsilon > 0$, form a basis for a topology $\topology$.
    This topology $\topology$ is called the \emph{metric topology} \emph{induced} by $\mu$.
\end{definition}

\begin{definition}[Cauchy sequence]
    Let $(X,\mu)$ be a metric space and $(x_i)_{i\in\omega}$ a sequence in $X$.
    Then $(x_i)_{i\in\omega}$ is called a \emph{Cauchy sequence} if for every positive real number $\epsilon > 0$ there exists $N \in \omega$ such that $(x_m,x_n)\mu < \epsilon$ for all $m,n \geq N$.
\end{definition}

\begin{definition}[Complete metric]
    Let $(X,\mu)$ be a metric space.
    Then we call $\mu$ a \emph{complete metric} and $(X,\mu)$ a \emph{complete space}, if every Cauchy sequence converges in $(X,\mu)$.
\end{definition}

If the topology $\topology$ is not explicitly used, we might simply refer to a topological space $(X,\topology)$ by the underlying set $X$ (similar to what we do for semigroups).
Before delving into the topic of Topological Algebra proper, we will first cover special type of topological spaces that are relevant to this thesis.

\subsection{Polish Spaces}

When studying topological spaces, one is typically only interested in spaces with certain properties.
There are various criteria by which one can judge whether a given topology is `interesting' or `nice'.
Among such criteria, the following definition describes a class of especially `nice' topologies.

\begin{definition}[Polish space] \label{def:polish_topology}
    Let $(X,\topology)$ be a topological space.
    We then refer to $\topology$ as a \emph{Polish topology}, if it satisfies the following conditions.
    \begin{itemize}
        \item[] \textbf{Separable}: $(X,\topology)$ contains a countable dense subset.
        That is, there exists a countable subset $A \subseteq X$ such that every non-empty open set in $\topology$ has non-empty intersection with $A$.

        \item[] \textbf{Completely metrisable}: There exists a metric $\mu$ on $X$ which induces the topology $\topology$ such that $(X,\mu)$ is a complete metric space.
    \end{itemize}
\end{definition}

Polish topologies have a lot of structure and are frequently studied in fields like Descriptive Set Theory \cite{kechris2012classical}.
Polish spaces are interesting as they closely resemble the real numbers $\R$ with the standard topology (clearly this is a Polish space, as $\Q$ is a countable dense subset of $\R$ and Euclidean distance is an example of a complete metric which induces the standard topology).
While Polish spaces are very useful, there is also a limit to which spaces can be Polish.
In particular, it follows from the Cantor-Bendixson Theorem that any Polish space can be at most size of the continuum.

\begin{theorem}[Cantor-Bendixson {\cite[Corollary 6.5]{kechris2012classical}}] \label{thm:cantor-bendixson}
    Any Polish space is either countable or has cardinality $2^{\aleph_0}$.
\end{theorem}

\subsection{Topological Semigroups}

The field of Topological Algebra is concerned with studying algebraic objects endowed with a `compatible' topology.
For semigroups this idea of compatibility is defined by having the semigroup operation be continuous.

\begin{definition}[Topological semigroup] \label{def:topological_semigroup}
    Let $S$ be a semigroup and $\topology$ a topology on $S$.
    Then $S$ is called a \emph{topological semigroup} and $\topology$ a \emph{semigroup topology} if the semigroup multiplication
	\begin{equation*}
		\mu \from S \times S \to S \text{ given by } (x,y) \mapsto xy
	\end{equation*}
	is continuous with respect to the product topology on $S \times S$.
\end{definition}

Similarly one can define topological groups and inverse semigroups.

\begin{definition}[Topological inverse semigroup] \label{def:topological_inverse_semigroup}
    A \emph{topological inverse semigroup} is a topological semigroup $S$ such that the inversion map
	\begin{equation*}
		\iota \from S \to S \text{ given by } x \mapsto \inv{x}
	\end{equation*}
	is well-defined and continuous with respect to the semigroup topology.
\end{definition}

\begin{definition}[Topological group] \label{def:topological_group}
    A \emph{topological group} is a topological inverse semigroup that happens to be a group.
\end{definition}

The aforementioned compatibility of the algebraic and topological structures of such topological semigroups are expressed in many ways.
An example of such is the following classic result in Topological Algebra, where certain algebraic and topological operations work in unison.

\begin{lemma} \label{lem:closure_subsemigroup}
    Let $S$ be a topological (inverse) semigroup.
    If $T$ is an (inverse) subsemigroup of $S$, then so is the closure of $T$.
\end{lemma}

\begin{proof}
    Let $x$ and $y$ be elements $\closure{T}$.
    We then wish to show that $xy$ is also in $\closure{T}$.
    If we let $\mu$ denote the binary operation associated with $S$, then $\mu$ is a continuous function from $S \times S$ with the product topology on $S$.
    This implies that the preimage of any open neighbourhood $V$ of $xy$ under $\mu$ can be written as a union of Cartesian products $\bigcup_{i \in I} U_i \times W_i$, where $U_i$ and $W_i$ are basic open subsets of $S$ for all $i \in I$.
    In particular there must exist a a Cartesian product $A \times B \subseteq V\inv{\mu}$, such that $A$ and $B$ are open neighbourhoods of $x$ and $y$ respectively.
    However, since $x$ and $y$ are in the closure of $T$, we must have that there exist elements $a \in A$ and $b \in B$ such that $a,b \in T$.
    Since $a$ and $b$ are both elements of the subsemigroup $T$, it follows that $ab \in T$.
    Furthermore, $(a,b)$ belongs to the basic open neighbourhood $A \times B$ in the preimage of the open neighbourhood $V$, so we can also conclude that $ab \in V$.
    Hence, the product of elements from $\closure{T}$ also belong to $\closure{T}$, meaning that $\closure{T}$ is itself a subsemigroup of $S$.

    If furthermore $S$ and $T$ are inverse semigroups and the inversion map $\iota \from S \to S$ is continuous, we will then show that $\inv{x} \in \closure{T}$ for all $x \in \closure{T}$.
    This follows a similar argument to the proof that $\closure{T}$ is a semigroup.
    Let $V$ denote an open neighbourhood of $\inv{x}$.
    Then $V\inv{\iota} = V\iota$ is an open neighbourhood of $x$, which implies that there exists an element $y$ in the intersection $T \cap V\iota$, since $x$ is in the closure of $T$.
    It follows that $\inv{y}$ is an element that belongs to both $V$ and $T$ (since $T = \inv{T}$), meaning that any open neighbourhood of an element from $\inv{\closure{T}}$ intersects $T$.
    Hence, $\inv{\closure{T}} = \closure{T}$.
\end{proof}

The discrete and trivial topologies are (inverse) semigroup topologies on every (inverse) semigroup, but the question of finding interesting (inverse) semigroup topologies has been a subject of much research. 
This is where Polish topologies become relevant, as demanding that the (inverse) semigroup topology be Polish is a way to guarantee a rich topological structure.
In fact, some topological semigroups have a unique Polish topology, making it in some sense the `canonical topology' on said semigroup.
The most famous of these unique Polish semigroup topologies is probably the \emph{pointwise topology} on $\Sym{X}$, under which the symmetric group is a topological group.
It follows from a result by Gaughan \cite{gaughan1967topological} and Kechris and Rosendal \cite{kechris2007turbulence} that the pointwise topology is the unique Polish group topology on $\Sym{X}$ whenever $X$ is a countably infinite set.
We will now define this pointwise topology, but rather than defining it on $\Sym{X}$ we will give the topology on $T_X$.

\begin{definition}[\emph{Pointwise topology}] \label{pointwise_top}
	Let $X$ be an infinite set.
    Then the \emph{pointwise} (or \emph{function} or \emph{finite}) topology on $T_X = X^X$ is the topology generated by the sets
    \begin{equation*}
        U_{x,y} = \set{f \in T_X \given (x,y) \in f}
    \end{equation*}
    for all $x,y \in X$.
	This is equivalent to endowing each component set $X$ with the discrete topology and then taking the product topology on $X^X$.	
	When we talk about the pointwise topology on $\Sym{X}$, the inferred meaning is then the subspace topology on this subsemigroup.
\end{definition}

In \cite{part0topological} Mesyan, Mitchell, and Péresse show that the pointwise topology is also the unique Polish semigroup topology on $T_X$ when $X$ is a countably infinite set.
The question of finding such `natural' or `canonical' Polish semigroup topologies for various semigroups is explored further by Elliott et al. in \cite{part1topological}. 
Before we introduce Polish topologies on $I_X$, we will first explore a generalisation of semigroup topologies.

\begin{definition}[Semitopological/shift-continuous] \label{def:semitopological}
    Let $S$ be a semigroup and $\topology$ a topology on $S$.
    Then $\mathcal{T}$ is \emph{left semitopological} for $S$ if for every $s \in S$ the function $\mathfrak{l}_s \from S \to S$ defined by $(x)\mathfrak{l}_s = sx$ is continuous. Analogously if every function $\mathfrak{r}_s \from S \to S$ defined by $(x)\mathfrak{r}_s = xs$ is continuous, then we say that $\mathcal{T}$ is \emph{right semitopological} for $S$. 
    If a topology $\mathcal{T}$ is both left and right semitopological on the semigroup $S$, then we say that $\mathcal{T}$ is \emph{semitopological} or \emph{shift-continuous} for $S$.
\end{definition}

Every topological semigroup is semitopological, but the converse is not true in general \cite{part0topological}. We refer to a semigroup with a left semitopological topology as a left semitopological semigroup. Analogous definitions can be made for right semitopological and semitopological/shift-continuous.

\subsection{Topologies on $I_X$}

For $I_X$ we can no longer rely on the pointwise topology.
In \cite{part1topological} the authors suggest two `natural' Polish topologies for $I_X$.

\begin{theorem}[{\cite[Theorem 5.12]{part1topological}}] \label{I1_theorem}
    Let $X$ be an infinite set and let $\mathcal{I}_1$ denote the topology on $I_X$ generated by the sets
    \begin{equation*}
	U_{x,y} = \set{f \in I_X \mid (x,y) \in f} 
	\quad \text{and } \quad
	V_{x,y} = \set{f \in I_X \mid (x,y) \notin f} 		
    \end{equation*}
    for all $x,y \in X$. 
    Then the following hold:
    \begin{enumerate}[\normalfont(i)]
        \item the topology $\mathcal{I}_1$ is compact, Hausdorff, and semitopological for $I_X$ and inversion is continuous;
        \item $\mathcal{I}_1$ is the least $T_1$ topology that is semitopological for $I_X$;
        \item if $X$ is countable, then $\mathcal{I}_1$ is Polish.
    \end{enumerate} 
\end{theorem}

$\mathcal{I}_1$ is an excellent topology with a lot of interesting properties, most important of these being that inversion is continuous.
It is however only semitopological for $I_X$ and not topological.
If this is considered to be an issue, then a larger topology on $I_X$ must be defined.

\begin{theorem}[{\cite[Theorem 5.15]{part1topological}}] \label{I4_theorem}
    Let $X$ be an infinite set and $\mathcal{I}_4$ the topology on $I_X$ generated by the sets
    \begin{equation*}
	U_{x,y} = \set{f \in I_X \mid (x,y) \in f},
    \end{equation*}
    \begin{equation*}
	W_x = \set{f \in I_X \mid x \notin \dom{f}}, 
	~\text{and}~\, W_{x}^{-1} = \set{f \in I_X \mid x \notin \im{f}}
    \end{equation*}
    for all $x,y \in X$.
    Then the following hold:
    \begin{enumerate}[\normalfont(i)]
        \item the topology $\mathcal{I}_4$ is a Hausdorff inverse semigroup topology for $I_X$;
        \item $\mathcal{I}_4$ is the unique $T_1$ inverse semigroup topology on $I_X$ inducing the pointwise topology on $\Sym{X}$;
        \item if $X$ is countable, then $\mathcal{I}_4$ is the unique $T_1$ second countable inverse semigroup topology on $I_X$;
        \item if $X$ is countable, then $\mathcal{I}_4$ is the unique Polish inverse semigroup topology on $I_X$.
    \end{enumerate} 
\end{theorem}

In \cite{part1topological} they list a number of other properties of $\mathcal{I}_4$ in addition to these, but we only mention the ones relevant to this thesis here.
So we are presented with two compelling topologies. 
$\mathcal{I}_1$ is very minimal but only shift-continuous whereas $\mathcal{I}_4$ is the unique Polish inverse semigroup topology for $I_X$ (when $X$ is countable).
However, we will later show that there is no need to distinguish between these two topologies for a certain type of subsets of $I_X$.
For a further study on the full classification of all semigroup topologies on $I_X$, see \cite{bardyla2024classifyingpolishsemigrouptopologies}.

\section{The Bergman-Shelah Preorder}

As mentioned earlier, Galvin showed in \cite{galvin1995generating} that the Sierpiński rank of $\Sym{X}$ is 2.
This in part constitutes a strengthening of an earlier result by Macpherson and Neumann.
\begin{theorem}[{\cite[Theorem 1.1]{subgroups_macpherson_neumann}}] \label{thm:symX_chain}
    Let $X$ be an infinite set.
    Then $\Sym{X}$ is not the union of a chain of $\leq \card{X}$ proper subgroups.
\end{theorem}
Theorem \ref{thm:symX_chain} implies that given a subgroup $G \subseteq \Sym{X}$, if there exists a subset $U \subseteq \Sym{X}$ of cardinality $\leq \card{X}$ such that $\genset{G,U} = \Sym{X}$, then there exists a finite subset $F \subseteq U$ which also satisfies $\genset{G,F} = \Sym{X}$.
We say that $G$ \emph{finitely generates} $\Sym{X}$.
To see that this implication holds let us assume that the $U$ defined above is infinite and of minimal size for generating $\Sym{X}$ (there exists no subset $V \subseteq \Sym{X}$ of cardinality $< \card{U}$ such that $\genset{G,V} = \Sym{X}$).
We then well-order the set $U = \set{u_i \given i \in \card{U}}$ (this is always possible since we assume the Axiom of Choice) and note that $\genset{G,U} = \bigcup_{i \in \card{U}} G_i$ where $G_i = \genset{G, \set{u_j \mid j<i}}$.
But then, since $\card{U}$ is assumed to be the minimal number of elements one must add to $G$ in order to generate $\Sym{X}$, the above defined $G_i$ form a chain of $\leq \card{X}$ proper subgroups whose union is all of $\Sym{X}$, which contradicts \cite[Theorem 1.1]{subgroups_macpherson_neumann}.
Hence $U$ cannot be minimal and must be be reducible to a finite set.

Combining this with Galvin's result that $\Sym{X}$ has Sierpiński rank 2, we get that for any subgroup $G$ of $\Sym{X}$ either $\rank{\Sym{X} : G} \leq 2$ or $\rank{\Sym{X} : G} > \card{X}$.
Improving on this further, Galvin showed in \cite{galvin1995generating} that if the relative rank of the subgroup $G$ is finite, then its relative rank as a group is at most 1.
The results by Galvin, Macpherson, and Neumann divide the subgroups of $\Sym{X}$ into two distinct classes. Those that generate all of $\Sym{X}$ by only adding a single element, and those where it is necessary to add more than $\card{X}$ elements to generate $\Sym{X}$.
The question of how to tell which of these two classes a given subgroup belongs to is studied by Bergman and Shelah in \cite{Bergman_2006}.

\subsection{Bergman-Shelah Preorder on \Sym{X}}

In fact, Bergman and Shelah go further than simply classifying which subgroups of the symmetric group finitely generate the entire group.
By doing what is in some sense a generalisation of the concept of relative rank, they define the following preorder.

\begin{definition}[Bergman-Shelah preorder] \label{def:BS_preorder}
	Let $T$ be a semigroup, $S_1, S_2$ subsemigroups of $T$, and $\kappa$ an infinite cardinal. 
    We shall then write $S_1 \bsleq[\kappa,T] S_2$ if there exists a subset $U \subseteq T$ of cardinality strictly less than $\kappa$ such that $S_1 \subseteq \genset{S_2 , U}$. 
    If $S_1 \bsleq[\kappa,T] S_2$ and $S_2 \bsleq[\kappa,T] S_1$, we shall write $S_1 \bsequal[\kappa,T] S_2$, while if $S_1 \bsleq[\kappa,T] S_2$ and $S_2 \not\bsleq[\kappa,T] S_1$, we shall write $S_1 \bsless[\kappa,T] S_2$. 
\end{definition}

Clearly $\bsleq[\kappa,T]$ is a preorder on subsemigroups of $T$ (when $\kappa$ is an infinite cardinal) and hence $\bsequal[\kappa,T]$ is an equivalence relation, equivalent to the assertion that there exists a subset $U \subseteq T$ of cardinality strictly less than $\kappa$ such that $\genset{S_1,U} = \genset{S_2,U}$.
Since we know from theorems \ref{thm:galvin_group_rank} and \ref{thm:sierpinksi_rank_IX} that every countable subset of $\Sym{X}$ and $I_X$ are contained in a subsemigroup generated by two elements, we can conclude that $\bsleq[3,T] = \bsleq[\aleph_0,T] = \bsleq[\aleph_1,T]$ (this is a slight abuse of notation, since $\kappa$ should always be infinite) when $T \in \set{\Sym{X}, I_X}$.
As such, going forward we will set $\kappa = \aleph_0$ in the Bergman-Shelah preorder and simply omit writing it (so $\bsleq[T] = \bsleq[\aleph_0,T]$).
For a similar ease of use we will also shorten $\bsleq[\Sym{X}]$ to $\bsleq[S]$ and $\bsleq[I_X]$ to $\bsleq[I]$.
In \cite{Bergman_2006} Bergman and Shelah study the equivalence classes of subgroups of $\Sym{X}$ induced by the preorder $\bsleq[S]$. 
The main result of the paper is stated in term of something called a \emph{pointwise stabiliser}.

\begin{definition}[Pointwise stabiliser]
Let $X$ be a set and $S$ a subsemigroup of  $P_X$.
Then the \emph{pointwise stabiliser subsemigroup} in $S$ of a subset $Y \subseteq X$, denoted $\pointStab{S}{Y}$, is the subsemigroup of $S$ which consists of all functions in $S$ that fix the elements of $Y$:
\begin{equation*}
    \pointStab{S}{Y} = \set{f \in S \given (\forall y \in Y)~ yf = y}
\end{equation*}
\end{definition}

We can now state the main theorem of \cite{Bergman_2006}, also referred to as the \emph{Bergman-Shelah Theorem}.

\begin{theorem}[Bergman \& Shelah \cite{Bergman_2006}] \label{BS_mainTheorem}
    Let $X$ be a countably infinite set. 
    Then the subgroups closed in the pointwise topology on $\Sym{X}$ lie in precisely four equivalence classes under the relation $\bsequal[S]$. 
    Which of these classes a closed subgroup $G$ belongs to depends on which of the following statements about pointwise stabiliser subgroups $G_{(\Gamma)}$ of finite subsets $\Gamma \subseteq X$ holds:
	\begin{enumerate}[\normalfont(i)]
		\item For every finite set $\Gamma$, the subgroup $G_{(\Gamma)}$ has at least one infinite orbit in $X$.
        
		\item There exist finite sets $\Gamma$ such that all orbits of $G_{(\Gamma)}$ are finite, but none such that the cardinalities of these orbits have a common finite bound.
        
		\item There exist finite sets $\Gamma$ such that the cardinalities of the orbits of $G_{(\Gamma)}$ have a common finite bound, but none such that $G_{(\Gamma)} = \set{1}$.
        
		\item There exist finite sets $\Gamma$ such that $G_{(\Gamma)} = \set{1}$.
	\end{enumerate}
    Furthermore, a closed subgroup satisfying a condition on the above list is $\bsgreat[S]$ any closed subgroups satisfying conditions lower on the list.
\end{theorem}

Theorem \ref{BS_mainTheorem} is extremely interesting, as it involves only a finite number of equivalence classes and it gives defining characteristics for all subgroups in each of these classes.
The result however only applies to closed subgroups in the pointwise topology, and the picture is much less clear for non-closed subgroups.
Theorem \ref{BS_mainTheorem} uses pointwise stabilisers to define each of the four equivalence classes, and the usefulness of such subgroups arises from the following lemma.

\begin{lemma}[{\cite[Lemma 2]{Bergman_2006}}] \label{BS_lemma2}
    Let $X$ be an infinite set, $G \subseteq \Sym{X}$ a subgroup, and $\Gamma \subseteq X$ a subset such that $\card{X}^{\card{\Gamma}} \leq \card{X}$ (e.g., a finite subset). 
    Then $G_{(\Gamma)} \bsequal[\card{X}^+,S] G$.
\end{lemma}

\begin{proof}
    Elements of $G$ in distinct right cosets of $G_{(\Gamma)}$ have distinct behaviours on $\Gamma$. 
    Hence, if $R$ is a set of representatives of these right cosets, then $\card{R} \leq \card{X}^{\card{\Gamma}} \leq \card{X}$.	
    Clearly, $\genset{G_{(\Gamma)}, R} = G = \genset{G , R}$, so $G_{(\Gamma)} \bsequal[\card{X}^+,S] G$, as claimed.
\end{proof}

Lemma \ref{BS_lemma2} is a powerful result that makes it very natural to work with stabiliser subgroups when defining the different Bergman-Shelah equivalence classes of $\Sym{X}$. 
Unfortunately however, the result does not generalise to subsemigroups of $T_X$ or $I_X$ (given a subsemigroup $S$ of $T_X$ or $I_X$ and a finite subset $\Gamma$ of $X$ it does not hold in general that $\pointStab{S}{\Gamma} \bsequal S$).
Theorem \ref{BS_mainTheorem} is the main result of \cite{Bergman_2006}, and it neatly divides the closed subgroups of $\Sym{X}$ into four classes using orbits on stabiliser subgroups as the criteria. 
Similar studies have been carried out for the full transformation monoid $T_X$ in \cite{Zach_2007,fulltransformationBS_2012}, but the results are not nearly as concise. 
It turns out that the closed subsemigroups of $T_X$ in the pointwise topology do not fall into a finite number of equivalence classes. 
There is however still some structure to be found on the subsemigroups with large stabilisers, which has been studied in \cite{Zach_2007}.

\subsection{Bergman-Shelah Preorder on $I_X$}

Theorem \ref{BS_mainTheorem} is a beautiful result and the hope would be that maybe something similar holds true for the symmetric inverse monoid.
The goal of this thesis is to study the subsemigroup structure of $I_X$ and in Chapter \ref{chap:BS_IX} we will classify some equivalence classes of $I_X$ under the Bergman-Shelah equivalence relation and their relative partial order.
Before we get to doing so in the thesis proper, we will cover some already well-know results about the relative ranks of certain subsemigroups of $I_X$.
First we'll show that the relative rank of $\Sym{X}$ in $I_X$ is at most 2.

\begin{lemma} \label{SymX=IX}
    Let $X$ be an infinite set and $S = \Sym{X}$.	
    Then there exists $g \in I_X$ such that $gSg^{-1} = I_X$.	
    In particular, $\Sym{X} \bsequal[I] I_X$.
\end{lemma}

\begin{proof}
    Let $g \in I_X$ be a bijection from $X$ to $Y$, where $Y$ is a moiety of $X$, and let $f \in I_X$ be an arbitrary chart.
    We can then define the chart $h=g^{-1}fg \in I_Y$ and denote by $h' \in \Sym{X}$ any extension of $h$ with the following properties.
    \begin{enumerate}[(a)]
        \item If $x \in \dom{f}$, then $(x)gh'=(x)gh=(x)fg$.
        This is a consequence of $h'$ being an extension of $h$.
        \item If $y \notin \dom{f}$, then we must demand that $y \notin \dom{gh'g^{-1}}$, so choose $(y)gh'$ outside $Y$. 
        This is always possible, since $Y$ is a moiety.
    \end{enumerate}
    Then $g(h')g^{-1} = g(h)g^{-1} = g(g^{-1}fg)g^{-1} = f$. 
    So $f \in \genset{S,g,g^{-1}}$.
\end{proof}	 

This result is easily extended to subsemigroups of $I_X$ that `contain' the symmetric group on some infinite  subset $Y \subseteq X$.

\begin{lemma} \label{SymYSub=IX}
    Let $X$ be an infinite set, $Y$ a subset of $X$ such that $\card{Y} = \card{X}$, and $S \subseteq I_X$ a subsemigroup.	
    If $\Sym{Y}$ is a subset of $\eval{S}_{Y} = \set{\eval{f}_{Y} \in I_Y \given f \in S}$, then there exists $g \in I_X$ such that $gSg^{-1} = I_X$.	
    In particular, $S \bsequal[I] I_X$.
\end{lemma}

\begin{proof}
    Let $g \in I_X$ be a bijection from $X$ to $Z$, where $Z$ is a moiety of $Y$, and let $f \in I_X$ be an arbitrary chart.	
    We can then define the chart $h=g^{-1}fg \in I_Z$ and denote by $h' \in \Sym{Y}$ any extension of $h$ with the following properties.
    \begin{enumerate}[\normalfont(i)]
        \item If $x \in \dom{f}$, then $(x)gh'=(x)gh=(x)fg$.
        This is a consequence of $h'$ being an extension of $h$.
        \item If $y \notin \dom{f}$, then we must demand that $y \notin \dom{gh'g^{-1}}$, so choose $(y)gh'$ in $Y\setminus Z$. 
        This is always possible since, $Z$ is a moiety of $Y$.
    \end{enumerate}
    Now let $h'' \in S$ be an extension of $h'$ such that $\eval{h''}_{Y} = \eval{h'}_{Y} = h'$.
    Then $g(h'')g^{-1} = g(h')g^{-1} = g(h)g^{-1} = g(g^{-1}fg)g^{-1} = f$. 
    So $f \in \genset{S,g,g^{-1}}$.
\end{proof}

Lemmas \ref{SymX=IX} and \ref{SymYSub=IX} tell us that any subgroups of the symmetric group that are considered `large' in the Bergman-Shelah preorder on $\Sym{X}$ are also to be considered `large' in $I_X$.

\begin{corollary} \label{cor:large_Hsubgroup}
	Let $X$ be an infinite set and $Y$ a subset of $X$ such that $\card{Y} = \card{X}$.
    If $G$ is a subgroup of $\Sym{Y}$ such that $G \bsequal[\Sym{Y}] \Sym{Y}$, then $G \bsequal[I] I_X$.
\end{corollary}

\begin{proof}
    This follows from Lemma \ref{SymYSub=IX} and transitivity of the Bergman-Shelah preorder.
\end{proof}

\chapter{Maximal Subsemigroups of $I_X$} \label{chap:maximal_subsemigroups}

A proper subsemigroup $T$ of a semigroup $S$ is called a \emph{maximal subsemigroup} if $T$ is maximal with respect to containment in the set of proper subsemigroups of $S$.
That is, $T \neq S$ and $\genset{T,s} = S$ for all $s \in S \setminus T$.
Analogously, a \emph{maximal inverse subsemigroup} $M$ and a \emph{maximal subgroup} $G$ are proper subsemigroups of a semigroup $S$, such that $M$ is a maximal element in the set of proper inverse subsemigroups of $S$ and $G$ is a maximal element in the set of proper subgroups of $S$ under containment.
In this chapter we will classify the maximal (inverse) subsemigroups of $I_X$ containing certain subgroups of $\Sym{X}$ when $X$ is a countably infinite set.

In the case where $X$ is finite, a complete classification of the maximal subsemigroups of $I_X$ already exists and was proved by Xiuliang \cite{xiuliang1999classification}.
We will state and prove Xiuliang's theorem in this text for convenience (Theorem \ref{thm:IX_finite}).
The result is dependent on the maximal subgroups of $\Sym{X}$ which, for finite $X$, are in some sense fully classified (see e.g. \cite{scott1980representations, aschbacher1985maximal, liebeck1988nan}).
It is in fact a general result that maximal subsemigroups of arbitrary finite semigroups are, in some sense, determined by their maximal subgroups \cite{graham1968maximal}.

For studies of the maximal subgroups of $\Sym{X}$ when $X$ is infinite, see e.g. \cite{ball1966maximal, ball1968indices, baumgartner1993maximal, biryukov2000set, brazil1994maximal, covington1996some, macpherson1993large, subgroups_macpherson_neumann, maximal_macpherson_preager, richman1967maximal}, and for maximal subsemigroups of the full transformation semigroup $T^X$, see e.g. \cite{gavrilov1965functional, pinsker2005maximal, maximal_east_mitchell_peresse}.
A full classification of the maximal subsemigroups of most infinite semigroups is probably infeasible, so in this thesis we take inspiration from East, Mitchell, and Péresse \cite{maximal_east_mitchell_peresse} by fully classifying the maximal (inverse) subsemigroups of $I_X$ containing certain subgroups of $\Sym{X}$.
The subgroups in question are: the symmetric group $\Sym{X}$ (Theorem \ref{thm:maximal_sym}), the pointwise stabiliser of a finite non-empty subset $\Sigma$ of $X$ (Theorem \ref{thm:maximal_pointwise}), the stabiliser of an ultrafilter $\filter$ on $X$ (Theorem \ref{thm:maximal_ultrafilter}), and the stabiliser of a finite partition $\partition$ of $X$ (Theorem \ref{thm:maximal_astab}).
We will in fact develop a method (Theorem \ref{thm:total->all}) for classifying all the maximal subsemigroups of $I_X$ containing certain `sufficiently large' inverse subsemigroups and then apply this method to the four subgroups mentioned above.

\section{Preliminary Results}

The following three concepts will be used extensively throughout the rest of this chapter and in the statement of the main results:
\begin{definition}[Rank]
The \emph{rank} of a chart $f \in I_X$ is the size of the image (and domain) of $f$. We denote the rank of $f$ by $r(f) = \card{\im{f}} \,(= \card{\dom{f}})$.
\end{definition}
\begin{definition}[Defect]
The \emph{defect} of a chart $f \in I_X$ is the size of the complement of the image of $f$. We denote the defect of $f$ by $d(f) = \card{X \setminus \im{f}}$.
\end{definition}
\begin{definition}[Collapse]
The \emph{collapse} of a chart $f \in I_X$ is the size of the complement of the domain of $f$. We denote the collapse of $f$ by $c(f) = \card{X \setminus \dom{f}}$.
\end{definition}
More generally one would say that the collapse of a partial function $f$ is the size of the complement of a transversal of $f$ with respect to $X$.
It is worth noting here that if $f$ is a partial bijection, then $c(f) = d(\inv{f})$ and $d(f) = c(\inv{f})$.
So one could in principle make do with only using either defects or collapses to state the results of this text, but we found it conceptually and notationally easier to use both.
We list below a number of useful properties of ranks, defects, and collapses.

\begin{lemma} \label{lem:defect_properties}
    Let $X$ be a set, $\mu$ an infinite cardinal such that $\mu \leq \card{X}$, and $f,g \in I_X$. 
    Then the following hold:
    \begin{enumerate}[~\normalfont(i)]
        \item $c(f) \leq c(fg) \leq c(f)+c(g)$;
        \label{lem:defect_properties/collapse_composition}
        
        \item if $f$ is \emph{surjective} (i.e. $d(f)=0$), then $c(fg) = c(f)+c(g)$;
        \label{lem:defect_properties/collapse_surjective}
        
        \item $d(g) \leq d(fg) \leq d(f)+d(g)$;
        \label{lem:defect_properties/defect_composition}
        
        \item if $g$ is \emph{total} (i.e. $c(g)=0$), then $d(fg) = d(f)+d(g)$;
        \label{lem:defect_properties/defect_total}
        
        \item $r(fg) \leq \min(r(f),r(g))$;
        \label{lem:defect_properties/rank_composition}

        \item if $d(f)=0$ and $c(g)=0$, then $r(fg) = \card{X}$;
        \label{lem:defect_properties/rank_permutation}
        
        \item if $c(g) < \mu \leq d(f)$, then $d(fg) \geq \mu$;
        \label{lem:defect_properties/collapse<defect}
        
        \item if $d(f) < \mu \leq c(g)$, then $c(fg) \geq \mu$;
        \label{lem:defect_properties/defect<collapse}

        \item if $d(f) + c(g) < \card{X}$, then $r(fg) \geq \card{X} - (d(f) + c(g))$.
        \label{lem:defect_properties/rank_bound}
    \end{enumerate}
\end{lemma}

\begin{proof}
    We prove each of the statements in order as listed above.
    \begin{enumerate}[\normalfont(i)]
        \item Since $\dom{fg} \subseteq \dom{f}$, it follows that $c(f) \leq c(fg)$.
        Also, $c(fg) = \card{X \setminus \dom{fg}} = \card{X \setminus \dom{f}} + \card{\dom{f} \setminus \dom{fg}}$ so it suffices to show that $\card{\dom{f} \setminus \dom{fg}} \leq c(g)$.
        For all $x \in \dom{f} \setminus \dom{fg}$ there exists $y \in X \setminus \dom{g}$ such that $(x)f = y$. 
        Hence, $\card{\dom{f} \setminus \dom{fg}} \leq \card{X \setminus \dom{g}} = c(g)$, as required.

        \item If $f$ is surjective, then $(\dom{f} \setminus \dom{fg})f = X \setminus \dom{g}$.
        Hence, $c(fg) = \card{X \setminus \dom{fg}} = \card{X \setminus \dom{f}} + \card{\dom{f} \setminus \dom{fg}} = \card{X \setminus \dom{f}} + \linebreak \card{X \setminus \dom{g}} = c(f) + c(g)$.

        \item This follows from the fact that $X \setminus \im{g}$ is a subset of $X \setminus \im{fg}$ which is again a subset of $(X \setminus \im{f})g \cup X \setminus \im{g}$.

        \item  If $g$ is total, then $X \setminus \im{fg} = (X \setminus \im{g}) \cup (\im{g} \setminus \im{fg}) = (X \setminus \im{g}) \cup \linebreak (X \!\setminus\! \im{f})g$ and $\card{(X \setminus \im{f})g} = \card{X \setminus \im{f}} = d(f)$.
        Hence, $d(fg) = \linebreak d(f)+d(g)$.
        
        \item This follows from the fact that $\dom{fg} \subseteq \dom{f}$ and $\im{fg} \subseteq \im{g}$.

        \item For any chart $h \in I_X$, if $d(h)=0$ or $c(h)=0$, then $r(h) = \card{X}$.
        Since $f$ is surjective and $g$ is total (hence $\im{f} = \dom{g}$), it follows that $\dom{fg} = \dom{f}$ and $\im{fg} = g$.
        Thus $r(fg) = r(f) = r(g) = \card{x}$.

        \item If $c(g) < \mu$ and $d(f) \geq \mu$, then $\card{\dom{g} \cap (X \setminus \im{f})} = \card{\dom{g} \setminus \im{f}} \geq \mu$.
        Hence, $\mu \leq \card{\im{g} \setminus (\im{f})g} = \card{\im{g} \setminus \im{fg}} \leq \card{X \setminus \im{fg}} = d(fg)$.

        \item If $d(f) < \mu$ and $c(g) \geq \mu$, then $\card{\im{f} \cap (X \setminus \dom{g})} = \card{\im{f} \setminus \dom{g}} \geq \mu$.
        Hence $\mu \leq \card{\dom{f} \setminus (\dom{g})\inv{f}} = \card{\dom{f} \setminus \dom{fg}} \leq \card{X \setminus \dom{fg}}$ where $\card{X \setminus \dom{fg}} = c(fg)$.

        \item The condition $d(f) + c(g) < \card{X}$ ensures that the subtraction $\card{X} - (d(f) + c(g))$ is well defined.
        Since $\card{X \setminus (\im{f} \cap \dom{g})} \leq d(f) + c(g)$, it follows that $r(fg) = \card{\im{f} \cap \dom{g}} \geq \card{X} - (d(f) + c(g))$.
    \end{enumerate}
\end{proof}

It is shown in \cite{baumgartner1993maximal} that when $X$ is infinite, and under certain set theoretic assumptions, there exists a subgroup of $\Sym{X}$ which is not contained in any maximal subgroup of $\Sym{X}$.
However, actually constructing such a subgroup is difficult and it holds in general that if a subgroup is sufficiently `small' or `large', then it will be contained in a maximal subgroup.
This question of determining whether a given sub(semi)group is contained in a maximal sub(semi)group or not has been considered for the symmetric group (see \cite{baumgartner1993maximal, subgroups_macpherson_neumann, maximal_macpherson_preager}) and the full transformation monoid (see \cite{maximal_east_mitchell_peresse}).
We will now consider the same question as well as other questions regarding the properties of maximal subsemigroups, but initially for semigroups in general and then afterwards showcase some choice results specifically for the symmetric inverse monoid.

\begin{proposition} \label{prop:contained_relative_rank}
    Let $S$ be a semigroup and $G$ a proper subsemigroup of $S$.
    If $G$ has finite relative rank in $S$, then $G$ is contained in a maximal subsemigroup of $S$.
\end{proposition}

\begin{proof}
    This is a consequence of Zorn's lemma.
    Let $\mathcal{A}$ be the partially ordered set of all proper subsemigroups of $S$ containing $G$, ordered by inclusion.
    Given a chain $(A_i)_{i \in I}$ in $\mathcal{A}$, the union $\bigcup_{i \in I} A_i$ is a subsemigroup of $S$ (since $(A_i)_{i \in I}$ is a chain under inclusion).
    If $\bigcup_{i \in I} A_i$ is a proper subsemigroup of $S$, then $(A_i)_{i \in I}$ is bounded by an element in $\mathcal{A}$, namely $\bigcup_{i \in I} A_i$.
    So the only way for a chain $(A_i)_{i \in I}$ in $\mathcal{A}$ to be unbounded, is if the union $\bigcup_{i \in I} A_i = S$.
    But since $G$ has finite relative rank in $S$, there exists a finite set $F \subseteq S$ such that $\genset{G,F} = S$.
    Since $(A_i)_{i \in I}$ is ordered by inclusion and $\bigcup_{i \in I} A_i = S$, there must exist a $j \in I$ such that $F \subseteq A_j$.
    Since $G \subseteq A_i$ for all $i \in I$ we get that $G \cup F \subseteq A_j$, which implies that $A_j = S$, contradicting the assumption that $(A_i)_{i \in I}$ is a chain in $\mathcal{A}$.
    Hence, we conclude that every chain in $\mathcal{A}$ is bounded from above, which by Zorn's lemma implies that $\mathcal{A}$ contains a maximal element.
\end{proof}

\begin{lemma} \label{lem:containing_large_semigroup}
    Let $S$ be a semigroup, $G$ a subsemigroup with finite relative rank in $S$, and $U$ any subset of $S$.
    Then $\genset{G,U} = S$ if and only if $U$ is not contained in any maximal subsemigroup of $S$ that contains $G$.
\end{lemma}

\begin{proof}
    If $G=S$, then the proof is trivial as there are no maximal subsemigroups of $S$ containing $G$ and $\genset{G,U} = S$ for all $U \subseteq S$.
    So assuming that $G$ is a proper subsemigroup of $S$, if $U$ is contained in a maximal subsemigroup of $S$ that contains $G$, then $\genset{G,U}$ is contained in that maximal subsemigroup and so $\genset{G,U} \neq S$.
    For the converse, by Proposition \ref{prop:contained_relative_rank} any proper subsemigroup of $S$ containing $G$ is contained in a maximal subsemigroup of $S$.
    Hence, if $U$ is not contained in any maximal subsemigroup containing $G$, then $\genset{G,U} = S$.
\end{proof}

\begin{lemma} \label{lem:exactly_maximal}
    Let $S$ be a semigroup, $G$ a subsemigroup of $S$, and $\set{M_i \given i \in I}$ a set of proper subsemigroups of $S$ such that $G \subseteq M_i$ for all $i \in I$.
    If $\genset{G,U}=S$ for all $U \subseteq S$ such that $U \nsubseteq M_i$ for all $i \in I$, then the maximal subsemigroups of $S$ containing $G$ are exactly the maximal elements of $\set{M_i \given i \in I}$ under containment.
\end{lemma}

\begin{proof}
    We first show by contradiction that all maximal subsemigroups of $S$ containing $G$ must be elements of $\set{M_i \given i \in I}$.
    Suppose that $N$ is a maximal subsemigroup of $S$ such that $G \subseteq N$ and $N \notin \set{M_i \given i \in I}$.
    Since $N$ is maximal in $S$ and the $M_i$ are proper subsemigroups of $S$, it follows that $N \nsubseteq M_i$ for all $i \in I$.
    Thus by the hypothesis $N = \genset{G,N} = S$, which contradicts $N$ being a maximal subsemigroup of $S$.
    Hence any maximal subsemigroup of $S$ containing $G$ must be in $\set{M_i \given i \in I}$. 
    
    Clearly non-maximal elements of $\set{M_i \given i \in I}$ cannot be maximal subsemigroups of $S$, so it only remains to show that every maximal element of $\set{M_i \given i \in I}$ is a maximal subsemigroup of $S$.
    Let $M_j$ be maximal in $\set{M_i \given i \in I}$ under containment and let $s \in S \setminus M_j$. Then $\genset{M_j,s} = S$, since $M_j \cup \set{s}$ is not contained in any of the elements of $\set{M_i \given i \in I}$ (because $M_j$ is maximal in $\set{M_i \given i \in I}$) but does contain $G$.
\end{proof}

The following result is an analogue to \cite[Lemma 6.9]{subgroups_macpherson_neumann} and \cite[Theorems 1.5 and 1.6]{maximal_macpherson_preager} as well as \cite[Propositions 3.1 and 3.3]{maximal_east_mitchell_peresse}, adapted to the symmetric inverse monoid.
Our proof hence follows a very similar argument.

\begin{lemma} \label{lem:contained_in_maximal}    
    Let $X$ be an infinite set and $S$ a proper subsemigroup of $I_X$ satisfying either of the following:
    \begin{enumerate}[\normalfont(i)]
        \item $S$ has finite relative rank in $I_X$;
        \label{lem:contained_in_maximal/S_fin_rank}
        \item $\card{S} \leq \card{X}$.
        \label{lem:contained_in_maximal/S<X}
    \end{enumerate}
    Then $S$ is contained in a maximal subsemigroup of $I_X$.
\end{lemma}

\begin{proof}
    We prove each statement in order as listed in the lemma above.
    \begin{enumerate}[\normalfont(i)]
        \item This is a direct corollary of Proposition \ref{prop:contained_relative_rank}.
    
        \item Let $V = \set{f \in S \given c(f) = 0}$ be the set of total charts in $S$.
        We then pick a point $y \in X$ and a subset $Y \subseteq X$ such that $y \in Y$, $(\forall f \in V)~ Y \cap \im{f} \neq \emptyset$, and $\card{X \setminus Y} = \card{X}$.
        Such a set $Y$ always exists since all elements of $V$ are total (and injective), meaning that $\card{\im{f}} = \card{X}$ for all $f \in V$.
        This implies that the `worst case scenario' is one in which all the images of elements of $V$ are pairwise disjoint moieties of $X$ (this is possible since $\card{V} \leq \card{S} \leq \card{X}$), in which case we can simply pick one point from each image that goes into $Y$ while still leaving $\card{X}$ many points unpicked.
        In any scenario where some of the images have non-empty intersection we can simply pick a single point from said intersection to count for both of these images.
        Furthermore, note that $Y$ is never empty due to the element $y$.
        Next, let $T = \set{f \in I_X \given Y \cap \dom{f} = \emptyset}$.
        Then $\genset{S,T}$ is a proper subsemigroup of $I_X$, since every total chart in $\genset{S,T}$ must be an element of $S$ and $\card{S} \leq \card{X}$.
        $T$ has finite relative rank in $I_X$ (let $g \in I_X$ be a total chart with $\im{g} = X \setminus Y$, then $\genset{T,g} = I_X$), so $\genset{S,T}$ also has finite relative rank in $I_X$.
        It follows from \eqref{lem:contained_in_maximal/S_fin_rank} that $\genset{S,T}$ is contained in a maximal subsemigroup of $I_X$, and since $S \subseteq \genset{S,T}$, then so is $S$.
    \end{enumerate}
\end{proof}

Lemma \ref{lem:contained_in_maximal} tells us that all sufficiently `large' and sufficiently `small' subsemigroups of $I_X$ will be contained in a maximal subsemigroup of $I_X$.
We will later expand this result by including topological conditions (Proposition \ref{prop:contained_in_maximal}), but for now this is sufficient.

As hinted in the introduction of this chapter, we will classify both maximal subsemigroups and maximal \emph{inverse} subsemigroups of $I_X$.
To this end we will prove a result (Proposition \ref{prop:inverse_intersection}), which states that the maximal inverse subsemigroups of any inverse semigroup can be derived from its maximal subsemigroups.
But first we must go over some basic properties of inverse semigroups.

\begin{lemma} \label{lem:inverse_maximal}
Let $S$ be an inverse semigroup and $M$ a maximal subsemigroup of $S$. 
Then the inverse $M^{-1} = \set{\inv{s} \in S \given s \in M}$ is also a maximal subsemigroup of $S$.
\end{lemma}

\begin{proof}
Inversion is an anti-isomorphism from $S$ to $S$.
\end{proof}

\begin{lemma} \label{lem:intersection_of_inverses}
    Let $S$ be an inverse semigroup and $M$ a subsemigroup of $S$.
    Then $M \cap \inv{M}$ is the largest (with respect to containment) inverse subsemigroup of $S$ which is contained in $M$.
\end{lemma}

\begin{proof}
    Being the intersection of two subsemigroups, $M \cap \inv{M}$ is itself a subsemigroup of $S$ and by the definition of intersection it is contained in $M$.
    Moreover, $\inv{(M \cap \inv{M})}=\inv{M}\cap M$ is closed under taking inverses.
    Finally, if $s \in M$ and $\inv{s} \in M$, then $s,\inv{s} \in M \cap \inv{M}$. 
    Hence all inverse subsemigroups of $M$ must be contained in $M \cap \inv{M}$.
\end{proof}

With this we are ready to state our main tool for classifying the maximal inverse subsemigroups of an inverse semigroup.

\begin{proposition} \label{prop:inverse_intersection}
Let $S$ be an inverse semigroup, $G$ an inverse subsemigroup with finite relative rank in $S$, and $\set{M_i \given i \in I}$ the set of all maximal subsemigroups of $S$ containing $G$.
Then the maximal inverse subsemigroups of $S$ containing $G$ are exactly the maximal elements of $\set{M_i \cap M_i^{-1} \given i \in I}$ under containment.
\end{proposition}

\begin{proof}
    By Lemma \ref{lem:intersection_of_inverses} the set $\mathcal{M}=\set{M_i \cap M_i^{-1} \given i \in I}$ consists of inverse subsemigroups of $S$. Since $G$ is a subset of $M_i$ for all $i \in I$, it is also a subset of $M_i \cap M_i^{-1}$ for all $i \in I$.
    To show that any maximal inverse subsemigroup $U$ of $S$ containing $G$ lies in $\mathcal{M}$, we note that $U$ must have finite relative rank in $S$ (it indeed has a relative rank of at most 2).
    It follows from Proposition \ref{prop:contained_relative_rank} that $U$ must be contained in some maximal subsemigroup $M_i$ of $S$ and so, by Lemma \ref{lem:intersection_of_inverses}, $U \subseteq M_i \cap \inv{M_i}$.
    In fact $U = M_i \cap \inv{M_i}$ since $U$ is a maximal inverse subsemigroup of $S$.
    Hence, $U$ is in $\mathcal{M}$ and clearly $U$ is a maximal element of $\mathcal{M}$.
    
    On the other hand, let $M_i \cap \inv{M_i}$ be a maximal element of $\mathcal{M}$ and suppose that $V$ is a proper inverse subsemigroup of $S$ containing $M_i \cap \inv{M_i}$.
    Since $G \subseteq M_i \cap \inv{M_i} \subseteq V$ and $G$ has finite relative rank in $S$, it follows from Proposition \ref{prop:contained_relative_rank} that $V$ must be contained in a maximal subsemigroup $M_j$ of $S$. 
    Then, by Lemma \ref{lem:intersection_of_inverses}, $M_i \cap \inv{M_i} \subseteq V \subseteq M_j \cap \inv{M_j} \in \mathcal{M}$. 
    Since we assumed $M_i \cap \inv{M_i}$ to be a maximal element of $\mathcal{M}$, it follows that $M_j \cap \inv{M_j} = M_i \cap \inv{M_i}$. In particular, $V=M_i \cap \inv{M_i}$, as required.
\end{proof}

We are now able to find the maximal inverse subsemigroups of any inverse semigroup given its maximal subsemigroups.
But this means that we have to find those maximal subsemigroups first.
One could do this using Lemma \ref{lem:exactly_maximal}, but we would like to to have a lemma that is more directly applicable to the symmetric inverse monoid.

\begin{lemma} \label{lem:in-and-out}
    Let $X$ be an infinite set, $Y$ a subset of $X$, $G$ a subset of $I_X$ such that $\inv{G} = G$, and $\set{M_i \given i \in I}$ a collection of subsets of $I_X$.
    Then the following are equivalent:
    \begin{enumerate}[~\normalfont(i)]
        \item For all subsemigroups $U$ of $I_X$ such that $G \subseteq U$ and $U \nsubseteq M_i$ for all $i \in I$, there exists $f \in U$ with $\dom{f} = X$ and $\im{f}$ a moiety of $Y$.
        \item For all subsemigroups $U$ of $I_X$ such that $G \subseteq U$ and $U \nsubseteq \inv{M_i}$ for all $i \in I$, there exists $g \in U$ with $\dom{g}$ a moiety of $Y$ and $\im{g} = X$.
    \end{enumerate}
\end{lemma}

\begin{proof}
    (i)$\implies$(ii): Assuming (i) to be true, we let $U$ be a subsemigroup of $I_X$ such that $G \subseteq U$ and $U \nsubseteq \inv{M_i}$ for all $i \in I$.
    Then $\inv{U}$ satisfies the hypothesis of (i), that is $G \subseteq \inv{U}$ and $\inv{U} \nsubseteq M_i$ for all $i \in I$.
    So by (i) there exists $f \in \inv{U}$ such that $\dom{f} = X$ and $\im{f}$ is a moiety of $Y$.
    Thus $\inv{f} \in U$ satisfies the conditions that $\dom{\inv{f}}$ is a moiety of $Y$ and $\im{\inv{f}} = X$.
    Let $g = \inv{f}$.

    (ii)$\implies$(i): Proof of this is entirely analogous to that of the converse.
\end{proof}

It is worth noting that Lemma \ref{lem:in-and-out} is a condition on the sets $Y$, $G$, and $\set{M_i \given i \in I}$ rather than on any given subset $U$.
We are now ready to state our main tool for classifying maximal subsemigroups of the symmetric inverse monoid.

\begin{theorem} \label{thm:total->all}
    Let $X$ be an infinite set, $G$ an inverse subsemigroup of $I_X$, and $\set{M_i \given i \in I}$ a collection of subsets of $I_X$.
   If the following six conditions hold, then $\set{M_i \given i \in I}$ is the set of all maximal subsemigroups of $I_X$ containing $G$.
    \begin{enumerate}[~\normalfont(i)]
        \item $G \subseteq M_i$ for all $i \in I$.
        \label{thm:total->all/G_subset}
        \item $M_i$ is a proper subsemigroup of $I_X$ for all $i \in I$.
        \label{thm:total->all/M_subsemigroup}
        \item $M_i \nsubseteq M_j$ whenever $i \neq j$.
        \label{thm:total->all/M_not_contained}
        \item For all $i \in I$ there exists $j \in I$ such that $M_j = \inv{M_i}$.
        \label{thm:total->all/M_inverses}
        \item There exists $Y \subseteq X$ such that $\card{Y} = \card{X}$ and $\Sym{Y} \subseteq \eval{G}_Y = \left\{\eval{g}_Y \;\middle|\; g \in G\right\}$.
        \label{thm:total->all/YSym}
        \item For all subsemigroups $U$ of $I_X$ such that $G \subseteq U$ and $U \nsubseteq M_i$ for all $i \in I$, there exists $f \in U$ such that $\dom{f} = X$ and $\im{f}$ is a moiety of $Y$ (where $Y$ is the same as in condition \eqref{thm:total->all/YSym}).
        \label{thm:total->all/in-out-condition}
    \end{enumerate}
\end{theorem}

\begin{proof}    
    Let $U$ be a subsemigroup of $I_X$ such that $G \subseteq U$ and $U \nsubseteq M_i$ for all $i \in I$.
    We will show that $U = I_X$.
    By condition \eqref{thm:total->all/in-out-condition} there exists $f \in U$ such that $\dom{f} = X$ and $\im{f}$ is a moiety of $Y$, and so, by Lemma \ref{lem:in-and-out}, there exists $g \in U$ with $\dom{g}$ a moiety of $Y$ and $\im{g} = X$ (since $\set{\inv{M_i} \given i \in I} = \set{M_i \given i \in I}$ by condition \eqref{thm:total->all/M_inverses}). \\
    Let $h \in I_X$ be any chart and $p=f^{-1}hg^{-1} \in I_X$.
    Then $\dom{p}=(\dom{hg^{-1}})f$ and $\dom{hg^{-1}}=\dom{h}$ since $g^{-1}$ is total.
    In particular, $\dom{p}=(\dom{h})f\subseteq \im{f}$ and $\im{p} \subseteq \dom{g}$. Since $\im{f}$ and $\dom{g}$ are moieties of $Y$, we may extend $p$ to an element $p'$ of $\Sym{Y}$ such that $(X \setminus \dom{h})f$ is mapped into $Y \setminus \dom{g}$ under $p'$. By condition \eqref{thm:total->all/YSym}, and since $G \subseteq U$, we may extend $p'$ to an element $p''$ of $U$. We now claim that $h=fp''g$.
    \begin{enumerate}[~\normalfont(a)]
        \item If $x \in \dom{h}$, then $(x)f \in \dom{p}$ and so $(x)fp''g = (x)f(f^{-1}h\inv{g})g=(x)h$;
        \item If $x \in X \setminus \dom{h}$, then $(x)fp''=(x)fp'$ lies in $Y \setminus \dom{g}$ and so $(x)fp''g=(x)fp'g$ is not defined. Hence, $x \notin \dom{fp''g}$.
    \end{enumerate}
    This confirms our claim that $h = fpg = fp'g = fp''g \in U$, and since $h$ was an arbitrary element of $I_X$, we conclude that $U = I_X$.
    We have shown that any subsemigroup $U$, which contains $G$ and is not contained in $M_i$ for all $i \in I$, is equal to $I_X$.
    Lemma \ref{lem:exactly_maximal} together with conditions \eqref{thm:total->all/G_subset}, \eqref{thm:total->all/M_subsemigroup}, and \eqref{thm:total->all/M_not_contained} then tells us that the set of maximal subsemigroups of $I_X$ containing $G$ is precisely $\set{M_i \given i\in I}$.
\end{proof}

Due to the symmetric nature of Lemma \ref{lem:in-and-out}, one could replace the chart $f$ in condition \eqref{thm:total->all/in-out-condition} of Theorem \ref{thm:total->all} with the later defined chart $g$ in the proof above.
As a final preliminary result we make the following observation.

\begin{lemma} \label{lem:small_charts}
    Let $X$ be a set, $\mu$ a cardinal such that $\mu \leq \card{X}$, and $\fin(X,\mu)$ the subsemigroup of all charts in $I_X$ with rank less than $\mu$.
    \begin{equation*}
        \fin(X,\mu) = \set{f \in I_X \given r(f) < \mu}
    \end{equation*}
    Then the following hold:
    \begin{enumerate}[~\normalfont(i)]
        \item $\fin(X,\mu)$ is an ideal of $I_X$;
        \label{lem:small_charts/F_ideal}
        
        \item if $S \subseteq I_X$ is an (inverse) subsemigroup, then so is $S \cup \fin(X,\mu)$;
        \label{lem:small_charts/S_cup_F_semigroup}
        
        \item if $X$ is infinite, then $\fin(X,\mu)$ is contained in every maximal (inverse) subsemigroup of $I_X$.
        \label{lem:small_charts/F_contained_maximal}
    \end{enumerate}
\end{lemma}

\begin{proof}
    We prove each statement in order as listed in the lemma above.
    \begin{enumerate}[~\normalfont(i)]
        \item This follows from Lemma \ref{lem:defect_properties} condition \eqref{lem:defect_properties/rank_composition}, which states that the rank of a composite chart is always less than or equal to the ranks of its factors.
        \item This follows from \eqref{lem:small_charts/F_ideal} and the fact that $\fin(X,\mu)$ is an inverse subsemigroup of $I_X$ ($\fin(X,\mu)$ is clearly closed under taking inverses, since for all $f \in I_X$, $r(f) = r(\inv{f})$).
        \item If $M$ is a maximal (inverse) subsemigroup of $I_X$, then, by \eqref{lem:small_charts/F_ideal} and \eqref{lem:small_charts/S_cup_F_semigroup} above, either $\fin(X,\mu) \subseteq M$ or $M \cup \fin(X,\mu) = I_X$.
        In the latter case $I_X \setminus \fin(X,\mu)$ is a subset of $M$. 
        But when $X$ is infinite, then $\genset{I_X \setminus \fin(X,\mu)} = I_X$ (any chart $f$ with $r(f) < \card{X}$ can be generated by composing charts $g,h \in I_X \setminus \fin(X,\mu)$ for which $\card{\im{g} \cap \dom{h}} = r(f)$), which contradicts the assumption that $M$ is maximal.
        Hence we are only left with the option that $\fin(X,\mu) \subseteq M$.
    \end{enumerate}
\end{proof}

Going forward we will use the shorthand $\fin_X = \fin(X,\card{X})$ to denote the ideal of all charts in $I_X$ with rank less than $\card{X}$ (also called the \emph{small charts} or \emph{finite charts} in the case of $X$ countably infinite).
With this we are finally ready to classify maximal (inverse) subsemigroups of $I_X$.
We begin with the finite case.

\section{The Finite Case}

In \cite{xiuliang1999classification} the maximal inverse subsemigroups of finite symmetric inverse monoids are classified.
In this section we will recount the result and show that it is in fact also a complete classification of the maximal (not necessarily inverse) subsemigroups of finite symmetric inverse monoids.

\begin{lemma}[{\cite[Theorem 3.1]{gomes1987ranks}}] \label{lem:n-1_generate}
    Let $n \geq 1$ be a natural number and $f \in I_n$ a chart such that $r(f) = n-1$.
    Then $\genset{\Sym{n},f} = I_n$.
\end{lemma}

\begin{proof}
    Let $g \in \fin_n$ be an arbitrary chart with $r(g) = m < n$.
    We will show that $g \in \genset{\Sym{n},f}$.
    We start by constructing a chain of permutations $(a_i)_{i \in n-m}$ where $a_i \in \Sym{n}$ for all $i$.
    Let $a_0 \in \Sym{n}$ be a permutation such that $(n \setminus \im{f}) a_0 \subseteq \dom{f}$.
    For the remaining $i \in n-m$ recursively define $a_i$ such that $(n \setminus \im{fa_0 \dots fa_{i-1}f}) a_i \subseteq \dom{f}$.
    Then the composite chart $\alpha = fa_0 \dots fa_{n-m-1} \in \genset{\Sym{n},f}$ has rank $m$.
    Finally, we pick permutations $b,c \in \Sym{n}$ such that $(\dom{g})b = \dom{\alpha}$ and $\inv{\alpha} \inv{b} g \subseteq c$.
    Then $g = b \alpha c$, as required.
\end{proof}

\begin{theorem}[{\cite[Theorem 2.3]{xiuliang1999classification}}] \label{thm:IX_finite}
    Let $n \geq 1$ be a natural number.
    Then the maximal subsemigroups of $I_n$ are:
    \begin{align*}
        S &= \Sym{n} \cup \fin(n,n-1) \\
        S_G &= G \cup \fin_n
    \end{align*}
    where $G$ is a maximal subgroup of $\Sym{n}$.
\end{theorem}

\begin{proof}
    It follows from Lemma \ref{lem:small_charts} condition \eqref{lem:small_charts/S_cup_F_semigroup} that $S$ and $S_G$ are indeed semigroups.
    $S$ being maximal follows from Lemma \ref{lem:n-1_generate} (since all elements of $I_n \setminus S$ have rank $n-1$), and $S_G$ being maximal follows from the fact that maximal subgroups of $\Sym{n}$ are also maximal subsemigroups (since $n$ is finite, so all elements of $\Sym{n}$ have finite order).
    All that remains is to show that there are no other maximal subsemigroups of $I_n$.

    Let $T$ be a maximal subsemigroup of $I_X$.
    If $\Sym{n} \subseteq T$, then $T$ cannot contain any charts with rank $n-1$, as otherwise we would get that $T = \Sym{n}$ by Lemma \ref{lem:n-1_generate}.
    As any semigroup satisfying these criteria is clearly a subset of $S$, we can conclude that $S$ is the unique maximal subsemigroup of $I_X$ containing $\Sym{X}$.
    If $\Sym{n} \nsubseteq T$, then there must exist a maximal sub(semi)group $G$ of $\Sym{n}$ such that $G \subseteq T$ (since $\fin_n$ is an ideal of $I_n$, meaning that no permutations can be generated using elements of $\fin_n$).
    Again, any semigroup satisfying these conditions must be a subset $S_G$, meaning that the maximal subsemigroups of $I_n$ not containing $\Sym{n}$ must be of the form $S_G = G \cup \fin_n$, as required.
\end{proof}

Theorem \ref{thm:IX_finite} above shows that the semigroups, which in \cite{xiuliang1999classification} were shown to be the maximal inverse subsemigroups of $I_n$, are actually also the maximal subsemigroups of $I_n$.
We will now show that a simple observation quickly recovers the result that they are the maximal inverse subsemigroups of $I_n$.

\begin{corollary}
    Let $n \geq 1$ be a natural number.
    Then the maximal inverse subsemigroups of $I_n$ are exactly the semigroups described in Theorem \ref{thm:IX_finite}.
\end{corollary}

\begin{proof}
    It should be clear that the semigroups described in Theorem \ref{thm:IX_finite} are indeed inverse semigroups (since groups are inverse semigroups and if $f \in I_n$, then $r(f) = r(\inv{f})$).
    With this information in mind, the proof of Theorem \ref{thm:IX_finite} then doubles as a proof that the described semigroups are the maximal inverse subsemigroups of $I_n$.
    Alternatively we can simply apply Proposition \ref{prop:inverse_intersection}, where we let $G = \emptyset$ (since finite semigroups are finitely generated).
\end{proof}

All of the previous results in this section have demanded that the natural number $n$ be no less than 1.
Not feeling like leaving out our favourite number, let us briefly discuss the case when $n=0$.
We will first make sense of the notion of a function on the empty set.
Clearly $\emptyset \subseteq \emptyset \times \emptyset$, and since $\emptyset$ has no elements, it vacuously satisfies any `for all' conditions on said elements.
Hence, $\emptyset$ is the one and only function from the empty set to itself.
This means that $\Sym{0} = I_0 = T_0 = 0^0 = 1 = \set{0} = \set{\emptyset}$ (remember that we define natural numbers to be ordinals and that $T_X = X^X$).

\begin{proposition}
    The maximal (inverse) subsemigroup of $I_0$ is $\emptyset$. 
\end{proposition}

\begin{proof}
    $I_0$ has exactly one element.
\end{proof}

\section{The Symmetric Group}

In this section we classify the maximal (inverse) subsemigroups of $I_X$ that contain $\Sym{X}$, when $X$ is a countably infinite set.

\begin{theorem} \label{thm:maximal_sym}
Let $X$ be a countably infinite set.
Then the maximal subsemigroups of $I_X$ which contain $\Sym{X}$ are:
\begin{align*}
    S &= \set{f \in I_X \given c(f)>0 ~\lor~ d(f)=0} \\
    \inv{S} &= \set{f \in I_X \given c(f)=0 ~\lor~ d(f)>0} \\
    T &= \set{f \in I_X \given c(f)=\aleph_0 ~\lor~ d(f)<\aleph_0} \\
    \inv{T} &= \set{f \in I_X \given c(f)<\aleph_0 ~\lor~ d(f)=\aleph_0}
\end{align*}
\end{theorem}

\begin{proof}
    We will prove this by showing that $S$, $\inv{S}$, $T$ and $\inv{T}$ satisfy the conditions of Theorem \ref{thm:total->all} with $G = \Sym{X}$.
    \begin{enumerate}[\normalfont(i)]
        \item Given any permutation $p \in \Sym{X}$ we have that $d(p)=0<\aleph_0$ and $c(p)=0<\aleph_0$, so $\Sym{X}$ is contained in all the semigroups $S$, $\inv{S}$, $T$ and $\inv{T}$.
        \item We must show that $S$, $\inv{S}$, $T$ and $\inv{T}$ are all proper subsemigroups of $I_X$.
        We define the following sets:
        \begin{align*}
            I &= \set{f \in I_X \given c(f) > 0}
            \\
            J &= \set{f \in I_X \given d(f) = 0}
            \\
            K &= \set{f \in I_X \given c(f) = \aleph_0}
            \\
            L &= \set{f \in I_X \given d(f) < \aleph_0}
        \end{align*}
        Using Lemma \ref{lem:defect_properties}, we see that $I$ and $K$ are right ideals by condition \eqref{lem:defect_properties/collapse_composition}, $J$ and $L$ are semigroups by condition \eqref{lem:defect_properties/defect_composition}, $JI \subseteq I$ by condition \eqref{lem:defect_properties/collapse_surjective}, and $KL \subseteq L$ by condition \eqref{lem:defect_properties/defect<collapse}.
        Thus $S = I \cup J$ and $T = K \cup L$ are semigroups.
        $\inv{S}$ and $\inv{T}$ being semigroups will follow from them being the inverses of $S$ and $T$, and all being proper subsemigroups will follow from none of them being contained in any of the other.
        \item We must show that none of the semigroups $S$, $\inv{S}$, $T$ or $\inv{T}$ are contained in any of the other.        
        To show that $S$ is not contained in any of other semigroups mentions above, choose $f,g,h \in I_X$ such that:
        \begin{itemize}
            \item $c(f)>0$ and $d(f)=0$;
            \item $0 < c(g) < \aleph_0$ and $d(g) = \aleph_0$;
            \item $c(h) = \aleph_0$ and $d(h) < \aleph_0$.
        \end{itemize}
        Then $f,g,h \in S$ but $f \notin \inv{S}$, $g \notin T$, and $h \notin \inv{T}$.        
        To show that $T$ is not contained in any of other semigroups mentions above, choose $f,g,h \in I_X$ such that:
        \begin{itemize}
            \item $d(f) < \aleph_0$ and $c(f) = \aleph_0$;
            \item $0 < d(g) < \aleph_0$ and $c(g)=0$;
            \item $d(h) = 0$ and $c(h) > 0$.
        \end{itemize}
        Then $f,g,h \in T$ but $f \notin \inv{T}$, $g \notin S$, and $h \notin \inv{S}$.        
        The rest will follow from $\inv{S}$ and $\inv{T}$ being the inverses of $S$ and $T$.
        \item $\inv{S}$ and $\inv{T}$ are indeed the inverses of $S$ and $T$ respectively, since $c(f) = d(\inv{f})$.
        So $f$ has simply been swapped for $\inv{f}$ in each of the definitions.
        \item Let $Y=X$, then clearly $\Sym{Y} \subseteq \eval{\Sym{X}}_Y$.
        \item Let $U$ be a subsemigroup of $I_X$ containing $\Sym{X}$, but which is itself not contained $S$ or $T$.
        Then there exists charts $\Bar{s},\Bar{t} \in U$ such that $\Bar{s} \notin S$ and $\Bar{t} \notin T$.
        That is:
        \begin{nalign} \label{notin_S_T}
            c(\Bar{s})=0 ~&\wedge~ d(\Bar{s})>0 \\
            c(\Bar{t})<\aleph_0 ~&\wedge~ d(\Bar{t})=\aleph_0
        \end{nalign}
        Since $d(\Bar{s})>0$ and $c(\Bar{t})<\aleph_0$, there exists $n \in \omega$ and $a \in \Sym{X}$ such that $\im{\Bar{s}^n a} \subseteq \dom{\Bar{t}}$.
        It follows from $c(\Bar{s})=0$ and $d(\Bar{t})=\aleph_0$ that $f = \Bar{s}^n a \Bar{t}$ is a total chart mapping $X$ into a moiety of itself (see Figure \ref{fig:max_SymX}, for an explanation of how to read the diagram see below).
        \label{proof:maximal_sym/in-out}
    \end{enumerate}
\end{proof}

\begin{figure}
    \centering
    \includegraphics[height=7cm]{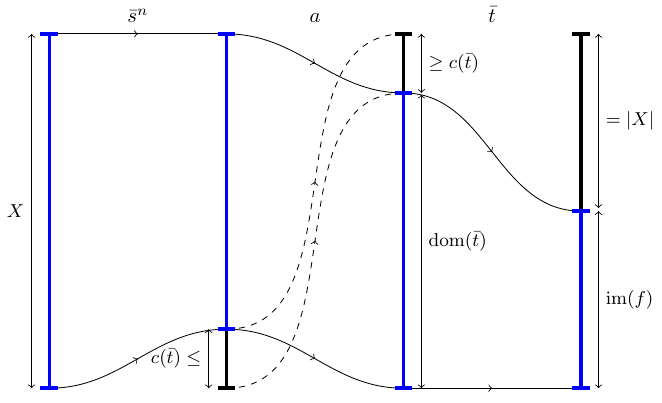}
    \caption{The chart $f = \Bar{s}^n a \Bar{t}$ in the proof of Theorem \ref{thm:maximal_sym} item \eqref{proof:maximal_sym/in-out}.}
    \label{fig:max_SymX}
\end{figure}

Diagrams like the one seen in Figure \ref{fig:max_SymX} will be used throughout this chapter as visual assistance to the proofs (they can help keep track of where points are mapped under composition of charts, but are not necessary to understand the proofs).
The following list explains the individual components that make up such a diamgram.
\begin{itemize}
    \item The coulombs (vertical bars) all represent the set $X$ and the points contained in it.

    \item The parts of the coulombs that are colour-coded \textcolor{blue}{blue} indicate the subset of $X$ which is included in the domain and image at each step of step of the function composition. Colours other than blue will be used later to indicate parts deserving of extra attention.

    \item The vertical lines with the arrowheads at each end are used to highlight certain subsets of $X$ and will typically be accompanied by a label relaying some information about the subset (such as it being the image/domain of some chart or the cardinality of the subset).

    \item The solid lines with an arrowhead in the middle indicate how a chart (the one written above the lines) maps a subset of $X$ to another subset of $X$.

    \item The dashed lines are used to either indicate where specific subsets are mapped or where a subset would have been mapped by the current chart, if it had been included in the image of the previous chart.
\end{itemize}
With the above list at hand, we will be using Figure \ref{fig:max_SymX} as an example of how to read the diagram.
\begin{enumerate}[~\normalfont(i)]
    \item The entirety of the first coulomb is coloured \textcolor{blue}{blue}, indicating that the composite chart $f$ is total.

    \item The first function section (the space between the first and the second coulomb) depicts the chart $\Bar{s}^n$ (indicated at the top of the section), which is total and has a defect of size $\geq c(\Bar{t})$ (indicated by the highlighted section at the bottom of the second coulomb).

    \item The second function section depicts how the permutation $a$ maps the image of $\Bar{s}^n$ into the domain of $\Bar{t}$. The dashed lines indicate how the permutation $a$ would have mapped the points outside the image of $\Bar{s}^n$ to the points outside the domain of $\Bar{t}$, had they been acted on by the permutation.

    \item The third and final function section depicts how the chart $\Bar{t}$ maps the image of $\Bar{s}^n a$ to a subset of $X$, whose complement has cardinality $\card{X}$. Since the composite chart $f = \Bar{s}^n a \Bar{t}$ has domain $X$, we can thus conclude that $f$ maps $X$ into a moiety of itself.
\end{enumerate}
With this minor intermission on how to read function diagrams out of the way, we will return to proving results about maximal subsemigroups of $I_X$.

\begin{corollary} \label{cor:maximal_inv_sym}
    Let $X$ be a countably infinite set.
    Then the maximal inverse subsemigroups of $I_X$ which contain $\Sym{X}$ are:
    \begin{align*}
        \mathcal{S} = S \cap \inv{S} &= \set{f \in I_X \given (c(f)>0 ~\land~ d(f)>0)} \cup \Sym{X} \\
        \mathcal{T} = T \cap \inv{T} &= \set{f \in I_X \given c(f)=d(f)=\aleph_0 ~\lor~ (c(f)<\aleph_0 \land d(f)<\aleph_0)}
    \end{align*}
\end{corollary}

\begin{proof}
    This follows from Proposition \ref{prop:inverse_intersection}, which states that the maximal inverse subsemigroups containing some inverse subsemigroup $G$ are given by the intersections of the maximal subsemigroups that contain $G$ with their inverses.
    So we let $G = \Sym{X}$ and use the semigroups described in Theorem \ref{thm:maximal_sym} as our set of maximal subsemigroups.
    All that remains is to show that neither the semigroup $\mathcal{S}$ nor $\mathcal{T}$ is contained in the other. Choose $f,g \in I_X$ such that:
    \begin{itemize}
        \item $c(f) = \aleph_0$ and $0 < d(f) < \aleph_0$;
        \item $c(g)=0$ and $0 < d(f) < \aleph_0$.
    \end{itemize}
    Then $f \in \mathcal{S}$ and $g \in \mathcal{T}$, but $f \notin \mathcal{T}$ and $g \notin \mathcal{S}$.
\end{proof}

The following is an interesting and noteworthy corollary to the above results.

\begin{corollary} \label{cor:contain_idempotents}
    Let $X$ be a countably infinite set.
    Then the semilattice of idempotents $E_X$ is contained in every maximal (inverse) subsemigroup of $I_X$.
    Said differently, every maximal (inverse) subsemigroup of $I_X$ is full.
\end{corollary}

\begin{proof}
    Suppose $M$ is a maximal (inverse) subsemigroup of $I_X$ such that $E_X \nsubseteq M$.
    It then follows from $M$ being maximal that $\genset{M,E_X} = I_X$ (idempotents are self-inverse, hence no extra assumptions are needed even if $M$ is only maximal as an inverse subsemigroup).
    This implies that every chart in $I_X$ can be written as a finite product of elements from $M$ and $E_X$.
    We should in particular be able to write any permutation $p \in \Sym{X}$ in this form, since $\Sym{X} \subseteq I_X$.
    That is, we can write any permutation $p \in \Sym{X}$ as $p = e_0m_0 \dots e_{n-1}m_{n-1}e_n$, where $n \in \omega$, $e_i \in E_X$, and $m_i \in M$ for all $i \leq n$ (we can write $p$ like this without loss of generality, as the identity is an element of $E_X$).    
    Furthermore, since $p$ is a permutation it is in particular a total chart.
    However, when composing charts $f,g \in I_X$, if the composite $fg$ is to be total, then we must have that $f$ is total and that $\im{f} \subseteq \dom{g}$.
    Hence, we can conclude that $e_0$ must be the identity and that for each $i \in n$, $\im{m_i} \subseteq \dom{e_{i+1}}$ (furthermore, we can conclude the $e_n$ must also be the identity since $p$ is surjective).
    Since the charts $e_i$ are all partial identities it then follows that they simply cancel in the expression for $p$ (that is, $m_{i}e_{i+1} = m_{i}$ for all $i \in n$ since $\im{m_i} \subseteq \dom{e_{i+1}}$ and $e_{i+1} \in E_X$).
    Hence, we can conclude that $p = e_0m_0 \dots e_{n-1}m_{n-1}e_n = m_0 \dots m_{n-1}$.
    Therefore, any product of elements from $M$ and $E_X$, which is equal to a permutation, can be reduced to a product of elements from $M$, meaning that the said permutation is already an element of $M$.
    We can thus conclude that $\Sym{X}$ is contained in $M$ (the same conclusion can be reached using a similar argument based on the fact that the permutation $p$ must be surjective).
    So for $M$ to be maximal, it must be one of the semigroups described in Theorem \ref{thm:maximal_sym} or Corollary \ref{cor:maximal_inv_sym}.
    It is however easy to check that $E_X$ is contained in all of these semigroups, hence all maximal (inverse) subsemigroups of $I_X$ must contain $E_X$.
\end{proof}

\section{The Pointwise Stabiliser}

In this section we classify the maximal (inverse) subsemigroups of $I_X$ containing the pointwise stabiliser of a finite non-empty subset of $X$, when $X$ is a countably infinite set.
For convenience we recount the definition a pointwise stabiliser subgroup here.

\begin{definition}[Pointwise Stabiliser]
The \emph{pointwise stabiliser} of a subset $Y \subseteq X$, denoted $\pointStab{\Sym{X}}{Y}$, is the group of permutations that fix all the elements of $Y$ while mapping $X \setminus Y$ to $X \setminus Y$.
\begin{equation*}
    \pointStab{\Sym{X}}{Y} = \set{f \in \Sym{X} \given (\forall y \in Y)~ yf = y}
\end{equation*}
\end{definition}

Similarly there is a larger subgroup called the setwise stabiliser.

\begin{definition}[Setwise Stabiliser]
The \emph{setwise stabiliser} of a subset $Y \subseteq X$, denoted $\setStab{\Sym{X}}{Y}$, is the group of permutations that map $Y$ to $Y$ and $X \setminus Y$ to $X \setminus Y$.
\begin{equation*}
    \setStab{\Sym{X}}{Y} = \set{f \in \Sym{X} \given Y f = Y}
\end{equation*}
\end{definition}

Setwise and pointwise stabiliser subgroups are of particular interest, as whenever $\Sigma$ is a non-empty finite subset of $X$, then $\setStab{\Sym{X}}{\Sigma}$ is a maximal subgroup of $\Sym{X}$ \cite{ball1966maximal}.

\begin{theorem} \label{thm:maximal_pointwise}
Let $X$ be a countably infinite set and $\Sigma$ a non-empty finite subset of $X$.
Then the maximal subsemigroups of $I_X$ which contain the pointwise stabiliser $\pointStab{\Sym{X}}{\Sigma}$ but not $\Sym{X}$ are:
\begin{align*}
    P_{\Gamma} &= \set{f \in I_X \given \Gamma \nsubseteq \dom{f} ~\lor~ \Gamma f = \Gamma} \cup \fin_X 
    \\
    \inv{P_{\Gamma}} &= \set{f \in I_X \given \Gamma \nsubseteq \im{f} ~\lor~ \Gamma f = \Gamma} \cup \fin_X 
    \\
    Q_{\Gamma} &= \set{f \in I_X \given \Gamma \nsubseteq \dom{f} ~\lor~ (\Gamma f = \Gamma \land d(f)<\aleph_0) ~\lor~ c(f)=\aleph_0} 
    \\
    \inv{Q_{\Gamma}} &= \set{f \in I_X \given \Gamma \nsubseteq \im{f} ~\lor~ (\Gamma f = \Gamma \land c(f)<\aleph_0) ~\lor~ d(f)=\aleph_0}
\end{align*}
where $\Gamma$ is a non-empty subset of $\Sigma$.
\end{theorem}

\begin{proof}
    We will prove this by showing that the above semigroups, together with the semigroups described in Theorem \ref{thm:maximal_sym}, satisfy the conditions of Theorem \ref{thm:total->all} with $G = \pointStab{\Sym{X}}{\Sigma}$.
    \begin{enumerate}[\normalfont(i)]
        \item $\pointStab{\Sym{X}}{\Sigma} \subseteq \set{f \in I_X \given \Gamma f = \Gamma ~\land~ c(f)=d(f)=0}$ whenever $\Gamma \subseteq \Sigma$.
        In fact, the intersection of any of the semigroups $P_{\Gamma}$, $\inv{P_{\Gamma}}$, $Q_{\Gamma}$, or $\inv{Q_{\Gamma}}$ with $\Sym{X}$ is exactly $\setStab{\Sym{X}}{\Gamma}$, which means that none of them contain all of $\Sym{X}$.
        \item We must show that $P_{\Gamma}$, $\inv{P_{\Gamma}}$, $Q_{\Gamma}$, and $\inv{Q_{\Gamma}}$ are indeed semigroups for all non-empty $\Gamma \subseteq \Sigma$.
        We define the following sets:
        \begin{align*}
            I &= \set{f \in I_X \given \Gamma \nsubseteq \dom{f}} 
            \\
            J &= \set{f \in I_X \given \Gamma f = \Gamma}
            \\
            K &= \set{f \in I_X \given \Gamma f = \Gamma ~\land~ d(f) < \aleph_0}
            \\
            L &= \set{f \in I_X \given c(f) = \aleph_0}
        \end{align*}
        Then $I$ and $L$ are right ideals, $J$ and $K$ are semigroups, $JI, KI \subseteq I$, $KL \subseteq L$ (by Lemma \ref{lem:defect_properties} \eqref{lem:defect_properties/defect<collapse}), and $\fin_X$ is an ideal (Lemma \ref{lem:small_charts}).
        Thus $P_{\Gamma} = I \cup J \cup \fin_X$ and $Q_{\Gamma} = I \cup K \cup L$ are semigroups. 
        $\inv{P_{\Gamma}}$ and $\inv{Q_{\Gamma}}$ being semigroups will follow from them being the inverses of $P_{\Gamma}$ and $Q_{\Gamma}$, and all being proper subsemigroups of $I_X$ will follow from none of them being contained in any of the other.
        \item We must show that, for all non-empty $\Gamma \subseteq \Sigma$, none of the semigroups $P_{\Gamma}$, $\inv{P_{\Gamma}}$, $Q_{\Gamma}$ or $\inv{Q_{\Gamma}}$ are contained in any of the other, nor in any of the semigroups described in Theorem \ref{thm:maximal_sym}.
        So let $\Gamma$ and $\Omega$ be arbitrary non-empty subsets of $\Sigma$.
        To show that $P_{\Gamma}$ is not contained in the other semigroups mentioned above, choose $f,g,h,i,j,k,l,m \in I_X \setminus \fin_X$ such that:
        \begin{itemize}
            \item $\Gamma f \subseteq \Gamma$, $\Omega \subseteq \dom{f}$, and $\Omega f \neq \Omega$ (when $\Omega \neq \Gamma$);
            \item $\Gamma g \subseteq \Gamma$, $\Omega \subseteq \im{g}$, and $\Omega g \neq \Omega$;
            \item $\Gamma h \subseteq \Gamma$, $\Omega \subseteq \dom{h}$, $d(h) = \aleph_0$, and $c(h) < \aleph_0$;
            \item $\Gamma i \subseteq \Gamma$, $\Omega \subseteq \im{i}$, $c(i) = \aleph_0$ and $d(i) < \aleph_0$;
            \item $\Gamma j \subseteq \Gamma$, $c(j) = 0$, and $d(j) > 0$;
            \item $\Gamma k \subseteq \Gamma$, $c(k) > 0$ and $d(k) = 0$;
            \item $\Gamma l \subseteq \Gamma$, $c(l) < \aleph_0$, and $d(l) = \aleph_0$;
            \item $\Gamma m \subseteq \Gamma$, $c(m) = \aleph_0$ and $d(m) < \aleph_0$. 
        \end{itemize}
        Then $f,g,h,i,j,k,l,m \in P_{\Gamma}$ but $f \notin P_{\Omega}$ (when $\Omega \neq \Gamma)$, $g \notin \inv{P_{\Omega}}$, $h \notin Q_{\Omega}$, $i \notin \inv{Q_{\Omega}}$, $j \notin S$, $k \notin \inv{S}$, $l \notin T$, and $m \notin \inv{T}$.
        To show that $Q_{\Gamma}$ is not contained in the other semigroups mentioned above, choose $f,g,h,i,j,k,l,m \in I_X \setminus \fin_X$ such that:
        \begin{itemize}
            \item $\Gamma f \subseteq \Gamma$, $d(f) < \aleph_0$, $\Omega \subseteq \dom{f}$, $\Omega f \neq \Omega$ (when $\Omega \neq \Gamma$), and $c(f) < \aleph_0$;
            \item $c(g) = \aleph_0$, $\Omega \subseteq \im{g}$, and $d(g) < \aleph_0$;
            \item $c(h) = \aleph_0$, $\Omega \subseteq \dom{h}$, and $\Omega h \neq \Omega$;
            \item $c(i) = \aleph_0$, $\Omega \subseteq \im{i}$, and $\Omega i \neq \Omega$;
            \item $\Gamma j = \Gamma$, $0 < d(j) < \aleph_0$ and $c(j) = 0$;
            \item $c(k) = \aleph_0$ and $d(k) = 0$;
            \item $\Gamma \nsubseteq \dom{l}$, $c(l) < \aleph_0$, and $d(l) = \aleph_0$;
            \item $c(m) = \aleph_0$ and $d(m) < \aleph_0$.
        \end{itemize}
        Then $f,g,h,i,j,k,l,m \in Q_{\Gamma}$ but $f \notin Q_{\Omega}$ (when $\Omega \neq \Gamma)$, $g \notin \inv{Q_{\Omega}}$, $h \notin P_{\Omega}$, $i \notin \inv{P_{\Omega}}$, $j \notin S$, $k \notin \inv{S}$, $l \notin T$, and $m \notin \inv{T}$.         
        The rest will follow from $\inv{P_{\Gamma}}$ and $\inv{Q_{\Gamma}}$ being the inverses of $P_{\Gamma}$ and $Q_{\Gamma}$.
        
        \item $\inv{P_{\Gamma}}$ and $\inv{Q_{\Gamma}}$ are indeed the inverses of $P_{\Gamma}$ and $Q_{\Gamma}$ respectively, since $c(f) = d(\inv{f})$, $\im{f} = \dom{\inv{f}}$, and if $\Gamma$ is finite then $\Gamma f = \Gamma \iff \Gamma \inv{f} = \Gamma$.
        
        \item Let $Y = X \setminus \Sigma$, then $\Sym{Y} \subseteq \eval{\pointStab{\Sym{X}}{\Sigma}}_{Y}$.
        
        \item Let $U$ be a subsemigroup of $I_X$ containing $\pointStab{\Sym{X}}{\Sigma}$, but which is itself not contained in $S$, $T$, $P_{\Gamma}$, or $Q_{\Gamma}$ for any non-empty $\Gamma \subseteq \Sigma$.
        Then there exists charts $\Bar{p}_{\Gamma},\Bar{q}_{\Gamma},\Bar{s},\Bar{t} \in U$ such that $\Bar{p}_{\Gamma} \notin P_{\Gamma}$, $\Bar{q}_{\Gamma} \notin Q_{\Gamma}$, $\Bar{s} \notin S$, and $\Bar{t} \notin T$.
        That is, $\Bar{s}$ and $\Bar{t}$ are the charts described in \eqref{notin_S_T} and:
        \begin{nalign} \label{notin_F_H}
            \Gamma \subseteq \dom{\Bar{p}_{\Gamma}} 
            ~&\land~ \Gamma \Bar{p}_{\Gamma} \neq \Gamma \\
            \Gamma \subseteq \dom{\Bar{q}_{\Gamma}} 
            ~&\land~ (\Gamma \Bar{q}_{\Gamma} \neq \Gamma 
            \lor d(\Bar{q}_{\Gamma}) = \aleph_0)
            ~\land~ c(\Bar{q}_{\Gamma}) < \aleph_0 
        \end{nalign}
        First we will show that there exists a total chart $f \in U$ such that $\im{f} \cap \Sigma = \emptyset$.
        Since $\Sigma$ is finite (and $U$ contains the total chart $\Bar{s}$), it suffices to show that for every total $f_0 \in U$ with $\im{f_0} \cap \Sigma \neq \emptyset$ there exists a total $f_1 \in U$ with $\im{f_1} \cap \Sigma \subsetneq \im{f_0} \cap \Sigma$.
        We shall denote $\im{f_0} \cap \Sigma$ by $\Gamma$.
        As an intermediate step we show that there exists a total $f_2 \in U$ such that $\im{f_2} \cap \Sigma \subseteq \Gamma$ and $\Gamma f_2 \neq \Gamma$.
        If $\Gamma f_0 \neq \Gamma$ then let $f_2 = f_0$, otherwise if $\Gamma \Bar{s} \neq \Gamma$ then we can set $f_2 = \Bar{s} f_0$.
        Thus the final case we need to consider is when $\Gamma f_0 = \Gamma \Bar{s} = \Gamma$.
        Since $c(\Bar{q}_{\Gamma}) < \aleph_0$ and $\Gamma \subseteq \dom{\Bar{q}_{\Gamma}}$, there exist $n \in \omega$ and $a \in \pointStab{\Sym{X}}{\Sigma}$ such that $\Bar{s}^n f_0 a \Bar{q}_{\Gamma}$ is total. 
        So if $\Gamma \Bar{q}_{\Gamma} \neq \Gamma$, then  we can set $f_2 = \Bar{s}^n f_0 a \Bar{q}_{\Gamma} f_0$ (see Figure \ref{fig:maximal_pointwise/qGamma_out}).
        \begin{figure}
            \centering
            \includegraphics[height=7cm]{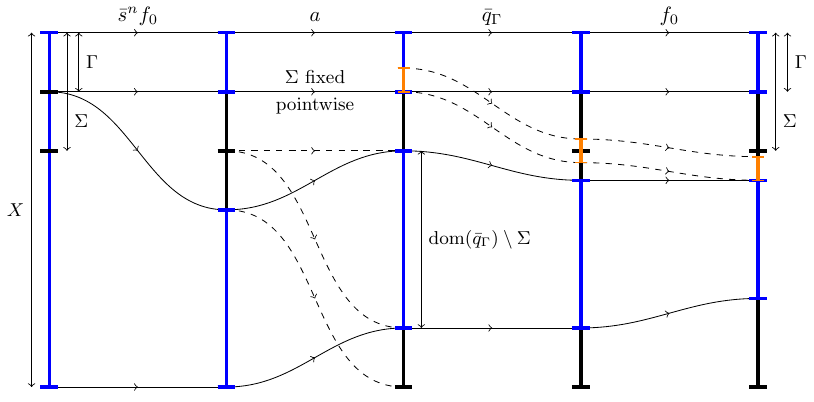}
            \caption{The composite $\Bar{s}^n f_0 a \Bar{q}_{\Gamma} f_0$ when $\Gamma \Bar{q}_{\Gamma} \neq \Gamma$ in the proof of Theorem \ref{thm:maximal_pointwise}.}
            \label{fig:maximal_pointwise/qGamma_out}
        \end{figure}
        Otherwise we know that $d(\Bar{q}_{\Gamma}) = \aleph_0$. Thus, since $\Gamma \subseteq \dom{\Bar{p}_{\Gamma}}$ and $\Bar{p}_{\Gamma} \notin \fin_X$, we can find $b \in \pointStab{\Sym{X}}{\Sigma}$ such that $\Bar{s}^n f_0 a \Bar{q}_{\Gamma} b \Bar{p}_{\Gamma}$ is total.
        Hence, since $\Gamma \Bar{p}_{\Gamma} \neq \Gamma$, we can set $f_2 = \Bar{s}^n f_0 a \Bar{q}_{\Gamma} b \Bar{p}_{\Gamma} f_0$ (see Figure \ref{fig:maximal_pointwise/qGamma_defect}).
        \begin{figure} 
            \centering
            \includegraphics[height=7cm]{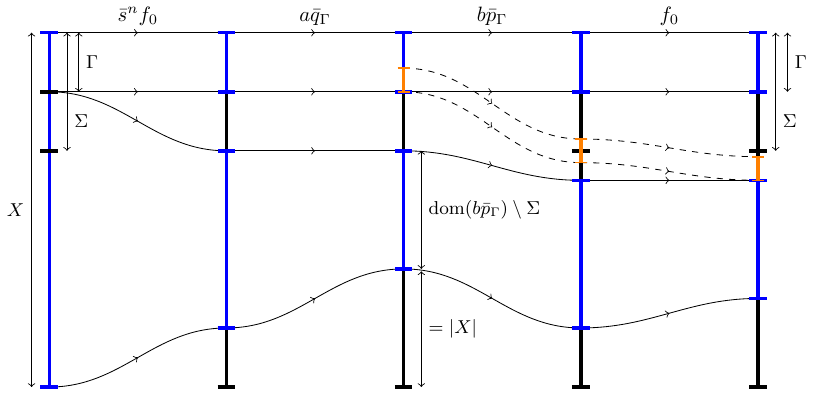}
            \caption{The composite $\Bar{s}^n f_0 a \Bar{q}_{\Gamma} b \Bar{p}_{\Gamma} f_0$ when $d(\Bar{q}_{\Gamma}) = \aleph_0$ in the proof of Theorem \ref{thm:maximal_pointwise}.}
            \label{fig:maximal_pointwise/qGamma_defect}
        \end{figure}
        
        We have shown that $U$ contains a total chart $f_2$ satisfying $\im{f_2} \cap \Sigma \subseteq \Gamma$ and $\Gamma f_2 \neq \Gamma$. We now use this to prove that there exists the promised total chart $f_1 \in U$ with $\im{f_1} \cap \Sigma \subsetneq \Gamma$.
        If $\im{f_2} \cap \Sigma \neq \Gamma$, then let $f_1 = f_2$.
        Otherwise, since $\Gamma f_2 \neq \Gamma$, it follows that $\Gamma \inv{f_2} \nsubseteq \Gamma = \im{f_2} \cap \Sigma$.
        From this there are two cases to consider: $\Gamma \inv{f_2} \nsubseteq \im{f_2}$ or $\Gamma \inv{f_2} \nsubseteq \Sigma$.
        If $\Gamma \inv{f_2} \nsubseteq \im{f_2}$, then $\im{{f_2}^2} \cap \Sigma \subsetneq \Gamma$, and we set $f_1 = {f_2}^2$.
        So suppose $\Gamma \inv{f_2} \nsubseteq \Sigma$.
        Let $x \in \Gamma \inv{f_2} \setminus \Sigma$ and $m \in \omega$ such that $d(\Bar{s}^m) > \card{\Sigma}$.
        Then there exists $y \in X \setminus \Sigma$ such that $y \notin \im{\Bar{s}^m}$ and there is $c \in \pointStab{\Sym{X}}{\Sigma}$ such that $(y)c = x$.
        So in this case we set $f_1 = \Bar{s}^m c f_2$.
        Then, since $y \notin \im{\Bar{s}^m}$ and $c f_2$ is total, it follows that $(x)f_2 = (y)c f_2 \notin \im{\Bar{s}^m c f_2} = \im{f_1}$. But $(x)f_2 \in \Gamma$, and so $\im{f_1} \cap \Sigma \subseteq \Gamma \setminus \set{(x)f_2} \subsetneq \Gamma$, as required (see Figure \ref{fig:maximal_pointwise/f1}).
        \begin{figure}
            \centering
            \includegraphics[height=7cm]{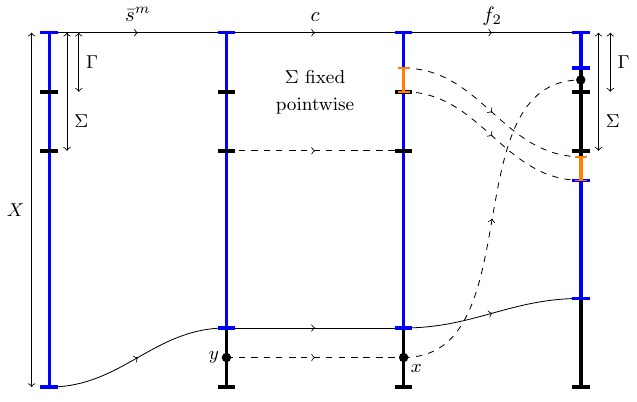}
            \caption{The composite $\Bar{s}^m c f_2$ when $\Gamma \inv{f_2} \nsubseteq \Sigma$ in the proof of Theorem \ref{thm:maximal_pointwise}.}
            \label{fig:maximal_pointwise/f1}
        \end{figure}
        
        It follows from the above that there exists a total chart $f \in U$ such that $\im{f} \cap \Sigma = \emptyset$.
        Lastly, we need to show that there exists a total chart $f' \in U$ whose image is a moiety of $X\setminus \Sigma$. 
        Since $d(f) \geq \card{\Sigma}$ and $c(\Bar{t}) < \aleph_0$, there exists $l \in \omega$ and $d \in \pointStab{\Sym{X}}{\Sigma}$ such that $\im{f^l d} \subseteq (X \setminus \Sigma) \inv{\Bar{t}}$.
        Let $f' = f^l d \Bar{t}$ (see Figure \ref{fig:maximal_pointwise/fprime}).
        \begin{figure}
            \centering
            \includegraphics[height=7cm]{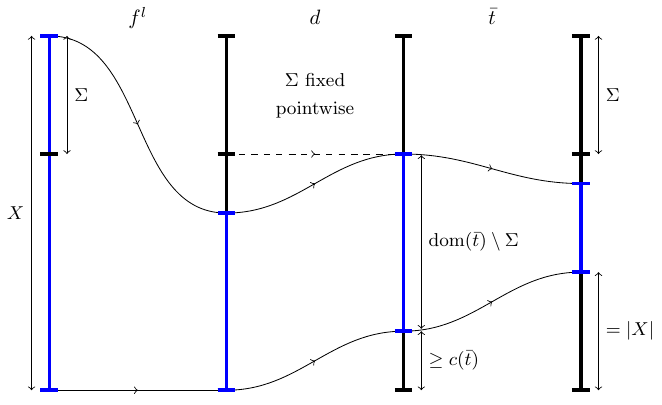}
            \caption{The chart $f' = f^l d \Bar{t}$ in the proof of Theorem \ref{thm:maximal_pointwise}.}
            \label{fig:maximal_pointwise/fprime}
        \end{figure}
    \end{enumerate}
\end{proof}

\begin{corollary} \label{cor:maximal_inv_pointwise}
Let $X$ be a countably infinite set and $\Sigma$ a non-empty finite subset of $X$.
Then the maximal inverse subsemigroups of $I_X$ which contain $\pointStab{\Sym{X}}{\Sigma}$ but not $\Sym{X}$ are:
\begin{align*}
    \mathcal{P}_\Gamma = P_{\Gamma} \cap \inv{P_{\Gamma}} = \lbrace\, f \in I_X \given\;& 
    (\Gamma \nsubseteq \dom{f} \land \Gamma \nsubseteq \im{f})
    ~\lor~ \Gamma f = \Gamma 
    \,\rbrace \cup \fin_X 
    \\
    \mathcal{Q}_\Gamma = Q_{\Gamma} \cap \inv{Q_{\Gamma}} = \lbrace\, f \in I_X \given\;&
    (\Gamma \nsubseteq \dom{f} \land \Gamma \nsubseteq \im{f})
    ~\lor~ (\Gamma \nsubseteq \dom{f} \land d(f)=\aleph_0) \\
    & \lor~ (\Gamma \nsubseteq \im{f} \land c(f)=\aleph_0) 
    ~\lor~ (\Gamma f = \Gamma \land c(f)+d(f) < \aleph_0) \\
    & \lor~ c(f)=d(f)=\aleph_0  
    \,\rbrace
\end{align*}
where $\Gamma$ is a non-empty subset of $\Sigma$.
\end{corollary}

\begin{proof}
    This follows from Proposition \ref{prop:inverse_intersection}, which states that the maximal inverse subsemigroups containing some inverse subsemigroup $G$ are given by the intersections of the maximal subsemigroups that contain $G$ with their inverses.
    So we let $G = \pointStab{\Sym{X}}{\Sigma}$ and use the semigroups described in Theorem \ref{thm:maximal_pointwise} as our set of maximal subsemigroups.
    All that remains is to show that for all non-empty subsets $\Gamma \subseteq \Sigma$ the semigroups $\mathcal{P}_\Gamma$ and $\mathcal{Q}_\Gamma$ are not contained in each other nor in the inverse semigroups described in Corollary \ref{cor:maximal_inv_sym}.
    So let $\Gamma$ and $\Omega$ be arbitrary non-empty subsets of $\Sigma$.
    The show that $\mathcal{P}_\Gamma$ is not contained in any of the other inverse semigroups mention above, choose $f,g,h,i \in I_X \setminus \fin_X$ such that:
    \begin{itemize}
        \item $\Gamma f \cup \Gamma \inv{f} \subseteq \Gamma$, 
        $\Omega \subseteq \dom{f} \cap \im{f}$,   
        and $\Omega f \neq \Omega$ (when $\Omega \neq \Gamma$);
        \item $\Gamma g = \Gamma$, 
        $\Omega \subseteq \dom{g}$, 
        $c(g) < \aleph_0$, 
        and $d(g) = \aleph_0$;
        \item $\Gamma h = \Gamma$,
        $c(h) = 0$,
        and $d(h) > 0$;
        \item $\Gamma i = \Gamma$,
        $c(i) < \aleph_0$,
        and $d(i) = \aleph_0$.
    \end{itemize}
    Then $f,g,h,i \in \mathcal{P}_{\Gamma}$ but $f \notin \mathcal{P}_{\Omega}$ (when $\Omega \neq \Gamma)$, $g \notin \mathcal{Q}_{\Omega}$, $h \notin \mathcal{S}$, and $i \notin \mathcal{T}$.
    The show that $\mathcal{Q}_\Gamma$ is not contained in any of the other inverse semigroups mentioned above, choose $f,g,h,i \in I_X \setminus \fin_X$ such that:
    \begin{itemize}
        \item $\Gamma f \cup \Gamma \inv{f} \subseteq \Gamma$, 
        $c(f) + d(f) < \aleph_0$, 
        $\Omega \subseteq \dom{f} \cap \im{f}$, 
        and $\Omega f \neq \Omega$ (when $\Omega \neq \Gamma$);
        \item $c(g) = d(g) = \aleph_0$, 
        $\Omega \subseteq \dom{g}$, 
        $\Omega g \neq \Omega$;
        \item $\Gamma h = \Gamma$, 
        $c(h) = 0$, 
        and $d(h) < \aleph_0$;
        \item $\Gamma \nsubseteq \dom{i}$, 
        $d(i) = \aleph_0$, 
        and $c(i) < \aleph_0$.
    \end{itemize}
    Then $f,g,h,i \in \mathcal{Q}_{\Gamma}$ but $f \notin \mathcal{Q}_{\Omega}$ (when $\Omega \neq \Gamma)$, $g \notin \mathcal{P}_{\Omega}$, $h \notin \mathcal{S}$, and $i \notin \mathcal{T}$.
\end{proof}

\section{Stabiliser of an Ultrafilter}

In this section we classify the maximal (inverse) subsemigroups of $I_X$ containing the stabiliser of an ultrafilter on $X$, when $X$ is a countably infinite set.
First, let us establish what a filter is.

\begin{definition}[Filter]
    Let $X$ be a set. A non-empty family of subsets $\filter$ of $X$ is called a \emph{filter} if it satisfies the following conditions:
    \begin{enumerate}[\normalfont(i)]
        \item $\emptyset \notin \filter$;
        \item If $\Sigma \in \filter$ and $\Sigma \subseteq \Gamma \subseteq X$, then $\Gamma \in \filter$;
        \item If $\Sigma, \Gamma \in \filter$, then $\Sigma \cap \Gamma \in \filter$.
    \end{enumerate}
\end{definition}

Some filters have properties that make them of particular interest.
We list a few such special types of filters here.

\begin{definition}[Ultrafilter]
    Let $\filter$ be a filter on a set $X$.
    Then $\filter$ is called an \emph{ultrafilter} if it is maximal with respect to containment among filters on $X$.
\end{definition}

\begin{definition}[Principal Ultrafilter]
    An ultrafilter $\filter$ on a set $X$ is called \emph{principal} if there exists an element $x \in X$ such that $\filter = \set{\Sigma \subseteq X \given x \in \Sigma}$.
\end{definition}

\begin{definition}[Uniform Ultrafilter]
    An ultrafilter $\filter$ on a set $X$ is called \emph{uniform} if $\card{\Sigma} = \card{X}$ for all $\Sigma \in \filter$.
\end{definition}

While the above definitions of different types of filters are conceptually good and well motivated, they do leave something to be desired in terms of practical applicability.
E.g. how can one tell whether a filter is maximal or not?
The following well-known result answers this question and gives us an alternative definition of ultrafilters (we include a proof of the statement here for convenience).

\begin{proposition} \label{prop:ultrafilter_definition}
    Let $\filter$ be a filter on a set $X$.
    Then $\filter$ is an ultrafilter if and only if for all $\Sigma \subseteq X$ either $\Sigma \in \filter$ or $X \setminus \Sigma \in \filter$.
\end{proposition}

\begin{proof}
    First we show that any filter $\filter$ satisfying the above condition is indeed maximal with respect to containment among filters on $X$.
    Taking any $\Sigma \subseteq X$ such that $\Sigma \notin \filter$, it follows that $X \setminus \Sigma \in \filter$.
    But then $\Sigma \cap (X \setminus \Sigma) = \emptyset$, hence $\filter \cup \set{\Sigma}$ does not generate a filter on $X$ (as filters must be closed under taking intersections, but filters also cannot contain the empty set).
    Next we show that any filter $\filter$, which does not satisfy the above condition, cannot be an ultrafilter on $X$.
    Let $\Sigma$ be any subset of $X$ such that both $\Sigma \notin \filter$ and $X \setminus \Sigma \notin \filter$ (we showed above why we cannot have both in the filter).
    Then neither are any subsets of $\Sigma$ nor $X \setminus \Sigma$ in $\filter$ (as this would otherwise contradict $\filter$ being closed under taking supersets).
    It thus follows that there are no elements in $\filter$, which have empty intersection with any superset of $\Sigma$.
    Hence, $\filter \cup \set{\Sigma}$ generates a filter on $X$, which strictly contains $\filter$.
\end{proof}

Henceforth we will be using Proposition \ref{prop:ultrafilter_definition} as our effective definition of an ultrafilter.
We will also give an alternative definition for principal ultrafilters, but first we go over some intermediate results, starting with the dual concept of a filter.

\begin{definition}[Set-Theoretic Ideal] \label{def:ideal}
    Let $X$ be a set. A non-empty family of subsets $\ideal$ of $X$ is called an \emph{ideal} if it satisfies the following conditions:
    \begin{enumerate}[\normalfont(i)]
        \item $\emptyset \in \ideal$; 
        \label{def:ideal_1}
        \item If $\Sigma \in \ideal$ and $\Gamma \subseteq \Sigma$, then $\Gamma \in \ideal$; 
        \label{def:ideal_2}
        \item If $\Sigma, \Gamma \in \ideal$, then $\Sigma \cup \Gamma \in \ideal$. 
        \label{def:ideal_3}
    \end{enumerate}
\end{definition}

\begin{proposition} \label{prop:ideal_of_filter}
    Let $\filter$ be a filter on a set $X$.
    Then $\ideal = \set{\Sigma \subseteq X \given X \setminus \Sigma \in \filter}$ is an ideal on $X$, and we will refer to it as \emph{the ideal of $\filter$}. Moreover, $\filter \cap \ideal = \emptyset$.
\end{proposition}

\begin{proof}
    It is immediate that $\ideal$ satisfies conditions \eqref{def:ideal_1} and \eqref{def:ideal_2} of Definition \ref{def:ideal}.
    To see that condition \eqref{def:ideal_3} holds, consider two elements $\Sigma, \Gamma \in \ideal$.
    It follows from definition that $X \setminus \Sigma, X \setminus \Gamma \in \filter$.
    Being in the filter is closed under finite intersections, so $(X \setminus \Sigma) \cap (X \setminus \Gamma) = X \setminus (\Sigma \cup \Gamma) \in \filter$.
    This in turn implies that $\Sigma \cup \Gamma \in \ideal$ as desired.

    Lastly, since $\Sigma \cap (X \setminus \Sigma) = \emptyset$, we get that if $\Sigma \in \ideal$, then $\Sigma \notin \filter$.
\end{proof}

\begin{corollary} \label{cor:ultra_ideal}
    If $\filter$ is an ultrafilter on a set $X$, then the ideal of $\filter$ consists of exactly those subsets of $X$ which are not in $\filter$.
\end{corollary}

\begin{lemma} \label{lem:ultra_subsets}
    Let $\filter$ be an ultrafilter on a set $X$ and $\ideal$ the ideal of $\filter$.
    If $\Sigma \in \filter$, then for every subset $\Gamma \subseteq \Sigma$ precisely one of $\Gamma$ and $\Sigma \setminus \Gamma$ is in $\filter$ and the other is in $\ideal$.
\end{lemma}

\begin{proof}
    Since $\Gamma \cup (\Sigma \setminus \Gamma) = \Sigma$, we cannot have both in $\ideal$ by Proposition \ref{prop:ideal_of_filter}.
    Since $\Gamma \cap (\Sigma \setminus \Gamma) = \emptyset$, we cannot have both in $\filter$ as $\emptyset \notin \filter$.
    Hence, one is in  $\filter$ and the other is in $\ideal$.
\end{proof}

Lemma \ref{lem:ultra_subsets} allows for the following alternative definition of a principal ultrafilter.

\begin{corollary} \label{cor:principal_alt}
    Let $\filter$ be an ultrafilter on a set $X$. 
    Then $\filter$ is principal if and only if there exists $\Sigma \in \filter$ such that $\Sigma$ is finite.
\end{corollary}

It is worth noting that when $X$ is a countably infinite set, then Corollary \ref{cor:principal_alt} implies that the uniform and non-principal ultrafilters on $X$ coincide.
We have covered the basic notions of what a filter is, and are now ready to define the subgroup mentioned in the introduction of this section.

\begin{definition}[Stabiliser of a Filter]
    The \emph{stabiliser of a filter} $\filter$ on a set $X$, denoted $\setStab{\Sym{X}}{\filter}$, is the group of permutations which preserve whether any given subset of $X$ is in the filter or not.
    \begin{equation*}
        \setStab{\Sym{X}}{\filter} = \set{f \in \Sym{X} \given (\forall \Sigma \subseteq X)~ \Sigma f \in \filter \iff \Sigma \in \filter}
    \end{equation*}
\end{definition}

It is shown in \cite[Theorem 6.4]{subgroups_macpherson_neumann} and \cite{richman1967maximal} that when $\filter$ is an ultrafilter, then $\setStab{\Sym{X}}{\filter}$ is a maximal subgroup of $\Sym{X}$.
Furthermore, \cite[Theorem 6.4]{subgroups_macpherson_neumann} also proves that the following useful alternative definition for the stabiliser of $\filter$ holds.

\begin{proposition}[{\cite[Theorem 6.4]{subgroups_macpherson_neumann}}] \label{prop:ultrafilter_alt}
    Let $\filter$ be an ultrafilter on a set $X$. 
    Then the stabiliser of $\filter$ is equal to the union of the pointwise stabilisers of the elements of $\filter$.
    \begin{equation} \label{eq:ultrafilter_alt}
        \setStab{\Sym{X}}{\filter} = \bigcup_{\Sigma \in \filter} \pointStab{\Sym{X}}{\Sigma}
    \end{equation}
\end{proposition}


The statement in equation \eqref{eq:ultrafilter_alt} above is written slightly different from the actual statement given in \cite[Theorem 6.4]{subgroups_macpherson_neumann}, which is that $\setStab{\Sym{X}}{\filter} = \pointStab{\Sym{X}}{\filter} := \set{f \in \Sym{X} \given \fix{f} \in \filter}$ (where $\fix{f} = \set{x \in X \given xf = x}$), but these two statements are equivalent since supersets of elements in a filter are also in the filter.
The following lemmas will prove useful for proving the main theorem of this section.

\begin{lemma} \label{lem:filter_map_in}
    Let $\filter$ be an ultrafilter on a set $X$, $\ideal$ the ideal of $\filter$, and $\Sigma, \Gamma \subseteq X$ subsets such that $\card{\Sigma} \leq \card{\Gamma}$ and $\card{X \setminus \Sigma} \geq \card{X \setminus \Gamma}$.
    If $\Sigma, \Gamma \in \filter$ or $\Sigma, \Gamma \in \ideal$, then there exists $f \in \setStab{\Sym{X}}{\filter}$ such that $\Sigma f \subseteq \Gamma$.
\end{lemma}

\begin{proof}
    We first consider the case where $\Sigma, \Gamma \in \filter$.
    It follows that the intersection $\Sigma \cap \Gamma$ is also in $\filter$.
    Thus, by Proposition \ref{prop:ultrafilter_alt}, we get that $\pointStab{\Sym{X}}{\Sigma \cap \Gamma} \subseteq \setStab{\Sym{X}}{\filter}$.
    Hence, since $\card{\Sigma} \leq \card{\Gamma}$ and $\card{X \setminus \Sigma} \geq \card{X \setminus \Gamma}$, we can pick a permutation $f \in \pointStab{\Sym{X}}{\Sigma \cap \Gamma}$ such that $\Sigma f \subseteq \Gamma$, as required.

    The case when $\Sigma, \Gamma \in \ideal$ follows directly from the above, since all elements in $\ideal$ are simply complements of elements in $\filter$.
\end{proof}

\begin{lemma} \label{lem:filter_transitive}
    Let $\filter$ be an ultrafilter on a set $X$, $\ideal$ the ideal of $\filter$, and $\Sigma, \Gamma \subseteq X$ subsets such that $\card{\Sigma} = \card{\Gamma}$ and $\card{X \setminus \Sigma} = \card{X \setminus \Gamma}$.
    If $\Sigma, \Gamma \in \filter$ or $\Sigma, \Gamma \in \ideal$, then there exists a permutation $f \in \setStab{\Sym{X}}{\filter}$ such that $\Sigma f = \Gamma$.
\end{lemma}

\begin{proof}
    We first consider the case where $\Sigma, \Gamma \in \filter$, from which it follows that $\Sigma \cap \Gamma \in \filter$.
    We now choose a subset $\Omega$ of $\Sigma \cap \Gamma$. 
    If $\Sigma \cap \Gamma$ is finite, simply let $\Omega=\Sigma\cap \Gamma$. 
    If $\Sigma \cap \Gamma$ is infinite, then, by Lemma \ref{lem:ultra_subsets}, there exists a moiety $\Omega$ of $\Sigma \cap \Gamma$ such that $\Omega \in \filter$. Note that in either case, 
    \begin{equation}\label{eq:perm_conditions}
        \Omega \in \filter \text{ and }\card{\Sigma \setminus \Omega} = \card{\Gamma \setminus \Omega}.
    \end{equation}    
    Since $\Omega \in \filter$, Proposition \ref{prop:ultrafilter_alt} tells us that $\pointStab{\Sym{X}}{\Omega} \subseteq \setStab{\Sym{X}}{\filter}$.
    Thus, by \eqref{eq:perm_conditions} and since $\card{\Sigma} = \card{\Gamma}$ and $\card{X \setminus \Sigma} = \card{X \setminus \Gamma}$, there exist a permutation $f \in \pointStab{\Sym{X}}{\Omega}$ such that $\Sigma f = \Gamma$, as required.

    As before, the case when $\Sigma, \Gamma \in \ideal$ follows directly from the above, since all elements in $\ideal$ are complements of elements in $\filter$.
\end{proof}

We are finally ready to prove the main result of this section.

\begin{theorem} \label{thm:maximal_ultrafilter}
Let $X$ be a countably infinite set, $\filter$ a non-principal ultrafilter on $X$, and $\ideal$ the ideal of $\filter$.
Then the maximal subsemigroups of $I_X$ which contain $\setStab{\Sym{X}}{\filter}$ but not $\Sym{X}$ are:
\begin{align*}
    V_{\filter} &= \set{f \in I_X \given \dom{f} \in \ideal ~\lor~ (\forall \Sigma \subseteq X)~ \Sigma f \in \filter \iff \Sigma \in \filter} \\
    \inv{V_{\filter}} &= \set{f \in I_X \given \im{f} \in \ideal ~\lor~ (\forall \Sigma \subseteq X)~ \Sigma f \in \filter \iff \Sigma \in \filter}
\end{align*}
\end{theorem}

\begin{figure}
    \centering
    \includegraphics[width=\linewidth]{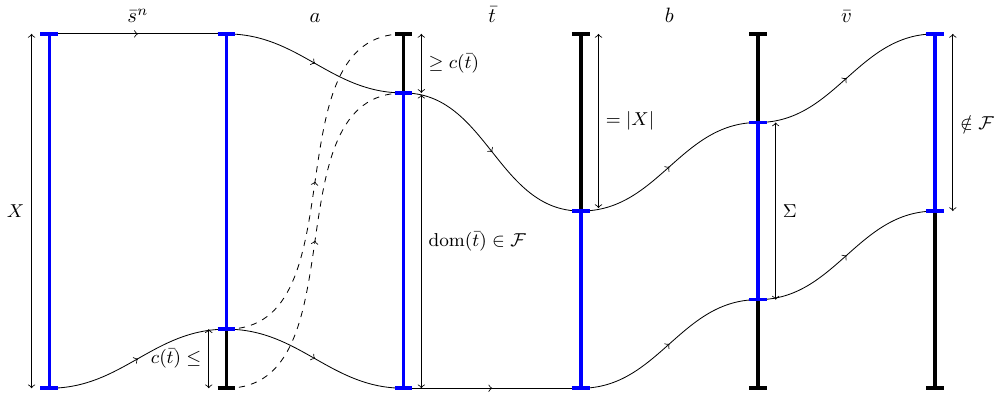}
    \caption{The chart $f = \Bar{s}^n a \Bar{t} b \Bar{v}$ in the proof of Theorem \ref{thm:maximal_ultrafilter}.}
    \label{fig:maximal_filter}
\end{figure}

\begin{proof}
    We will prove this by showing that $V_{\filter}$ and $\inv{V_{\filter}}$, together with the semigroups described in Theorem \ref{thm:maximal_sym}, satisfy the conditions of Theorem \ref{thm:total->all} with $G=\setStab{\Sym{X}}{\filter}$.
    \begin{enumerate}[\normalfont(i)]
        \item $\setStab{\Sym{X}}{\filter} \subseteq \set{f \in I_X \given (\forall \Sigma \subseteq X)~ \Sigma f \in \filter \iff \Sigma \in \filter}$, which is contained in both $V_{\filter}$ and $\inv{V_{\filter}}$.

        \item We must show that $V_{\filter}$ and $\inv{V_{\filter}}$ are proper subsemigroups of $I_X$. 
        We define the following sets:
        \begin{align*}
            I &= \set{f \in I_X \given \dom{f} \in \ideal}
            \\
            J &= \set{f \in I_X \given (\forall \Sigma \subseteq X)~ \Sigma f \in \filter \iff \Sigma \in \filter}
        \end{align*}
        Then $I$ is a right ideal, $J$ is a semigroup, and $JI \subseteq I$.
        Thus $V_{\filter} = I \cup J$ is a semigroup, and $\inv{V_{\filter}}$ being a semigroup will then follow from it being the inverse of $V_{\filter}$.
        Both being proper subsemigroups of $I_X$ will follow from neither being contained in the other.
        
        \item We must show that none of the semigroups $V_{\filter}$ or $\inv{V_{\filter}}$ are contained in the other, nor in any of the semigroups described in Theorem \ref{thm:maximal_sym}.
        To show that $V_{\filter}$ is not contained in any of the other semigroups mentioned above, choose $f,g,h,i,j \in I_X$ such that:
        \begin{itemize}
            \item $\dom{f} \in \ideal$ and $\im{f} \in \filter$;

            \item $(\forall \Sigma \subseteq X)(\Sigma \in \filter \iff \Sigma f \in \filter)$, $c(g) = 0$, and $d(g) > 0$;

            \item $\dom{h} \in \ideal$, $c(h)> 0$, and $d(h) = 0$;

            \item $(\forall \Sigma \subseteq X)(\Sigma \in \filter \iff \Sigma f \in \filter)$, $c(i) < \aleph_0$, and $d(i) = \aleph_0$;

            \item $\dom{j} \in \ideal$, $c(j) = \aleph_0$, and $d(j) < \aleph_0$.
        \end{itemize}
        Then $f,g,h,i,j \in V_{\filter}$ but $f \notin \inv{V_{\filter}}$, $g \notin S$, $h \notin \inv{S}$, $i \notin T$, and $j \notin \inv{T}$.
        $\inv{V_{\filter}}$ not being contained in any of the other semigroups will follow from it being the inverse of $V_{\filter}$.
        
        \item $\inv{V_{\filter}}$ is indeed the inverse of $V_{\filter}$, since $\im{f} = \dom{\inv{f}}$ and $\{\, f \in I_X \given \allowbreak (\forall \Sigma \subseteq X)~ \Sigma f \in \filter \iff \Sigma \in \filter \,\}$ is self inverse.
        
        \item Let $Y$ be a subset of $X$ such that $\card{Y} = \card{X}$ and $Y \in \ideal$. Then $\Sym{Y} \subseteq \eval{\setStab{\Sym{X}}{\filter}}_Y$, since $\pointStab{\Sym{X}}{X \setminus Y} \subseteq \setStab{\Sym{X}}{\filter}$.
        
        \item Let $U$ be a subsemigroup of $I_X$ containing $\setStab{\Sym{X}}{\filter}$, but which is itself not contained in $V_{\filter}$, $S$, or $T$.
        Then there exists charts $\Bar{v},\Bar{s},\Bar{t} \in U$ such that $\Bar{v} \notin V_{\filter}$, $\Bar{s} \notin S$, and $\Bar{t} \notin T$.
        That is, $\Bar{s}$ and $\Bar{t}$ are the charts described in \eqref{notin_S_T} and:
        \begin{equation*}
            \dom{\Bar{v}} \in \filter ~\land~ \paa{(\exists \Sigma \subseteq X)~ (\Sigma \Bar{v} \in \filter ~\land~ \Sigma \in \ideal) ~\lor~ (\Sigma \Bar{v} \in \ideal ~\land~ \Sigma \in \filter)}
        \end{equation*}
        However, if there exists an subset $\Sigma$ in $\ideal$ such that $\Sigma \Bar{v}$ is in $\filter$, then $\Gamma = \dom{\Bar{v}} \setminus \Sigma$ is an element of $\filter$ (by Lemma \ref{lem:ultra_subsets}) and $\Gamma \Bar{v}$ must be in $\ideal$ (since $\Bar{v}$ is injective, $\Gamma \Bar{v}$ must be the relative complement of $\Sigma \Bar{v}$ in $\im{\Bar{v}}$).
        Hence, we can reduce the reduce the definition of $\Bar{v}$ to the following:
        \begin{equation} \label{notin_V}
            \dom{\Bar{v}} \in \filter ~\land~ (\exists \Sigma \in \filter)~ \Sigma \Bar{v} \in \ideal
        \end{equation}
        Our first goal is to show that there exists a total chart $f \in U$ such that $\im{f} \in \ideal$.
        If there exists $m \in \omega$ such that $\im{\Bar{s}^m} \in \ideal$, then simply let $f = \Bar{s}^m$.
        Otherwise, since $\dom{\Bar{t}} \in \filter$ (this follows from $c(\Bar{t}) < \aleph_0$ and the fact that $\filter$ is a non-principal ultrafilter), we can find $n \in \omega$ and, by Lemma \ref{lem:filter_map_in}, a permutation $a \in \setStab{\Sym{X}}{\filter}$ such that $\im{\Bar{s}^n a} \subseteq \dom{\Bar{t}}$.
        If $\im{\Bar{t}} \in \ideal$, then let $f = \Bar{s}^n a \Bar{t}$.
        Otherwise, we use that $\Bar{v}$ maps an element in $\filter$, call it $\Sigma$, to an element in $\ideal$.
        Using Lemma \ref{lem:filter_map_in} again, we can choose a permutation $b \in \setStab{\Sym{X}}{\filter}$ such that $\im{\Bar{s}^n a \Bar{t} b} \subseteq \Sigma$.
        Then we let $f = \Bar{s}^n a \Bar{t} b \Bar{v}$ (see Figure \ref{fig:maximal_filter}).
        We note here that $\im{f}$ is a moiety of $X$.

        Finally, we use Lemma \ref{lem:filter_transitive}, together with the fact that subsets of elements in an ideal are also in the ideal, to pick a permutation $c \in \setStab{\Sym{X}}{\filter}$, such that $f^* = fc$ is a total chart and its image is a moiety of $Y$, as required.
    \end{enumerate}
\end{proof}

\begin{corollary} \label{cor:maximal_inv_ultrafilter}
Let $X$ be a countably infinite set, $\filter$ a non-principal ultrafilter on $X$, and $\ideal$ the ideal of $\filter$.
Then the maximal inverse subsemigroup of $I_X$ which contains $\setStab{\Sym{X}}{\filter}$ but not $\Sym{X}$ is:
\begin{equation*}
    \mathcal{V}_\filter = V_{\filter} \cap \inv{V_{\filter}} = \set{f \in I_X \given \dom{f}, \im{f} \in \ideal ~\lor~ (\forall \Sigma \subseteq X)~ \Sigma f \in \filter \iff \Sigma \in \filter}
\end{equation*}
\end{corollary}

\begin{proof}
    This follows from Proposition \ref{prop:inverse_intersection}, which states that the maximal inverse subsemigroups containing some inverse subsemigroup $G$ are given by the intersections of the maximal subsemigroups that contain $G$ with their inverses.
    So we let $G = \setStab{\Sym{X}}{\filter}$ and use the semigroups described in Theorem \ref{thm:maximal_ultrafilter} as our set of maximal subsemigroups.
    All that remains is to show that $\mathcal{V}_\filter$ is not contained in any of the inverse semigroups described in Corollary \ref{cor:maximal_inv_sym}.
    So choose $f,g \in I_X$ such that:
    \begin{itemize}
        \item $(\forall \Sigma \subseteq X)(\Sigma \in \filter \iff \Sigma f \in \filter)$, $c(f) > 0$, and $d(f) = 0$;

        \item $(\forall \Sigma \subseteq X)(\Sigma \in \filter \iff \Sigma g \in \filter)$, $c(g) < \aleph_0$, and $d(g) = \aleph_0$.
    \end{itemize}
    Then $f,g \in \mathcal{V}_\filter$ but $f \notin \mathcal{S}$ and $g \notin \mathcal{T}$.
\end{proof}

The main results of this section lead us to an interesting conclusion about the subsemigroup structure of $I_X$.
Pospi\u sil's Theorem \cite[Theorem 7.6]{Jech2003} states that in \zfc, the cardinality of the set of uniform ultrafilters on any infinite set $X$ is $2^{2^{\vert\! X \!\vert}}$, which is the same as the cardinality of the powerset of $I_X$.
It therefore follows from Theorem \ref{thm:maximal_ultrafilter} and Corollary \ref{cor:maximal_inv_ultrafilter} that, for $X$ countably infinite, there are $\beth_2 = 2^{2^{\aleph_0}}$ distinct maximal (inverse) subsemigroups of $I_X$, which is as many as there are subsemigroups of $I_X$ in total.

\section{Stabiliser of a Partition}

In this section we classify the maximal (inverse) subsemigroups of $I_X$ containing the stabiliser of a finite partition of $X$, when $X$ is a countably infinite set.
Throughout this section, when we talk about a finite partition of a set $X$, we specifically mean a partition $\partition = \set{\Sigma_0, \dots, \Sigma_{n-1}}$ of $X$ into $n \geq 2$ parts such that for all $i \in n$, $\card{\Sigma_i} = \card{X}$.

\begin{definition}[Stabiliser of a Finite Partition]
    The \emph{stabiliser of a finite partition} $\partition$ of a set $X$, denoted $\Stab{\partition}$, is the group of permutations which maps all parts of the partition to other parts.
    \begin{equation*}
        \Stab{\partition} = \set{f \in \Sym{X} \given (\forall i \in n) (\exists j \in n)~ \Sigma_i f = \Sigma_j}
    \end{equation*}
\end{definition}

Very closely related to the stabiliser of a finite partition is the more general notion of an almost stabiliser.

\begin{definition}[Almost Stabiliser of a Finite Partition]
    The \emph{almost stabiliser of a finite partition} $\partition$ of a set $X$, denoted $\Stab{\partition}$, is the group of permutations defined as follows:
    \begin{equation*}
        \AStab{\partition} = \set{f \in \Sym{X} \given (\forall i \in n) (\exists j \in n)~ \card{\Sigma_i f \setminus \Sigma_j} +\card{\Sigma_j \setminus \Sigma_i f} < \card{X}}
    \end{equation*}
\end{definition}

It is shown in \cite{richman1967maximal} that the almost stabiliser of a finite partition of $X$ is a maximal subgroup of $\Sym{X}$.
Keeping track of where each part of the partition $\partition$ is mapped will be essential throughout this section.
To assist in this endeavour, we define the following binary relation.

\begin{definition}
    Let $\partition$ be a finite partition of an infinite set $X$ and $f \in I_X$.
    We then define the binary relation $\rho_f \subseteq n \times n$ by
    \begin{equation*}
    \rho_f = \set{(i,j) \in n \times n \given \card{\Sigma_i f \cap \Sigma_j} = \card{X}}
    \end{equation*}
\end{definition}

In the following lemma we will show that the binary relation defined above will still `keep track' of where parts of the partition $\partition$ are being mapped under composition of charts, modulo composition with an element from \Stab{\partition}.

\begin{lemma} \label{lem:compose_relations}
    Let $\partition$ be a finite partition of an infinite set $X$ and $f,g \in I_X$. Then there exists $a \in \Stab{\partition}$ such that $\rho_{fag} = \rho_f \rho_g$.
\end{lemma}

\begin{proof}
    Let $i \in n$ be arbitrary.
    If $j \in (i)\inv{\rho_f}$, then $\card{\Sigma_j f \cap \Sigma_i} = \card{X}$, and so $\Sigma_j f \cap \Sigma_i$ can be partitioned into $\card{(i)\rho_g} + 1$ moieties.
    If $k \in (i)\rho_g$, then $\card{\Sigma_k \inv{g} \cap \Sigma_i} = \card{X}$.
    Hence, $\Sigma_k \inv{g} \cap \Sigma_i$ can be partitioned into $\card{(i)\inv{\rho_f}} + 1$ moieties.
    Let $a_i \in \Stab{\partition}$ be any element mapping one of the moieties partitioning $\Sigma_j f \cap \Sigma_i$ to one of the moieties partitioning $\Sigma_k \inv{g} \cap \Sigma_i$ for all $j \in (i)\inv{\rho_f}$ and for all $k \in (i)\rho_g$, while fixing everything else.
    The required $a \in \Stab{\partition}$ is then just the composite $a_0 \dots a_{n-1}$.
\end{proof}

The following are Lemmas 9.3 and 9.4 in \cite{maximal_east_mitchell_peresse}, which will be needed for proving the main theorem of this section (Theorem \ref{thm:maximal_astab}).

\begin{lemma}[\cite{maximal_east_mitchell_peresse}, Lemma 9.3] \label{lem:nxn}
    Let $n$ be a natural number and $\rho, \sigma \subseteq n \times n$ binary relations such that $\rho$ and $\inv{\sigma}$ are total but $\rho, \sigma \notin \Sym{n}$.
    Then the semigroup $\genset{\Sym{n}, \rho, \sigma}$ contains the relation $n \times n$.
\end{lemma}

\begin{lemma}[{\cite[Lemma 9.4]{maximal_east_mitchell_peresse}}] \label{lem:missing_all_part}
    Let $\partition$ be a finite partition of an infinite set $X$ and $f \in I_X$ a total non-surjective chart (that is $c(f) = 0 < d(f)$).
    Then there exists a total $f^* \in \genset{\Stab{\partition}, f}$ such that $\card{\Sigma_i \setminus \im{f^*}} \geq d(f)$ for all $i \in n$. 
    If $d(f)$ is infinite, then $\card{\Sigma_i \setminus \im{f^*}} = d(f)$ for all $i \in n$.
\end{lemma}

We can now prove the main result of this section.

\begin{theorem} \label{thm:maximal_astab}
Let $X$ be a countably infinite set and $\partition = \set{\Sigma_0, \dots, \Sigma_{n-1}}$ a finite partition of $X$ into $n \geq 2$ part such that for all $i \in n$, $\card{\Sigma_i} = \card{X}$.
Then the maximal subsemigroups of $I_X$ which contain $\Stab{\partition}$ but not $\Sym{X}$ are:
\begin{align*}
    A_{\partition} &= \set{f \in I_X \given \rho_f \in \Sym{n} ~\lor~ \dom{\rho_f} \neq n} \\
    \inv{A_{\partition}} &= \set{f \in I_X \given \rho_f \in \Sym{n} ~\lor~ \im{\rho_f} \neq n}
\end{align*}
\end{theorem}

\begin{figure} 
    \centering
    \includegraphics[height=7cm]{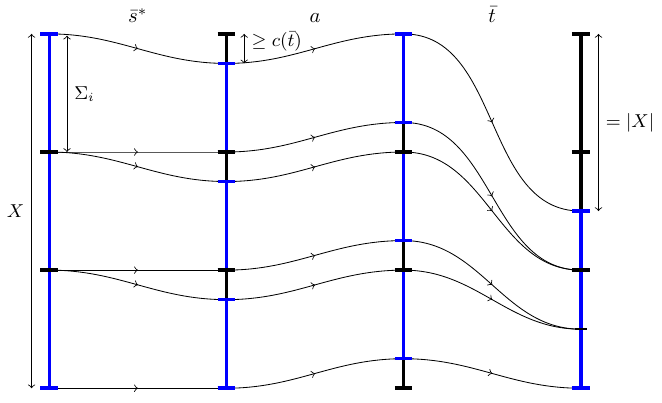}
    \caption{The chart $g = \Bar{s}^* a \Bar{t}$ in the proof of Theorem \ref{thm:maximal_pointwise}.}
    \label{fig:maximal_partition/g}
\end{figure}

\begin{figure} 
    \centering
    \includegraphics[height=7cm]{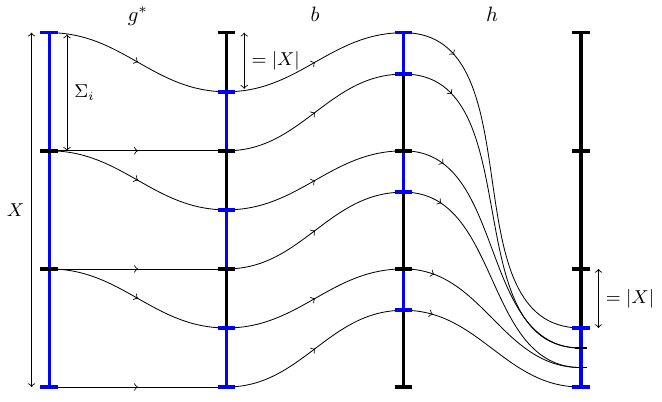}
    \caption{The chart $f = g^* b h$ in the proof of Theorem \ref{thm:maximal_astab}.}
    \label{fig:maximal_partition/f}
\end{figure}

\begin{proof}
    We will prove this by showing that $A_{\partition}$ and $\inv{A_{\partition}}$, together with the semigroups described in Theorem \ref{thm:maximal_sym}, satisfy the conditions of Theorem \ref{thm:total->all} with $G=\Stab{\partition}$.
    \begin{enumerate}[\normalfont(i)]
        \item If $f \in \Stab{\partition}$, then $\rho_f \in \Sym{n}$. Hence, $\Stab{\partition}$ is contained in both $A_{\partition}$ and $\inv{A_{\partition}}$. In fact, the intersection of either $A_{\partition}$ or $\inv{A_{\partition}}$ with $\Sym{X}$ is exactly $\AStab{\partition}$.

        \item We must show that $A_{\partition}$ and $\inv{A_{\partition}}$ are proper subsemigroups of $I_X$.
        We define the following sets:
        \begin{align*}
            I &= \set{f \in I_X \given \rho_f \in \Sym{n}}
            \\
            J &= \set{f \in I_X \given \dom{\rho_f} \neq n}
        \end{align*}
        Then $I$ is a semigroup, $J$ is a right ideal, and $IJ \subseteq J$.
        Thus $A_{\partition} = I \cup J$ is a semigroup and $\inv{A_{\partition}}$ being a semigroup will follow from it being the inverse of $A_{\partition}$.
        Both being proper subsemigroups of $I_X$ will follow from neither being contained in the other.
        
        \item We must show that neither of the semigroups $A_{\partition}$ or $\inv{A_{\partition}}$ are contained in the other, nor in any of the semigroups described in Theorem \ref{thm:maximal_sym}.
        To show that $A_{\partition}$ is not contained in any of the semigroups mentioned above, choose $f,g,h,i,j \in I_X$ such that:
        \begin{itemize}
            \item $\dom{\rho_f} \neq n$ and $\im{\rho_f} = n$;

            \item $\rho_g \in \Sym{n}$, $c(g) = 0$, and $d(g) > 0$;

            \item $\dom{\rho_h} \neq n$, $c(h)>0$, and $d(h)=0$;

            \item $\rho_i \in \Sym{n}$, $c(i) < \aleph_0$, and $d(i) = \aleph_0$, 

            \item $\dom{\rho_j} \neq n$, $c(j) = \aleph_0$, and $d(j) < \aleph_0$.
        \end{itemize}
        Then $f,g,h,i,j \in A_{\partition}$ but $f \notin \inv{A_{\partition}}$, $g \notin S$, $h \notin \inv{S}$, $i \notin T$, and $j \notin \inv{T}$.
        $\inv{A_{\partition}}$ not being contained in any of the other semigroups will follow from it being the inverse of $A_{\partition}$.
        
        \item $\inv{A_{\partition}}$ is indeed the inverse of $A_{\partition}$, since $\inv{\rho_f} = \rho_{\inv{f}}$ (since if $\card{\Sigma_i f \cap \Sigma_j} = \aleph_0$, then $\card{\Sigma_j \inv{f} \cap \Sigma_i} = \aleph_0$).
        
        \item Let $Y$ be an element of $\partition$. Then $\Sym{Y} \subseteq \eval{\Stab{\partition}}_Y$.
        
        \item Let $U$ be a subsemigroup of $I_X$ containing $\Stab{\partition}$, but which is itself not contained $A_{\partition}$, $\inv{A_{\partition}}$, $S$, or $T$.
        Then there exists charts $\Bar{a}, \Bar{b}, \Bar{s}, \Bar{t} \in U$ such that $\Bar{a} \notin A_{\partition}$, $\Bar{b} \notin \inv{A_{\partition}}$, $\Bar{s} \notin S$, and $\Bar{t} \notin T$.
        That is, $\Bar{s}$ and $\Bar{t}$ are the charts described in \eqref{notin_S_T} and:
        \begin{nalign} \label{notin_A}
        \rho_{\Bar{a}} \notin \Sym{n} &~\land~ \dom{\rho_a} = n \\
        \rho_{\Bar{b}} \notin \Sym{n} &~\land~ \im{\rho_b} = n
        \end{nalign}
        As in the proofs of all the previous theorems there exists a natural number $m \in \omega$ such that $\Bar{s}^m$ is a total chart with $d(\Bar{s}^m) \geq c(\Bar{t})$.
        It then follows from Lemma \ref{lem:missing_all_part} that there exists a chart $\Bar{s}^* \in U$ such that $\card{\Sigma_i \setminus \im{\Bar{s}^*}} \geq d(\Bar{t})$ for all $i \in n$.
        We can then find a permutation $a \in \Stab{\partition}$ such that $g = \Bar{s}^* a \Bar{t}$ is a total chart whose image is a moiety of $X$ (see Figure \ref{fig:maximal_partition/g}).
        We then apply Lemma \ref{lem:missing_all_part} again to show that there exists a chart $g^* \in U$ such that $\card{\Sigma_i \setminus \im{g^*}} = \card{X}$ for all $i \in n$.
        Next we note that the charts $\Bar{a}, \Bar{b} \in U$ satisfy the conditions of Lemma \ref{lem:nxn}, so it follows from this and Lemma \ref{lem:compose_relations} that there exists a chart $h \in U$ such that $\rho_h = n \times n$. 
        Then $(Y)\inv{h} \cap \Sigma_i$ is a moiety of $\Sigma_i$ for all $i \in n$.
        Since $\im{g^*} \cap \Sigma_i$ has infinite complement in $\Sigma_i$ for every $i\in \n$, there exists $b \in \Stab{\partition}$ mapping $\im{g^*} \cap \Sigma_i$ into a moiety of $(Y)\inv{h} \cap \Sigma_i$ for every $i\in n$. 
        Then $f = g^* b h$ is a total chart whose image is a moiety of $Y$, as required (see Figure \ref{fig:maximal_partition/f}).
    \end{enumerate}
\end{proof}

\begin{corollary} \label{cor:maximal_inv_astab}
Let $X$ be a countably infinite set and $\partition$ a finite partition of $X$.
Then there exists a unique maximal inverse subsemigroup of $I_X$ which contains $\Stab{\partition}$ but not $\Sym{X}$, namely:
\begin{equation*}
    \mathcal{A}_\partition = A_{\partition} \cap \inv{A_{\partition}} = \set{f \in I_X \given \rho_f \in \Sym{n} ~\lor~ \dom{\rho_f} \neq n \neq \im{\rho_f}}
\end{equation*}
\end{corollary}

\begin{proof}
    This follows from Proposition \ref{prop:inverse_intersection}, which states that the maximal inverse subsemigroups containing some inverse subsemigroup $G$ are given by the intersections of the maximal subsemigroups that contain $G$ with their inverses.
    So we let $G = \Stab{\partition}$ and use the semigroups described in Theorem \ref{thm:maximal_astab} as our set of maximal subsemigroups.
    All that remains is to show that $\mathcal{A}_\partition$ is not contained in any of the inverse semigroups described in Corollary \ref{cor:maximal_inv_sym}.
    So choose $f,g \in I_X$ such that:
    \begin{itemize}
        \item $\rho_f \in \Sym{n}$, $c(f) > 0$, and $d(f) = 0$;

        \item $\rho_g \in \Sym{n}$, $c(g) < \aleph_0$, and $d(g) = \aleph_0$.
    \end{itemize}
    Then $f,g \in \mathcal{A}_\partition$ but $f \notin \mathcal{S}$ and $g \notin \mathcal{T}$.
\end{proof}

Throughout this entire chapter, we have stated our lemmas and preliminary results in terms of a general infinite set $X$, but in all of the main theorems we have demanded that $X$ be countable.
At the time of writing, these are the only results that we have, however we believe that generalising to an arbitrarily large infinite set $X$ is feasible.
However, this complicates the results by introducing more maximal subsemigroups (at least one for each infinite cardinal smaller than $\card{X}$), similar to the main results of \cite{maximal_east_mitchell_peresse} on maximal subsemigroups of $T_X$.
Furthermore, an uncountable $X$ would require doing a somewhat complicated transfinite induction over all infinite cardinals smaller than $\card{X}$ in each of our proofs.
We will attempt to address these issues and generalise the results in an upcoming paper.

\chapter{Topology and Relations} \label{chap:relational_structures}

In this chapter we study semigroups of homomorphisms and isomorphisms between induced substructures of relational structures on an infinite set.
In the first section we show how such semigroups relate to the topologies $\topI_1$ and $\topI_4$, and in the second section we find conditions under which such semigroups must have finite relative rank in $I_X$.

\section{Closed Submonoids of $I_X$}

It is a well-known result that a subgroup $G$ of $\Sym{X}$ is closed in the pointwise topology if and only if $G$ is the group $\Aut{\relStructure}$ of automorphisms of some relational structure $\relStructure$ on $X$ \cite{cameron2009oligomorphic}.
Similarly, a submonoid $S$ of $T_X$ is closed in the pointwise topology if and only if $S$ is the monoid $\Endo{\relStructure}$ of endomorphisms of some relational structure $\relStructure$ on $X$ \cite[Proposition 6.1]{cameron2006homomorphism}.
In this section, we will derive the analogous result for $I_X$.
In fact, we will do so for both subsemigroups and inverse subsemigroups of $I_X$.
We begin by defining the analogous objects to automorphism groups and endomorphism semigroups.

\begin{definition}[Injective partial endomorphism] \label{def:ipEnd}
Let $X$ be a set, $f \in I_X$, and $R$ a relation on $X$.
Then $f$ is called an \emph{injective partial endomorphism} on $R$, if it is an injective homomorphism between induced subrelations of $R$. 
That is, $(\forall \alpha \in \dom{f}^{\ary{R}})~ \alpha \in R \implies \alpha f \in R$. 
We then say that $f$ \emph{respects} $R$ as an endomorphism.

Similarly, if $\relStructure = (X,\relations)$ is a relational structure on $X$, we say that $f$ is an \emph{injective partial endomorphism} on $\relStructure$, if it is an injective homomorphism between induced substructures.
That is, $f$ is an injective partial endomorphism on all relations $R \in \relations$.
We then say that $f$ \emph{respects} $\relStructure$ as an endomorphism.
\end{definition}

\begin{definition}[Partial automorphism] \label{def:pAut}
Let $X$ be a set, $f \in I_X$, and $R$ a relation on $X$.
Then $f$ is called a \emph{partial automorphism} on $R$, if it is an isomorphism between induced subrelations of $R$. 
That is, $(\forall \alpha \in \dom{f}^{\ary{R}})~ \alpha f \in R \Leftrightarrow a \in R$. 
We then say that $f$ \emph{respects} $R$.

Similarly, if $\relStructure = (X,\relations)$ is a relational structure on $X$, we say that $f$ is a \emph{partial automorphism} on $\relStructure$, if it is an isomorphism between induced substructures.
That is, $f$ is a partial automorphism on all relations $R \in \relations$.
We then say that $f$ \emph{respects} $\relStructure$.
\end{definition}

In this chapter, semigroups of injective partial endomorphisms and partial automorphisms will play the role of `interesting' submonoids and inverse submonoids of $I_X$ respectively.
We start by studying monoids of injective partial endomorphisms.

\begin{lemma} \label{lem:ipEnd_properties}
    Let $\relStructure = (X,\relations)$ be a relational structure.
    Then the set of injective partial endomorphisms on $\relStructure$
    \begin{equation*}
        \ipEnd{\relStructure} = \set{f \in I_X \given (\forall R \in \relations) (\forall \alpha \in \dom{f}^{\ary{R}}) ~ \alpha \in R \implies \alpha f \in R}
    \end{equation*}
    is a full submonoid of $I_X$.
\end{lemma}

\begin{proof}
    We first show that $\ipEnd{\relStructure}$ is a subsemigroup of $I_X$.
    Let $f,g \in \ipEnd{\relStructure}$, $R \in \relations$ and $\alpha \in \dom{fg}^{\ary{R}}$.
    If $\alpha \in R$, then $\alpha f \in R \cap \dom{g}^{\ary{R}}$ and hence $\alpha fg \in R$.

    Next, we show that $\ipEnd{\relStructure}$ is full (i.e. $E_X \subseteq \ipEnd{\relStructure}$).
    If $e \in E_X$, then $(\forall n \in \omega) (\forall \alpha \in \dom{e}^n)~ \alpha e = \alpha$. Hence, any relation $R$ on $X$ is respected by $e$.
\end{proof}

Lemma \ref{lem:ipEnd_properties} immediately confronts us with an issue, as it shows that these partial endomorphism monoids are all \emph{full} submonoids of $I_X$.
This means that they cannot make up the entire set of closed subsemigroups of $I_X$ in either of the topologies $\topI_1$ or $\topI_4$ (we will see examples of closed subsemigroups of $I_X$ that are not full in Chapter \ref{chap:BS_IX}).
Instead we must restrict our focus to the set of subsemigroups of $I_X$ that are both \emph{full} and \emph{closed}.
However, considering that the idempotents of $I_X$ consist exclusively of partial identities, this might not be too much of a sacrifice to make (it is in some sense similar to the fact that the endomorphism monoids of relational structures make up the set of all closed submonoids of $T_X$, not the set of all closed subsemigroups).
Additionally, looking at the set of full submonoids of $I_X$ we get an interesting result regarding the closed sets in the topologies $\topI_1$ and $\topI_4$.

\begin{lemma} \label{lem:full+closed->I1=I4}
    Let $X$ be a set and $S$ a full subsemigroup of $I_X$.
    Then $S$ is closed in the topology $\topI_1$ if and only if it is closed in $\topI_4$.
\end{lemma}

\begin{proof}
    We first recall some definitions.
    The topology $\topI_1$ is generated by the sets:
    \begin{align} \label{I_1 subbase}
        U_{x,y} = \set{h \in I_X \given (x,y) \in h},     \qquad
        V_{x,y} = \set{h \in I_X \given (x,y) \notin h}
    \end{align}
    The topology $\topI_4$ is generated by the sets:
    \begin{equation} \label{I_4 subbase}
        U_{x,y} = \set{h \in I_X \given (x,y) \in h},
    \end{equation}
    \begin{equation*}
        W_{x} = \set{h \in I_X \given x \notin \dom{h}},      \quad
        W_{x}^{-1} = \set{h \in I_X \given x \notin \im{h}}
    \end{equation*}
    The subsemigroup $S$ of $I_X$ being \emph{full} means that $E_X \subseteq S$, where $E_X$ is the semilattice of idempotents.
    We will actually not be using the full power of this condition, but rather the weaker condition that $(\forall f \in S) (\forall e \in E_X) (ef,fe \in S)$.
    We will refer to this weaker condition as $S$ being \emph{semi-full} in $I_X$.

    It follows from Theorems \ref{I1_theorem} and \ref{I4_theorem} that any closed set in $\topI_1$ is also closed in $\topI_4$.
    In fact, we can easily show that $\topI_1 \subseteq \topI_4$ by noting that the $\topI_1$ basic open sets $V_{x,y} = \set{h \in I_X \given (x,y) \notin h}$ can be written as:
    \begin{align*}
        V_{x,y} = \paa{\bigcup_{z \neq x} U_{z,y}} \cup W_x = \paa{\bigcup_{z \neq y} U_{x,z}} \cup W_y^{-1}
    \end{align*}
    Hence, we only need to show that $S$ being semi-full and closed in $\topI_4$ implies that it is also closed in $\topI_1$.

    $S$ being closed in $\topI_4$ means that for all $f \in I_X \setminus S$, there exists a basic open neighbourhood $N$ of $f$ in $\topI_4$ such that $S \cap N = \emptyset$.
    Borrowing the notation used in \eqref{I_1 subbase} and \eqref{I_4 subbase}, the basic open neighbourhood $N$ can in general be written as 
    \begin{equation} \label{eq:open_I4}
        N = A \cap B \cap C,    \quad \text{ where }
        A = \bigcap_{i \in n} U_{x_i,y_i},~ 
        B = \bigcap_{i \in m} W_{z_i}, \text{ and }
        C = \bigcap_{i \in l} \inv{W_{w_i}}
    \end{equation}
    with $n,m,l \in \omega$ and $x_i,y_i,z_i,w_i \in X$.
    Here we use that the intersection over an empty index set is equal to $I_X$, which occurs when $n$, $m$, or $l$ is equal to zero. 
    The sets $A$, $B$, and $C$ are also open neighbourhoods of $f$, even basic open neighbourhoods in the case $n, m, l \neq 0$. 
    We now claim that if $N = A \cap B \cap C$ has empty intersection with $S$, then so has $A$.
    We will prove this via contradiction.
    So we assume there exists some chart $g \in S \cap A$.
    Letting $z_i$ and $w_i$ denote the elements by the same names in \eqref{eq:open_I4}, we can then pick idempotents $e,d \in E_X$ such that $(\forall i \in m) (z_i \notin \dom{e})$ and $(\forall i \in l) (w_i \notin \im{d})$.
    Then $g' = egd \in A \cap B \cap C = N$. 
    But it would then follow from $S$ being semi-full that $g' \in S$, which contradicts the hypothesis that $S \cap N = \emptyset$.
    Hence, $A$ has empty intersection with $S$, and since $A$ is also a basic open neighbourhood of $f$ in $\topI_1$, we conclude that $S$ is closed in $\topI_1$ as well.
\end{proof}

Lastly, we must show that partial endomorphism monoids are indeed themselves closed in some interesting topology on $I_X$.

\begin{lemma} \label{lem:ipEnd_closed}
    Let $\relStructure = (X,\relations)$ be a relational structure. 
    Then $\ipEnd{\relStructure}$ is closed in every $T_1$ shift-continuous topology on $I_X$.
\end{lemma}

\begin{proof}
    By Theorem \ref{I1_theorem} we only have to check that $\ipEnd{\relStructure}$ is closed in the topology $\topI_1$.
    So either $\ipEnd{\relStructure} = I_X$, in which case the conclusion is trivial, or there exists some chart $f \in I_X \setminus \ipEnd{\relStructure}$. 
    If such an $f$ exists, then there exists some $R \in \relations$ and $\alpha \in R \cap \dom{f}^{\ary{R}}$ such that $\alpha f \notin R$.
    Then $U_{\alpha,\alpha f} = \set{h \in I_X \given \alpha h = \alpha f}$ forms a basic open neighbourhood of $f$ in $\topI_1$, which is entirely disjoint from $\ipEnd{\relStructure}$.
    Hence, $I_X \setminus \ipEnd{\relStructure}$ is open in the topology $\topI_1$, and we can conclude that $\ipEnd{\relStructure}$ is closed.
\end{proof}

So in particular the partial endomorphism monoids are indeed closed in both of the topologies $\topI_1$ and $\topI_4$.
Finally, we are ready to prove that any full and closed submonoid of $I_X$ is the monoid of injective partial endomorphisms on some relational structure.

\begin{theorem} \label{thm:full+closed=ipEnd}
    Let $X$ be a set and $S$ a full submonoid of $I_X$. 
    Then the following are equivalent:
    \begin{enumerate}[\normalfont(i)]
        \item $S$ is closed in the topology $\topI_1$.
        \item $S$ is closed in the topology $\topI_4$.
        \item There exists a relational structure $\relStructure = (X,\relations)$ such that $S = \ipEnd{\relStructure}$.
    \end{enumerate}
\end{theorem}

\begin{proof}
    We prove this by constructing the implication chain (i)$\implies$(ii)$\implies$(iii)$\implies$(i).
    
    (i)$\implies$(ii): This follows from Lemma \ref{lem:full+closed->I1=I4}.
    
    (ii)$\implies$(iii): Let $A_n = \set{\alpha = (a_0,\dots,a_{n-1})  \in X^n \given a_i \neq a_j \text{ when } i \neq j}$ be the set of all $n$-tuples with distinct elements in $X$.
    We then define for each $n \in \omega$ and each $\alpha \in A_n$ an $n$-ary relation $R_\alpha$ by
    \begin{align*}
	    R_\alpha = \set{\alpha f \in X^n \given f \in S \,\wedge\, (\forall i \in n) (a_i \in \dom{f})}
    \end{align*}
    where $\alpha f$ is the $n$-tuple $(a_0 f,\dots,a_{n-1} f)$.
    Since $f$ is injective, it follows that $\alpha f$ is also an $n$-tuple of distinct elements.
    Let $\relations = \set{R_\alpha \given \alpha \in \cup_{n \in \omega} A_n}$ and $\relStructure = (X,\relations)$.
    We then claim that $S = \ipEnd{\relStructure}$. 

    First we will show that $S \subseteq \ipEnd{\relStructure}$.	
    That is, we wish to show that for all $f \in S$, $f$ respects all of the relations in  $\relStructure$.
    Given any tuple $\alpha \in \cup_{n \in \omega} A_n$, then $\alpha f$ is an element of $R_\alpha$ by the definition of $R_\alpha$, and if $\alpha \in R_\beta$ for some other tuple $\beta \neq \alpha$, then there exist charts $g,h \in S$ such that $\beta g = \alpha$ and $gf = h$.
    This means that $\beta h = \beta gf = \alpha f \in R_\beta$.
    Hence, we conclude that $S \subseteq \ipEnd{\relStructure}$.
	
    Next we show that $\ipEnd{\relStructure} \subseteq S$.
    We will do this by showing that every element of $\ipEnd{\relStructure}$ is in the closure of $S$.
    Let $h \in \ipEnd{\relStructure}$ and $\alpha$ a tuple of distinct elements from $\dom{h}$. 
    Then, by the definition of $R_\alpha$, there exists some $g \in S$ such that $\alpha g = \alpha h$, and since $S$ is full it follows that $g' = \eval{g}_{\dom{h}}^{\im{h}}$ is also in $S$. 
    Then:
    \begin{equation*}
        g' \in U_{\alpha, \alpha h} \cap \paa{\bigcap_{x \notin \dom{h}} W_x} \cap \paa{\bigcap_{y \notin \im{h}} \inv{W}_y}
    \end{equation*}
    where the sets $U_{\alpha, \alpha h} = \set{f \in I_X \given \alpha f = \alpha h}$, $W_{x} = \set{f \in I_X \given x \notin \dom{f}}$, and $W_{y}^{-1} = \set{f \in I_X \given y \notin \im{f}}$ are the basic open neighbourhoods of $h$ in $\topI_4$.
    Since neither $h$ nor $\alpha$ were specified, this means that given any open neighbourhood $N$ around an element of $\ipEnd{\relStructure}$, we can always find an element of $S$ which also belongs to $N$.
    Since $S$ is closed in $\topI_4$, we can then conclude that $\ipEnd{\relStructure} \subseteq S$. 
    Hence, $S = \ipEnd{\relStructure}$.
	
    (iii)$\implies$(i): This is Lemma \ref{lem:ipEnd_closed}
\end{proof}

With this we have proven that the partial endomorphism monoids make up the set of all full closed submonoids of $I_X$ in our favourite topologies.
We are not done however, as $I_X$ is an inverse semigroup (perhaps the canonical example of such) and we would therefore like to classify the closed inverse submonoids as well.
For this we shift our focus to monoids of partial automorphisms on relational structures.
First, we show how these relate to partial endomorphism monoids.

\begin{lemma} \label{lem:pAut=ipEndInverses}
    Let $\relStructure = (X,\relations)$ be a relational structure.
    Then the intersection $\ipEnd{\relStructure} \cap \inv{\ipEnd{\relStructure}}$ is equal to
    \begin{equation*}
        \pAut{\relStructure} = \set{f \in I_X \given (\forall R \in \relations) (\forall \alpha \in \dom{f}^{\ary{R}}) ~ \alpha f \in R \Leftrightarrow \alpha \in R},
    \end{equation*}
    the set of all partial automorphisms on $\relStructure$.
\end{lemma}

\begin{proof}
    We start by writing out $\ipEnd{\relStructure}$ and its inverse $\inv{\ipEnd{\relStructure}}$:
    \begin{align*}
        \ipEnd{\relStructure} &= \set{f \in I_X \given (\forall R \in \relations) (\forall \alpha \in \dom{f}^{\ary{R}})~ \alpha \in R \implies \alpha f \in R}
        \\
        \inv{\ipEnd{\relStructure}} &= \set{f \in I_X \given (\forall R \in \relations) (\forall \alpha \in \dom{f}^{\ary{R}})~ \alpha f \in R \implies \alpha \in R}
    \end{align*}
    Then the intersection $\ipEnd{\relStructure} \cap \inv{\ipEnd{\relStructure}}$ consist of those charts $f \in I_X$ such that for every relation $R \in \relations$ and tuple $\alpha \in \dom{f}^{\ary{R}}$, the image and preimage of $\alpha$ under $f$ is also in $R$.
    This is precisely the set of all partial automorphisms of $\relStructure$, as claimed.
\end{proof}

Lemma \ref{lem:pAut=ipEndInverses} will allow us to easily derive results for partial automorphism monoids analogous to those that we have above for partial endomorphism monoids.
So that is what we shall do.

\begin{lemma}
    Let $\relStructure = (X,\relations)$ be a relational structure.
    Then $\pAut{\relStructure}$ is a full inverse submonoid of $I_X$.
\end{lemma}

\begin{proof}


    %
    %
    %
    This follows from the fact that $\pAut{\relStructure} = \ipEnd{\relStructure} \cap \inv{\ipEnd{\relStructure}}$ together with Lemma \ref{lem:ipEnd_properties}, since $E_X$ is its own inverse and the intersection of a semigroup with its inverse is an inverse semigroup.
\end{proof}

\begin{lemma} \label{lem:pAut_closed}
    Let $\relStructure = (X,\relations)$ be a relational structure. 
    Then $\pAut{\relStructure}$ is closed in every $T_1$ shift-continuous topology on $I_X$.
\end{lemma}

\begin{proof}
    %
    %
    This again follows from the fact $\pAut{\relStructure} = \ipEnd{\relStructure} \cap \inv{\ipEnd{\relStructure}}$ together with Lemma \ref{lem:ipEnd_closed}, since inversion is continuous under $\topI_1$ and $\pAut{\relStructure}$ is therefore the intersection of two closed sets.
\end{proof}

And finally, we also have a classification theorem of the full inverse submonoids of $I_X$ that are closed in $\topI_1$ and $\topI_4$.

\begin{theorem} \label{thm:full+inverse+closed=pAut}
Let $X$ be a set and $M$ a full inverse submonoid of $I_X$. 
Then the following are equivalent:
\begin{enumerate}[\normalfont(i)]
	\item $M$ is closed in the topology $\topI_1$.
	\item $M$ is closed in the topology $\topI_4$.
	\item There exists a relational structure $\relStructure = (X,\relations)$ such that $M = \pAut{\relStructure}$.
\end{enumerate}
\end{theorem}

\begin{proof}
    This follows from Theorem \ref{thm:full+closed=ipEnd}, since $M$ also satisfies the hypothesis of this theorem.
    So if $M$ satisfies condition (i) or (ii), then there exists a relational structure $\relStructure'$ on $X$ such that $M = \ipEnd{\relStructure'}$.
    However, since $M$ is an inverse monoid it would follow that $\ipEnd{\relStructure'}$ is also an inverse monoid.
    This means that $\ipEnd{\relStructure'} = \inv{\ipEnd{\relStructure'}}$ and therefore by Lemma \ref{lem:pAut=ipEndInverses} we get that $\pAut{\relStructure'} = \ipEnd{\relStructure'} \cap \inv{\ipEnd{\relStructure'}} = \ipEnd{\relStructure'} = M$.
    So let $\relStructure = \relStructure'$.
    Since $M$ is the monoid of injective partial endomorphisms of some relational structure, it then also follows from Theorem \ref{thm:full+closed=ipEnd} that condition (iii) implies conditions (i) and (ii).
\end{proof}

Theorems \ref{thm:full+closed=ipEnd} and \ref{thm:full+inverse+closed=pAut} are an indication that full submonoids are of particular interest when studying the closed subsemigroups of $I_X$.
The following result further motivates this connection between closed and full subsemigroups of $I_X$.

\begin{proposition} \label{prop:closed_I1->full_closed}
    Let $X$ be a countably infinite set and $S$ a closed subsemigroup in the topology $\mathcal{I}_1$ on $I_X$.
    Then $S' = \genset{S,E_X}$ is also closed in $\mathcal{I}_1$.
\end{proposition}

\begin{proof}
    We will prove that $S'$ is closed by showing that it is its own closure.
    Let $g \in I_X$ be an element the closure of $S'$. 
    Enumerate the domain of $g$ as $\dom{g} = \set{x_i \mid i \in \omega}$ and for each $n \in \omega$ we define $s_n \in I_X$ to be the restriction of $g$ to the first $n$ elements of $\dom{g}$.
    We will show that these $s_n$ are all elements of $S'$.
    Since $g$ is in the closure of $S'$ and for each $n \in \omega$ the open neighbourhood $U = \set{h \in I_X \given (\forall i \in n)~ x_i h = x_i g}$ of $g$ contains an element $f \in S'$. 
    Since $E_X \subseteq S'$ it follows that $s_n = \eval{g}_n = \eval{f}_n$ is indeed an element of $S'$.
    And the sequence $(s_n)_{n\in\omega}$ converges to $g$, since for any finite tuple $a$ of elements in $X$, there exists $m \in \omega$ such that for all $n > m$, $(a)s_n = (a)g$. 
	
    Being elements of $S'$, each $s_n$ extends to an element $t_n \in S$.
    So we can `extend' the entire sequence $(s_n)_{n\in\omega}$ in $S'$ to a sequence $(t_n)_{n\in\omega}$ in $S$.
    Here we encounter a problem, as this sequence does not necessarily converge.
    To get around this we will use the fact that $\topI_1$ is a compact Polish topology (Theorem \ref{I1_theorem}). Being Polish implies that the topology is metrisable, which together with compactness implies that $\topI_1$ is sequentially compact (assuming the Axiom of Choice) \cite[Theorem 28.2]{munkrestopology}.
    This means that every sequence has a converging subsequence. 
    Since $S$ is closed, we can conclude that there exists a subsequence $(t_{n_k})_{k\in\omega}$ of $(t_n)_{n\in\omega}$, which converges to a chart $t \in S$.
    Since $(s_n)_{n\in\omega}$ converges to $g$, we can conclude that the subsequence $(s_{n_k})_{k\in\omega}$, defined as the subsequence which `extends' to $(t_{n_k})_{k\in\omega}$, also converges to $g$.
    As we let $k$ increase, $s_{n_k}$ agrees with $g$ on more and more elements in $\dom{g}$, and by extension so does $t_{n_k}$.
    Hence the limit $t$ must agree with $g$ on the entire domain of $g$, which means that $t$ is an extension of $g$. 
    It follows that $g$ must be an element of $S'$, completing the proof.
\end{proof}

\section{Partial Automorphism Monoids}

In this section we will show that for most of the relational structures $\relStructure$ that one would usually consider, the submonoid $\pAut{\relStructure}$ has finite relative rank in $I_X$.
That is, $\pAut{\relStructure} \bsequal[I] I_X$ whenever $\relStructure$ is sufficiently simple and $\card{X}$ is in some sense `similar' to $\aleph_0$.
This work is partially inspired by \cite{higgins2003generating}, where the relative ranks of semigroups of order preserving maps in the full transformation monoid are studied.
First, we introduce the notion of a \emph{local permutation}.

\begin{definition}[Local permutation]
    Let $\alpha$ be an ordinal. 
    We call a permutation $g \in \Sym{\alpha}$ \emph{local} if for all $i \in \alpha$ there exists $j > i$ in $\alpha$ such that $g$ maps the initial segment $s(j) = \set{k \in \alpha \given k < j}$ to itself.
\end{definition}

Such local permutations will be immensely useful throughout the rest of the thesis.
This is due to the following lemma, which is a generalisation of \cite[Lemma 9]{Bergman_2006}.

\begin{lemma} \label{lem:local_generate}
    Let $\kappa$ be a regular cardinal.
    Then the local permutations on $\kappa$ form a generating set for $\Sym{\kappa}$.
\end{lemma}

\begin{proof}
    It is well known that every permutation is the product of two involutions \cite[Lemma 2.2]{galvin1995generating}.
    It thus suffices to show that the local permutations generate the set of all involutions on $\kappa$.
    So given an involution $f \in \Sym{\kappa}$, we will recursively construct a sequence of ordinals $(a_i)_{i \in \kappa}$.
    Let $a_i$ be any ordinal in $\kappa$ such that for all $j < i$, $a_i > a_j$ and $s(a_j)f \,\cup\, s(a_j)\inv{f} \subseteq s(a_i)$.
    This is always possible since $\kappa$ is regular.
    From this we define a partition $\partition = \set{\Sigma_i \given i \in \kappa}$ of $\kappa$, where $\Sigma_i = \set{k \in \kappa \given a_i \leq k < a_{i+1}}$.
    If $i = 0$ or $i$ is a limit ordinal we denote $\Sigma_{i-1} = \emptyset$.
    The partition $\partition$ then has the property for all $i \in \kappa$ that $\card{\Sigma_i} < \kappa$ and $\Sigma_i f \subseteq \Sigma_{i-1} \cup \Sigma_i \cup \Sigma_{i+1}$.
    Since $f$ is an involution it follows that $(\forall i \in \kappa) (\Sigma_i f \cap \Sigma_{i+1} = \Sigma_{i+1} f \cap \Sigma_i)$.
    That is, $f$ `carries' equally many elements from $\Sigma_i$ to $\Sigma_{i+1}$ as it does from $\Sigma_{i+1}$ to $\Sigma_i$ (this holds true in general for any two subsets when $f$ is an involution).
    Our goal now is to construct a permutation $g$ such that for all $i \in \kappa$, $g$ maps $\Sigma_{2i} \cup \Sigma_{2i+1}$ to itself and $\inv{g}f$ maps $\Sigma_{2i-1} \cup \Sigma_{2i}$ to itself.
    Then $g$ and $\inv{g}f$ are local permutations and their product is $f$, as required.

    For each $i \in \kappa$ let $g$ `act as $f$' on the elements that $f$ maps from $\Sigma_{2i}$ to $\Sigma_{2i+1}$ and those it maps from $\Sigma_{2i+1}$ to $\Sigma_{2i}$.
    Let $g$ fix all remaining elements.
    Then $g$ preserves $\Sigma_{2i} \cup \Sigma_{2i+1}$ and for all $i \in \kappa$, $ig$ lies in the same $\Sigma_{2i-1} \cup \Sigma_{2i}$ interval as $if$ does.
    Hence, $\inv{g}f$ preserves the interval $\Sigma_{2i-1} \cup \Sigma_{2i}$, as required.
\end{proof}

A special type of relational structures, which we are very familiar with, is ordered sets.
As such, we introduce special notation and nomenclature for partial automorphism monoids of order relations.

\begin{definition}[Order preserving chart]
    Let $(X,\leq)$ be an ordered set.
    We shall then use the shorthand notation;
    \begin{equation*}
        \order{(X,\leq)} = \set{f \in I_{X} \given (\forall x,y \in \dom{f})~ x \leq y \Leftrightarrow (x)f \leq (y)f}
    \end{equation*}
    for the monoid $\pAut{\leq}$, or simply $\order{X}$ if the order is clear from context.
    Furthermore, we will refer to $f \in \order{X}$ as an \emph{order preserving chart} on $X$.
\end{definition}

We are introducing special notation for monoids of order preserving charts, as these will play a special role in the coming results.
Among the different type of order relations, well-orders are particularly `nice.'

\begin{lemma} \label{lem:order_preserving}
    Let $X$ be a well-ordered set of regular cardinality $\kappa$ such that for all $\lambda < \kappa$ we have $2^\lambda \leq \kappa$.
    Then $\order{X} \bsequal[I] I_{X}$.
\end{lemma}

\begin{proof}
    We may assume that $\kappa$ is a suborder of $X$ and since $\Sym{X} \approx I_X$, it suffices by Lemma \ref{lem:local_generate} to show that $\order{X}$ finitely generates the local permutations on $\kappa$.
    Let $\partition$ be a partition of $\kappa$ such that for all $\lambda \in \kappa$ there are $\kappa$ many elements in $\partition$ of size $\lambda$.
    Let $p \in \Sym{\kappa}$ be a permutation such that for all $\lambda \in \kappa$ and all $q \in \Sym{\lambda}$, there exist $\kappa$ many $A \in \partition$ such that $q = \eval{p}_A$.
    It is possible to construct such a $p$, since we assume that for all $\lambda < \kappa$ we have $2^\lambda \leq \kappa$.    
    Finally, let $g \in \Sym{\kappa}$ be a local permutation.
    We will then construct an order preserving chart $f \in \order{X}$ such that $g = fp\inv{f}$.

    Since $g$ is a local permutation, it defines a partition $\set{\Sigma_i \given i \in \kappa}$ of $\kappa$ into subsets all of size $< \kappa$, where $g$ maps each $\Sigma_i$ to itself.
    For each $i \in \kappa$ let $f$ map $\Sigma_i$ to some $A \in \partition$ such that $\Sigma_i fp\inv{f} = \Sigma_i g$.
    We can choose such an $f$ to be order preserving, since there are $\kappa$ many such $A \in \partition$, meaning that they form an unbounded subset of $\kappa$.
    Then $fp\inv{f} = g$ with $f, \inv{f} \in \order{X}$, completing the proof.
\end{proof}

What Lemma \ref{lem:order_preserving} tells us is that if $\card{X}$ is regular and not a counterexample to the generalised continuum hypothesis (the set theoretic assumption that for all infinite cardinals $\kappa$, $2^\kappa = \kappa^+$), then $\order{X,\leq} \bsequal[I] I_X$ for all well-orders $\leq$ on $X$.
To progress from here however, we will need to introduce some new concepts.
The field of \emph{Ramsey Theory} concerns itself with combinatorics of partitions of sets and finding `large subsets' with particular properties.
We first go over some definitions and then present some well-known results from the field, which we will be using.

\begin{definition}
    Let $X$ be a set and $n$ a cardinal.
    We then denote by
    \begin{equation*}
        [X]^n = \set{Y \subseteq X \given \card{Y} = n}
    \end{equation*}
    the set of all $n$-element subsets of $X$.
\end{definition}

\begin{definition}
    Let $X$ be a set, $n$ a cardinal, $\partition$ a partition of $[X]^n$, and $H \subseteq X$. We then say that $H$ is homogeneous with respect to $\partition$ if there exists $A \in \partition$ such that $[H]^n \subseteq A$.
\end{definition}

\begin{definition}[Arrow notation]
    Let $\kappa$, $\lambda$, and $\mu$ be cardinals and $n \in \omega$.
    If every partition of $[\kappa]^n$ into $\mu$ parts has a homogeneous set of size $\lambda$, we write;
    \begin{equation*}
        \kappa \arrows{\lambda}{n}{\mu}
    \end{equation*}
    which reads as `$\kappa$ \emph{arrows} $\lambda$'.
\end{definition}

With the basic definition out of the way, we can now state the main theorem of Ramsey Theory that we will be using.
This theorem is also know as \emph{Ramsey's (infinite) Theorem}.

\begin{theorem}[{Ramsey's Theorem \cite[Theorem A]{ramsey1987problem}}] \label{RamseyTheorem}
    Let $X$ be an infinite set and $n, m \in \omega$.
    Then every partition of $[X]^n$ into $m$ parts admits an infinite homogeneous set $H \subseteq X$.
    That is, $\aleph_0 \arrows{\aleph_0}{n}{m}$.
\end{theorem}

The following definition and lemmas are also staples of Ramsey Theory, which we will need.

\begin{definition}
    An infinite cardinal $\kappa$ is called a \emph{strong limit cardinal}, if it is regular and for all $\lambda < \kappa$ we have $2^\lambda < \kappa$.
\end{definition}

\begin{lemma}[{\cite[Lemma 9.9]{Jech2003}}] \label{lem:Ramsey_strong_limit}
    If $\kappa \arrows{\kappa}{2}{2}$, then $\kappa$ is a strong limit cardinal.
\end{lemma}

\begin{lemma}[{\cite[Chapter 7, Theorem 3.5]{drake1974set}}] \label{lem:Ramsey_equivalent}
    Let $\kappa$ be an infinite cardinal. 
    Then the following are equivalent.
    \begin{enumerate}[~\normalfont(i)]
        \item $\kappa \arrows{\kappa}{2}{2}$;
        \item $\kappa \arrows{\kappa}{n}{\mu}$ for all $\mu < \kappa$ and $n \in \omega$.
    \end{enumerate}
\end{lemma}

With this we are finally able to proceed.
In the following results we will be using the well-known fact, that if two well-ordered sets $X$ and $Y$ are isomorphic as ordered sets, then the order isomorphism between them is unique \cite[Corollary 2.6]{Jech2003}.
Also, given a relational structure $\relStructure = (X, \relations)$ and a relation $R$ on $X$, we will be using the notation $\relStructure \land R$ to indicate the extended relational structure $(X, \relations \cup \set{R})$.

\begin{lemma} \label{lem:hom->contain_order}
    Let $\relStructure = (X, \relations)$ be a relational structure, $n = \sup(\ary{R} \given R \in \relations)$, and $\leq$ a well-ordering of $X$.
    If a subset $H \subseteq X$ is homogeneous with respect to the partition of $[X]^n$ into orbits under the action of $\pAut{\relStructure ~\land \leq}$, then $\order{(H,\leq)} \subseteq \pAut{\relStructure}$.
\end{lemma}

\begin{proof}
    The set $H$ being homogeneous with respect to the orbits of $[X]^n$ under $\pAut{\relStructure ~\land \leq}$ means that for all $A,B \in [H]^n$, there exists an order preserving chart (with respect to $\leq$) $g \in \pAut{\relStructure}$ such that $Ag = B$.
    Let $f$ be an order preserving chart on $H$ with respect to $\leq$.
    We will show that $f \in \pAut{\relStructure}$.
    If $\card{\dom{f}} < n$, then there exists an extension $f'$ of $f$ with $\dom{f'} \geq n$ in $\order{(H,\leq)}$. 
    So if $f' \in \pAut{\relStructure}$, then so is $f$ since the restriction of a partial automorphism is also a partial automorphism.
    Hence we can assume without loss of generality that $\card{\dom{f}} \geq n$.
    Since $(X,\leq)$ is a well-order, any subset of $X$ is also well-ordered by $\leq$, hence the chart $f$ must necessarily be the unique order isomorphism from $\dom{f}$ to $\im{f}$.
    Consider an arbitrary $m$-tuple $a \in \dom{f}^m$ with $m \leq n$.
    Choose an $A \in [\dom{f}]^n$ such that $a \in A^m$.
    Since $A,Af \in [H]^n$ there necessarily exists $g \in \pAut{\relStructure ~\land \leq}$ such that $Ag = Af$. 
    Hence, by the uniqueness of order isomorphisms of well-orders, $\eval{g}_A = \eval{f}_A$.
    So if $R$ is an arbitrary $m$-ary relation in $\relStructure$, then $(a)f \in R \Leftrightarrow (a)g \in R \Leftrightarrow a \in R$.
    Therefore $f \in \pAut{\relStructure}$.
\end{proof}

\begin{corollary} \label{cor:hom->pAut=I_X}
    Let $\relStructure = (X, \relations)$ be a relational structure on an infinite set $X$ of regular cardinality $\kappa$ such that for all $\lambda < \kappa$ we have $2^\lambda \leq \kappa$.
    If there exists a well-ordering $\leq$ of $X$ such that the partition of $[X]^n$ into orbits under the action of $\pAut{\relStructure ~\land \leq}$ admits a homogeneous subset $H \subseteq X$ of size $\kappa$, then $\pAut{\relStructure} \bsequal[I] I_X$.
\end{corollary}

\begin{proof}
    This follows from Lemmas \ref{lem:order_preserving} and \ref{lem:hom->contain_order}.
\end{proof}

We can now state the main theorem of this section.

\begin{theorem} \label{thm:finite_rel=I_X}
    Let $X$ be a set of cardinality $\kappa$ with $\kappa \arrows{\kappa}{2}{2}$.
    If $\relStructure = (X, \relations)$ is a relational structure on $X$ such that $\card{\relations} < \kappa$ and the set $\set{\ary{R} \given R \in \relations}$ is bounded in $\omega$,
    then $\pAut{\relStructure} \bsequal[I] I_X$.
\end{theorem}

\begin{proof}
    Lemma \ref{lem:Ramsey_strong_limit} states that $\kappa$ must be a strong limit cardinal and by Lemma \ref{lem:Ramsey_equivalent} we have that $\kappa \arrows{\kappa}{n}{\mu}$ for all $n \in \omega$ and $\mu < \kappa$.
    Let $m = \max(\ary{R} \given R \in \relations)$ and choose a well-order $\leq$ on $X$.    
    For any two $m$-element subsets $A,B \in [X]^m$ to be in different orbits under $\pAut{\relStructure ~\land \leq}$ they would have to have `different' induced substructures under $\relStructure$ when viewed as labelled by their order under $\leq$.
    But the number of possible induced substructures on an $m$-element subset of $X$ given a relational structure with $\card{\relations}$ many relations is bounded by $2^{m^m} \card{\relations}$.
    Hence the number of orbits under $\pAut{\relStructure ~\land \leq}$ is also bounded by $2^{m^m} \card{\relations}$.
    Since $m$ is finite and $\card{\relations} < \kappa$, it follows that $2^{m^m} \card{\relations} < \kappa$. 
    If we denote the actual number of orbits on $[X]^m$ under $\pAut{\relStructure ~\land \leq}$ by $\lambda$, then $\kappa \arrows{\kappa}{m}{\lambda}$.
    Since $\kappa$ is a strong limit cardinal it satisfies the conditions of Corollary \ref{cor:hom->pAut=I_X} from which the result then follows.
\end{proof}

Since the countably infinite cardinal $\aleph_0$ has the property that $\aleph_0 \arrows{\aleph_0}{2}{2}$, we can state the following corollary to Theorem \ref{thm:finite_rel=I_X}.

\begin{corollary} \label{cor:finite_rel=IX}
    Let $\relStructure = (X, \relations)$ be a relational structure on a countably infinite set $X$.
    If $\card{\relations}$ is finite, then $\pAut{\relStructure} \bsequal[I] I_X$.
\end{corollary}

This means that for far most relational structures $\relStructure$ that we are familiar with, such as graphs or partial orders on a countable set, we get that $\pAut{\relStructure} \bsequal[I] I_X$.
This is a surprising result, as it is vastly different from both the $\Sym{X}$ and $T_X$ cases!
Namely, there exist graphs with no non-trivial automorphisms or endomorphisms (see e.g. \cite[Theorem 3]{hell1973groups}).
Furthermore, Corollary \ref{cor:finite_rel=IX} allows us to expand on a previous result, namely Lemma \ref{lem:contained_in_maximal}, concerning when a given subgroup of $I_X$ is guaranteed to be contained in a maximal subsemigroup.

\begin{proposition} \label{prop:contained_in_maximal}    
    Let $X$ be an infinite set and $S$ a proper subsemigroup of $I_X$ satisfying any of the following:
    \begin{enumerate}[\normalfont(i)]
        \item $S$ has finite relative rank in $I_X$;
        \label{prop:contained_in_maximal/S_fin_rank}
        \item $\card{S} \leq \card{X}$;
        \label{prop:contained_in_maximal/S<X}
        \item $X$ is countable, $S$ is inverse, and there exists $n \in \omega$ such that $S$ has more than one orbit on $X^n$; or
        \label{prop:contained_in_maximal/orbits>1}
        \item $X$ is countable, $S$ is inverse, and the closure of $S$ in the topology $\topI_4$ is a proper subset of $I_X$.
        \label{prop:contained_in_maximal/closure}
    \end{enumerate}
    Then $S$ is contained in a maximal subsemigroup of $I_X$.
\end{proposition}

\begin{proof}
    We prove each statement in order as listed above.
    \begin{enumerate}[\normalfont(i)]
        \item This given by Lemma \ref{lem:contained_in_maximal}.
    
        \item This given by Lemma \ref{lem:contained_in_maximal}.
    
        \item Given $n \in \omega$, we can construct a relational structure $\relStructure = (X, \relations)$ where the relations $R \in \relations$ are exactly the orbits of $S$ in $X^n$. 
        Then $S$ is a subsemigroup of $\pAut{\relStructure}$.
        Furthermore, if $\relStructure' = (X,\relations'$ is another relational structure on $X$ such that for all $R' \in \relations'$ there exists a subset $A$ of $\relations$ such that $R' = \bigcup A$, then $\pAut{\relStructure} \subseteq \pAut{\relStructure'}$.
        Hence, if there exists $n \in \omega$ such that $S$ has strictly more than one orbit on $X^n$, we can construct a relational structure $\relStructure' = (X,\set{R_1, R_2})$ with exactly two distinct relations such that $S \subseteq \pAut{\relStructure'}$, by letting $R_1$ be an orbit of $S$ on $X^n$ and $R_2$ the union of the remaining orbits of $S$ on $X^n$.
        Since $\relStructure'$ contains only finitely many relations, Theorem \ref{thm:finite_rel=I_X} states that $\pAut{\relStructure'}$ has finite relative rank in $I_X$.
        Also, since $\relStructure'$ contains two distinct relations, it follows that $\pAut{\relStructure'} \neq I_X$.
        It thus follows from \eqref{prop:contained_in_maximal/S_fin_rank} that $\pAut{\relStructure'}$, and hence $S$, is contained in a maximal subsemigroup of $I_X$.
    
        \item We only need to show that this holds in the case where $S$ does not already satisfy condition \eqref{prop:contained_in_maximal/orbits>1}.
        That is, when $S$ has exactly one orbit on $X^n$ for all $n \in \omega$.
        We say that $S$ is \emph{transitive on tuples} in $X$.
        We then claim that in this case, $\Sym{X}$ is contained in the closure of $S$.
        Given any permutation $g \in \Sym{X}$, basic open neighbourhoods of $g$ in the topology $\topI_4$ take the form $U_a = \set{f \in I_X \given (a)f = (a)g}$, where $a$ is any tuple of elements in $X$ (we do not need to consider any other basic open sets in $\topI_4$, since $c(g) = d(g) = 0$).
        But since $S$ is transitive on tuples, we know that there exists $h \in S$ such that $(a)h = (a)g$.
        We can conclude that $g$ is in the closure of $S$, and $\Sym{X}$ is hence contained in the closure of $S$.
        Thus, since $\Sym{X}$ has finite relative rank in $I_X$, we get that the closure of $S$ also has finite relative rank in $I_X$.
        We know from Lemma \ref{lem:closure_subsemigroup} that the closure of a subsemigroup of any topological semigroup is again a subsemigroup, so it follows from \eqref{prop:contained_in_maximal/S_fin_rank} that the closure of $S$, and hence also $S$, is contained in a maximal subsemigroup of $I_X$ (since the closure of $S$ was assumed to be a proper subset of $I_X$).
    \end{enumerate}
\end{proof}

From Proposition \ref{prop:contained_in_maximal} we see that there are very strict requirements for an inverse subsemigroup $S$ of $I_X$ to not be contained in a maximal subsemigroup of $I_X$ when $X$ is countable.
If such a subsemigroup $S$ exists, then $S$ must be a dense subset $I_X$ in the topology $\topI_4$ and it must be transitive on tuples in $X$.

\chapter{Bergman-Shelah Preorder on $I_X$} \label{chap:BS_IX}

In this chapter we study the equivalence classes of subsemigroups of $I_X$ under the Bergman-Shelah equivalence relation.
As mentioned in the introduction, it was discovered by George Bergman and Saharon Shelah in their 2006 paper \cite{Bergman_2006} that the closed subgroups in the pointwise topology on $\Sym{X}$ fall into exactly four distinct equivalence classes, which form a chain under the Bergman-Shelah preorder (Definition \ref{def:BS_preorder}), when $X$ is a countably infinite set.
In this chapter we will attempt to carry out the same analysis for $I_X$ in hopes of obtaining a similar result.
In the first section, we extend some of the results achieved by Bergman and Shelah for $\Sym{X}$ to $I_X$, and we show that the four equivalence classes found in \cite{Bergman_2006} are still distinct in this more general setting.
In the second section, we study inverse subsemigroups of $I_X$ consisting of bijections between paths in rooted trees, as they provide interesting case studies, and we fully classify all such subsemigroups under the Bergman-Shelah preorder.
Finally in the third section, we show examples of subsemigroups of $I_X$ that lie outside the original four equivalence classes introduced by Bergman and Shelah, and we use these examples to provide a list of necessary conditions to formulate a new theorem similar to the main theorem of \cite{Bergman_2006} (Theorem \ref{BS_mainTheorem}).
At the end, we compile these conditions into a conjecture, which states what we believe such an analogous theorem for $I_X$ could look like.

\section{Generalising Group Results}

Since the symmetric group \Sym{X} is a sub(semi)group of the symmetric inverse monoid $I_X$ and we have Lemma \ref{SymX=IX}, which states that $\Sym{X} \bsequal[I] I_X$, a natural approach to classifying the Bergman-Shelah equivalence classes of $I_X$ would be to check whether any results from the \Sym{X} case follow through for $I_X$.
One of the results needed to prove the main theorem by Bergman and Shelah in \cite{Bergman_2006} is that the `relative size' of subgroups of \Sym{X} is `limited' by the size of their orbits.
We will make this more clear later, but first we need to introduce some concepts.

\begin{definition}[Partwise stabiliser]
    The \emph{partwise stabiliser} of a partition $\partition$ of a set $X$, denoted $\PStab{\partition}$, is the group of permutations, which maps each part of the partition to itself.
    \begin{equation*}
        \PStab{\partition} = \set{f \in \Sym{X} \given (\forall \Sigma \in \partition)(\Sigma f = \Sigma)}
    \end{equation*}
\end{definition}

The partwise stabiliser is to the stabiliser of a partition as the pointwise stabiliser is to the setwise stabiliser, and as such can be seen as a generalisation of this.
The following are the main types of partitions we will consider throughout this chapter.

\begin{definition}[Uniform partition]
    Let $\partition$ be a partition of an infinite set $X$ and $\kappa \leq \card{X}$ a cardinal. We shall call $\partition$ a \emph{$\kappa$-uniform partition} of $X$ if all parts of $\partition$ have cardinality $\kappa$.
\end{definition}

\begin{definition}[Unbounded partition]
    Let $\partition$ be a partition of an infinite set $X$ and $\kappa \leq \card{X}$ a limit cardinal. We shall call $\partition$ a \emph{$\kappa$-unbounded partition} of $X$ if all parts of $\partition$ have cardinalities strictly smaller than $\kappa$, but there is no common bound $\lambda < \kappa$ on those cardinalities.
\end{definition}

\begin{definition}[Bounded partition]
    Let $\partition$ be a partition of an infinite set $X$ and $\kappa \leq \card{X}$ a cardinal. We shall call $\partition$ a \emph{$\kappa$-bounded partition} of $X$ if there exists a cardinal $\lambda < \kappa$ such that all parts of $\partition$ have cardinalities at most $\lambda$.
\end{definition}

Partwise stabilisers of partitions with any of the properties listed above make up a set of particularly `nice' subgroups of $\Sym{X}$ and consequently $I_X$.
If the parts are small enough, then the partwise stabilisers are also `nice' topologically.

\begin{lemma} \label{lem:PStab_finite_closed}
    Let $X$ be a set and $\partition$ a partition of $X$.
    If all parts in $\partition$ are finite, then $\PStab{\partition}$ is closed in every $T_1$ shift-continuous topology on $I_X$.
\end{lemma}

\begin{proof}
    Theorem \ref{I1_theorem} states that $\topI_1$ is minimal among all $T_1$ shift-continuous topologies on $I_X$, hence we only need to check that this holds for $\topI_1$.
    The defining feature of the elements of $\PStab{\partition}$ is that they maps every point in $X$ to the part which the point itself belongs to in $\partition$. 
    We can therefore conclude that for any chart $f \in I_X \setminus \PStab{\partition}$ there must exist a point $x$, which belongs to some part $A$ in \partition, such that either $x \notin \dom{f}$ or $xf \in X \setminus A$.
    We thus get that the following set;
    \begin{equation*}
        U = \bigcap_{y \in A} V_{x,y} = \bigcap_{y \in A} \set{h \in I_X \given (x,y) \notin h}
    \end{equation*}
    is a basic open neighbourhood of $f$ in $\topI_1$, which has empty intersection with $\PStab{\partition}$.
\end{proof}

In \cite{Bergman_2006} Bergman and Shelah prove the following theorem, which states how the Bergman-Shelah ordering of subgroups of \Sym{X} depends on their orbit sizes.
The theorem has here been restated to only include the parts relevant to this thesis.

\begin{theorem}[{\cite[Theorem 6]{Bergman_2006}}] \label{BS_Theorem6}
    Let $X$ be an infinite set, $\kappa$ an infinite regular cardinal such that $\kappa \leq \card{X}$, and $A,B,C$ partitions of $X$ with the following properties:
    \begin{enumerate}[\normalfont ~(a)]
        \item Some part in $A$ has cardinality at least $\kappa$.
        \label{BS_Theorem6:A}

        \item $B$ is $\kappa$-unbounded (this is only possible when $\kappa$ is a limit cardinal).
        \label{BS_Theorem6:B}

        \item $C$ is $\kappa$-bounded.
        \label{BS_Theorem6:C}
    \end{enumerate}
    Then $\PStab{A} \not\bsleq[S] \PStab{B} \not\bsleq[S] \PStab{C}$ and $\PStab{A} \not\bsleq[S] \PStab{C}$.
\end{theorem}

It follows from Theorem \ref{BS_Theorem6} that given any subgroup $G \subseteq \Sym{X}$, if the partition of $X$ into the orbits of $G$, which we denote by $\orbit{G}$, is $\kappa$-unbounded, then $G \not\bsgeq[S] \PStab{A}$ where $A$ is any partition of $X$ with at least one part of size $\kappa$.
This holds true, since $G$ is a subset of $\PStab{\orbit{G}}$ and $\orbit{G}$ is $\kappa$-unbounded.
Similarly, if $\orbit{G}$ is $\kappa$-bounded, then $G \not\bsgeq[S] \PStab{B}$ where $B$ is any $\kappa$-unbounded partition of $X$ or contains a part of cardinality at least $\kappa$.
We will now prove a series of lemmas aiming to generalise this result to encompass partial permutation semigroups (and singular cardinals).

\begin{lemma} \label{lem:iso_partitions_bsequal}
    Let $X$ be a set, $\partition$ and $\mathcal{Q}$ isomorphic partitions of subsets of $X$, and $S$ a subsemigroup of $I_X$ such that every element of $\PStab{\mathcal{Q}}$ extends to an element of $S$.
    Then $\PStab{\partition} \bsequal[I] \PStab{\mathcal{Q}} \bsleq[I] S$.
\end{lemma}

\begin{proof}
    This lemma has two parts to it.

    The first part is the statement $\PStab{\partition} \bsequal[I] \PStab{\mathcal{Q}}$, which follows from $\partition$ and $\mathcal{Q}$ being isomorphic as partitions.
    This means that there exists a chart $f \in I_X$ such that $f \PStab{\partition} \inv{f} = \PStab{\mathcal{Q}}$ and $\inv{f} \PStab{\mathcal{Q}} f = \PStab{\partition}$.
    
    The second part is the statement $\PStab{\mathcal{Q}} \bsleq[I] S$, which follows from the fact that every element of $\PStab{\mathcal{Q}}$ extends to an element of $S$.
    Let $e \in I_X$ be the partial identity on the subset of $X$ on which $\mathcal{Q}$ is a partition, then $\PStab{\mathcal{Q}} \subseteq eS$.
\end{proof}

\begin{lemma} \label{lem:orbits_leq->bsleq_uniform}
    Let $X$ be an infinite set, $\kappa$ a cardinal with $2 \leq \kappa \leq \card{X}$, $\partition_\kappa$ a $\kappa$-uniform partition of $X$, and $M$ an inverse subsemigroup of $I_X$.
    If all orbits of $M$ have cardinalities no greater than $\kappa$, then $M \bsleq[I] \PStab{\partition_\kappa}$.
\end{lemma}

\begin{proof}
    We split the proof into three cases depending on $\kappa$.

    \emph{Case 1: $\kappa = \card{X}$.} 
    In this case $\PStab{\partition_\kappa} \bsequal[I] I_X$ due to Lemma \ref{SymYSub=IX}, since $\PStab{\partition_\kappa}$ contains the symmetric group on each part in $\partition_\kappa$.
    The result is then immediate, since for all $U \subseteq I_X$ it follows that $U \bsleq[I] I_X$.

    \emph{Case 2: $\kappa$ is infinite, but $\kappa < \card{X}$.}
    This necessitates that $\card{\partition_\kappa} = \card{X}$, since $\kappa$ is a bound in \card{X} on the sizes of the parts of $\partition_\kappa$.
    So let $f \in I_X$ be a chart which induces an injection from \orbit{M} into $\partition_\kappa$ such that for each $\Sigma \in \partition_\kappa$ there exists a subset $\Gamma \subseteq \Sigma$ with $\card{\Gamma} = \kappa$ and $\Gamma \cap \im{f} = \emptyset$.
    We will then show that $M \subseteq \genset{\PStab{\partition_\kappa}, f, \inv{f}}$.
    Let $h \in M$ be arbitrary and $a \in \PStab{\partition_\kappa}$ a permutation with the following properties:
    \begin{enumerate}[(i)]
        \item If $x \in \dom{h}$, then $xfa\inv{f} = xh$.
        This is always possible, since the orbit $xM$ is mapped into a part in $\partition_\kappa$ under $f$.

        \item If $y \in X \setminus \dom{h}$, then $yfa \notin \dom{\inv{f}} = \im{f}$. 
        This is always possible, since each part in $\partition_\kappa$ contains a subset $\Gamma$ of cardinality $\kappa$ which is not hit by anything under $f$.
    \end{enumerate}
    Then $h = fa\inv{f}$, and since $h$ was arbitrary it follows that $M \subseteq \genset{\PStab{\partition_\kappa}, f, \inv{f}}$.

    \emph{Case 3: $\kappa$ is finite.}
    Again we know that $\card{\partition_\kappa} = \card{X}$.
    So let $f,g \in I_X$ be charts such that $f$ induces an injection from $X$ into $\partition_\kappa$ and $g$ induces an injection from \orbit{M} into $\partition_\kappa$.
    We will then show that $M \subseteq \genset{\PStab{\partition_\kappa}, f, \inv{f}, g, \inv{g}}$.
    Let $h \in M$ be arbitrary and $a,b \in \PStab{\partition_\kappa}$ be permutations with the following properties:
    \begin{enumerate}[(i)]
        \item If $x \in \dom{h}$, then $xfa\inv{f} = x$ and $xgb\inv{g} = xh$.
        This is always possible, since you can always fix any subset of $X$ under actions of $\PStab{\partition_\kappa}$ and the orbit $xM$ is mapped into a part in $\partition_\kappa$ under $f$.

        \item If $y \in X \setminus \dom{h}$, then $yfa \notin \dom{\inv{f}} = \im{f}$. This is always possible, since each part in $\partition_\kappa$ has cardinality at least 2.
    \end{enumerate}
    Then $h = fa\inv{f}gb\inv{g}$, since each orbit under $M$ is finite.
    Because $h$ was arbitrary, it then follows that $M \subseteq \genset{\PStab{\partition_\kappa}, f, \inv{f}, g, \inv{g}}$.
\end{proof}

\begin{lemma} \label{lem:cones_less->not_bsgeq_uniform}
    Let $X$ be an infinite set, $\kappa \leq \card{X}$ an infinite cardinal, $\partition$ a partition of $X$ such that some part in $\partition$ has cardinality at least $\kappa$, and $S$ a subsemigroup of $I_X$.
    If all cones of $S$ have cardinalities strictly smaller than $\kappa$, then $S \not\bsgeq[I] \PStab{\partition}$.
\end{lemma}

\begin{proof}
    We will show this using a diagonal argument.
    So we assume the contrary, that is, we assume there exists a finite subset $U \subseteq I_X$ such that $\PStab{\partition} \subseteq \genset{S,U}$.
    We can assume without loss of generality that this subset $U$ contains the identity element.
    We then define a chain $(S_n)_{n\in\omega}$ of subsets in $I_X$ such that for each $n \in \omega$,
    \begin{equation} \label{eq:Sn_chain_genset}
        S_n = \set{s_0u_0 \dots s_{n-1}u_{n-1} \in I_X \given (\forall i \in n) (s_i \in S^1 ~\land~ u_i \in U)}.
    \end{equation}
    That is, $S_n$ is the set of all words of length $2n$ using alternating elements from $S^1$ and $U$.
    It follows that $\genset{S,U} = \bigcup_{n \in \omega} S_n$.
    We note however, since $U$ is finite, that all cones of $S_n$ are still strictly smaller than $\kappa$ for all $n \in \omega$ (even though $\genset{S,U}$ might very well have orbits of size $\kappa$).
    We will use this to construct a permutation $f \in \PStab{\partition}$ such that $f \notin \genset{S,U}$.
    First, we pick a part $\Sigma \in \partition$ with $\card{\Sigma} \geq \kappa$ and partition $\Sigma$ into three parts $\set{A,B,C}$ such that $\card{A} = \aleph_0$ and $\card{B} = \card{C} = \card{\Sigma}$.
    Next, let $(a_n)_{n\in\omega}$ be an enumeration of the elements of $A$.
    We will now construct the permutation $f$, starting with defining how it acts on $A$.
    For each $n \in \omega$, we choose the image of $a_n$ under $f$ such that $a_n f \in B \setminus \paa{a_n S_n \cup \set{a_i \given i \in n} f}$.
    This is always possible, since $\card{a_n S_n} < \kappa$ as mentioned above and $\set{a_i \given i \in n}$ is finite for each $n \in \omega$.
    It already follows from the above that $f \notin \genset{S,U}$, since for each $n \in \omega$ there exists a point $x \in X$ (namely $a_n$) such that $xf \neq xs$ for all $s \in S_n$, meaning that $f \notin \bigcup_{n\in\omega} S_n = \genset{S,U}$.
    All that remains is to define $f$ on the remaining points of $X$ such that $f \in \PStab{\partition}$.
    On the remaining elements of $\Sigma$ simply let $f$ act as a permutation, which is always possible due to the thus far unused subset $C$ (e.g. map $B$ to a moiety of $C$ and map $C$ onto $\Sigma \setminus (B \cup \set{a_n \given n \in \omega}) f$).
    It should then be clear that $\Sigma f = \Sigma$.
    For the final step, simply let $f$ act in any way on the remaining points of $X$ such that it partwise stabilises the partition $\partition$ (e.g. let $f$ fix the remaining points outside $\Sigma$).
    Then $f \in \PStab{\partition}$ and $f \notin \genset{S,U}$, as required.
    Since the above construction is possible for any finite subset $U \subseteq I_X$ it follows that $S \not\bsgeq[I] \PStab{\partition}$.
\end{proof}

\begin{corollary} \label{cor:obits_less->bsless_uniform}
    Let $X$ be an infinite set, $\kappa \leq \card{X}$ an infinite cardinal, $\partition_\kappa$ a $\kappa$-uniform partition of $X$, and $M$ an inverse subsemigroup of $I_X$.
    If all orbits of $M$ have cardinalities strictly smaller than $\kappa$, then $M \bsless[I] \PStab{\partition_\kappa} \bsleq[I] I_X$.
\end{corollary}

\begin{proof}
    This follows from Lemmas \ref{lem:orbits_leq->bsleq_uniform} and \ref{lem:cones_less->not_bsgeq_uniform}, which state that $M \bsleq[I] \PStab{\partition_\kappa}$ and $M \not\bsgeq[I] \PStab{\partition_\kappa}$.
    For the final part it follows from the fact that $\PStab{\partition_\kappa} \subseteq I_X$ that $\PStab{\partition_\kappa} \bsleq[I] I_X$ and hence $M \bsless[I] I_X$ by transitivity.
\end{proof}

\begin{lemma} \label{lem:orbits_less->bsleq_unbounded}
    Let $X$ be an infinite set, $\kappa \leq \card{X}$ a limit cardinal, $\partition_\kappa$ a $\kappa$-unbounded partition of $X$, and $M$ an inverse subsemigroup of $I_X$.
    If all orbits of $M$ have cardinalities strictly smaller than $\kappa$, then $M \bsleq[I] \PStab{\partition_\kappa}$.
\end{lemma}

\begin{proof}
    We will first show that we can assume without loss of generality that $\card{\partition_\kappa} = \card{X}$.
    Notably, this is only a potential issue when \card{X} is singular.
    Let $\mathcal{Q}_\kappa$ be a refinement of $\partition_\kappa$ such that for each infinite $\Sigma \in \partition_\kappa$, $\mathcal{Q}_\kappa$ contains $\card{\Sigma}$ many moieties of $\Sigma$.
    It then follows from $\bigcup \partition_\kappa = X$ that $\card{\mathcal{Q}_\kappa} = \card{X}$ and $\mathcal{Q}_\kappa$ is $\kappa$-unbounded because $\partition_\kappa$ is (since for each part in $\partition_\kappa$ there exists a part in $\mathcal{Q}_\kappa$ with the same cardinality).
    Since $\mathcal{Q}_\kappa$ is a refinement of $\partition_\kappa$ it follows that $\PStab{\mathcal{Q}_\kappa} \subseteq \PStab{\partition_\kappa}$, and hence $\PStab{\mathcal{Q}_\kappa} \bsleq[I] \PStab{\partition_\kappa}$.

    So, assuming $\card{\partition_\kappa} = \card{X}$, we pick a chart $f \in I_X$ which which induces an injection from \orbit{M} into $\partition_\kappa$ such that each orbit is mapped to a strictly larger part in $\partition_\kappa$ (this is always possible, since $\partition_\kappa$ is $\kappa$-unbounded).
    If an orbit is finite, let $f$ map it to a part that is at least twice the size of the orbit.
    We will then show that $M \subseteq \genset{\PStab{\partition_\kappa}, f, \inv{f}}$.
    Let $h \in M$ be arbitrary and $a \in \PStab{\partition_\kappa}$ a permutation with the following properties:
    \begin{enumerate}[(i)]
        \item If $x \in \dom{h}$, then $xfa\inv{f} = xh$.
        This is always possible, since the orbit $xM$ is mapped into a part in $\partition_\kappa$ under $f$.

        \item If $y \in X \setminus \dom{h}$, then $xfa \notin \dom{\inv{f}} = \im{f}$.
        This is always possible, since $f$ maps each orbit into a strictly larger part in $\partition_\kappa$.
    \end{enumerate}
    Then $h = fa\inv{f}$, and since $h$ was arbitrary it follows that $M \subseteq \genset{\PStab{\partition_\kappa}, f, \inv{f}}$.
\end{proof}

\begin{lemma} \label{lem:cones_bounded->not_bsgeq_unbounded}
    Let $X$ be an infinite set, $\kappa \leq \card{X}$ a limit cardinal, $\partition_\kappa$ a $\kappa$-unbounded partition of $X$, and $S$ a subsemigroup of $I_X$.
    If there exists a cardinal $\lambda < \kappa$ such that all cones of $S$ have cardinality no greater than $\lambda$, then $S \not\bsgeq[I] \PStab{\partition_\kappa}$.
\end{lemma}

\begin{proof}
    We will use the same setup as in the proof of Lemma \ref{lem:cones_less->not_bsgeq_uniform}.
    So, let $U \subseteq I_X$ be any finite subset (containing the identity element) and $(S_n)_{n\in\omega}$ a chain of subsets of $I_X$ as defined in equation \eqref{eq:Sn_chain_genset}.
    It again follows that $\genset{S,U} = \bigcup_{n\in\omega} S_n$ and that for all $n \in \omega$ the cones of $S_n$ have their sizes bounded by some $\lambda_n < \kappa$ (but there is not necessarily a bound in $\kappa$ on how large $\lambda_n$ can grow as $n$ increases).
    Let $(\lambda_n)_{n\in\omega}$ be a sequence of such cardinals in $\kappa$, where each $\lambda_n$ is an upper bound on the sizes of the cones under $S_n$.
    We now wish to construct a chart $f$ such that $f$ is in $\PStab{\partition_\kappa}$ but not in $\genset{S,U}$.
    Let $(\Sigma_n)_{n\in\omega} \in {\partition_\kappa}^\omega$ and $(x_n)_{n\in\omega} \in X^\omega$ be sequences such that for all $n \in \omega$, $\card{\Sigma_n} > \lambda_n$ (this is always possible since $\partition_\kappa$ is $\kappa$-unbounded) and $x_n \in \Sigma_n$.
    We then construct the desired chart $f$ by first defining how it acts on $(x_n)_{n\in\omega}$.
    For each $n \in \omega$, let $x_n f$ be in $\Sigma_n$ such that $x_n f \neq x_n s$ for all $s \in S_n$ (this is always possible since $\card{\Sigma_n} > \lambda_n$).
    On the remaining points in $X$, simply let $f$ act a permutation which partwise stabilises $\partition_\kappa$.
    Then $f \in \PStab{\partition_\kappa}$ but $f \notin \bigcup_{n\in\omega} S_n = \genset{S,U}$.
    Hence, $S \not\bsgeq[I] \PStab{\partition_\kappa}$.
\end{proof}

\begin{corollary} \label{cor:obits_bounded->bsless_unbounded}
    Let $X$ be an infinite set, $\kappa \leq \card{X}$ a limit cardinal, $\partition_\kappa$ a $\kappa$-unbounded partition of $X$, and $M$ an inverse subsemigroup of $I_X$.
    If there exists a cardinal $\lambda < \kappa$ such that all orbits of $M$ have cardinality no greater than $\lambda$, then $M \bsless[I] \PStab{\partition_\kappa} \bsless[I] I_X$.
\end{corollary}

\begin{proof}
    This follows from Lemmas \ref{lem:orbits_less->bsleq_unbounded} and \ref{lem:cones_bounded->not_bsgeq_unbounded}, which state that $M \bsleq[I] \PStab{\partition_\kappa}$ and $M \not\bsgeq[I] \PStab{\partition_\kappa}$.
    The final part follows from Corollary \ref{cor:obits_less->bsless_uniform} since all orbits of $\PStab{\partition_\kappa}$ have cardinalities strictly less than $\kappa \leq \card{X}$.
\end{proof}

\begin{lemma} \label{lem:orbits_bounded->bsleq_bounded}
    Let $X$ be an infinite set, $n \geq 2$ a natural number, $\partition_n$ an $n$-uniform partition of $X$, and $M$ an inverse subsemigroup of $I_X$.
    If there exists a natural number $m \in \omega$ such that all orbits of $M$ have cardinalities no greater than $m$, then $M \bsleq[I] \PStab{\partition_n}$.
\end{lemma}

\begin{proof}
    We first prove an intermediate result.
    Let $\partition_{2n}$ be a $2n$-uniform partition of $X$.
    We will then prove that $\PStab{\partition_{2n}} \bsleq[I] \PStab{\partition_{n}}$.
    We assume without loss of generality that for all $\Sigma \in \partition_{2n}$ there exists a pair of subsets $\Gamma_0, \Gamma_1 \in \partition_n$ such that $\Sigma = \Gamma_0 \cup \Gamma_1$ (we can assume this, since all $2n$-uniform partitions of $X$ are isomorphic).
    Let $g \in \Sym{X}$ be an involution which for each $\Sigma = \Gamma_0 \cup \Gamma_1 \in \partition_{2n}$ (with $\Gamma_0, \Gamma_1 \in \partition_n$) swaps exactly one element in $\Gamma_0$ with an element in $\Gamma_1$ while fixing everything else.
    We then claim that $\PStab{\partition_{2n}} \subseteq \genset{\PStab{\partition_n}, g}$.
    To see this, let us consider what happens on a single part $\Sigma = \Gamma_0 \cup \Gamma_1$ in $\partition_{2n}$ under a permutation $h \in \PStab{\partition_{2n}}$.
    Here an important quantity to consider is $\card{\Gamma_0 h \cap \Gamma_1}$, which is the number of elements $h$ carries from $\Gamma_0$ to $\Gamma_1$ and vice versa.
    Given a permutation $f \in \PStab{\partition_{2n}}$ let us refer to the quantity $\card{\Gamma_0 f \cap \Gamma_1}$ as the `crossing number' of $f$.
    Let us further denote $\PStab{\partition_n}$ by $P$ and consider an as of yet unfixed sequence of permutations $(f_i)_{i \in n}$ with $f_i \in PgP$ for all $i \in n$.
    It should then be clear that $\card{\Gamma_0 f_i \cap \Gamma_1} = 1$ for all $i \in n$ no matter how we end up choosing  $(f_i)_{i \in n}$.
    However, if we consider the composite $f_0f_1$, then the crossing number of this chart can assume three possible values, namely 0, 1, or 2 depending on how $f_1$ acts on $\Sigma$.
    Let us go through how each of these cases may come about.
    To achieve a crossing number of 0, we can let $f_1 = \inv{f_0}$ so that $f_1$ undoes what $f_0$ did.
    To achieve a crossing number of 1, we can let $f_1$ swap an element in $\Gamma_0$ with the element that $f_0$ mapped from $\Gamma_0$ to $\Gamma_1$ (then the crossing number remains unchanged).
    To achieve a crossing number of 2, we can let $f_1$ swap any two elements from $\Gamma_0$ and $\Gamma_1$ which have not previously been swapped over.
    Using the same argument as above, we see that the composite $f_0f_1f_2$ (if such exists) will have have crossing number which differs from that of $f_0f_1$ by a value of at most 1 (that is, $\card{\Gamma_0 f_0f_1 \cap \Gamma_1} - \card{\Gamma_0 f_0f_1f_2 \cap \Gamma_1} \in \set{-1, 0, 1}$).
    So by induction on $n$ we see that the crossing number of any composite $f_0 \dots f_{i-1} f_i$ differs by a value of at most 1 from the crossing number of $f_0 \dots f_{i-1}$ for all $i \in n$, and we can always choose $f_i$ such that this difference is a desired value.
    It follows that we can choose the sequence $(f_i)_{i \in n}$ such that the composite $f_0 \dots f_{n-1}$ has the same crossing number as the permutation $h$.
    But since any such composite consists of pre- and post-multiplication by elements of $P$, it follows that the composite can be chosen to act exactly as $h$ on $\Sigma$.
    Furthermore, since each part in $\partition_{2n}$ is acted on independently, we can choose the sequence $(f_i)_{i \in n}$ such that $f_0 \dots f_{n-1} = h$.
    Hence $\PStab{\partition_{2n}} \subseteq \genset{\PStab{\partition_n}, g}$, as we wanted to show.

    It then follows from transitivity of $\bsleq$ and induction on the proof above, that $\PStab{\partition_{2kn}} \bsleq[I] \PStab{\partition_{n}}$ for any $k \in \omega$ and $2kn$-uniform partition $\partition_{2kn}$.
    So choose $k \in \omega$ such that $2kn > m$, it then follows from Lemma \ref{lem:orbits_leq->bsleq_uniform} that $M \bsleq[I] \PStab{\partition_{2kn}} \bsleq[I] \PStab{\partition_n}$.
\end{proof}

\begin{corollary} \label{cor:n_partition=m_partition}
    Let $X$ be an infinite set, $n,m \geq 2$ natural numbers, $\partition_n$ an $n$-uniform partition, and $\partition_m$ an $m$-uniform partition of $X$.
    Then $\PStab{\partition_n} \bsequal[I] \PStab{\partition_m}$.
\end{corollary}

\begin{proof}
    It follows from Lemma \ref{lem:orbits_bounded->bsleq_bounded} that $\PStab{\partition_n} \bsleq[I] \PStab{\partition_m}$ and $\PStab{\partition_n} \bsgeq[I] \PStab{\partition_m}$.
\end{proof}

\begin{lemma} \label{lem:countable->bsequal}
    Let $X$ be an infinite set and $U$ a subset of $I_X$.
    Then $\genset{U} \bsequal[I] \emptyset$ if and only if $U$ is countable.
\end{lemma}

\begin{proof}
    The generated semigroup $\genset{U}$ is the union of all finite products of elements from $U$.
    So if $U$ is countable, then $\genset{U}$ can be at most countable.
    Hence, $\genset{U} \subseteq \genset{\emptyset,U}$.
    If $U$ is not countable, then no countable set can generate $U$ and hence $\emptyset \not\bsgeq[I] \genset{U}$.
\end{proof}

\begin{corollary} \label{cor:countable->bsleq}
    Let $X$ be an infinite set, $S$ a subsemigroup of $I_X$, and $U$ a countable subset of $I_X$.
    Then $\genset{U} \bsleq[I] S$.
\end{corollary}

\begin{proof}
    It follows from Lemma \ref{lem:countable->bsequal} that $\genset{U} \bsequal[I] \emptyset$ and $\emptyset \subseteq S$.
\end{proof}

We combine the above results to give us the following theorem.

\begin{theorem} \label{thm:partition_eqclasses_distinct}
    Let $X$ be an infinite set, $\kappa \leq \card{X}$ an infinite cardinal, and $A,B,C$ partitions of $X$ with the following properties:
    \begin{enumerate}[\normalfont ~(a)]
        \item $A$ is $\kappa$-uniform.

        \item $B$ is $\kappa$-unbounded (this is only possible when $\kappa$ is a limit cardinal).

        \item $C$ is $\kappa$-bounded.
    \end{enumerate}
    Then $\PStab{A} \bsgreat[I] \PStab{B} \bsgreat[I] \PStab{C}$.
\end{theorem}

\begin{proof}
    It follows from Lemmas \ref{lem:orbits_leq->bsleq_uniform} and \ref{lem:orbits_less->bsleq_unbounded} that $\PStab{A} \bsgeq[I] \PStab{B} \bsgeq[I] \PStab{C}$, while Lemmas \ref{lem:cones_less->not_bsgeq_uniform} and \ref{lem:cones_bounded->not_bsgeq_unbounded} show that $\PStab{A} \not\bsleq[I] \PStab{B} \not\bsleq[I] \PStab{C}$.
\end{proof}

From here on out we will mainly be concerned with the case where $X$ is countably infinite.
So for the remainder of this chapter fix $\partition_2$ and $\partition_\omega$ to be your favourite 2-uniform and $\aleph_0$-unbounded partitions respectively of the countably infinite set $X$.
We can then restate the main result of Bergman and Shelah (Theorem \ref{BS_mainTheorem}) in the following way.

\begin{theorem} \label{thm:BS_alt}
    Let $X$ be a countably infinite set and $G$ a closed subgroup of $\Sym{X}$ in the pointwise topology.
    Then exactly one of the following holds:
    \begin{enumerate}[\normalfont(i)]
        \item $G \bsequal[I] \Sym{X}$;
        \item $G \bsequal[I] \PStab{\partition_\omega}$;
        \item $G \bsequal[I] \PStab{\partition_2}$; or
        \item $G \bsequal[I] \set{1}$.
    \end{enumerate}
    Moreover, $\set{1} \bsless[I] \PStab{\partition_2} \bsless[I] \PStab{\partition_\omega} \bsless[I] I_X$.
\end{theorem}

This restatement of the theorem does three things.
First, it generalises the result by allowing elements of $I_X$ for generating sets.
Second, it uses specific groups as representatives of each of the four equivalence classes.
Third, it reduces the scope of the result by not describing any properties of the group $G$ in terms of stabiliser subgroups.
The reason why the Theorem \ref{BS_mainTheorem} version of the Bergman-Shelah theorem can be stated in terms of stabiliser subgroups is mainly due to Lemma \ref{BS_lemma2}, which states that for any subgroup $G \subseteq \Sym{X}$ and finite subset $\Gamma \subseteq X$, $\pointStab{G}{\Gamma} \bsequal[S] G$.
This result does not generalise to subsemigroups, or even inverse subsemigroups, of $I_X$ (for a counterexample, see Proposition \ref{prop:stab_neq_S(T)}) and as such we will not be considering stabiliser subsemigroups for classification of the Bergman-Shelah equivalence classes of $I_X$.

There is a final set of results that we wish to carry over from the paper by Bergman and Shelah.
The following three lemmas are analogues of Lemmas 10, 12, and 14 in \cite{Bergman_2006}.
The lemmas here are in fact so similar to their group counterparts that the proof carry directly over with only minor adjustments.
As such, we have only stated the proofs of the following three lemmas in terms of what would have to be changed from the proofs given in \cite{Bergman_2006}.

\begin{lemma} \label{BS_lemma10}
    Let $X$ be a countably infinite set and $S$ a subsemigroup of $I_X$.
    Suppose there exist a sequence $(\alpha_i)_{i\in\omega} \in X^\omega$ of distinct elements and a sequence of non-empty subsets $D_i \subseteq X^i$ $(i\in\omega)$ such that
    \begin{enumerate}[\normalfont(i)]
        \item for each $i \in \omega$ and each $(\beta_0,\dots,\beta_i) \in D_{i+1}$, we have $(\beta_0,\dots,\beta_{i-1}) \in D_i$;
        \label{BS_lemma10:(i)}
        \item for each $i \in \omega$ and each $(\beta_0,\dots,\beta_{i-1}) \in D_i$, there exist infinitely many elements $\beta \in X$ such that $(\beta_0,\dots,\beta_{i-1}, \beta) \in D_{i+1}$; and
        \label{BS_lemma10:(ii)}
        \item if $(\beta_i)_{i\in\omega} \in X^\omega$ has the property that $(\beta_0,\dots,\beta_{i-1}) \in D_i$ for each $i \in \omega$, then there exists $g,g' \in S$ such that $(\alpha_ig)_{i\in\omega} = (\beta_i)_{i\in\omega}$ and $(\alpha_i)_{i\in\omega} = (\beta_ig')_{i\in\omega}$.
        \label{BS_lemma10:(iii)}
    \end{enumerate}
    
    Then $S \bsequal[I] I_X$.
\end{lemma}

\begin{proof}
    The proof of this lemma can be obtained from the proof of \cite[Lemma 10]{Bergman_2006} by making the following adjustments throughout:
    \begin{itemize}
        \item Replace $\Omega$ by $X$.
        \item Replace $S = \Sym{\Omega}$ by $I_X$.
        \item Replace the group $G$ by the semigroup $S$.
        \item Replace the word `group' with `semigroup'.
        \item Replace the inverses $\inv{g}$ and $\inv{h}$ by $g'$ and $h'$ respectively as in accordance with condition \eqref{BS_lemma10:(iii)}.
        \item Replace statement (8) by Lemma \ref{SymYSub=IX}.
        \item Replace the permutation $x \in \Sym{\Omega}$ by charts $f, \inv{f} \in I_X$.
    \end{itemize}
\end{proof}

\begin{lemma} \label{BS_lemma12}
    Let $X$ be a countably infinite set and $S$ a subsemigroup of $I_X$.
    Suppose there exist a sequence $(\alpha_i)_{i\in\omega} \in X^\omega$ of distinct elements, an unbounded	sequence of positive integers $(N_i)_{i\in\omega}$, and a sequence of non-empty subsets $D_i \subseteq X^i$ $(i\in\omega)$ such that
    \begin{enumerate}[\normalfont(i)]
        \item for each $i \in \omega$ and each $(\beta_0,\dots,\beta_i) \in D_{i+1}$, we have $(\beta_0,\dots,\beta_{i-1}) \in D_i$;
        \label{BS_lemma12:(i)}
        \item for each $i \in \omega$ and each $(\beta_0,\dots,\beta_{i-1}) \in D_i$, there exist at least $N_i$ elements $\beta \in X$ such that $(\beta_0,\dots,\beta_{i-1}, \beta) \in D_{i+1}$; and
        \label{BS_lemma12:(ii)}
        \item if $(\beta_i)_{i\in\omega} \in X^\omega$ has the property that $(\beta_0,\dots,\beta_{i-1}) \in D_i$ for each $i \in \omega$, then there exists $g,g' \in S$ such that $(\alpha_ig)_{i\in\omega} = (\beta_i)_{i\in\omega}$ and $(\alpha_i)_{i\in\omega} = (\beta_ig')_{i\in\omega}$.
        \label{BS_lemma12:(iii)}
    \end{enumerate}
    
    Then $S \bsgeq[I] \PStab{\partition_\omega}$.
\end{lemma}

\begin{proof}
    The proof of this lemma can be obtained from the proof of \cite[Lemma 12]{Bergman_2006} by making the following adjustments throughout:
    \begin{itemize}
        \item Replace $\Omega$ by $X$.
        \item Replace $S = \Sym{\Omega}$ by $I_X$.
        \item Replace the group $G$ by the semigroup $S$.
        \item Replace the notation $\Sym{\Delta}_{(A)}$ for $\PStab{A}$.
        \item Replace the inverse $\inv{g}$ by $g'$ as in accordance with condition \eqref{BS_lemma12:(iii)}.
        \item Replace statement (27) by Lemma \ref{lem:iso_partitions_bsequal}.
        \item Replace the partition $A'$ by $\partition_\omega$.
        \item Replace $\bsgeq$ by $\bsgeq[I]$.
    \end{itemize}
\end{proof}

\begin{lemma} \label{BS_lemma14}
    Let $X$ be a countably infinite set and $S$ a subsemigroup of $I_X$.
    Suppose there exist two disjoint sequences of distinct elements $(\alpha_i)_{i \in \omega}, (\beta_i)_{i \in \omega} \in X^\omega$ such that for every sequence $(\gamma_i)_{i \in \omega} \in \prod_{i \in \omega} \set{\alpha_i, \beta_i} \subseteq X^\omega$, there exist $g,g' \in S$ such that $(\alpha_i g)_{i \in \omega} = (\gamma_i)_{i \in \omega}$ and $(\alpha_i)_{i \in \omega} = (\gamma_i g')_{i \in \omega}$.

    Then $S \bsgeq[I] \PStab{\partition_2}$.
\end{lemma}

\begin{proof}
    The proof of this lemma can be obtained from the proof of \cite[Lemma 14]{Bergman_2006} by making the following adjustments throughout:
    \begin{itemize}
        \item Replace $\Omega$ by $X$.
        \item Replace $S = \Sym{\Omega}$ by $I_X$.
        \item Replace the group $G$ with the semigroup $S$.
        \item Replace the notation $\Sym{\Delta}_{(A)}$ for $\PStab{A}$.
        \item Replace the inverse $\inv{g}$ by $g'$ as in accordance with the hypothesis of Lemma \ref{BS_lemma14}.
        \item Replace statement (38) by Lemma \ref{lem:iso_partitions_bsequal}.
        \item Replace the partition $A_0$ by $\partition_2$.
        \item Replace $\bsleq$ by $\bsleq[I]$.
    \end{itemize}
\end{proof}

The group counterparts of Lemmas \ref{BS_lemma10}, \ref{BS_lemma12}, and \ref{BS_lemma14} play a major role in proving the main theorem of \cite{Bergman_2006}.
This is due to results by Bergman and Shelah, which show that for closed subgroups of \Sym{X} in the pointwise topology satisfying certain conditions on their stabiliser subgroups, the converse of these lemmas hold.

\section{Rooted Trees}

In this section we will study a specific class of inverse subsemigroups of $I_X$ and the Bergman-Shelah preorder on these subsemigroups.
The semigroups in question are defined on rooted trees and will provide important counterexamples to various claims one could make about the subsemigroup structure of $I_X$.

\begin{definition}[Rooted tree] \label{def:rooted_tree}
    A \emph{rooted tree} $T$ is a connected graph on a non-empty vertex set $X$ such that $T$ contains no cycles (that is, for any two vertices $x,y \in X$ there exists a unique path from $x$ to $y$).
    In addition, a single element $r \in X$ is designated as the \emph{root} of the tree.
\end{definition}

\begin{figure}[h]
    \centering
    \includegraphics[height=0.3\textheight]{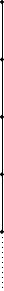}
    \hspace{1.5cm}
    \includegraphics[width=0.6\linewidth,height=0.3\textheight]{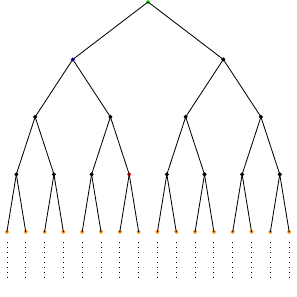}
    \caption{The unary tree \regtree{1} and the binary rooted tree \regtree{2} on countably infinite vertex sets. The \textcolor{green}{green} vertex at the top is the root, the \textcolor{blue}{blue} vertex just below is a child of the root (in turn the root is the parent of the \textcolor{blue}{blue} vertex), the \textcolor{red}{red} vertex further below is a descendant of the \textcolor{blue}{blue} vertex (in turn the \textcolor{blue}{blue} vertex is an ancestor of the \textcolor{red}{red} vertex), and the \textcolor{orange}{orange} vertices make up the $4$th level of the tree. The dotted lines at the bottom indicate to repeat the above demonstrated pattern ad infinitum.}
    \label{fig:Tree_bi}
\end{figure}

We will quickly cover some common terminology associated with rooted trees.
Given distinct vertices $x,y$ of a rooted tree, if $x$ lies on the path from $y$ to the root, we then say that $x$ is an \emph{ancestor} of $y$ and $y$ is a \emph{descendant} of $x$.
Furthermore, if $x$ and $y$ are also adjacent, we then call $x$ the \emph{parent} (or \emph{mother}) of $y$ and $y$ the \emph{child} (or \emph{daughter}) of $x$.
The set of all vertices with the same distance to the root in a rooted tree is referred to as a \emph{level} (or \emph{generation}) of the tree, and these levels are usually labelled by their distance to the root (so level 0 will consist solely of the root, while level $n \in \omega$ of a rooted tree will consist of all vertices with distance $n$ to the root).
Two rooted trees are considered to be \emph{isomorphic as rooted trees} if there exists a graph isomorphism between them, which maps the root of one tree to the root of the other.
For an illustration of a rooted tree, see Figure \ref{fig:Tree_bi}.

\begin{definition}[Subtree]
    Let $T$ be a rooted tree. A subgraph $T'$ of $T$ is then called a \emph{subtree} if it is itself a rooted tree and it shares its root with $T$.
\end{definition}

\begin{figure}[h]
    \centering
    \includegraphics[width=0.95\linewidth]{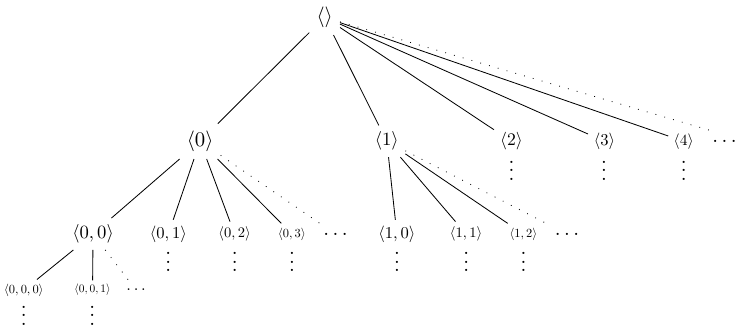}
    \caption{The infinitary rooted tree $T^{(\aleph_0)}$. The dotted lines indicate to continue the pattern}
    \label{fig:Tree_inf}
\end{figure}

\begin{definition}[Regular rooted tree]
    A rooted tree $T$ in which all vertices have the same number of children is called a \emph{regular rooted tree}.
    For any cardinal $\kappa$, we will refer to the regular rooted tree in which every vertex has $\kappa$ many children as the $\kappa$-ary rooted tree and denote it by $T^{(\kappa)}$.
    Examples of such trees are the trivial rooted tree \regtree{0} (the trivial graph with one vertex and no edges), the unary tree \regtree{1} (see Figure \ref{fig:Tree_bi}), and the binary rooted tree \regtree{2} (see Figure \ref{fig:Tree_bi}).
    Additionally, we will usually refer to the $\aleph_0$-ary tree $\regtree{\aleph_0}$ (see Figure \ref{fig:Tree_inf}) as simply the \emph{infinitary rooted tree}.
\end{definition}

Figure \ref{fig:Tree_inf} makes use of a useful way to label the vertices of a rooted tree.
In this convention, the vertices of a rooted tree consist of finite sequences over some set.
The $n$th level of the tree consists of sequences of length $n$ and two vertices have an edge between them if and only if one is an extension of the other by a single element (that is, for any vertex $x$ its initial segment of length $\card{x}-1$ is the parent of $x$).
For some examples, in this convention the infinitary rooted tree $T^{(\aleph_0)}$ can be represented as the set of all finite sequences in $\omega$ (as in Figure \ref{fig:Tree_inf}) and the binary rooted tree $T^{(2)}$ can be represented as the set of all finite sequences in $\set{0,1} = 2$.
This labelling convention is useful when describing infinite rooted trees, as it is not always easy to visually depict what they look like.

With this, we are now ready to introduce the type of semigroups we will be working with in this section.
Given a path $p$ we denote by $\verset{p}$ the set of all vertices on $p$.
Furthermore, recall that given a sequence $p = (x_i)_{i \in \lambda}$ and a chart $f$, then $pf = (x_i f)_{i \in \lambda}$. 
Given another chart $g$, the resulting sequence $pfg$ is simply a sequence over a subset of the ordinal $\lambda$, depending on the overlap between $pf$ and the domain of $g$ (in the usual sense of chart composition, where we compose `wherever possible').

\begin{definition}[Semigroup on rooted tree] \label{def:semigroup_on_tree}
    Let $T$ be a rooted tree on a vertex set $X$ and $P$ the set of all paths in $T$ starting at the root.
    Then the \emph{semigroup of path-bijections on $T$}, denoted $\tree{T}$, is the inverse semigroup of isomorphisms between rooted paths in $T$.
    \begin{align*}
        \tree{T} = \set{f \in I_X \given (\exists p,q \in P)~ \dom{f}=\verset{p} ~\land~ pf=q}
    \end{align*}
\end{definition}

\begin{figure} [h]
    \centering
    \includegraphics[width=0.6\linewidth,height=0.3\textheight]{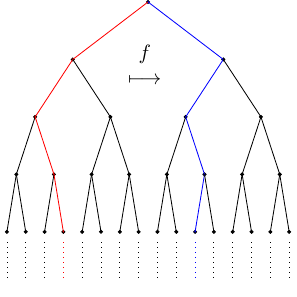}
    \caption{The chart $f$ maps the \textcolor{red}{red} path to the \textcolor{blue}{blue} path while preserving the levels of the rooted tree.}
    \label{fig:Tree_func}
\end{figure}

It is important to note, that path-bijections map paths in a rooted to tree to other paths of the same length (whether finite or infinite) such that levels are preserved (the ordering of a path as a sequence is preserved when mapped to another path).
Also note that given a rooted tree $T$, the path-bijections on $T$ make up a subset of $\pAut{T}$.
Furthermore, definition \ref{def:semigroup_on_tree} claims that the described $\tree{T}$ is an \emph{inverse subsemigroup} of $I_X$.
To verify this claim, let $f,g \in \tree{T}$ be path-bijections on a rooted tree $T$ and $p,q$ the rooted paths for which $\dom{f} = \verset{p}$ and $\im{g} = \verset{q}$.
Then the composite $fg$ is a path-bijection from a subpath of $p$ to a subpath of $q$ (depending on the overlap between $\im{f}$ and $\dom{g}$) and the inverse $\inv{f}$ is a path-bijection from $pf$ to $p$.
Note that whenever $\im{f} \neq \dom{g}$, then $fg$ is a path-bijection between finite rooted paths in $T$.
For an illustration of how path-bijections act on a rooted tree, see Figure \ref{fig:Tree_func}.

These semigroups of path-bijections make up an interesting class of subsemigroups of the symmetric inverse monoid.
They are in some sense very minimal, but we will also show that they can be very large in terms of the Bergman-Shelah preorder.
The first thing we will look at is their topological properties.

\begin{lemma} \label{lem:trees_closed}
    Let $T$ be a rooted tree on a vertex set $X$.
    Then $\tree{T}$ is closed in every $T_1$ shift-continuous topology on $I_X$.
\end{lemma}

\begin{proof}
    By Theorem \ref{I1_theorem} we only have to check that $\tree{T}$ is closed in the topology $\topI_1$.
    Let $P$ be the set of all paths in $T$ starting at the root $r \in X$ and $f \in I_X$ a chart such that $f \notin \tree{T}$.
    Then one of the following must hold:
    \begin{enumerate}[\normalfont ~(i)]
        \item There exists no path $p \in P$ such that $\dom{f} = \verset{p}$; or
        \label{lem:trees_closed/domain}
        
        \item There exists a path $p \in P$ such that $\dom{f} = \verset{p}$, but $pf \notin P$.
        \label{lem:trees_closed/path_bijection}
    \end{enumerate}
    
    For condition \eqref{lem:trees_closed/domain}, there are three ways in which the domain of $f$ might not form a rooted path in $T$.        
    First, if $r \notin \dom{f}$.
    In this case, $V_{r,r} = \set{h \in I_X \given (r,r) \notin h}$ is a basic open neighbourhood of $f$ in $\topI_1$, which has empty intersection with $\tree{T}$.
    Second, if there exist vertices $x,y \in \dom{f}$ such that no path starting at $r$ in $T$ contains both $x$ and $y$.
    In this case, $U_{x,xf} \cap U_{y,xf} = \set{h \in I_X \given (x,xf) \in h \,\land\, (y,yf) \in h}$ is a basic open neighbourhood of $f$, which is disjoint from $\tree{T}$.
    Third, if there exists a vertex $x \in \dom{f}$ such that the parent of $x$ is not in $\dom{f}$.
    In this case, let $m \in X^X$ be a function which maps non-root vertices in $T$ to their parent vertex and maps the root to itself.
    Then $U_{x,xf} \cap V_{xm,xfm} = \set{h \in I_X \given (x,xf) \in h \,\land\, (xm,xfm) \notin h}$ is a basic open neighbourhood of $f$, which has empty intersection with $\tree{T}$.

    For condition \eqref{lem:trees_closed/path_bijection}, we need to find an open neighbourhood around $f$ which has empty intersection with $\tree{T}$, given that $pf$ is not a rooted path in $T$.
    Then either $(r,r) \notin f$ or there exist adjacent vertices $x,y \in \verset{p}$ such that $xf$ and $yf$ are not adjacent in $T$.
    But then either $V_{r,r} = \set{h \in I_X \given (r,r) \notin h}$ or $U_{x,xf} \cap U_{y,yf} = \set{h \in I_X \given (x,xf) \in h \,\land\, (y,yf) \in h}$ are basic open neighbourhoods of $f$, which are disjoint from $\tree{T}$ in $\topI_1$.
    
    We have thus shown that that the complement of $\tree{T}$ in $I_X$ is open, which concludes the proof that $\tree{T}$ is closed in $\topI_1$.
\end{proof}

Lemma \ref{lem:trees_closed} above is entirely independent of the tree in question, meaning that rooted trees provide us with a rich source of closed inverse subsemigroups of $I_X$.
Another property of these semigroups is that they have minimal group substructures.

\begin{lemma} \label{lem:trees_Hsubgroups}
    Let $T$ be a rooted tree on a vertex set $X$.
    Then all subgroups of $\tree{T}$ are trivial.
\end{lemma}

\begin{proof}
    All groups contain exactly one idempotent, which acts as the identity on all elements in the group.
    The idempotents of $I_X$ are the partial identities, so it follows that the idempotents of $\tree{T}$ must be partial identities on paths starting at the root in $T$.
    However, such charts can only act as a two-sided identity for themselves.
    To see this, let $e \in E_X$ be the identity on some path $p$ in $T$ and let $f \in \tree{T}$ be a chart such that $f \notin E_X$ (that is $\dom{f} \neq \im{f}$).
    Then either $ef \neq f$ or $fe \neq f$.
    Hence $\set{e}$ is a maximal subgroup of $\tree{T}$ and the same holds for all other idempotents in $\tree{T}$.
\end{proof}

From Lemma \ref{lem:trees_Hsubgroups} we see that semigroups of path-bijections have very minimal group substructures.
However, this does not mean that they are necessarily `small' subsemigroups of $I_X$ in the Bergman-Shelah sense.

\begin{lemma} \label{lem:tree_inf_bsequal_IX}
    Let $\regtree{\aleph_0}$ be the infinitary rooted tree on a countable vertex set $X$.
    Then $\tree{\regtree{\aleph_0}} \bsequal[I] I_X$.
\end{lemma}

\begin{proof}
    We prove this using Lemma \ref{BS_lemma10}.
    Let $P$ be the set of paths in $\regtree{\aleph_0}$ starting at the root.
    We can then choose $(\alpha_i)_{i \in \omega}$ to be any infinite path in $P$ and $D_{i+1}$ to be the set of all paths in $P$ of length $i$ for each $i \in \omega$.
    Then $\tree{\regtree{\aleph_0}}$ together with the sequence $(\alpha_i)_{i \in \omega}$ and the family of non-empty subsets $\set{D_i \subseteq X^i \given i \in \omega}$ satisfy the hypothesis of Lemma \ref{BS_lemma10}.
\end{proof}

\begin{lemma} \label{lem:tree_bi_bsequal_PStab(omega)}
    Let \regtree{2} be the binary rooted tree on a countable vertex set $X$.
    Then $\tree{\regtree{2}} \bsequal[I] \PStab{\partition_\omega}$.
\end{lemma}

\begin{proof}
    We use Lemma \ref{BS_lemma12} to show that $\tree{\regtree{2}} \bsgeq[I] \PStab{\partition_\omega}$.
    Let $P$ be the set of paths in \regtree{2} starting at the root.
    We can then pick the following sequences and families of sequences:
    \begin{enumerate}[\normalfont (i)]
        \item Let $(\alpha_i)_{i \in \omega}$ to be a subsequence of  any path in $P$ such that $\alpha_i$ belongs to level $i^2$ of \regtree{2} for each $i \in \omega$.

        \item Let $(N_i)_{i \in \omega}$ to be the unbounded sequence $(2^{2i+1})_{i \in \omega}$ of positive integers.

        \item For each $i \in \omega$ let $D_{i+1}$ to be the set of all subsequences of paths in $P$ with length $i^2$ such that for all $(\beta_j)_{j \in i}$ the element $\beta_j$ belongs to level $j^2$ of \regtree{2}.
    \end{enumerate}
    Then $\tree{\regtree{2}}$ together with the sequences $(\alpha_i)_{i \in \omega}$ and $(N_i)_{i \in \omega}$ as well as the family of non-empty subsets $\set{D_i \subseteq X^i \given i \in \omega}$ satisfy the hypothesis of Lemma \ref{BS_lemma12}.

    Finally we show that $\tree{\regtree{2}} \bsleq[I] \PStab{\partition_\omega}$.
    This follows from Lemma \ref{lem:orbits_less->bsleq_unbounded}, since \tree{\regtree{2}} has all finite orbits (therefore all orbits are strictly smaller than $\card{X}$) and $\partition_\omega$ is a $\card{X}$-unbounded partition of $X$.
\end{proof}

Lemmas \ref{lem:tree_inf_bsequal_IX} and \ref{lem:tree_bi_bsequal_PStab(omega)} only hold when $X$ is countably infinite due to how we have defined rooted trees (since all charts in $\tree{T}$ have cardinality at most $\aleph_0$).
One could imagine using a more general definition of rooted trees in which vertices can have infinite distance to the root, but we will not consider such objects here.
Instead we will choose to focus on the case where $X$ is countably infinite throughout the rest of this section.
With Lemmas \ref{lem:tree_inf_bsequal_IX} and \ref{lem:tree_bi_bsequal_PStab(omega)} at hand, we can also prove a previously claimed statement about stabiliser subsemigroups.

\begin{proposition} \label{prop:stab_neq_S(T)}
    Let $T$ be a non-trivial rooted tree on an infinite vertex set $X$.
    Then there exists a finite subset $\Gamma \subseteq X$ such that $\pointStab{\tree{T}}{\Gamma} = \emptyset$.
    Hence, there exists a rooted tree $T$ and a finite subset $\Gamma \subseteq X$ such that $\tree{T} \not\bsequal[I] \pointStab{\tree{T}}{\Gamma}$.
\end{proposition}

\begin{proof}
    Since $T$ is non-trivial, it follows that there exists two distinct points $x,y \in X$ such that $x$ and $y$ do not belong to the same rooted path in $T$.
    Then setting $\Gamma = \set{x,y}$ we get that $\pointStab{\tree{T}}{\Gamma} = \emptyset$, since all elements of $\tree{T}$ must contain a single path in their domain and only finitely many elements from other rooted paths.
    Furthermore, by Lemmas \ref{lem:tree_inf_bsequal_IX} and \ref{lem:tree_bi_bsequal_PStab(omega)} we know that there exist rooted trees (namely \regtree{2} and $\regtree{\aleph_0})$ such that $\tree{T}$ is not equivalent to the empty set under the Bergman-Shelah preorder.
\end{proof}

Proposition \ref{prop:stab_neq_S(T)} tells us that pointwise stabiliser subsemigroups cannot be used for classification of the Bergman-Shelah preorder on $I_X$ similar to how they were used in \cite{Bergman_2006} (not even for closed and inverse subsemigroups specifically, as semigroups on rooted trees satisfy these conditions).
Returning to Lemmas \ref{lem:tree_inf_bsequal_IX} and \ref{lem:tree_bi_bsequal_PStab(omega)}, these are also great examples of applications of Lemmas \ref{BS_lemma10} and \ref{BS_lemma12}.
This naturally raises the question of whether there exists any rooted trees, which can be used for applications of Lemma \ref{BS_lemma14}?
It turns out that there are actually no such trees.

\begin{proposition} \label{prop:no_rooted_tree_BSlemma14}
    Let $T$ be any rooted tree on a countably infinite vertex set $X$.
    Then \tree{T} does not satisfy the hypothesis of Lemma \ref{BS_lemma14}.
\end{proposition}

\begin{proof}
    Let $(\alpha_i)_{i \in \omega}, (\beta_i)_{i \in \omega} \in X^\omega$ be two disjoint sequences of distinct elements.
    If there exists $g \in \tree{T}$ such that $(\alpha_ig)_{i \in \omega} = (\beta_i)_{i \in \omega}$ (which the hypothesis of Lemma \ref{BS_lemma14} demands), then $(\alpha_i)_{i \in \omega}$ and $(\beta_i)_{i \in \omega}$ must be subsequences of distinct paths starting at the root in $T$.
    However, if $(\alpha_i)_{i \in \omega}$ and $(\beta_i)_{i \in \omega}$ are subsequences of distinct paths, then there exists no $f \in \tree{T}$ such that $(\alpha_0 f) = (\alpha_0)$, $(\alpha_1 f) = (\beta_1)$, and $(\alpha_2 f) = (\alpha_2)$, since any two distinct paths can only agree on an initial segment (and $\beta_1 \neq \alpha_1$ by assumption).
    Since Lemma \ref{BS_lemma14} requires the existence of such an $f$, it follows that $\tree{T}$ does not satisfy the lemma.
\end{proof}

Proposition \ref{prop:no_rooted_tree_BSlemma14} shows us that $\tree{\regtree{2}}$ and $\tree{\regtree{\aleph_0}}$ are counterexamples to the converse of Lemma \ref{BS_lemma14}.
Not that we necessarily expected the converse of Lemma \ref{BS_lemma14} to hold, but this confirms that it doesn't even hold for closed inverse subsemigroups of $I_X$ in the topology $\topI_1$.
Perhaps more interestingly, we will later use another rooted tree to produce a similar counterexample for the converse of Lemma \ref{BS_lemma10}.
But first we will prove another important result for determining the Bergman-Shelah ordering of semigroups of path-bijections.
The proof of Lemma \ref{lem:tree_bi_bsequal_PStab(omega)} above highlights an interesting property of semigroups on rooted trees.
Rather than considering the tree as a whole, the proof is only concerned with a substructure of the rooted tree which is not quite a subtree but still very `tree-like'.

\begin{figure}
    \centering
    \includegraphics[width=0.7\linewidth, height=0.3\textheight]{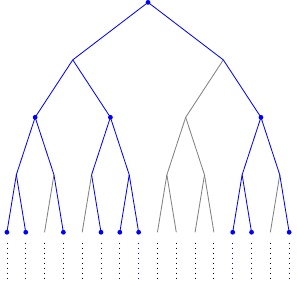}
    \caption{The ternary rooted tree (colour-coded \textcolor{blue}{blue}) as a quasi-subtree of the binary rooted tree. Note that only the highlighted vertices are included in the ternary quasi-subtree.}
     \label{fig:quasi_subtree}
\end{figure}

\begin{definition}[Quasi-subtree]
    Let $T$ be a rooted tree on a vertex set $X$ and $Q$ a rooted tree on some subset $Y \subseteq X$.
    If $Q$ preserves the ancestry/descendancy relations and levels of $T$, then $Q$ is called a \emph{quasi-subtree} of $T$.
    That is, for all vertices $x,y \in Y$, $x$ is an ancestor/descendant of $y$ in $Q$ if and only it is in $T$ and $x$ belongs to the same level as $y$ in $Q$ if and only if it does in $T$.
    For an illustrated example of a quasi-subtree, see Figure \ref{fig:quasi_subtree}.
\end{definition}

These quasi-subtrees will play an important role in determining the Bergman-Shelah ordering of semigroups of path-bijections.
It is worth noting that all subtrees are quasi-subtrees, but not the other way around, and that a quasi-subtree of a quasi-subtree is again a quasi-subtree of the original tree.
Further note that any path in a quasi-subtree is a subsequence of a path in the original tree.
Also, the root of any quasi-subtree $Q$ is simply the unique vertex in $Q$ with the shortest distance to the root in the original tree.
The method, by which we will typically construct quasi-subtrees from a given rooted tree $T$, is by considering substructures generated by any combination of the following two processes: restricting to a subtree of $T$ and `skipping levels' (deleting all the vertices belonging to any number of levels).
It should be clear the process described above results in a quasi-subtree of the rooted tree $T$.

\begin{lemma} \label{lem:quasi_subtree->bsgeq}
    Let $T$ and $R$ be rooted trees and $X$ the union of the vertex sets of $T$ and $R$.
    If there exists a quasi subtree $Q$ in $T$ such that $Q$ and $R$ are isomorphic as rooted trees, then $\tree{T} \bsgeq[I] \tree{R}$.
\end{lemma}

\begin{proof}
    We will start by proving that $\tree{T} \bsgeq[I] \tree{Q}$.
    Let $Y \subseteq X$ denote the vertex set of $Q$.
    Since every path in $Q$ is a subsequence of some path in $T$ and any set of vertices belong to the same level in $Q$ if and only they belong to the same level in $T$, it follows that every chart in $\tree{Q}$ can be extended to a chart in $\tree{T}$.
    So let $e \in I_X$ be the partial identity on $Y$, then $\tree{Q} \subseteq \genset{\tree{T},e}$.

    We then complete the proof by showing that $\tree{Q} \bsgeq[I] \tree{R}$.
    This follows from $Q$ and $R$ being isomorphic as rooted trees.
    Let $f \in I_X$ be the isomorphism from $R$ to $Q$ which maps the root of $R$ to the root of $Q$.
    Then for any two paths $(x_i)_{i \in I}$ and $(y_i)_{i \in I}$ starting at the root in $R$, $(x_if)_{i \in I}$ and $(y_if)_{i \in I}$ are paths starting at the root in $Q$.
    Hence, there exists some $q \in \tree{Q}$ such that $(x_ifq)_{i \in I} = (y_if)_{i \in I}$ and thus $fq\inv{f}$ is the path-bijection from $(x_i)_{i \in I}$ to $(y_i)_{i \in I}$ in $R$.
    Since the above holds for any pair of paths in $R$, we get that $\tree{R} \subseteq \genset{\tree{Q},f,\inv{f}}$.
\end{proof}

Lemma \ref{lem:quasi_subtree->bsgeq} will be our main tool for determining the Bergman-Shelah ordering of semigroups of path-bijections going forward.
We have one last object to introduce before stating the main theorem of this section.
We will construct a rooted tree which is in some sense `sparse' while still allowing for enough path-bijections that the associated semigroup is `large' in the Bergman-Shelah sense.

\begin{figure}[h] 
    \includegraphics[width=1\linewidth]{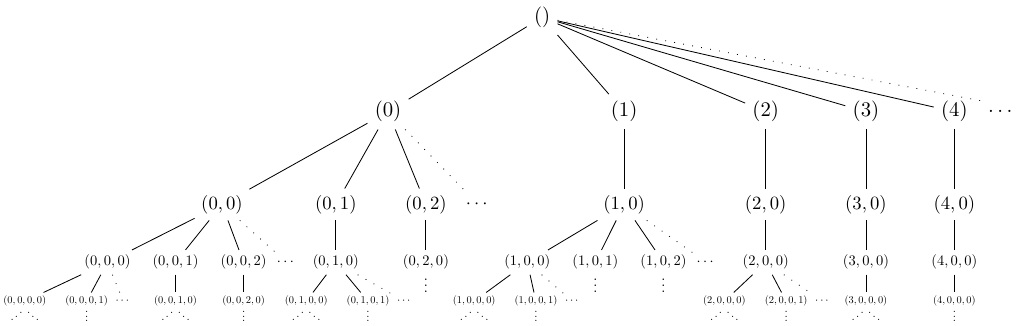}
    \caption{The recursive tree $\regtree{\text{rec}}$. The dotted lines indicate to continue the pattern. Vertical dotted lines specifically mean that the above vertex only has a single child, while split dotted lines indicate that the above vertex has infinitely many children.}
    \label{fig:tree_rec}
\end{figure}

\begin{definition}[Recursive tree $\regtree{\text{rec}}$]
    The \emph{recursive tree} $\regtree{\text{rec}}$ is a rooted tree on a countably infinite vertex set $X$ which can be defined as follows, using the convention of labelling vertices by finite sequences of natural numbers:
    \begin{equation*}
        X = \bigcup_{k \in \omega} \set{(n_i)_{i \in k} \in \omega^k \given (\forall i,j \in k) ~ i<j\leq i+n_i \implies n_j=0}.
    \end{equation*}
    Said in words, the recursive tree $\regtree{\text{rec}}$ is a rooted tree on a countably infinite vertex set $X$ such that the root has infinitely many children and for each $n \in \omega$ there is exactly one child of the root which has $n$ nearest descendants (including the vertex itself) which have exactly one child. Any vertex which is not limited to one child is then the root of recursive subtree (i.e. has infinitely many children which again have the same properties as described above).
    For an illustration of $\regtree{\text{rec}}$, see Figure \ref{fig:tree_rec}.
\end{definition}

\begin{lemma} \label{lem:tree_rec_bsequal_IX}
Let \regtree{\text{rec}} be the recursive tree on a countable vertex set $X$. Then $\tree{\regtree{\text{rec}}} \bsequal[I] I_X$.
\end{lemma}

\begin{proof}
    Let $\partition$ be the partition of $X$ which consists of the `straight segments' of \regtree{\text{rec}}.
    That is, two distinct vertices $x,y \in X$ belong to the same part in $\partition$ if and only if $x$ or $y$ is the only child of the other.
    We will then show that $\tree{\regtree{\text{rec}}} \bsgeq[I] \PStab{\partition}$ and that $\genset{\tree{\regtree{\text{rec}}}, \PStab{\partition}} \bsequal[I] I_X$.
    This will then imply that $\tree{\regtree{\text{rec}}} \bsequal[I] I_X$, since $\tree{\regtree{\text{rec}}} \bsgeq[I] \PStab{\partition}$ implies that there exists a finite subset $U \subseteq I_X$ such that $\PStab{\partition} \subseteq \genset{\tree{\regtree{\text{rec}}}, U}$ and $\genset{\tree{\regtree{\text{rec}}}, \PStab{\partition}} \bsequal[I] I_X$ implies that there exists a finite subset $V \subseteq I_X$ such that $\genset{\tree{\regtree{\text{rec}}}, \PStab{\partition}, V} = I_X$, thus $U \cup V$ is also finite and $\genset{\tree{\regtree{\text{rec}}}, U, V} = I_X$.

    We start by proving that $\tree{\regtree{\text{rec}}} \bsgeq[I] \PStab{\partition}$.
    Consider the subtree $T$ of $\regtree{\text{rec}}$ which consist of those vertices whose label only consist 0s and 1s (so e.g. (), (0), and (0,1,0) but not (0,0,0,2)).
    Now consider the quasi-subtree $Q$ which consist of only the even numbered levels in $T$.
    We will then show that every vertex in $Q$ has at least 2 children.
    The root () has 3 children in $Q$, namely (0,0), (0,1), and (1,0), and given any vertex $(...,n)$ of $Q$ with $n \in \set{0,1}$, either of the following will hold:
    \begin{itemize}
        \item If $n = 0$, then $(...,n,0,0)$, $(...,n,0,1)$, and $(...,n,1,0)$ are all children of $(...,n)$.

        \item If $n = 1$, then $(...,n,0,0)$, and $(...,n,0,1)$ are both children of $(...,n)$.
    \end{itemize}
    Since every vertex in $Q$ has at least 2 children, we can pick a subtree $Q'$ of $Q$ in which every vertex has exactly 2 children.
    That is, $Q'$ is a binary rooted quasi-subtree of $\regtree{\text{rec}}$.
    Since $Q'$ is isomorphic as rooted trees to any binary rooted tree $\regtree{2}$ on $X$, it follows from Lemma \ref{lem:quasi_subtree->bsgeq} that $\tree{\regtree{\text{rec}}} \bsgeq[I] \tree{\regtree{2}}$.
    And since all parts in the partition $\partition$ are finite (hence smaller than $\card{X}$), we get that $\tree{\regtree{2}} \bsgeq[I] \PStab{\partition}$ by Lemmas \ref{lem:tree_bi_bsequal_PStab(omega)} and \ref{lem:orbits_less->bsleq_unbounded}.
    Thus by transitivity of the Bergman-Shelah preorder we get that $\tree{\regtree{\text{rec}}} \bsgeq[I] \PStab{\partition}$.
    Alternatively one could prove this by showing that $\tree{\regtree{\text{rec}}}$ satisfies the conditions of Lemma \ref{BS_lemma12}. We leave this as an exercise to the reader.

    We now prove that $\genset{\tree{\regtree{\text{rec}}}, \PStab{\partition}} \bsequal[I] I_X$.
    Let $(\alpha_i)_{i\in\omega}$ be any infinite path in $\regtree{\text{rec}}$ starting at the root and let $(\beta_i)_{i\in\omega} = (\alpha_{i+1})_{i\in\omega}$ (that is, $(\beta_i)_{i\in\omega}$ is the subpath of $(\alpha_i)_{i\in\omega}$ which only excludes the root).
    We will then show that, in some sense, the symmetric group on $\set{\beta_i \given i \in \omega}$ is contained in $\genset{\tree{\regtree{\text{rec}}}, \PStab{\partition}}$.
    We will prove this by showing that $\genset{\tree{\regtree{\text{rec}}}, \PStab{\partition}}$ contains charts that act as the set of all local permutations on $\set{\beta_i \given i \in \omega}$, where we have ordered the set according to the index $i \in \omega$.
    Let $p \in I_X$ be any local permutation on $\set{\beta_i \given i \in \omega}$.
    Recall that $p$ being local means that for all $i \in \omega$ there exists $j > i$ in $\omega$ such that $p$ maps the initial segment $\set{\beta_k \given k < j}$ to itself.
    An alternative definition would be that there exists a partition $\mathcal{Q}$ of $\set{\beta_i \given i \in \omega}$ which consists of finite intervals in the Well-order given by the index $i \in \omega$, and $p$ then belongs to $\PStab{\mathcal{Q}}$.
    By construction, the recursive tree $\regtree{\text{rec}}$ contains `straight segments' (recall from the definition of the partition $\partition$ above) of every length $n \in \omega$ directly below every vertex of infinite degree.
    So if we exclude the root, there must exist a subpath $(\gamma_i)_{i\in\omega}$ in $\regtree{\text{rec}}$ where the `straight segments' exactly correspond to the intervals in $\mathcal{Q}$.
    So let $f \in \tree{\regtree{\text{rec}}}$ be a path-bijection, which maps $(\beta_i)_{i\in\omega}$ to $(\gamma_i)_{i\in\omega}$.
    Since the partition $\partition$ contains all `straight segments', it follows that there exists a permutation $g \in \PStab{\partition}$ such that $fg\inv{f}$ acts as $p$ on $\set{\beta_i \given i \in \omega}$.
    Since the above construction holds for any local permutation on $\set{\beta_i \given i \in \omega}$, we have shown that $\genset{\tree{\regtree{\text{rec}}}, \PStab{\partition}}$ induces all the local permutations on $\set{\beta_i \given i \in \omega}$.
    It then follows from Lemma \ref{lem:local_generate} that $\genset{\tree{\regtree{\text{rec}}}, \PStab{\partition}}$ induces the symmetric group on $\set{\beta_i \given i \in \omega}$ and it hence follows from Lemma \ref{SymYSub=IX} that $\genset{\tree{\regtree{\text{rec}}}, \PStab{\partition}} \bsequal[I] I_X$, as required.
\end{proof}

Lemma \ref{lem:tree_rec_bsequal_IX} shows us that the recursive tree $\regtree{\text{rec}}$ is an interesting example of a rooted tree.
It is in a sense very `sparse' because it has relatively few vertices of infinite degree (each level only has finitely many), but it is somehow still `large' in that its semigroup of path-bijection finitely generates the entire symmetric inverse monoid on its vertex set.
Furthermore, the reader might have noticed that we did not use Lemma \ref{BS_lemma10} in the proof of Lemma \ref{lem:tree_rec_bsequal_IX}.
That is because $\tree{\regtree{\text{rec}}}$ actually does not satisfy the conditions of the lemma!

\begin{proposition}
    Let $\regtree{\text{rec}}$ be the recursive tree on a countable vertex set $X$.
    Then $\tree{\regtree{\text{rec}}}$ does not satisfy the hypothesis of Lemma \ref{BS_lemma10}.
\end{proposition}

\begin{proof}
    Let $(\alpha_i)_{i\in\omega}$ and $(\beta_i)_i\in\omega$ be sequences in $X$.
    If there exists a chart $f \in \tree{\regtree{\text{rec}}}$ such $(\alpha_if)_{i\in\omega} = (\beta_i)_i\in\omega$, then $(\alpha_i)_{i\in\omega}$ and $(\beta_i)_i\in\omega$ must be subsequences of paths starting at the root in $\regtree{\text{rec}}$ such that $\alpha_i$ and $\beta_i$ belong to the same level for all $i \in \omega$.
    However, each level of the recursive tree $\regtree{\text{rec}}$ only has finitely many vertices of infinite degree.
    This means that it is impossible to find infinitely many $\beta$ and $\gamma_\beta$ in $X$ such that there exists some $g \in \tree{\regtree{\text{rec}}}$ which maps $(\alpha_0,\dots,\alpha_i,\alpha_{i+1},\alpha_{i+2})$ to $(\beta_0,\dots,\beta_i,\beta,\gamma_\beta)$ for any $i \in \omega$.
    We have thus shown that $\tree{\regtree{\text{rec}}}$ cannot satisfy both condition \eqref{BS_lemma10:(ii)} and condition \eqref{BS_lemma10:(iii)} at the same time for any set of sequences in $X$.
\end{proof}

This makes $\tree{\regtree{\text{rec}}}$ a counterexample to the converse of Lemma \ref{BS_lemma10}.
This is an interesting result, as the proof of \cite[Theorem 11]{Bergman_2006} shows that the converse of Lemma \ref{BS_lemma10} holds for closed subgroups of $\Sym{X}$ in the pointwise topology.
However, here we have an example of a closed inverse subsemigroup of $I_X$ in the topology $\topI_1$, for which the converse does not hold.
Furthermore, we will use the property of containing a recursive quasi-subtree (a quasi-subtree isomorphic to the recursive tree $\regtree{\text{rec}}$) as our benchmark for determining, whether a given rooted tree induces a semigroup of path-bijections which finitely generates $I_X$ or not.
We have one last fact about the recursive rooted tree that will prove useful later.

\begin{lemma} \label{lem:rec_tree_growth}
    Let $\regtree{\text{rec}}$ be the recursive tree on a countable vertex set $X$ and $n \geq 1$ a natural number. 
    Then level $n$ of $\regtree{\text{rec}}$ has exactly $2^{n-1}$ many vertices of infinite degree.
\end{lemma}

\begin{proof}
    We will prove this by induction.
    First, we note that the base case of $n = 1$ does indeed hold (level 1 has exactly $2^{1-1} = 1$ vertex of infinite degree).
    So for $n \geq 2$ we assume that for all $k < n$, level $k$ of $\regtree{\text{rec}}$ has exactly $2^{k-1}$ many vertices of infinite degree.
    But from this it follows that there are $2^{n-2}$ vertices of infinite degree on level $n-1$ of $\regtree{\text{rec}}$.
    Each of these infinite degree vertices then have one child of infinite degree on level $n$.
    Hence, there is in some sense a `contribution' of $2^{n-2}$ many vertices of infinite degree from level $n-2$.
    Similarly, there will be $2^{n-3}$ vertices of infinite degree on level $n-2$ (if this level exists), and these vertices have one child whose child (grandchild?) has infinite degree on level $n$.
    Hence, we can say there is a contribution of $2^{n-3}$ vertices of infinite degree from level $n-2$.
    This pattern continues, giving us that level $n-m$ contributes $2^{n-m-1}$ vertices of infinite degree to level $n$ with the exception of level 0, which contributes 1 vertex.
    Thus, we get that the total number of vertices of infinite degree on level $n$ is:
    \begin{equation*}
        1 + \sum_{m=1}^{n-1} 2^{n-m-1} = 1 + (2^{n-1} - 1) = 2^{n-1}
    \end{equation*}
\end{proof}

We are now ready to state the main theorem of this section.

\begin{theorem} \label{thm:tree_trichotomy}
    Let $T$ be a rooted tree on a countably infinite vertex set $X$ and $P$ the set of all paths in $T$ starting at the root. Then one of the following holds:		
        \begin{enumerate}[\normalfont (i)]
            \item $T$ contains a subtree for which every vertex has a descendant of infinite degree and $\tree{T} \bsequal[I] I_X$.
            \label{thm:tree_trichotomy/inf_tree}
			
            \item $T$ does not satisfy \eqref{thm:tree_trichotomy/inf_tree} but contains a subtree for which every vertex has a descendant of degree at least 3 and $\tree{T} \bsequal[I] \PStab{\partition_\omega}$.
            \label{thm:tree_trichotomy/2_tree}
			
            \item $P$ is countable and $\tree{T} \bsequal[I] \emptyset$.
            \label{thm:tree_trichotomy/countable_tree}
        \end{enumerate}
\end{theorem}

We will prove Theorem \ref{thm:tree_trichotomy} by proving a series of intermediate results.
The recursive and binary trees are examples that represent cases \eqref{thm:tree_trichotomy/inf_tree} and \eqref{thm:tree_trichotomy/2_tree} of the theorem respectively.
We will now prove case \eqref{thm:tree_trichotomy/countable_tree}.

\begin{lemma} \label{lem:tree_countable}
Let $T$ be a rooted tree and $P$ the set of all paths in $T$ starting at the root. Then $\tree{T} \bsequal[I] \emptyset$ if and only if $P$ is countable.
\end{lemma}

\begin{proof}
    For any path $p \in P$ all vertices on $p$ have distinct distances to the root. 
    This means that any chart $f \in \tree{T}$ is uniquely determined by the combination of its domain and image.
    Thus $|\tree{T}| \leq |P \times P|$ and $|\tree{T}| \geq |P|$.
    So $\tree{T}$ is countable if and only if $P$ is countable, and it is exactly the countable subsets of $I_X$ that are equivalent to $\emptyset$ under the equivalence relation $\bsequal[I]$.
\end{proof}

Lemmas \ref{lem:tree_inf_bsequal_IX}, \ref{lem:tree_bi_bsequal_PStab(omega)}, and \ref{lem:tree_countable} show us that there indeed exist rooted trees for which the associated semigroups of path-bijections belong to each of the equivalence classes described in Theorem \ref{thm:tree_trichotomy}. That is, we know that these equivalence classes are not devoid of path-bijection semigroups (and we know from Theorem \ref{thm:partition_eqclasses_distinct} that these three classes are distinct).
To show that these three are the only equivalence classes of semigroups of path-bijections, we will introduce the concept of pruning.

\begin{definition}[Pruning process]
    Let $\lambda$ be an ordinal and $(T_\alpha)_{\alpha \in \lambda}$ a $\lambda$-sequence (a sequence of length $\lambda$) of rooted trees.
    If for all $\beta, \gamma \in \lambda$, $\beta > \gamma$ implies that $T_\beta$ is a subtree of $T_\gamma$, then we call $(T_\alpha)_{\alpha \in \lambda}$ a \emph{pruning process} on the initial rooted tree $T_0$.
\end{definition}

The concept of a subtree can alternatively be seen as the act of deleting vertices in the original rooted tree in a way such that all descendants of any deleted vertex are also deleted.
In this line of thinking, a pruning process can be thought of as an iterative process in which some amount of vertices are deleted at each step of the process.

\begin{lemma} \label{lem:pruning_countable}
    Let $T$ be a rooted tree on a countable vertex set $X$ and $(T_\alpha)_{\alpha \in \omega_1}$ a pruning process on $T$ of length $\omega_1$.
    Then there exists some countable $\lambda \in \omega_1$ such that for all $\kappa \geq \lambda$, $T_\kappa = T_\lambda$.
\end{lemma}

\begin{proof}
    Let $\phi$ be a partial function from $X$ to $\omega_1$, which maps each vertex $x \in X$ to the least ordinal $(x)\phi$ such that $x$ is not a vertex of $T_{(x)\phi}$ (or simply deletes $x$ if no such ordinal exists).
    This is well-defined since all ordinals are Well-ordered.
    Since $\dom{\phi} \subseteq X$ is countable, then the image of $\phi$ is also countable.
    But $|\omega_1| = \aleph_1$ is a regular cardinal \cite[Corollary 5.3]{Jech2003}, which means that any countable subset of $\omega_1$ is bounded.
    Hence, $\im{\phi}$ is bounded in $\omega_1$.
    This implies that there exists some countable ordinal $\lambda \in \omega_1$ such that for all $\kappa \geq \lambda$, $T_\kappa$ and $T_\lambda$ share the same vertex set.
    Since $T_\kappa$ is a subtree of $T_\lambda$, this in turn implies that $T_\kappa = T_\lambda$.
\end{proof}

Lemma \ref{lem:pruning_countable} tells us that any pruning process of length $\omega_1$ eventually becomes constant after some countable number of steps (when applied to a tree with countable vertex set).
That is, there exists no rooted tree on a countable vertex set which admits an uncountable down-chain of proper subtrees!
This will be useful shortly.

\begin{definition}[Weak pruning process] \label{def:weak_pruning}
    Let $T$ be a rooted tree on a countable vertex set $X$ and $(T_\alpha)_{\alpha \in \omega_1}$ a pruning process on $T$ of length $\omega_1$ satisfying the following:
    \begin{enumerate}[\normalfont (i)]
        \item For all successor ordinals $\alpha \in \omega_1$, $T_\alpha$ is obtained from $T_{\alpha-1}$ by deleting all non-root vertices with the property that they and all their descendants have degree at most 2 (that is, straight segments where all vertices and their descendants have at most one child).
        
        \item For all limit ordinals $\beta \in \omega_1$, $T_\beta = \bigcap_{\gamma < \beta} T_\gamma$.
    \end{enumerate}
    We then call $(T_\alpha)_{\alpha \in \omega_1}$ the \emph{weak pruning process} on $T$ and denote the eventually constant tree under this process by $\prune_w(T)$.
\end{definition}

Definition \ref{def:weak_pruning} holds the implicit assumption that the intersection of rooted subtrees is itself a subtree.
This is luckily true, since for each vertex $x$ in the intersection all ancestors of $x$ must also be vertices of all the intersected subtrees (as otherwise $x$ could not be a vertex of some of these rooted trees).
And we know that the intersection is not empty, as all of the intersected trees share the same root (hence the root must be in the intersection).
It thus follows from Lemma \ref{lem:pruning_countable} that the resulting `weakly pruned subtree' $\prune_w(T)$ is well-defined.

\begin{lemma} \label{lem:weak_pruning_bsequal}
    Let $T$ be a rooted tree on a countable vertex set $X$.
    Then $\tree{\prune_w(T)} \bsequal[I] \tree{T}$.
\end{lemma}

\begin{proof}
    Since $\prune_w(T)$ is a subtree of $T$, it follows from Lemma \ref{lem:quasi_subtree->bsgeq} that $\tree{\prune_w(T)} \bsleq[I] \tree{T}$.
    So we only need to show that $\tree{\prune_w(T)} \bsgeq[I] \tree{T}$.
    Let $(T_\alpha)_{\alpha \in \omega_1}$ be the weak pruning process on $T$.
    We established in the proof of Lemma \ref{lem:tree_countable} that the number of path-bijections on a rooted tree depends on the number of paths starting at the root.
    It thus follows that if two rooted trees only differ by a countable number of rooted paths, then their associated semigroups also differ by a countable number of path-bijections (since there is only a countable number of path bijections between a countable set of paths and then we only need to add a single path-bijection from one path in the set of differing paths to a path in the collection of shared paths as well as its inverse).
    So we will show that for all $\alpha \in \omega_1$, $T_\alpha$ only differs from $T = T_0$ by a countable number of rooted paths.
    We will do this by transfinite induction.

    Let $\alpha \in \omega_1$ be a countable ordinal.
    Assume that for all $\beta < \alpha$, $T_\beta$ only differs from $T_0$ by a countable number of rooted paths.
    We then consider the following cases.
    \begin{enumerate}[\normalfont (i)]
        \item If $\alpha = 0$, then the statement is trivially true.
        
        \item If $\alpha$ is a successor ordinal, then $T_\alpha$ is obtained from $T_{\alpha-1}$ be deleting all vertices with the property that their descendants are all of degree at most 2.
        But a non-root vertex being of degree 2 means that it is a `straight segment', that is, it only has one child.
        If this is true for all descendants of some vertex $x$, then there is at most countably many paths starting at the root that run through $x$ (since $x$ can have at most countably many children when $X$ is countable).
        So the set of removed paths from $T_{\alpha-1}$ to $T_\alpha$ is a countable union of countable many paths (since there are only countably many vertices).
        Hence $T_\alpha$ only differs by a countable number of paths from $T_{\alpha-1}$, which again only differs by a countable number of paths from $T_0$.

        \item If $\alpha$ is a limit cardinal, then $T_\alpha = \bigcap_{\beta < \alpha} T_\beta$.
        Since each $T_\beta$ only differs from $T_0$ by a countable number of paths, it follows that the set of paths deleted between $T_0$ and $T_\alpha$ is a countable union of countably many paths.
        Hence, $T_\alpha$ only varies from $T_0$ by countably many paths.
    \end{enumerate}
    So we can conclude that $\prune_w(T)$ only differs from $T$ by a countable number of paths.
    Hence, there exists some countable set of path-bijections $U \subseteq \tree{T}$ such that $\tree{T} \subseteq \genset{\tree{\prune_w(T)}, U}$.
\end{proof}

We then get the following dichotomy.

\begin{lemma} \label{lem:weak_dichotomy}
    Let $T$ be a rooted tree on a countable vertex set $X$.
    Then one of the following holds:
    \begin{enumerate}[\normalfont~(i)]
        \item Every vertex in $\prune_w(T)$ has a descendant of degree at least 3 and $\tree{T} \bsgeq[I] \PStab{\partition_\omega}$.
        \label{lem:weak_dichotomy/bsgeq_PStab}

        \item $\prune_w(T) = \regtree{0}$ and $\tree{T} \bsequal[I] \emptyset$.
        \label{lem:weak_dichitomy/bsequal_empty}
    \end{enumerate}
\end{lemma}

\begin{proof}
    We start by showing that each of the statements \eqref{lem:weak_dichotomy/bsgeq_PStab} and \eqref{lem:weak_dichitomy/bsequal_empty} are true on their own.
    \begin{enumerate}[\normalfont (i)]
        \item Assume that every vertex in $\prune_w(T)$ has a descendant of degree at least 3.
        We will then recursively construct a binary rooted quasi-subtree $Q$ of $\prune_w(T)$.
        Let the root of $Q$ be any vertex $x$ of $\prune_w(T)$.
        By our assumption we know that $x$ has a descendant of degree at least 3 (that is, a descendant with at least two children).
        Pick any two children $y_0, y_1$ of this descendant of $x$ be the children of $x$ in $Q$.
        For each of the vertices $y_0$ and $y_1$ there then exist descendants $d_0$ and $d_1$, which each have at least two children.
        However, $d_0$ and $d_1$ might not belong to the same level in $\prune_w(T)$.
        This is no problem however, since all vertices in $\prune_w(T)$ have descendants, so simply pick out two descendants of both $d_0$ and $d_1$ from the maximal level of which the children of each belong to and let these vertices make up level 2 of $Q$
        Repeat this process ad infinitum, where we always pick two children for each vertex.
        Then $Q$ is a binary rooted tree, and it follows from Lemmas \ref{lem:tree_bi_bsequal_PStab(omega)}, \ref{lem:quasi_subtree->bsgeq}, and \ref{lem:weak_pruning_bsequal} that $\tree{T} \bsequal[I] \tree{\prune_w(T)} \bsgeq[I] \PStab{\partition_\omega}$.

        \item Assume that $\prune_w(T)$ is the trivial rooted tree \regtree{0}.
        Then $\prune_w(T)$ contains exactly one path, which is the one starting and ending at the root.
        It then follows from Lemmas \ref{lem:tree_countable} and \ref{lem:weak_pruning_bsequal} that $\tree{T} \bsequal[I] \emptyset$.
    \end{enumerate}
    All that remains is to show that statements \eqref{lem:weak_dichotomy/bsgeq_PStab} and \eqref{lem:weak_dichitomy/bsequal_empty} are indeed a dichotomy.
    Let $T_\alpha$ be some rooted tree in the weak pruning process on $T$ such that not every vertex in $T_\alpha$ has a descendant of degree at least 3 and $T_\alpha$ is not the trivial tree.
    But then $T_\alpha$ would contain some non-root vertex which has only descendants of degree at most 2.
    Then $T_{\alpha+1}$ must be distinct from $T_\alpha$ by the definition of the weak pruning process and hence $T_\alpha \neq \prune_w(T)$.
\end{proof}

Combining Lemmas \ref{lem:tree_countable} and \ref{lem:weak_dichotomy} we see that for any rooted tree $T$, either $T$ has countably many rooted paths and $\tree{T} \bsequal[I] \emptyset$ or $T$ has uncountably many rooted paths and $\tree{T} \bsgeq[I] \PStab{\partition_\omega}$.
Next step would be to prove a dichotomy of rooted trees with uncountably many rooted paths, but we will first take a moment to reflect on some of the results achieved above.
To be exact, we will show that Lemma \ref{lem:pruning_countable} is sharp for the weak pruning process (that is, we answer whether it is truly necessary for us to define this pruning process over all countable ordinals).

\begin{proposition} \label{prop:weak_all_countable_sharp}
    For every countable ordinal $\lambda \in \omega_1$ there exists a rooted tree $T$ on a countable vertex set $X$ such that in the weak pruning process $(T_\alpha)_{\alpha \in \omega_1}$ on $T$, $T_\lambda \neq \prune_w(T)$.
\end{proposition}

\begin{figure}
    \centering
    \includegraphics[scale=1.3]{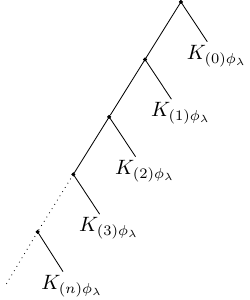}
    \caption{The rooted tree $K_\lambda$ is constructed by `grafting' the rooted trees $\set{K_\kappa \given \kappa \in \lambda}$ onto the unary rooted tree \regtree{1}.}
    \label{fig:Tree_Kalpha}
\end{figure}

\begin{proof}
    We will prove this by transfinite induction.
    So assume that for all ordinals $\kappa < \lambda$ there exists a rooted tree $K_\kappa$, such that the weak pruning process on $K_\kappa$ is not yet constant at step $\kappa$.
    It should be clear that this is true for $\kappa = 0$, as we can then simply let $K_0$ be any rooted tree other than the trivial tree \regtree{0}.
    We will then construct a rooted tree $K_\lambda$ such that the weak pruning process on $K_\lambda$ is not yet constant at step $\lambda$.
    Let $\phi_\lambda$ be a surjection from $\omega$ to $\lambda$ with the property that for all $n \in \omega$ there exists an $m > n$ in $\omega$ such that $(m)\phi_\lambda \geq (n)\phi_\lambda$.
    We then define $K_\lambda$ to be the rooted tree generated by taking the union of the unary rooted tree \regtree{1} with $K_\kappa$ for all $\kappa < \lambda$ and for each $n \in \omega$ we add an edge from the vertex on level $n$ of \regtree{1} to the root of $K_{(n)\phi_\lambda}$ (see Figure \ref{fig:Tree_Kalpha}).
    If we still consider the root of the original unary tree \regtree{1} to be the root of this supertree, then after $\lambda$ steps of the weak pruning process on $K_\lambda$ we get that at the very least the unary tree \regtree{1} still remains, which can then be reduced to the trivial rooted tree \regtree{0} after another step of the process (or there is some $K_\kappa$ which is still not constant after $\lambda$ steps of the weak pruning process).
\end{proof}

Proposition \ref{prop:weak_all_countable_sharp} shows us that iterating the weak pruning process over all countable ordinals is indeed necessary to get the desired results like Lemma \ref{lem:weak_dichotomy}.
The method used in the proof of Proposition \ref{prop:weak_all_countable_sharp} to generate supertrees by taking the union of some base tree with other rooted trees and adding downwards edges we will refer to as `grafting' (as opposed to pruning) and we will make use of this method again later.
We now define our second pruning process.

\begin{definition}[Strong pruning process] \label{def:strong_pruning}
    Let $T$ be a rooted tree on a countable vertex set $X$ and $(T_\alpha)_{\alpha \in \omega_1}$ a pruning process on $T$ of length $\omega_1$ satisfying the following:
    \begin{enumerate}[\normalfont (i)]
        \item For all successor ordinals $\alpha \in \omega_1$, $T_\alpha$ is obtained from $T_{\alpha-1}$ by deleting all non-root vertices with the property that they and all their descendants have finite degree.
        
        \item For all limit ordinals $\beta \in \omega_1$, $T_\beta = \bigcap_{\gamma < \beta} T_\gamma$.
    \end{enumerate}
    We then call $(T_\alpha)_{\alpha \in \omega_1}$ the \emph{strong pruning process} on $T$ and denote the eventually constant tree under this process by $\prune_s(T)$.
\end{definition}

This strong pruning process does not necessarily preserve the Bergman-Shelah equivalence class of the semigroups induced by the rooted trees at each step, however it does satisfy a similar albeit weaker statement.

\begin{lemma} \label{lem:strong_pruning_finite_orbits}
    Let $T$ be a rooted tree on a countable vertex set $X$.
    Then there exists a subset $U \subseteq I_X$ such that $U \bsleq[I] \PStab{\partition_\omega}$ and $\tree{T} \subseteq \genset{\tree{\prune_s(T)}, U}$.
\end{lemma}

\begin{proof}
    As mentioned previously, each step of a pruning process can be thought of as deleting some number rooted paths from the previous tree.
    We can group these `deleted paths' together to form `deleted subtrees' by considering the sets of all rooted paths which share a common deleted vertex as subtrees.
    So consider a rooted tree $R$ and a subtree $R'$ such that $R$ and $R'$ differ by a countable union of deleted subtrees $\bigcup_{\i \in I} R_\i$ with the property that $\tree{R_i} \bsleq[I] \PStab{\partition_\omega}$ for all $i \in I$ (it is worth noting, that such a division of deleted subtrees is not necessarily unique, so we only require that a countable one exists).
    We then claim that there exists a subset $U \subseteq I_X$ such that $U \bsleq[I] \PStab{\partition_\omega}$ and $\tree{R} \subseteq \genset{\tree{R'}, U}$.
    Since $X$ is countable, any finite path in $R$ will only result in a countable number path-bijections, so we will without loss of generality assume that $R_i$ contains an infinite rooted path for all $i \in I$.
    So let $p$ be an infinite rooted path in $R$ and for each $i \in I$ let $f_i$ be a path-bijection from some infinite rooted path in $R_i$ to $p$.
    Then $U = \bigcup_{i \in I} \paa{\tree{R_i} \cup \set{f_i, \inv{f_i}}}$ is a subset of $I_X$ such that $U \bsleq[I] \PStab{\partition_\omega}$, since $I$ is a countable index and $\tree{R_i} \bsleq[I] \PStab{\partition_\omega}$ for all $i \in I$.
    And since $R$ and $R'$ differ exactly by the set of deleted subtrees $\bigcup_{i \in I}$ and the charts $\set{f_i, \inv{f_i} \given i \in I}$ allow for mappings between all of these deleted subtrees, it follows that $\tree{R} \subseteq \genset{\tree{R'}, U}$, as required.
    
    With this we are now ready to prove that there exists a subset $U \subseteq I_X$ such that $U \bsleq[I] \PStab{\partition_\omega}$ and $\tree{T} \subseteq \genset{\tree{\prune_s(T)}, U}$.
    We will do so by using transfinite induction to prove that for all $\alpha \in \omega_1$, $T_\alpha$ differs from $T = T_0$ by a countably union of deleted subtrees $U_\alpha$ such that for all deleted subtrees $R \in U_\alpha$, $\tree{R} \bsleq[I] \PStab{\partition_\omega}$.
    So let $\alpha$ be a countable ordinal and assume that for all $\beta < \alpha$, $T_\beta$ only differs from $T$ by a countable union of deleted subtrees $U_\beta$ such that for all deleted subtrees $R \in U_\beta$, $\tree{R} \bsleq[I] \PStab{\partition_\omega}$.
    We will then show that the desired conclusion holds for $T_\alpha$ by considering the following cases.
    \begin{enumerate}[\normalfont (i)]
        \item If $\alpha = 0$, then the statement is trivially true.

        \item If $\alpha$ is a successor ordinal, then $T_\alpha$ differs from $T_{\alpha-1}$ by the deletion of all vertices with the property that they and all their descendants have finite degree.
        Since all deleted vertices have finite degree it follows that all the semigroups of path-bijections on the deleted subtrees have only finite orbits. 
        So by Lemma \ref{lem:orbits_less->bsleq_unbounded} we get that all the deleted subtrees are $\bsleq[I] \PStab{\partition_\omega}$.
        And the total number of deleted subtrees can at must be countable, since $X$ is countable and the set of deleted vertices are unique to each deleted subtree (if they shared any deleted vertices, they would by definition be the same deleted subtree).
        So $T_\alpha$ only differs from $T_{\alpha-1}$ by a countable set of deleted subtrees $V_\alpha$ and for all $R \in V_\alpha$, $\tree{R} \bsleq[I] \PStab{\partition_\omega}$.
        Hence, $T_\alpha$ differs from $T$ by the countable union of deleted subtrees $U_\alpha = U_\beta \cup V_\alpha$ and for all $R \in U_\alpha$, $\tree{R} \bsleq[I] \PStab{\partition_\omega}$, as required.

        \item If $\alpha$ is a limit ordinal, let $U_\alpha = \bigcup_{\beta \in \alpha} U_\beta$.
        Then $U_\alpha$ is a countable union of countable sets of deleted subtrees that are all $\bsleq[I] \PStab{\partition_\omega}$ and $T_\alpha$ differs from $T$ by this countable set $U_\alpha$, as required.
    \end{enumerate}
    So we can conclude that $\prune_s(T)$ differs from $T$ by a countably union of deleted subtrees $V$ such that for all deleted subtrees $R \in V$, $\tree{R} \bsleq[I] \PStab{\partition_\omega}$.
    From this it follows that there exists a subset $U \subseteq I_X$ such that $U \bsleq[I] \PStab{\partition_\omega}$ and $\tree{T} \subseteq \genset{\tree{\prune_s(T)}, U}$.
\end{proof}

With this we can now use the strong pruning process to get another dichotomy.

\begin{lemma} \label{lem:strong_dichotomy}
    Let $T$ be a rooted tree on a countable vertex set $X$.
    Then one of the following holds:
    \begin{enumerate}[\normalfont~(i)]
        \item Every vertex in $\prune_s(T)$ has a descendant of infinite degree and $\tree{T} \bsequal[I] I_X$.
        \label{lem:strong_dichotomy/bsequal_IX}

        \item $\prune_s(T) = \regtree{0}$ and $\tree{T} \bsleq[I] \PStab{\partition_\omega}$.
        \label{lem:strong_dichitomy/bsleq_PStab(omega)}
    \end{enumerate}
\end{lemma}

\begin{proof}
    We start by showing that each of the statements \eqref{lem:strong_dichotomy/bsequal_IX} and \eqref{lem:strong_dichitomy/bsleq_PStab(omega)} are true on their own.
    \begin{enumerate}[\normalfont (i)]
        \item Assume that every vertex in $\prune_s(T)$ has a descendant of infinite degree.
        We will then show that $\prune_s(T)$ contains a quasi-subtree $Q$, which is isomorphic to the recursive rooted tree $\regtree{\text{rec}}$.
        We will recursively construct each level of the quasi-subtree $Q$.
        We start the recursion by letting the root of $Q$ be any vertex $x$ of $\prune_s(T)$ and level $1$ consist of any infinite set of descendants of $x$ which all belong to the same level in $\prune_s(T)$ (we know that this is possible since all vertices in $\prune_s(T)$ have a descendant of infinite degree).
        Assuming levels 0 though $n-1$ given, we choose the vertices of $\prune_s(T)$ that shall go into the $n$th level of $Q$.
        Since we have already chosen the root and the first level, we will assume that $n \geq 2$.
        To match the recursive tree there must be $2^{n-2}$ vertices on level $n-1$ (Lemma \ref{lem:rec_tree_growth}) which we have to find infinitely many children for. 
        But since all these $2^{n-2}$ many vertices all have a descendant of infinite degree, we can choose the supremum of the level numbers in $\prune_s(T)$ on which all the $2^{n-2}$ vertices have infinitely many descendants and choose infinitely many children for each vertex there.
        For the remaining vertices of level $n-1$ in $Q$ simply chose one child for each from the same level as above.
        Then $Q$ is a recursive rooted tree, and it follows from Lemmas \ref{lem:tree_rec_bsequal_IX} and \ref{lem:quasi_subtree->bsgeq} that $\tree{T} \bsgeq[I] \tree{\prune_s(T)} \bsequal[I] I_X$.

        \item Assume that $\prune_s(T)$ is the trivial rooted tree \regtree{0}.
        It then follows from Lemma \ref{lem:strong_pruning_finite_orbits} that there exists a subset $U \subseteq I_X$ such that $U \bsleq[I] \PStab{\partition_\omega}$ and $\tree{T} \subseteq \genset{\tree{\regtree{0}}, U}$.
        But since \tree{\regtree{0}} is finite, it follows that $\genset{\tree{\regtree{0}}, U} \bsleq[I] \PStab{\partition_\omega}$.
    \end{enumerate}
    All that remains is to show that statements \eqref{lem:strong_dichotomy/bsequal_IX} and \eqref{lem:strong_dichitomy/bsleq_PStab(omega)} are indeed a dichotomy.
    Let $T_\alpha$ be some rooted tree in the weak pruning process on $T$ such that not every vertex in $T_\alpha$ has a descendant of infinite degree and $T_\alpha$ is not the trivial tree.
    Then $T_\alpha$ would contain some non-root vertex, whose descendants all have finite degree.
    But then $T_{\alpha+1}$ must be distinct from $T_\alpha$ by the definition of the strong pruning process and hence $T_\alpha \neq \prune_s(T)$.
\end{proof}

We will now do yet another short intermission to reflect on the above result.
This time we will show that Lemma \ref{lem:pruning_countable} is also sharp for the strong pruning process.

\begin{proposition} \label{prop:strong_all_countable_sharp}
    For every countable ordinal $\lambda \in \omega_1$ there exists a rooted tree $T$ on a countable vertex set $X$ such that in the strong pruning process $(T_\alpha)_{\alpha \in \omega_1}$ on $T$, $T_\lambda \neq \prune_s(T)$.
\end{proposition}

\begin{figure}
    \centering
    \includegraphics[scale=1.3]{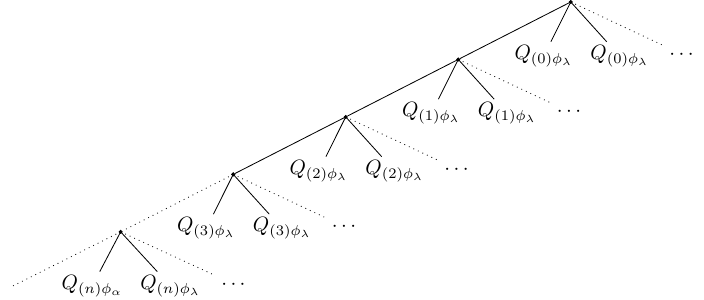}
    \caption{The rooted tree $Q_\lambda$ is constructed by `grafting' infinitely many copies of the rooted trees $\set{Q_\kappa \given \kappa \in \lambda}$ onto each vertex of the unary rooted tree \regtree{1}.}
    \label{fig:Tree_Qalpha}
\end{figure}

\begin{proof}
    We will prove this by transfinite induction.
    So assume that for all ordinals $\kappa < \lambda$ there exists a rooted tree $Q_\kappa$, such that the strong pruning process on $Q_\kappa$ is not yet constant at step $\kappa$.
    It should be clear that this is true for $\kappa = 0$, as we can then simply let $Q_0$ be any non-trivial rooted tree.
    We will then construct a rooted tree $Q_\lambda$ such that the strong pruning process on $Q_\lambda$ is not yet constant at step $\lambda$.
    Let $\phi_\lambda$ be a surjection from $\omega$ to $\lambda$ with the property that for all $n \in \omega$ there exists an $m > n$ in $\omega$ such that $(m)\phi_\lambda \geq (n)\phi_\lambda$.
    We then define $Q_\lambda$ to be the rooted tree generated by taking the unary rooted tree \regtree{1} and for each $n \in \omega$ we graft $\aleph_0$ many copies of $Q_{(n)\phi_\lambda}$ onto the vertex on level $n$ of \regtree{1} (see Figure \ref{fig:Tree_Qalpha}).
    Then after $\lambda$ steps of the strong pruning process on $Q_\lambda$ we get that at the very least the unary tree \regtree{1} still remains, which can then be reduced to the trivial rooted tree \regtree{0} after another step of the process (or there is some $Q_\kappa$ which is still not constant after $\lambda$ steps of the strong pruning process).
\end{proof}

Finally we can give a short proof of Theorem \ref{thm:tree_trichotomy}.

\begin{proof}[Proof of Theorem \ref{thm:tree_trichotomy}]
    This follows from Lemmas \ref{lem:weak_dichotomy} and \ref{lem:strong_dichotomy}.
    We get the following three cases.
    \begin{enumerate}[\normalfont (i)]
        \item $\prune_w(T) = \regtree{0}$ and $\prune_s(T) = \regtree{0}$, so $\tree{T} \bsequal[I] \emptyset$.

        \item $\prune_w(T) \neq \regtree{0}$ and $\prune_s(T) = \regtree{0}$, so $\tree{T} \bsgeq[I] \PStab{\partition_\omega}$ and $\tree{T} \bsleq[I] \PStab{\partition_\omega}$.
        Hence, $\tree{T} \bsequal[I] \PStab{\partition_\omega}$.

        \item $\prune_w(T) \neq \regtree{0}$ and $\prune_s(T) \neq \regtree{0}$, so $\tree{T} \bsequal[I] I_X$.
    \end{enumerate}
    It should be clear that there is no $\prune_w(T) = \regtree{0}$ and $\prune_s(T) \neq \regtree{0}$ case, since everything that gets deleted by the weak pruning process also gets deleted by the strong pruning process.
\end{proof}

\section{New Equivalence Classes}

So far in this chapter we have shown that the four Bergman-Shelah equivalence classes of closed subgroups of \Sym{X} found in \cite{Bergman_2006} are still distinct when extending the setting to $I_X$, and we have only provided examples of subsemigroups of $I_X$ which belong to one of the four equivalence classes (the classes in question being those containing either $\emptyset$, $\PStab{\partition_2}$, $\PStab{\partition_\omega}$, or $I_X$).
This raises a natural question: are these perhaps the only equivalence classes of subsemigroups of $I_X$ under the Bergman-Shelah equivalence relation?
The answer to this question is no, as we will show shortly.
Before that however, we will need the following two lemmas, which pertain to subsemigroups of an inverse semigroup that are themselves not necessarily inverse.

\begin{lemma} \label{lem:inverse_same_relative}
    Let $M$ be an inverse semigroup, $S$ and $T$ subsemigroups of $M$, and $\kappa$ a cardinal.
    Then $S \bsleq[\kappa,M] T$ if and only if $\inv{S} \bsleq[\kappa,M] \inv{T}$.
\end{lemma}

\begin{proof}
    This follows from the fact that inversion is an anti-automorphism on $M$.
    That is, if there exists a subset $U$ of $M$ with $\card{U} < \kappa$ such that $T \subseteq \genset{S, U}$, then it follows that $\inv{T} \subseteq \genset{\inv{S}, \inv{U}}$ (where $\card{\inv{U}} = \card{U}$).
    The same argument holds for the inverse and therefore it also follows that if no such subset $U$ exists for $S$, then neither does one for $\inv{S}$.
\end{proof}

It follows immediately from Lemma \ref{lem:inverse_same_relative} that the analogous statements, in which one replaces $\bsleq$ by either $\bsequal$ or $\bsless$, also hold.
Furthermore, we get the following interesting dichotomy.

\begin{lemma} \label{lem:inverse_dichotomy}
    Let $M$ be an inverse semigroup, $S$ a subsemigroup of $M$, and $\kappa$ an infinite cardinal.
    Then precisely one of the following holds:
    \begin{enumerate}[\normalfont (i)]
        \item $S \bsequal[\kappa,M] \inv{S} \bsequal[\kappa,M] \genset{S, \inv{S}}$.
        \item $S$ and $\inv{S}$ are incomparable under the Bergman-Shelah preorder on $M$ (over the cardinal $\kappa$) and neither is equivalent to any inverse subsemigroups of $M$.
    \end{enumerate}
\end{lemma}

\begin{proof}
    It follows from Lemma \ref{lem:inverse_same_relative} that if $S \bsleq[\kappa,M] \inv{S}$, then $\inv{S} \bsleq[\kappa,M] S$.
    Hence, $S$ and $\inv{S}$ must either be equivalent or incomparable.
    If $S \bsequal[\kappa,M] \inv{S}$, then there exists a subset $U$ with $\card{U} < \kappa$ such that $S \cup \inv{S} \subseteq \genset{S, U}$ and hence $S \bsequal[\kappa,M] \inv{S} \bsequal[\kappa,M] \genset{S, \inv{S}}$ (which also means that $S$ and $\inv{S}$ are equivalent to an inverse semigroup).
    Lastly, if there exists some inverse subsemigroup $T$ of $M$ such that $S \bsequal[\kappa,M] T$, then it again follows from Lemma \ref{lem:inverse_same_relative} that $\inv{S} \bsequal[\kappa,M] \inv{T} = T$ and thus $S \bsequal[\kappa,M] \inv{S}$.
    Hence, if $S$ and $\inv{S}$ are incomparable, then neither can be equivalent to any inverse subsemigroups of $M$.
\end{proof}

The remainder of this chapter will be spent finding examples of subsemigroups that lie outside the four familiar Bergman-Shelah equivalence classes, as well as formulating a new set of conditions that can go into an analogous result to the Bergman-Shelah theorem (Theorem \ref{thm:BS_alt}) for $I_X$ when $X$ is countably infinite.

\begin{proposition} \label{prop:eq_class_EX}
    Let $X$ be a countably infinite set and $E_X$ the semilattice of idempotents in $I_X$.
    Then $\emptyset \bsless[I] E_X \bsless[I] \PStab{\partition_2}$.
\end{proposition}

\begin{proof}
     It follows immediately from the fact that $\card{E_X} = 2^{\card{X}}$ and $\emptyset$ is countable, that $\emptyset \bsless[I] E_X$.
    To show that $E_X \bsless[I] \PStab{\partition_2}$ we must prove that both $E_X \bsleq[I] \PStab{\partition_2}$ and $E_X \not\bsgeq[I] \PStab{\partition_2}$.
    Let $e \in E_X$ be a partial identity and $f \in I_X$ a chart which acts as a bijection from $X$ to $\partition_2$ (that is, for each part $A \in \partition_2$ there exists a unique $x \in X$ such that $xf \in A$).
    We can then pick a permutation $a \in \PStab{\partition_2}$ using the following criteria:
    \begin{enumerate}[\normalfont~(i)]
        \item If $x \in \dom{e} = \im{e}$, then $xfa = xf$.
        \item If $y \in X \setminus \dom{e}$, then $yfa \neq yf$.
    \end{enumerate}
    This is always possible, since $f$ only maps one element into each part in $\partition_2$ and $\partition_2$ is a 2-uniform partition.
    Then $fa\inv{f} = e$, and since $e$ was chosen arbitrarily it follows that $E_X \subseteq \genset{\PStab{\partition_2}, f, \inv{f}}$.

    Finally, to show that $E_X \not\bsgeq[I] \PStab{\partition_2}$, we use the fact that all elements of $\PStab{\partition_2}$ are permutations of $X$.
    Meanwhile, composition with any element of $E_X$ will only result in a restriction to some subset of $X$.
    So for the subset $U \subseteq I_X$, all permutations in $\genset{E_X, U}$ other than the identity must be elements of $\genset{U}$.
    Since $\PStab{\partition_2}$ is uncountable, it therefore follows that there exists no countable subset $U \subseteq I_X$ such that $\PStab{\partition_2} \subseteq \genset{E_X, U}$.
\end{proof}

So we find that the semilattice of idempotents $E_X$ does not belong to any of the four previously described equivalence classes (since both of the two remaining equivalence classes are $\bsgeq[I] \PStab{\partition_2}$).
This on its own would not be a huge issue, as adding just one more equivalence class would not disrupt the existing picture too much, but it turns out that we can find even more equivalence classes.

\begin{proposition} \label{prop:eq_class_E(T)}
    Let $T$ be a rooted tree on a countably infinite vertex set $X$ and $E(T) = \tree{T} \cap E_X$ the semilattice of idempotents in \tree{T}.
    If $T$ has uncountably many rooted paths, then $\emptyset \bsless[I] E(T) \bsless[I] E_X$.
\end{proposition}

\begin{proof}
    It follows from the fact that $T$ has uncountably many rooted paths that $E(T)$ must also be uncountable, and hence $\emptyset \bsless[I] E(T)$.
    Since $E(T)$ is a subsemigroup of $E_X$ it also clearly holds that $E(T) \bsleq[I] E_X$, so all that remains is to prove that $E(T) \not\bsgeq[I] E_X$.
    
    We will do so using a diagonal argument.
    So let $U$ be any countable subset of $I_X$.
    We will then construct a partial identity $e \in E_X$ such that $e \notin \genset{E(T), U}$. 
    We will assume without loss of generality that the identity of $I_X$ is contained in $U$.
    Let $\mathbf{U} = \set{(u_i)_{i \in n} \in \genset{U}^n \given n \in \omega}$ denote the set of all finite tuples of elements in $\genset{U}$ and let $\partition = \set{X_{\mathbf{u}} \given \mathbf{u} \in \mathbf{U}}$ be a partition of $X$ into \card{\mathbf{U}} many moieties (we note that $\mathbf{U}$ is still a countable set).
    To show that $e \notin \genset{E(T), U}$ it will thus suffice to show that for each $\mathbf{u} = (u_0, u_1, \dots, u_n) \in \mathbf{U}$ we can define $e$ on $X_{\mathbf{u}}$ in such a way that $e \notin u_0E(T) u_1E(T) \dots E(T)u_n$, since
    \begin{equation*}
        \genset{E(T), U} = \bigcup_{\mathbf{u} = (u_0, u_1, \dots, u_n) \in \mathbf{U}} u_0E(T) u_1E(T) \dots E(T)u_n.
    \end{equation*}
    The proof then falls into two cases.

    \emph{Case 1: Either $X_\mathbf{u}$ is not contained in the domain of $u_0u_1 \dots u_n$ or there exists a $k \in n$ such that $(X_\mathbf{u}) u_0 u_1 \dots u_k$ is not contained in any rooted path in $T$.}
    We then claim that we can let $e$ act as the identity on all of $X_\mathbf{u}$.
    Since the elements of $E(T)$ can only restrict the domain of the product further, it is clear that $e$ cannot be in $u_0E(T) u_1E(T) \dots E(T)u_n$ if $X_\mathbf{u}$ is not contained in the domain of $u_0u_1 \dots u_n$.
    If there exists a $k \in n$ such that $(X_\mathbf{u}) u_0 u_1 \dots u_k$ is not contained in any rooted path in $T$, then there exist vertices $x,y \in X_\mathbf{u}$ such that no rooted path in $T$ contains both $x u_0 u_1 \dots u_k$ and $y u_0 u_1 \dots u_k$.
    Since all elements of $E(T)$ only contain rooted paths in their domains, it follows that either $x$ or $y$ must be excluded from the domain of all elements in $u_0E(T) u_1E(T) \dots E(T)u_n$. 
    Hence, some points in $X_\mathbf{u}$ will be deleted under the product $u_0E(T) u_1E(T) \dots E(T)u_n$, and letting $e$ act as the identity on $X_\mathbf{u}$ will exclude it from this product.

    \emph{Case 2: $X_\mathbf{u}$ is fully contained in the domain of $u_0u_1 \dots u_n$ and $(X_\mathbf{u}) u_0 u_1 \dots u_k$ is contained in some rooted path of $T$ for all $k \in n$.}
    Since all elements of $E(T)$ act as the identity on rooted paths of $T$, it then follows that for any chart $g \in u_0E(T) u_1E(T) \dots E(T)u_n$ either $g$ has finite domain on $X_\mathbf{u}$ or $\eval{g}_{X_\mathbf{u}} = \eval{u_0 u_1 \dots u_n}_{X_\mathbf{u}}$. That is, either $X_\mathbf{u} \cap \dom{g}$ is finite or $X_\mathbf{u} \subseteq \dom{g}$.
    So let $e$ act as the identity on a moiety of $X_\mathbf{u}$ while deleting the rest of $X_\mathbf{u}$.

    So we can always define $e$ so that it disagrees with all elements of the product set $u_0E(T) u_1E(T) \dots E(T)u_n$ on the subset $X_\mathbf{u}$ for all $\mathbf{u} \in \mathbf{U}$, as required.
\end{proof}

We have now found two subsemigroups of $I_X$ which belong to equivalence classes sitting in between those containing $\emptyset$ and $\PStab{\partition_2}$ respectively.
However, similar issues also occur at the other end of the spectrum.

\begin{proposition} \label{prop:eq_class_I<}
    Let $(X,\leq)$ be a partially ordered set and define the following subsemigroups of $I_X$:
    \begin{itemize}
        \item $I_\leq = \set{f \in I_X \given (\forall x \in \dom{f})~ xf \leq x}$;
        \item $I_\geq = \inv{I_\leq} = \set{f \in I_X \given (\forall x \in \dom{f})~ xf \geq x}$.
    \end{itemize}
    If $(X,\leq)$ is order isomorphic to $\omega$ (with the usual ordering), then the following hold:
    \begin{enumerate}[\normalfont (i)]
        \item If $T \in \set{I_\leq, I_\geq}$, then $\PStab{\partition_\omega} \bsless[I] T \bsless[I] I_X$.
        \label{prop:eq_class_I</above_below}
        
        \item $I_\leq$ and $I_\geq$ are incomparable under the Bergman-Shelah preorder on $I_X$ and neither are equivalent to any inverse subsemigroups of $I_X$.
        \label{prop:eq_class_I</incomparable}

        \item $\genset{I_\leq, I_\geq} = I_X$.
        \label{prop:eq_class_I</generate_IX}

        \item If $S \subseteq I_X$ has only finite cones, then $S \bsleq[I] I_\leq$ and $\inv{S} \bsleq[I] I_\geq$.
        \label{prop:eq_class_I</finite_cones_bsleq}

        \item  If $S \subseteq I_X$ has only finite cones, then $S \not\bsgeq[I] I_\geq$ and $\inv{S} \not\bsgeq[I] I_\leq$.
        \label{prop:eq_class_I</finite_cones_not_bsgeq}
    \end{enumerate}
\end{proposition}

\begin{proof}
    It should be clear that $X$ must be a countable set (since it is order isomorphic to $\omega$) and that $I_\leq$ and $I_\geq$ are indeed subsemigroups of $I_X$ (if $xf \leq x$ and $xfg \leq xf$, then $xfg \leq x$, and the same applies for the $\geq$ relation).
    It should also be clear that $I_\geq$ is indeed the inverse of $I_\leq$, since if the image of any point $x$ under a chart $f$ is greater than or equal to $x$, then the pre-image of any point $y$ in the image of $f$ must be less than or equal to $y$.
    For the structure of this proof we will first prove statements \eqref{prop:eq_class_I</generate_IX}, \eqref{prop:eq_class_I</finite_cones_bsleq}, and \eqref{prop:eq_class_I</finite_cones_not_bsgeq}, then afterwards we will use the previous results to prove statements \eqref{prop:eq_class_I</above_below} and \eqref{prop:eq_class_I</incomparable}.

    \emph{Proof of \eqref{prop:eq_class_I</generate_IX}.}
    Let $f \in I_X$ be any chart.
    We can then pick $g,h \in I_\geq$ such that $\dom{g} = \dom{f}$ and for all $x \in \dom{f}$, $xg \geq xf$ and $xgh = xf$ (this is always possible since $(X,\leq)$ is an infinite up-chain and $xf \leq xg$).
    That is, $f = gh$, and since this holds for any $f \in I_X$ we get that $I_\leq \cup I_\geq$ generates all of $I_X$.

    \emph{Proof of \eqref{prop:eq_class_I</finite_cones_bsleq}.}
    Due to Lemma \ref{lem:inverse_same_relative} we only need to prove that $S \bsleq[I] I_\leq$. 
    Let $f \in I_X$ be a chart such that for all $x \in X$ and all $y \in xS$, $xf > y$ (this is possible for all $x$, since all cones of $S$ are finite and $(X,\leq)$ is an infinite up-chain).
    Then for any $s \in S$ we can find a $g \in I_\leq$ such that $fg = s$, since $f$ maps all elements of $X$ above the maximal element in their image under $S$.
    That is, $S \subseteq \genset{I_\leq, f}$.

    \emph{Proof of \eqref{prop:eq_class_I</finite_cones_not_bsgeq}.}
    Again we only need to prove that $S \not\bsgeq[I] I_\geq$ due to Lemma \ref{lem:inverse_same_relative}.
    We prove this using a diagonal argument.
    So let $U$ be a finite subset of $I_X$ (we will assume without loos of generality that $U$ contains the identity of $I_X$).
    Similar to the proof of Lemma \ref{lem:cones_less->not_bsgeq_uniform} we then define for each $n \in \omega$ a subset $S_n$ of $I_X$ in the following way:
    \begin{equation}
        S_n = \set{s_0u_0 \dots s_{n-1}u_{n-1} \given (\forall i \in n) (s_i \in S ~\land~ u_i \in U)}
    \end{equation}
    Then, as in the proof of Lemma \ref{lem:cones_less->not_bsgeq_uniform}, we get that $S_n$ only has finite cones for all $n \in \omega$.
    So we enumerate the set $X = \set{x_n \given n \in \omega}$ (there is a natural choice for such an enumeration, since $(X,\leq)$ is order isomorphic to $\omega$) and define a chart $f \in I_\geq$ such that for each $n \in \omega$ and all $y \in x_nS_n$, $(x_n)f > y$ (this is always possible, since $S_n$ only has finite cones for all $n \in \omega$).
    Such a chart $f \in I_\geq$ exists, since $(X,\leq)$ is an infinite up-chain, and $f \notin \bigcup_{n \in \omega} S_n = \genset{S, U}$, as required.

    \emph{Proof of \eqref{prop:eq_class_I</incomparable}.}
    It follows from statement \eqref{prop:eq_class_I</finite_cones_not_bsgeq} that $I_\leq \not\bsequal[I] I_\geq$ (since $I_\leq$ only has finite cones).
    Since $I_\leq$ and $I_\geq$ are inverses and not equivalent under the Bergman-Shelah preorder, Lemma \ref{lem:inverse_dichotomy} then tells us that $I_\leq$ and $I_\geq$ are incomparable and neither are equivalent to any inverse subsemigroups of $I_X$. 

    \emph{Proof of \eqref{prop:eq_class_I</above_below}.}
    Since $\PStab{\partition_\omega}$ and $I_X$ are inverse semigroups, it follows from \eqref{prop:eq_class_I</incomparable} that both are are not equivalent to neither $I_\leq$ nor $I_\geq$ (and thus $T \bsless[I] I_X$ for $T \in \set{I_\leq, I_\geq}$ since $T \subseteq I_X$).
    And it follows from statements \eqref{prop:eq_class_I</finite_cones_bsleq} and \eqref{prop:eq_class_I</finite_cones_not_bsgeq} that $\PStab{\partition_\omega} \bsleq[I] I_\leq$ and $\PStab{\partition_\omega} \not\bsgeq[I] I_\geq$ respectively (since $\PStab{\partition_\omega}$ only has finite orbits).
    It then follows from Lemma \ref{lem:inverse_same_relative} that $\PStab{\partition_\omega} \bsless[I] T$ for $T \in \set{I_\leq, I_\geq}$, since $\PStab{\partition_\omega}$ is an inverse semigroup.
\end{proof}

Propositions \ref{prop:eq_class_EX}, \ref{prop:eq_class_E(T)}, and \ref{prop:eq_class_I<} have shown us that there are at least four more equivalence classes of subsemigroups of $I_X$ under the Bergman-Shelah equivalence relation on top of the four classes inherited from \Sym{X}.
And we cannot easily discard these results, as it turns out that these semigroups are closed in the topology $\topI_1$.

\begin{lemma} \label{lem:EX_E(T)_I<_closed}
    Let $X$ be a set, $T$ a rooted tree on $X$, and $\leq$ a partial order on $X$.
    Then the semigroups $E_X$, $E(T)$, $I_\leq$, and $I_\geq$ are all closed in every $T_1$ shift-continuous topology on $I_X$
\end{lemma}

\begin{proof}
    Due to Theorem \ref{I1_theorem} we only have to check that the above-mentioned semigroups are closed in the topology $\topI_1$.
    
    If $f \in I_X$ is a chart such that $f \notin E_X$, then there exists some $x \in \dom{f}$ such that $xf \neq x$.
    Then $U_{x,xf} = \set{h \in I_X \given (x,xf) \in h}$ is a basic open neighbourhood of $f$ in $\topI_1$, which has empty intersection with $E_X$.

    $E(T)$ is the intersection of $\tree{T}$ and $E_X$, which are both closed in $\topI_1$, and hence $E(T)$ is also closed.

    If $f \in I_X$ is a chart such that $f \notin I_\leq$, then there exists some $x \in \dom{f}$ such that $xf \not\leq x$.
    Then $U_{x,xf} = \set{h \in I_X \given (x,xf) \in h}$ is a basic open neighbourhood of $f$ in $\topI_1$, which has empty intersection with $I_\leq$.

    If $f \in I_X$ is a chart such that $f \notin I_\geq$, then there exists some $x \in \dom{f}$ such that $xf \not\geq x$.
    Then $U_{x,xf} = \set{h \in I_X \given (x,xf) \in h}$ is a basic open neighbourhood of $f$ in $\topI_1$, which has empty intersection with $I_\geq$.
\end{proof}

So the semigroups $E_X$, $E(T)$, $I_\leq$, $I_\geq$ give us at least four examples of closed subsemigroups of $I_X$, which do not belong to any of the equivalence classes inherited from \Sym{X}.
There is no indication however, that these should be the only new equivalence classes.
In fact, given two rooted trees $T$ and $R$ on a common vertex set $X$, we do not even know how $E(T)$ and $E(R)$ compare under the Bergman-Shelah preorder on $I_X$.
This makes it seem unlikely that one can find a simple finite list of the equivalence classes of the closed subsemigroups of $I_X$ similar to the main result in \cite{Bergman_2006} for closed subgroups of \Sym{X}.
However, looking at all the closed subsemigroups of $I_X$ might not be the correct thing to do, when making a comparison to the closed subgroups of \Sym{X}.

Recall that the closed subgroups of \Sym{X} in the pointwise topology are exactly the groups of automorphisms of relational structures (see \cite{cameron2009oligomorphic}).
This means that in the statement of the main theorem by Bergman and Shelah in \cite{Bergman_2006} (Theorem \ref{BS_mainTheorem} or Theorem \ref{thm:BS_alt}), one could replace `closed subgroup of \Sym{X}' by `automorphism group of some relational structure on $X$.'
Meanwhile, we have a result which states that the full inverse submonoids of $I_X$ which are closed in the topology $\topI_1$ (equivalently $\topI_4$) are exactly the monoids of partial automorphisms of relational structures (Theorem \ref{thm:full+inverse+closed=pAut}).
As such, perhaps the correct analogue to the Bergman-Shelah Theorem is to consider the submonoids of $I_X$ which are partial automorphism monoids of relational structures on $X$.

This excludes the semigroups $E(T)$, $I_\leq$, and $I_\geq$, as $E(T)$ is not full and $I_\leq$ as well as $I_\geq$ are not inverse semigroups.
This new classification does not exclude $E_X$, but instead we get that $E_X$ belongs to the new `bottom' equivalence class in the Bergman-Shelah preorder on $I_X$, as all full subsemigroups of $I_X$ contain $E_X$ (this is the definition of being full).
And making this demand that we only consider full submonoids of $I_X$ is not a huge issue, since we have Proposition \ref{prop:closed_I1->full_closed} which states that if a subsemigroup $S$ of $I_X$ is closed in $\topI_1$, then so is $\genset{S, E_X}$.
And if $S \subseteq I_X$ already has the property that $S \bsgeq[I] E_X$, then $S \bsequal[I] \genset{S, E_X}$.
All of the above makes for an argument that we can make a new analogous statement for $I_X$ to that of the Bergman-Shelah Theorem for \Sym{X}.
So we have dealt with the `anomalies' which we have introduced in this chapter, but as of yet we do not have any proof that there are no other equivalence classes other than the ones mentioned previously.
As such, we close this thesis with the following conjecture.

\begin{conjecture} \label{conj:BS_analogous}
    Let $X$ be a countably infinite set and $M$ a full inverse submonoid of $I_X$ which is closed in the topology $\topI_1$.
    Then one of the following holds:
    \begin{enumerate}[\normalfont (i)]
        \item $M \bsequal[I] I_X$;
        \item $M \bsequal[I] \PStab{\partition_\omega}$;
        \item $M \bsequal[I] \PStab{\partition_2}$; or
        \item $M \bsequal[I] E_X$.
    \end{enumerate}
\end{conjecture}

	\clearpage 
    \addcontentsline{toc}{chapter}{Bibliography}
	\bibliographystyle{naturemag}
	\bibliography{Kilder}

\begin{thebibliography}{10}
\expandafter\ifx\csname url\endcsname\relax
  \def\url#1{\texttt{#1}}\fi
\expandafter\ifx\csname urlprefix\endcsname\relax\def\urlprefix{URL }\fi
\providecommand{\bibinfo}[2]{#2}
\providecommand{\eprint}[2][]{\url{#2}}

\bibitem{artin2016geometric}
\bibinfo{author}{Artin, E.}
\newblock \emph{\bibinfo{title}{Geometric algebra}}
  (\bibinfo{publisher}{Courier Dover Publications}, \bibinfo{year}{2016}).

\bibitem{weinstein1996groupoids}
\bibinfo{author}{Weinstein, A.}
\newblock \bibinfo{title}{Groupoids: unifying internal and external symmetry}.
\newblock \emph{\bibinfo{journal}{Notices of the AMS}}
  \textbf{\bibinfo{volume}{43}}, \bibinfo{pages}{744--752}
  (\bibinfo{year}{1996}).

\bibitem{wagner1952groups}
\bibinfo{author}{Wagner, V.~V.}
\newblock \bibinfo{title}{Generalized groups}.
\newblock \emph{\bibinfo{journal}{Doklady Akademit Nauk SSSR}}
  \textbf{\bibinfo{volume}{84}} (\bibinfo{year}{1952}).

\bibitem{preston1954inverse}
\bibinfo{author}{Preston, G.~B.}
\newblock \bibinfo{title}{Inverse semi-groups}.
\newblock \emph{\bibinfo{journal}{Journal of the London Mathematical Society}}
  \textbf{\bibinfo{volume}{1}}, \bibinfo{pages}{396--403}
  (\bibinfo{year}{1954}).

\bibitem{preston1954ideals}
\bibinfo{author}{Preston, G.}
\newblock \bibinfo{title}{Inverse semi-groups with minimal right ideals}.
\newblock \emph{\bibinfo{journal}{Journal of the London Mathematical Society}}
  \textbf{\bibinfo{volume}{1}}, \bibinfo{pages}{404--411}
  (\bibinfo{year}{1954}).

\bibitem{preston1954representations}
\bibinfo{author}{Preston, G.~B.}
\newblock \bibinfo{title}{Representations of inverse semi-groups}.
\newblock \emph{\bibinfo{journal}{Journal of the London Mathematical Society}}
  \textbf{\bibinfo{volume}{1}}, \bibinfo{pages}{411--419}
  (\bibinfo{year}{1954}).

\bibitem{paterson2012groupoids}
\bibinfo{author}{Paterson, A.}
\newblock \emph{\bibinfo{title}{Groupoids, inverse semigroups, and their
  operator algebras}}, vol. \bibinfo{volume}{170} (\bibinfo{publisher}{Springer
  Science \& Business Media}, \bibinfo{year}{2012}).

\bibitem{lawson_inverse_semigroups}
\bibinfo{author}{Lawson, M.~V.}
\newblock \emph{\bibinfo{title}{Inverse semigroups}} (\bibinfo{publisher}{World
  Scientific Publishing Co., Inc., River Edge, NJ}, \bibinfo{year}{1998}).
\newblock \urlprefix\url{https://doi.org/10.1142/9789812816689}.
\newblock \bibinfo{note}{The theory of partial symmetries}.

\bibitem{ball1966maximal}
\bibinfo{author}{Ball, R.~W.}
\newblock \bibinfo{title}{Maximal subgroups of symmetric groups}.
\newblock \emph{\bibinfo{journal}{Transactions of the American Mathematical
  Society}} \textbf{\bibinfo{volume}{121}}, \bibinfo{pages}{393--407}
  (\bibinfo{year}{1966}).

\bibitem{ball1968indices}
\bibinfo{author}{Ball, R.~W.}
\newblock \bibinfo{title}{Indices of maximal subgroups of infinite symmetric
  groups}.
\newblock \emph{\bibinfo{journal}{Proceedings of the American Mathematical
  Society}} \textbf{\bibinfo{volume}{19}}, \bibinfo{pages}{948--950}
  (\bibinfo{year}{1968}).

\bibitem{baumgartner1993maximal}
\bibinfo{author}{Baumgartner, J.~E.}, \bibinfo{author}{Shelah, S.} \&
  \bibinfo{author}{Thomas, S.}
\newblock \bibinfo{title}{Maximal subgroups of infinite symmetric groups}.
\newblock \emph{\bibinfo{journal}{Notre Dame J. Formal Logic}}
  \textbf{\bibinfo{volume}{34}}, \bibinfo{pages}{1--11} (\bibinfo{year}{1993}).

\bibitem{biryukov2000set}
\bibinfo{author}{Biryukov, P.} \& \bibinfo{author}{Mishkin, V.}
\newblock \bibinfo{title}{Set ideals with maximal symmetry group and minimal
  dynamical systems}.
\newblock \emph{\bibinfo{journal}{Bulletin of the London Mathematical Society}}
  \textbf{\bibinfo{volume}{32}}, \bibinfo{pages}{39--46}
  (\bibinfo{year}{2000}).

\bibitem{brazil1994maximal}
\bibinfo{author}{Brazil, M.}, \bibinfo{author}{Covington, J.},
  \bibinfo{author}{Penttila, T.}, \bibinfo{author}{Praeger, C.~E.} \&
  \bibinfo{author}{Woods, A.~R.}
\newblock \bibinfo{title}{Maximal subgroups of infinite symmetric groups}.
\newblock \emph{\bibinfo{journal}{Proceedings of the London Mathematical
  Society}} \textbf{\bibinfo{volume}{3}}, \bibinfo{pages}{77--111}
  (\bibinfo{year}{1994}).

\bibitem{covington1996some}
\bibinfo{author}{Covington, J.}, \bibinfo{author}{Macpherson, D.} \&
  \bibinfo{author}{Mekler, A.}
\newblock \bibinfo{title}{Some maximal subgroups of infinite symmetric groups}.
\newblock \emph{\bibinfo{journal}{Quarterly Journal of Mathematics}}
  \textbf{\bibinfo{volume}{47}}, \bibinfo{pages}{297--311}
  (\bibinfo{year}{1996}).

\bibitem{macpherson1993large}
\bibinfo{author}{Macpherson, D.}
\newblock \bibinfo{title}{Large subgroups of infinite symmetric groups}.
\newblock In \emph{\bibinfo{booktitle}{Finite and Infinite Combinatorics in
  Sets and Logic}}, \bibinfo{pages}{249--278} (\bibinfo{publisher}{Springer},
  \bibinfo{year}{1993}).

\bibitem{subgroups_macpherson_neumann}
\bibinfo{author}{Macpherson, H.~D.} \& \bibinfo{author}{Neumann, P.~M.}
\newblock \bibinfo{title}{Subgroups of infinite symmetric groups}.
\newblock \emph{\bibinfo{journal}{Journal of the London Mathematical Society}}
  \textbf{\bibinfo{volume}{s2-42}}, \bibinfo{pages}{64--84}
  (\bibinfo{year}{1990}).
\newblock \urlprefix\url{https://doi.org/10.1112/jlms/s2-42.1.64}.

\bibitem{maximal_macpherson_preager}
\bibinfo{author}{Macpherson, H.~D.} \& \bibinfo{author}{Praeger, C.~E.}
\newblock \bibinfo{title}{Maximal subgroups of infinite symmetric groups}.
\newblock \emph{\bibinfo{journal}{Journal of the London Mathematical Society}}
  \textbf{\bibinfo{volume}{s2-42}}, \bibinfo{pages}{85--92}
  (\bibinfo{year}{1990}).
\newblock \urlprefix\url{https://doi.org/10.1112/jlms/s2-42.1.85}.

\bibitem{richman1967maximal}
\bibinfo{author}{Richman, F.}
\newblock \bibinfo{title}{Maximal subgroups of infinite symmetric groups}.
\newblock \emph{\bibinfo{journal}{Canadian Mathematical Bulletin}}
  \textbf{\bibinfo{volume}{10}}, \bibinfo{pages}{375--381}
  (\bibinfo{year}{1967}).

\bibitem{xiuliang1999classification}
\bibinfo{author}{Xiuliang, Y.}
\newblock \bibinfo{title}{A classification of maximal inverse subsemigroups of
  the finite symmetric inverse semigroups}.
\newblock \emph{\bibinfo{journal}{Communications in Algebra}}
  \textbf{\bibinfo{volume}{27}}, \bibinfo{pages}{4089--4096}
  (\bibinfo{year}{1999}).

\bibitem{Bergman_2006}
\bibinfo{author}{Bergman, G.~M.} \& \bibinfo{author}{Shelah, S.}
\newblock \bibinfo{title}{Closed subgroups of the infinite symmetric group}.
\newblock \emph{\bibinfo{journal}{Algebra universalis}}
  \textbf{\bibinfo{volume}{55}}, \bibinfo{pages}{137–173}
  (\bibinfo{year}{2006}).
\newblock \urlprefix\url{http://dx.doi.org/10.1007/s00012-006-1959-z}.

\bibitem{part1topological}
\bibinfo{author}{Elliott, L.} \emph{et~al.}
\newblock \bibinfo{title}{Automatic continuity, unique polish topologies, and
  zariski topologies on monoids and clones}.
\newblock \emph{\bibinfo{journal}{Transactions of the American Mathematical
  Society}}  (\bibinfo{year}{2023}).
\newblock \urlprefix\url{http://dx.doi.org/10.1090/tran/8987}.

\bibitem{halmos1960naive}
\bibinfo{author}{Halmos, P.~R.}
\newblock \emph{\bibinfo{title}{Naive set theory}} (\bibinfo{publisher}{van
  Nostrand}, \bibinfo{year}{1960}).

\bibitem{Jech2003}
\bibinfo{author}{Jech, T.}
\newblock \emph{\bibinfo{title}{Set Theory: The Third Millennium Edition}}
  (\bibinfo{publisher}{Springer}, \bibinfo{year}{2003}).

\bibitem{part0topological}
\bibinfo{author}{Mesyan, Z.}, \bibinfo{author}{Mitchell, J.~D.} \&
  \bibinfo{author}{Péresse, Y.~H.}
\newblock \bibinfo{title}{Topological transformation monoids}
  (\bibinfo{year}{2018}).
\newblock \eprint{1809.04590}.

\bibitem{howie1995fundamentals}
\bibinfo{author}{Howie, J.~M.}
\newblock \emph{\bibinfo{title}{Fundamentals of semigroup theory}}
  (\bibinfo{publisher}{oxford university Press}, \bibinfo{year}{1995}).

\bibitem{hungerford2013abstract}
\bibinfo{author}{Hungerford, T.~W.}
\newblock \emph{\bibinfo{title}{Abstract algebra: an introduction}}
  (\bibinfo{publisher}{Cengage Learning}, \bibinfo{year}{2013}).

\bibitem{aluffi2021algebra}
\bibinfo{author}{Aluffi, P.}
\newblock \emph{\bibinfo{title}{Algebra: Notes from the Underground}}
  (\bibinfo{publisher}{Cambridge University Press}, \bibinfo{year}{2021}).

\bibitem{lipscomb1996symmetric}
\bibinfo{author}{Lipscomb, S.}
\newblock \emph{\bibinfo{title}{Symmetric inverse semigroups}},
  vol.~\bibinfo{volume}{46} (\bibinfo{publisher}{American Mathematical Soc.},
  \bibinfo{year}{1996}).

\bibitem{howie1998relative}
\bibinfo{author}{Howie, J.~M.}, \bibinfo{author}{Ru{\v{s}}kuc, N.} \&
  \bibinfo{author}{Higgins, P.~M.}
\newblock \bibinfo{title}{On relative ranks of full transformation semigroups}.
\newblock \emph{\bibinfo{journal}{Communications in Algebra}}
  \textbf{\bibinfo{volume}{26}}, \bibinfo{pages}{733--748}
  (\bibinfo{year}{1998}).

\bibitem{higgins2003countable}
\bibinfo{author}{Higgins, P.~M.}, \bibinfo{author}{Howie, J.~M.},
  \bibinfo{author}{Mitchell, J.~D.} \& \bibinfo{author}{Ru{\v{s}}kuc, N.}
\newblock \bibinfo{title}{Countable versus uncountable ranks in infinite
  semigroups of transformations and relations}.
\newblock \emph{\bibinfo{journal}{Proceedings of the Edinburgh Mathematical
  Society}} \textbf{\bibinfo{volume}{46}}, \bibinfo{pages}{531--544}
  (\bibinfo{year}{2003}).

\bibitem{higgins2003generating}
\bibinfo{author}{Higgins, P.}, \bibinfo{author}{Mitchell, J.~D.} \&
  \bibinfo{author}{Ru{\v{s}}kuc, N.}
\newblock \bibinfo{title}{Generating the full transformation semigroup using
  order preserving mappings}.
\newblock \emph{\bibinfo{journal}{Glasgow Mathematical Journal}}
  \textbf{\bibinfo{volume}{45}}, \bibinfo{pages}{557--566}
  (\bibinfo{year}{2003}).

\bibitem{sierpinski1935suites}
\bibinfo{author}{Sierpi{\'n}ski, W.}
\newblock \bibinfo{title}{Sur les suites infinies de fonctions d{\'e}finies
  dans les ensembles quelconques}.
\newblock \emph{\bibinfo{journal}{Fundamenta mathematicae}}
  \textbf{\bibinfo{volume}{24}}, \bibinfo{pages}{209--212}
  (\bibinfo{year}{1935}).

\bibitem{galvin1995generating}
\bibinfo{author}{Galvin, F.}
\newblock \bibinfo{title}{Generating countable sets of permutations}.
\newblock \emph{\bibinfo{journal}{Journal of the London Mathematical Society}}
  \textbf{\bibinfo{volume}{51}}, \bibinfo{pages}{230--242}
  (\bibinfo{year}{1995}).

\bibitem{hyde2012sierpi}
\bibinfo{author}{Hyde, J.~T.} \& \bibinfo{author}{P{\'e}resse, Y.}
\newblock \bibinfo{title}{Sierpi$\backslash$'nski rank of the symmetric inverse
  semigroup}.
\newblock \emph{\bibinfo{journal}{arXiv preprint arXiv:1211.6284}}
  (\bibinfo{year}{2012}).

\bibitem{munkrestopology}
\bibinfo{author}{Munkres, J.}
\newblock \emph{\bibinfo{title}{Topology: Second Edition}}
  (\bibinfo{publisher}{Pearson}, \bibinfo{year}{2003}).

\bibitem{kechris2012classical}
\bibinfo{author}{Kechris, A.}
\newblock \emph{\bibinfo{title}{Classical descriptive set theory}}, vol.
  \bibinfo{volume}{156} (\bibinfo{publisher}{Springer-Verlag},
  \bibinfo{year}{1995}).

\bibitem{gaughan1967topological}
\bibinfo{author}{Gaughan, E.~D.}
\newblock \bibinfo{title}{Topological group structures of infinite symmetric
  groups}.
\newblock \emph{\bibinfo{journal}{Proceedings of the National Academy of
  Sciences}} \textbf{\bibinfo{volume}{58}}, \bibinfo{pages}{907--910}
  (\bibinfo{year}{1967}).

\bibitem{kechris2007turbulence}
\bibinfo{author}{Kechris, A.~S.} \& \bibinfo{author}{Rosendal, C.}
\newblock \bibinfo{title}{Turbulence, amalgamation, and generic automorphisms
  of homogeneous structures}.
\newblock \emph{\bibinfo{journal}{Proceedings of the London Mathematical
  Society}} \textbf{\bibinfo{volume}{94}}, \bibinfo{pages}{302--350}
  (\bibinfo{year}{2007}).

\bibitem{bardyla2024classifyingpolishsemigrouptopologies}
\bibinfo{author}{Bardyla, S.}, \bibinfo{author}{Elliott, L.},
  \bibinfo{author}{Mitchell, J.} \& \bibinfo{author}{Péresse, Y.}
\newblock \bibinfo{title}{Classifying the polish semigroup topologies on the
  symmetric inverse monoid} (\bibinfo{year}{2024}).
\newblock \urlprefix\url{https://arxiv.org/abs/2405.20134}.
\newblock \eprint{2405.20134}.

\bibitem{Zach_2007}
\bibinfo{author}{Mesyan, Z.}
\newblock \bibinfo{title}{Generating self-map monoids of infinite sets}.
\newblock \emph{\bibinfo{journal}{Semigroup Forum}}
  \textbf{\bibinfo{volume}{75}}, \bibinfo{pages}{648–675}
  (\bibinfo{year}{2007}).
\newblock \urlprefix\url{http://dx.doi.org/10.1007/s00233-007-0731-9}.

\bibitem{fulltransformationBS_2012}
\bibinfo{author}{Mesyan, Z.}, \bibinfo{author}{Mitchell, J.~D.},
  \bibinfo{author}{Morayne, M.} \& \bibinfo{author}{Péresse, Y.~H.}
\newblock \bibinfo{title}{The bergman-shelah preorder on transformation
  semigroups}.
\newblock \emph{\bibinfo{journal}{Mathematical Logic Quarterly}}
  \textbf{\bibinfo{volume}{58}}, \bibinfo{pages}{424–433}
  (\bibinfo{year}{2012}).
\newblock \urlprefix\url{http://dx.doi.org/10.1002/malq.201200002}.

\bibitem{scott1980representations}
\bibinfo{author}{Scott, L.~L.}
\newblock \bibinfo{title}{Representations in characteristic p}.
\newblock In \emph{\bibinfo{booktitle}{The Santa Cruz Conference on Finite
  Groups (Univ. California, Santa Cruz, Calif., 1979)}},
  vol.~\bibinfo{volume}{37}, \bibinfo{pages}{319--331}
  (\bibinfo{organization}{Amer. Math. Soc. Providence, RI},
  \bibinfo{year}{1980}).

\bibitem{aschbacher1985maximal}
\bibinfo{author}{Aschbacher, M.} \& \bibinfo{author}{Scott, L.}
\newblock \bibinfo{title}{Maximal subgroups of finite groups}.
\newblock \emph{\bibinfo{journal}{Journal of Algebra}}
  \textbf{\bibinfo{volume}{92}}, \bibinfo{pages}{44--80}
  (\bibinfo{year}{1985}).

\bibitem{liebeck1988nan}
\bibinfo{author}{Liebeck, M.~W.}, \bibinfo{author}{Praeger, C.~E.} \&
  \bibinfo{author}{Saxl, J.}
\newblock \bibinfo{title}{On the o'nan-scott theorem for finite primitive
  permutation groups}.
\newblock \emph{\bibinfo{journal}{Journal of the Australian Mathematical
  Society}} \textbf{\bibinfo{volume}{44}}, \bibinfo{pages}{389--396}
  (\bibinfo{year}{1988}).

\bibitem{graham1968maximal}
\bibinfo{author}{Graham, N.}, \bibinfo{author}{Graham, R.} \&
  \bibinfo{author}{Rhodes, J.}
\newblock \bibinfo{title}{Maximal subsemigroups of finite semigroups}.
\newblock \emph{\bibinfo{journal}{Journal of Combinatorial Theory}}
  \textbf{\bibinfo{volume}{4}}, \bibinfo{pages}{203--209}
  (\bibinfo{year}{1968}).

\bibitem{gavrilov1965functional}
\bibinfo{author}{Gavrilov, G.}
\newblock \bibinfo{title}{On functional completeness in countable-valued
  logic}.
\newblock \emph{\bibinfo{journal}{Problemy Kibernetiki}}
  \textbf{\bibinfo{volume}{15}}, \bibinfo{pages}{5--64} (\bibinfo{year}{1965}).

\bibitem{pinsker2005maximal}
\bibinfo{author}{Pinsker, M.}
\newblock \bibinfo{title}{Maximal clones on uncountable sets that include all
  permutations}.
\newblock \emph{\bibinfo{journal}{Algebra Universalis}}
  \textbf{\bibinfo{volume}{54}}, \bibinfo{pages}{129--148}
  (\bibinfo{year}{2005}).

\bibitem{maximal_east_mitchell_peresse}
\bibinfo{author}{East, J.}, \bibinfo{author}{Mitchell, J.~D.} \&
  \bibinfo{author}{Péresse, Y.}
\newblock \bibinfo{title}{Maximal subsemigroups of the semigroup of all
  mappings on an infinite set}.
\newblock \emph{\bibinfo{journal}{Trans. Amer. Math. Soc. 367 (2015),
  1911-1944}}  (\bibinfo{year}{2014}).
\newblock \urlprefix\url{https://doi.org/10.1090/S0002-9947-2014-06110-2}.

\bibitem{gomes1987ranks}
\bibinfo{author}{Gomes, G.~M.} \& \bibinfo{author}{Howie, J.~M.}
\newblock \bibinfo{title}{On the ranks of certain finite semigroups of
  transformations}.
\newblock In \emph{\bibinfo{booktitle}{Mathematical Proceedings of the
  Cambridge Philosophical Society}}, vol. \bibinfo{volume}{101},
  \bibinfo{pages}{395--403} (\bibinfo{organization}{Cambridge University
  Press}, \bibinfo{year}{1987}).

\bibitem{cameron2009oligomorphic}
\bibinfo{author}{Cameron, P.~J.}
\newblock \bibinfo{title}{Oligomorphic permutation groups}.
\newblock In \emph{\bibinfo{booktitle}{Perspectives in Mathematical Sciences
  II: Pure Mathematics}}, \bibinfo{pages}{37--61} (\bibinfo{publisher}{World
  Scientific}, \bibinfo{year}{2009}).

\bibitem{cameron2006homomorphism}
\bibinfo{author}{Cameron, P.~J.} \& \bibinfo{author}{Ne{\v{s}}et{\v{r}}il, J.}
\newblock \bibinfo{title}{Homomorphism-homogeneous relational structures}.
\newblock \emph{\bibinfo{journal}{Combinatorics, probability and computing}}
  \textbf{\bibinfo{volume}{15}}, \bibinfo{pages}{91--103}
  (\bibinfo{year}{2006}).

\bibitem{ramsey1987problem}
\bibinfo{author}{Ramsey, F.~P.}
\newblock \bibinfo{title}{On a problem of formal logic}.
\newblock In \emph{\bibinfo{booktitle}{Classic Papers in Combinatorics}},
  \bibinfo{pages}{1--24} (\bibinfo{publisher}{Springer}, \bibinfo{year}{1987}).

\bibitem{drake1974set}
\bibinfo{author}{Drake, F.}
\newblock \emph{\bibinfo{title}{Set Theory: An Introduction to Large
  Cardinals}}.
\newblock Studies in logic and the foundations of mathematics
  (\bibinfo{publisher}{North-Holland Publishing Company},
  \bibinfo{year}{1974}).

\bibitem{hell1973groups}
\bibinfo{author}{Hell, P.} \& \bibinfo{author}{Ne{\v{s}}et{\v{r}}il, J.}
\newblock \bibinfo{title}{Groups and monoids of regular graphs (and of graphs
  with bounded degrees)}.
\newblock \emph{\bibinfo{journal}{Canadian Journal of Mathematics}}
  \textbf{\bibinfo{volume}{25}}, \bibinfo{pages}{239--251}
  (\bibinfo{year}{1973}).

\end{thebibliography}


\end{document}